Tsymbaliuk Alexander

# Shuffle approach towards quantum affine and toroidal algebras

— Lecture Notes —

# Preface

These are detailed lecture notes of the crash-course on shuffle algebras delivered by the author at Tokyo University of Marine Science and Technology during the second week of March 2019. Historically, the shuffle approach to Drinfeld-Jimbo quantum groups $U_q(\mathfrak{g})$ (embedding the "positive" subalgebras $U_q^>(\mathfrak{g})$ into $q$-deformed shuffle algebras) was first developed independently in the 1990s by J. Green, M. Rosso, and P. Schauenburg. Motivated by similar ideas, B. Feigin and A. Odesskii proposed a shuffle approach to elliptic quantum groups around the same time.

The shuffle algebras featuring in the present notes may be viewed as trigonometric degenerations of the aforementioned Feigin-Odesskii elliptic algebras. They provide combinatorial models for the "positive" subalgebras of various quantum affinized algebras (such as affine or toroidal). Conceptually, they make it possible to work conveniently with the elements of quantum affinized algebras that are given by rather complicated non-commutative polynomials in the original loop generators.

Starting from our joint work with B. Feigin (where an action of the Heisenberg algebra on the equivariant $K$-theory of the Hilbert schemes of points on a plane was constructed, generalizing the classical result of H. Nakajima in the cohomological setup), and being immensely generalized in a profound series of work by A. Neguț (in his study of affine Laumon spaces, Nakajima quiver varieties associated to a cyclic quiver, and $W$-algebras associated to surfaces), the shuffle approach has found major applications in the geometric representation theory within the last decade. Explicitly, while the action of the loop generators is classically given by the simplest Hecke correspondences equipped with tautological line bundles (following Nakajima's geometric construction of representations for quantum affine algebras), the geometry of the corresponding moduli spaces (such as Nakajima quiver varieties or Laumon spaces) allows for a much wider family of natural geometric correspondences that are best handled using the shuffle realization. Yet another interesting algebro-geometric application of the shuffle algebras has recently been presented in joint work with M. Finkelberg in the study of quantized $K$-theoretic Coulomb branches. In the simplest rank one case, the latter also provides a link to type $A$ quantum $Q$-systems, as first noticed by P. Di Francesco and R. Kedem.

These notes consist of three chapters, providing a separate treatment for: (1) the quantum loop algebras of $\mathfrak{sl}_n$ and their generalizations; (2) the quantum toroidal algebras of $\mathfrak{gl}_1$; and (3) the quantum toroidal algebras of $\mathfrak{sl}_n$. The shuffle realization of the corresponding "positive" subalgebras is established as well as of the commutative subalgebras and some combinatorial representations for the toroidal algebras. One of the key techniques involved is that of "specialization maps". Each chapter emphasizes a different aspect of the theory: in the first chapter, shuffle algebras are used to construct a family of new PBWD bases for type $A$ quantum loop algebras and their integral forms; in the second chapter, a geometric interpretation of the Fock modules is provided and the shuffle description of a commutative subalgebra is used to construct an action of the Heisenberg algebra on the equivariant $K$-theory of the Hilbert schemes of points on a plane; in the last chapter, vertex and combinatorial representations of the quantum toroidal algebras of $\mathfrak{sl}_n$ are related using Miki's isomorphism, and shuffle realization is used to explicitly compute commutative Bethe subalgebras and their limits. The latter construction is inspired by B. Enriquez's work relating shuffle algebras to the correlation functions of quantum affinized algebras.

Given size restrictions of the present volume, the focus is mostly on the algebraic aspects of the above theory (with a single exception made in Chapter 2 to provide a flavor of the geometry involved). To the reader interested in the beautiful geometric applications, we highly recommend a series of papers by A. Neguţ. While most of the presented results have previously appeared in the literature, the aim here is to unify the exposition by simplifying some arguments as well as adding new details.

These notes are intended as an introduction to the subject of shuffle algebras with an indication of the key algebraic techniques involved and are appropriate for any interested researcher, including both graduate and advanced undergraduate students.

**Acknowledgments**

I am indebted to Pavel Etingof, Boris Feigin, Michael Finkelberg, and Andrei Neguţ for numerous stimulating discussions and collaboration over the years; to Luan Bezerra and Evgeny Mukhin for their correspondence on quantum affine superalgebras; to Naihuan Jing for correspondence on two-parameter quantum algebras; to Joshua Wen for a discussion of a beautiful application of the results from the last chapter to the study of wreath Macdonald polynomials and operators; and to the anonymous referees for useful suggestions that improved the overall exposition.

Special thanks are due to Hitoshi Konno for inviting me to deliver this crash-course on shuffle algebras at Tokyo University of Marine Science and Technology (Japan) in March 2019; and to Hitoshi Konno, Atsuo Kuniba, Tomoki Nakanishi, and Masato Okado for inviting me to speak at the workshop "Infinite Analysis 2019: Quantum symmetries in integrable systems", held at the University of Tokyo, a week prior to the course. I am also grateful to the audience for their interest and questions.

I gratefully acknowledge NSF grants DMS-1821185 and DMS-2037602.

West Lafayette, Indiana, USA                                               *Alexander Tsymbaliuk*
May 2022

# Contents

# Chapter 1
# Quantum loop $\mathfrak{sl}_n$, its two integral forms, and generalizations

**Abstract** In this chapter, we establish the shuffle algebra realization in the simplest case of type $A$ quantum loop algebra $U_v^>(L\mathfrak{sl}_n)$. As an important application we construct a family of PBWD (Poincaré-Birkhoff-Witt-Drinfeld) bases for the quantum loop algebra $U_v(L\mathfrak{sl}_n)$ in the new Drinfeld realization. The shuffle approach also allows to strengthen this by constructing a family of PBWD bases for the RTT form (arising naturally from a different, historically the first, realization [9] of $U_v(L\mathfrak{sl}_n)$ shown in [4] to be equivalent to the one presented below) and Lusztig form (introduced in [15]) of $U_v(L\mathfrak{sl}_n)$, both defined over the ring $\mathbb{Z}[v, v^{-1}]$. We also discuss the generalizations of these results to two-parameter and super-cases. Let us emphasize that the previously known PBW result (see [1]) for the quantum loop algebras $U_v(L\mathfrak{g})$ (more generally, for quantum affine $U_v(\widehat{\mathfrak{g}})$) was established in the Drinfeld-Jimbo realization of the latter, and thus it is compatible with the Drinfeld-Jimbo triangular decomposition, while ours is compatible with the new Drinfeld triangular decomposition. The exposition of this chapter closely follows the author's paper [19].

## 1.1 Algebras $U_v(L\mathfrak{sl}_n)$, $\mathfrak{U}_v(L\mathfrak{sl}_n)$, $\mathsf{U}_v(L\mathfrak{sl}_n)$ and their bases

In this section, we recall the definition of the quantum loop algebra $U_v(L\mathfrak{sl}_n)$, defined over $\mathbb{C}(v)$, as well as of its two integral forms, defined over $\mathbb{Z}[v, v^{-1}]$. Evoking the triangular decompositions of these algebras, we provide a family of PBWD-type bases for their "positive" subalgebras (the proofs are postponed till Section 1.2).

### 1.1.1 Quantum loop algebra $U_v(L\mathfrak{sl}_n)$ and its PBWD bases

For an integer $n \geq 2$, let $I = \{1, \ldots, n-1\}$, $(c_{ij})_{i,j \in I}$ be the Cartan matrix of $\mathfrak{sl}_n$ (explicitly: $c_{ii} = 2$, $c_{i,i\pm1} = -1$, and $c_{ij} = 0$ if $|i - j| > 1$), and $v$ be a formal variable. Following [3], we define the quantum loop algebra of $\mathfrak{sl}_n$ (a.k.a. the quantum affine algebra with a trivial central charge), denoted by $U_{\boldsymbol{v}}(L\mathfrak{sl}_n)$, to be the associative $\mathbb{C}(\boldsymbol{v})$-algebra generated by $\{e_{i,r}, f_{i,r}, \psi_{i,\pm s}^{\pm}\}_{i\in I}^{r\in\mathbb{Z},s\in\mathbb{N}}$ with the following defining relations:

$$[\psi_i^{\epsilon}(z),\psi_j^{\epsilon'}(w)]=0, \quad \psi_{i,0}^{\pm}\cdot\psi_{i,0}^{\mp}=1, \tag{1.1}$$

$$(z-\boldsymbol{v}^{c_{ij}}w)e_i(z)e_j(w)=(\boldsymbol{v}^{c_{ij}}z-w)e_j(w)e_i(z), \tag{1.2}$$

$$(\boldsymbol{v}^{c_{ij}}z-w)f_i(z)f_j(w)=(z-\boldsymbol{v}^{c_{ij}}w)f_j(w)f_i(z), \tag{1.3}$$

$$(z-\boldsymbol{v}^{c_{ij}}w)\psi_i^{\epsilon}(z)e_j(w)=(\boldsymbol{v}^{c_{ij}}z-w)e_j(w)\psi_i^{\epsilon}(z), \tag{1.4}$$

$$(\boldsymbol{v}^{c_{ij}}z-w)\psi_i^{\epsilon}(z)f_j(w)=(z-\boldsymbol{v}^{c_{ij}}w)f_j(w)\psi_i^{\epsilon}(z), \tag{1.5}$$

$$[e_i(z),f_j(w)]=\frac{\delta_{ij}}{\boldsymbol{v}-\boldsymbol{v}^{-1}}\delta\left(\frac{z}{w}\right)\left(\psi_i^{+}(z)-\psi_i^{-}(z)\right), \tag{1.6}$$

$$e_i(z)e_j(w)=e_j(w)e_i(z) \text{ if } c_{ij}=0,$$
$$[e_i(z_1),[e_i(z_2),e_j(w)]_{\boldsymbol{v}^{-1}}]_{\boldsymbol{v}}+[e_i(z_2),[e_i(z_1),e_j(w)]_{\boldsymbol{v}^{-1}}]_{\boldsymbol{v}}=0 \text{ if } c_{ij}=-1, \tag{1.7}$$

$$f_i(z)f_j(w)=f_j(w)f_i(z) \text{ if } c_{ij}=0,$$
$$[f_i(z_1),[f_i(z_2),f_j(w)]_{\boldsymbol{v}^{-1}}]_{\boldsymbol{v}}+[f_i(z_2),[f_i(z_1),f_j(w)]_{\boldsymbol{v}^{-1}}]_{\boldsymbol{v}}=0 \text{ if } c_{ij}=-1, \tag{1.8}$$

where $[a,b]_x:=ab-x\cdot ba$ and the generating series are defined as follows:

$$e_i(z):=\sum_{r\in\mathbb{Z}}e_{i,r}z^{-r}, \quad f_i(z):=\sum_{r\in\mathbb{Z}}f_{i,r}z^{-r}, \quad \psi_i^{\pm}(z):=\sum_{s\geq 0}\psi_{i,\pm s}^{\pm}z^{\mp s}, \quad \delta(z):=\sum_{r\in\mathbb{Z}}z^r.$$

The last two relations (1.7, 1.8) are usually referred to as quantum Serre relations. It is convenient to use the generators $\{h_{i,\pm r}\}_{r>0}$ instead of $\{\psi_{i,\pm s}^{\pm}\}_{s>0}$, defined via:

$$\exp\left(\pm(\boldsymbol{v}-\boldsymbol{v}^{-1})\sum_{r>0}h_{i,\pm r}z^{\mp r}\right)=\bar{\psi}_i^{\pm}(z):=(\psi_{i,0}^{\pm})^{-1}\psi_i^{\pm}(z) \tag{1.9}$$

with $h_{i,\pm r}\in\mathbb{C}[\psi_{i,0}^{\mp},\psi_{i,\pm 1}^{\pm},\psi_{i,\pm 2}^{\pm},\ldots]$. Then, the relations (1.5, 1.6) are equivalent to:

$$\psi_{i,0}e_{j,l}=\boldsymbol{v}^{c_{ij}}e_{j,l}\psi_{i,0}, \quad [h_{i,k},e_{j,l}]=\frac{[kc_{ij}]_{\boldsymbol{v}}}{k}e_{j,l+k} \quad (k\neq 0), \tag{1.10}$$

$$\psi_{i,0}f_{j,l}=\boldsymbol{v}^{-c_{ij}}f_{j,l}\psi_{i,0}, \quad [h_{i,k},f_{j,l}]=-\frac{[kc_{ij}]_{\boldsymbol{v}}}{k}f_{j,l+k} \quad (k\neq 0), \tag{1.11}$$

where $[s]_{\boldsymbol{v}}=\frac{\boldsymbol{v}^s-\boldsymbol{v}^{-s}}{\boldsymbol{v}-\boldsymbol{v}^{-1}}$ for $s\in\mathbb{Z}$.

Let $U_{\boldsymbol{v}}^{<}(L\mathfrak{sl}_n)$, $U_{\boldsymbol{v}}^{>}(L\mathfrak{sl}_n)$, $U_{\boldsymbol{v}}^{0}(L\mathfrak{sl}_n)$ be the $\mathbb{C}(\boldsymbol{v})$-subalgebras of $U_{\boldsymbol{v}}(L\mathfrak{sl}_n)$ generated respectively by $\{f_{i,r}\}_{i\in I}^{r\in\mathbb{Z}}$, $\{e_{i,r}\}_{i\in I}^{r\in\mathbb{Z}}$, and $\{\psi_{i,\pm s}^{\pm}\}_{i\in I}^{s\in\mathbb{N}}$. The following result is standard:

**Proposition 1.1** *(a) (Triangular decomposition of $U_{\boldsymbol{v}}(L\mathfrak{sl}_n)$) The multiplication map*

$$m: U_{\boldsymbol{v}}^{<}(L\mathfrak{sl}_n)\otimes_{\mathbb{C}(\boldsymbol{v})}U_{\boldsymbol{v}}^{0}(L\mathfrak{sl}_n)\otimes_{\mathbb{C}(\boldsymbol{v})}U_{\boldsymbol{v}}^{>}(L\mathfrak{sl}_n)\longrightarrow U_{\boldsymbol{v}}(L\mathfrak{sl}_n)$$

*is an isomorphism of $\mathbb{C}(\boldsymbol{v})$-vector spaces.*

*(b) The algebra $U_{\boldsymbol{v}}^{>}(L\mathfrak{sl}_n)$ (resp. $U_{\boldsymbol{v}}^{<}(L\mathfrak{sl}_n)$ and $U_{\boldsymbol{v}}^{0}(L\mathfrak{sl}_n)$) is isomorphic to the associative $\mathbb{C}(\boldsymbol{v})$-algebra generated by $\{e_{i,r}\}_{i\in I}^{r\in\mathbb{Z}}$ (resp. $\{f_{i,r}\}_{i\in I}^{r\in\mathbb{Z}}$ and $\{\psi_{i,\pm s}^{\pm}\}_{i\in I}^{s\in\mathbb{N}}$) with the defining relations (1.2, 1.7) (resp. (1.3, 1.8) and (1.1)).*

For completeness of our exposition let us sketch the proof (cf. [12, Theorem 2]).

*Proof* Let $\widetilde{U}_{\boldsymbol{v}}(L\mathfrak{sl}_n)$ be the $\mathbb{C}(\boldsymbol{v})$-algebra generated by $\{e_{i,r}, f_{i,r}, \psi_{i,\pm s}^{\pm}\}_{i\in I}^{r\in\mathbb{Z}, s\in\mathbb{N}}$ with the defining relations (1.1, 1.4–1.6), and define its subalgebras $\widetilde{U}_{\boldsymbol{v}}^{<}, \widetilde{U}_{\boldsymbol{v}}^{0}, \widetilde{U}_{\boldsymbol{v}}^{>}$ as before. Let $I^{>}$ and $J^{>}$ (respectively, $I^{<}$ and $J^{<}$) be the 2-sided ideals of $\widetilde{U}_{\boldsymbol{v}}(L\mathfrak{sl}_n)$ and $\widetilde{U}_{\boldsymbol{v}}^{>}$ (respectively, $\widetilde{U}_{\boldsymbol{v}}(L\mathfrak{sl}_n)$ and $\widetilde{U}_{\boldsymbol{v}}^{<}$) generated by the following families of elements:

$$A_{i,j}^{r,s} = e_{i,r+1}e_{j,s} - \boldsymbol{v}^{c_{ij}}e_{i,r}e_{j,s+1} - \boldsymbol{v}^{c_{ij}}e_{j,s}e_{i,r+1} + e_{j,s+1}e_{i,r} \ \text{ if } \ c_{ij}\neq 0, \quad (1.12)$$

$$B_{i,j}^{r,s} = e_{i,r}e_{j,s} - e_{j,s}e_{i,r} \ \text{ if } \ c_{ij}=0, \quad (1.13)$$

$$C_{i,j}^{r_1,r_2,s} = \underset{r_1,r_2}{\mathrm{Sym}}\left\{e_{i,r_1}e_{i,r_2}e_{j,s} - (\boldsymbol{v}+\boldsymbol{v}^{-1})e_{i,r_1}e_{j,s}e_{i,r_2} + e_{j,s}e_{i,r_1}e_{i,r_2}\right\} \ \text{ if } \ c_{ij}=-1, \quad (1.14)$$

(respectively, by their $f$-analogues) which encode the relations (1.2, 1.7). We have:

**Lemma 1.1** *(a) (Triangular decomposition of $\widetilde{U}_{\boldsymbol{v}}(L\mathfrak{sl}_n)$) The multiplication map $m\colon \widetilde{U}_{\boldsymbol{v}}^{>} \otimes_{\mathbb{C}(\boldsymbol{v})} \widetilde{U}_{\boldsymbol{v}}^{0} \otimes_{\mathbb{C}(\boldsymbol{v})} \widetilde{U}_{\boldsymbol{v}}^{>} \to \widetilde{U}_{\boldsymbol{v}}(L\mathfrak{sl}_n)$ is an isomorphism of $\mathbb{C}(\boldsymbol{v})$-vector spaces.*

*(b) $\widetilde{U}_{\boldsymbol{v}}^{>}$ and $\widetilde{U}_{\boldsymbol{v}}^{<}$ are free associative algebras in $\{e_{i,r}\}_{i\in I}^{r\in\mathbb{Z}}$ and $\{f_{i,r}\}_{i\in I}^{r\in\mathbb{Z}}$, respectively, while $\widetilde{U}_{\boldsymbol{v}}^{0}$ is the algebra generated by $\{\psi_{i,\pm s}^{\pm}\}_{i\in I}^{s\in\mathbb{N}}$ subject only to (1.1).*

*(c) We have $I^{>} = m(\widetilde{U}_{\boldsymbol{v}}^{<} \otimes \widetilde{U}_{\boldsymbol{v}}^{0} \otimes J^{>})$ and $I^{<} = m(J^{<} \otimes \widetilde{U}_{\boldsymbol{v}}^{0} \otimes \widetilde{U}_{\boldsymbol{v}}^{>})$ for $m$ from (a).*  $\square$

*Proof* Parts (a, b) are standard and we leave their proofs to the interested reader. The first equality in (c) is equivalent to $\widetilde{U}_{\boldsymbol{v}}^{<} \cdot \widetilde{U}_{\boldsymbol{v}}^{0} \cdot J^{>}$ being a 2-sided ideal of $\widetilde{U}_{\boldsymbol{v}}(L\mathfrak{sl}_n)$). To this end, it suffices to verify that any $X$ from (1.12)–(1.14) satisfies the following:

$$[X, h_{\iota,l}] \in \widetilde{U}_{\boldsymbol{v}}^{0} \cdot J^{>} \quad \text{for any} \quad \iota \in I, \ l\neq 0, \quad (1.15)$$

$$[X, f_{\iota,k}] \in \widetilde{U}_{\boldsymbol{v}}^{0} \cdot J^{>} \quad \text{for any} \quad \iota \in I, \ k\in\mathbb{Z}. \quad (1.16)$$

The first of those follows from (1.10), e.g. $[h_{\iota,l}, A_{i,j}^{r,s}] = \frac{[lc_{\iota,i}]_{\boldsymbol{v}}}{l}A_{i,j}^{r+l,s} + \frac{[lc_{\iota,j}]_{\boldsymbol{v}}}{l}A_{i,j}^{r,s+l}$. Let us now verify (1.16) for $X = C_{i,j}^{r_1,r_2,s}$ by working with the generating series (the cases $X = A_{i,j}^{r,s}$ or $X = B_{i,j}^{r,s}$ are similar but simpler). According to (1.6) and using $\psi_j(z) = \psi_j^{+}(z) - \psi_j^{-}(z) = \sum_{r\in\mathbb{Z}}\psi_{j,r}z^{-r}$, we get:

$$\big[e_i(z_1)e_i(z_2)e_j(w) - (\boldsymbol{v}+\boldsymbol{v}^{-1})e_i(z_1)e_j(w)e_i(z_2) + e_j(w)e_i(z_1)e_i(z_2), f_i(y)\big] =$$

$$\frac{\delta_{j,\iota}\delta(w/y)}{\boldsymbol{v}-\boldsymbol{v}^{-1}}\Big(e_i(z_1)e_i(z_2)\psi_j(w) - (\boldsymbol{v}+\boldsymbol{v}^{-1})e_i(z_1)\psi_j(w)e_i(z_2) + \psi_j(w)e_i(z_1)e_i(z_2)\Big) +$$

$$\frac{\delta_{i,\iota}\delta(z_2/y)}{\boldsymbol{v}-\boldsymbol{v}^{-1}}\Big(e_i(z_1)\psi_i(z_2)e_j(w) - (\boldsymbol{v}+\boldsymbol{v}^{-1})e_i(z_1)e_j(w)\psi_i(z_2) + e_j(w)e_i(z_1)\psi_i(z_2)\Big) +$$

$$\frac{\delta_{i,\iota}\delta(z_1/y)}{\boldsymbol{v}-\boldsymbol{v}^{-1}}\Big(\psi_i(z_1)e_i(z_2)e_j(w) - (\boldsymbol{v}+\boldsymbol{v}^{-1})\psi_i(z_1)e_j(w)e_i(z_2) + e_j(w)\psi_i(z_1)e_i(z_2)\Big) .$$

Symmetrizing with respect to $z_1 \leftrightarrow z_2$ and pulling the $\psi_*$-factors to the leftmost part by using (1.4), we will end up with the series whose coefficients are linear combinations of $\psi_{j,m}A_{i,i}^{r,s}$ (if $\iota = j$) or $\psi_{i,m}A_{i,j}^{r,s}$ (if $\iota = i$), and thus are in $\widetilde{U}_{\nu}^0 \cdot J^>$. $\square$

This lemma implies Proposition 1.1 since $\widetilde{U}_{\nu}(L\mathfrak{sl}_n)/(I^> + I^<) \simeq U_{\nu}(L\mathfrak{sl}_n)$. $\square$

We now proceed to the construction of a family of bases for $U_{\nu}(L\mathfrak{sl}_n)$. Let $\{\alpha_i\}_{i=1}^{n-1}$ be the standard simple positive roots of $\mathfrak{sl}_n$, and $\Delta^+$ be the set of positive roots: $\Delta^+ = \{\alpha_j + \alpha_{j+1} + \cdots + \alpha_i\}_{1 \le j \le i < n}$. Consider the following total order $\le$ on $\Delta^+$:

$$\alpha_j + \alpha_{j+1} + \cdots + \alpha_i \le \alpha_{j'} + \alpha_{j'+1} + \cdots + \alpha_{i'} \quad \text{iff} \quad j < j' \text{ or } j = j', i \le i' . \quad (1.17)$$

We also pick a total order $\preceq_\beta$ on $\mathbb{Z}$ for any $\beta \in \Delta^+$. This gives rise to the total order $\le$ on $\Delta^+ \times \mathbb{Z}$:

$$(\beta, r) \le (\beta', r') \quad \text{iff} \quad \beta < \beta' \text{ or } \beta = \beta', r \preceq_\beta r' . \quad (1.18)$$

For every pair $(\beta, r) \in \Delta^+ \times \mathbb{Z}$, we choose:

(1) a decomposition $\beta = \alpha_{i_1} + \cdots + \alpha_{i_p}$ such that $[\cdots [e_{\alpha_{i_1}}, e_{\alpha_{i_2}}], \cdots, e_{\alpha_{i_p}}]$ is a nonzero root vector $e_\beta$ of $\mathfrak{sl}_n$ ($e_{\alpha_i}$ being the standard Chevalley generator of $\mathfrak{sl}_n$),
(2) a decomposition $r = r_1 + \cdots + r_p$ with $r_k \in \mathbb{Z}$,
(3) a sequence $(\lambda_1, \ldots, \lambda_{p-1}) \in \{\nu, \nu^{-1}\}^{p-1}$,

and define the *PBWD basis elements* $e_\beta(r) \in U_\nu^>(L\mathfrak{sl}_n)$ and $f_\beta(r) \in U_\nu^<(L\mathfrak{sl}_n)$ via:

$$\begin{aligned} e_\beta(r) &:= [\cdots [[e_{i_1,r_1}, e_{i_2,r_2}]_{\lambda_1}, e_{i_3,r_3}]_{\lambda_2}, \cdots, e_{i_p,r_p}]_{\lambda_{p-1}} , \\ f_\beta(r) &:= [\cdots [[f_{i_1,r_1}, f_{i_2,r_2}]_{\lambda_1}, f_{i_3,r_3}]_{\lambda_2}, \cdots, f_{i_p,r_p}]_{\lambda_{p-1}} . \end{aligned} \quad (1.19)$$

In particular, we have $e_{\alpha_i}(r) = e_{i,r}$ and $f_{\alpha_i}(r) = f_{i,r}$.

*Remark 1.1* (a) We note that $e_\beta(r)$, $f_\beta(r)$ degenerate to the corresponding root generators $e_\beta \otimes t^r$, $f_\beta \otimes t^r$ of $\mathfrak{sl}_n[t, t^{-1}] = \mathfrak{sl}_n \otimes_\mathbb{C} \mathbb{C}[t, t^{-1}]$ as $\nu \to 1$, hence, the terminology.

(b) The following particular choice features manifestly in [10]:

$$\begin{aligned} e_{\alpha_j + \alpha_{j+1} + \cdots + \alpha_i}(r) &:= [\cdots [[e_{j,r}, e_{j+1,0}]_\nu, e_{j+2,0}]_\nu, \cdots, e_{i,0}]_\nu , \\ f_{\alpha_j + \alpha_{j+1} + \cdots + \alpha_i}(r) &:= [\cdots [[f_{j,r}, f_{j+1,0}]_\nu, f_{j+2,0}]_\nu, \cdots, f_{i,0}]_\nu . \end{aligned} \quad (1.20)$$

Let $\mathsf{H}$ be the set of all functions $h: \Delta^+ \times \mathbb{Z} \to \mathbb{N}$ with finite support. The monomials

$$e_h := \overrightarrow{\prod_{(\beta,r) \in \Delta^+ \times \mathbb{Z}}} e_\beta(r)^{h(\beta,r)} \quad \text{and} \quad f_h := \overleftarrow{\prod_{(\beta,r) \in \Delta^+ \times \mathbb{Z}}} f_\beta(r)^{h(\beta,r)} \quad \text{for all} \ \ h \in \mathsf{H} \quad (1.21)$$

will be called the *ordered PBWD monomials* of $U_\nu^>(L\mathfrak{sl}_n)$ and $U_\nu^<(L\mathfrak{sl}_n)$, respectively. Here, the arrows $\to$ and $\leftarrow$ over the product signs indicate that the products are ordered with respect to (1.18) or its opposite, respectively.

Our first main result establishes the PBWD property of $U_\nu^>(L\mathfrak{sl}_n)$ and $U_\nu^<(L\mathfrak{sl}_n)$ (the proof is presented in Section 1.2 and is based on the shuffle approach):

**Theorem 1.1** *The ordered PBWD monomials $\{e_h\}_{h\in\mathsf{H}}$ and $\{f_h\}_{h\in\mathsf{H}}$ form $\mathbb{C}(\boldsymbol{v})$-bases of the "positive" and "negative" subalgebras $U_{\boldsymbol{v}}^>(L\mathfrak{sl}_n)$ and $U_{\boldsymbol{v}}^<(L\mathfrak{sl}_n)$, respectively.*

Let us relabel the Cartan generators via $\{\psi_{i,r}\}_{i\in I}^{r\in\mathbb{Z}}$ with $\psi_{i,r}=\psi_{i,r}^+$ for $r\geq 0$ and $\psi_{i,r}=\psi_{i,r}^-$ for $r<0$, so that $(\psi_{i,0})^{-1}=\psi_{i,0}^-$. Let $\mathsf{H}_0$ denote the set of all functions $h_0\colon I\times\mathbb{Z}\to\mathbb{Z}$ with finite support and such that $h_0(i,r)\geq 0$ for $r\neq 0$. The monomials

$$\psi_{h_0}:=\prod_{(i,r)\in I\times\mathbb{Z}}\psi_{i,r}^{h_0(i,r)}\quad\text{for all}\quad h_0\in\mathsf{H}_0\tag{1.22}$$

will be called the *PBWD monomials* of $U_{\boldsymbol{v}}^0(L\mathfrak{sl}_n)$; we note that these products are unordered due to (1.1). They form a $\mathbb{C}(\boldsymbol{v})$-basis of $U_{\boldsymbol{v}}^0(L\mathfrak{sl}_n)$, due to Proposition 1.1(b).

Combining this with Theorem 1.1 and Proposition 1.1(a), we finally obtain:

**Theorem 1.2** *For any choices (1)–(3) made prior to (1.19), the elements*

$$\left\{f_{h_-}\cdot\psi_{h_0}\cdot e_{h_+}\;\middle|\;h_-\in\mathsf{H}\,,h_0\in\mathsf{H}_0\,,h_+\in\mathsf{H}\right\}$$

*form a $\mathbb{C}(\boldsymbol{v})$-basis of the quantum loop algebra $U_{\boldsymbol{v}}(L\mathfrak{sl}_n)$.*

### 1.1.2 Integral form $\mathfrak{U}_{\boldsymbol{v}}(L\mathfrak{sl}_n)$ and its PBWD bases

Following the above notations, we define $\widetilde{e}_\beta(r)\in U_{\boldsymbol{v}}^>(L\mathfrak{sl}_n)$ and $\widetilde{f}_\beta(r)\in U_{\boldsymbol{v}}^<(L\mathfrak{sl}_n)$ via:

$$\widetilde{e}_\beta(r):=(\boldsymbol{v}-\boldsymbol{v}^{-1})e_\beta(r)\ \text{and}\ \widetilde{f}_\beta(r):=(\boldsymbol{v}-\boldsymbol{v}^{-1})f_\beta(r)\ \text{for any}\ (\beta,r)\in\Delta^+\times\mathbb{Z}\,.\tag{1.23}$$

We also define $\widetilde{e}_h,\widetilde{f}_h$ via (1.21) but using $\widetilde{e}_\beta(r),\widetilde{f}_\beta(r)$ instead of $e_\beta(r),f_\beta(r)$. Finally, let us define integral forms $\mathfrak{U}_{\boldsymbol{v}}^>(L\mathfrak{sl}_n)$ and $\mathfrak{U}_{\boldsymbol{v}}^<(L\mathfrak{sl}_n)$ as the $\mathbb{Z}[\boldsymbol{v},\boldsymbol{v}^{-1}]$-subalgebras of $U_{\boldsymbol{v}}^>(L\mathfrak{sl}_n)$ and $U_{\boldsymbol{v}}^<(L\mathfrak{sl}_n)$ generated by $\{\widetilde{e}_\beta(r)\}_{\beta\in\Delta^+}^{r\in\mathbb{Z}}$ and $\{\widetilde{f}_\beta(r)\}_{\beta\in\Delta^+}^{r\in\mathbb{Z}}$ of (1.23), respectively.

The above definition of $\mathfrak{U}_{\boldsymbol{v}}^>(L\mathfrak{sl}_n),\mathfrak{U}_{\boldsymbol{v}}^<(L\mathfrak{sl}_n)$ seems to have an obvious flaw, due to all the choices (1)–(3) used in the construction (1.19) of $e_\beta(r),f_\beta(r)$ and hence $\widetilde{e}_\beta(r),\widetilde{f}_\beta(r)$ in (1.23). The next result shows that $\mathfrak{U}_{\boldsymbol{v}}^>(L\mathfrak{sl}_n),\mathfrak{U}_{\boldsymbol{v}}^<(L\mathfrak{sl}_n)$ are actually independent of these choices and posses PBWD bases similar to those from Theorem 1.1 (the proof is presented in Section 1.2 and is based on the shuffle approach):

**Theorem 1.3** *(a) $\mathfrak{U}_{\boldsymbol{v}}^>(L\mathfrak{sl}_n)$, $\mathfrak{U}_{\boldsymbol{v}}^<(L\mathfrak{sl}_n)$ are independent of all such choices (1)–(3).*

*(b) The ordered PBWD monomials $\{\widetilde{e}_h\}_{h\in\mathsf{H}}$ and $\{\widetilde{f}_h\}_{h\in\mathsf{H}}$ form bases of the free $\mathbb{Z}[\boldsymbol{v},\boldsymbol{v}^{-1}]$-modules $\mathfrak{U}_{\boldsymbol{v}}^>(L\mathfrak{sl}_n)$ and $\mathfrak{U}_{\boldsymbol{v}}^<(L\mathfrak{sl}_n)$, respectively.*

Let $\mathfrak{U}_{\boldsymbol{v}}^0(L\mathfrak{sl}_n)$ be the $\mathbb{Z}[\boldsymbol{v},\boldsymbol{v}^{-1}]$-subalgebra of $U_{\boldsymbol{v}}^0(L\mathfrak{sl}_n)$ generated by $\{\psi_{i,\pm s}^{\pm}\}_{i\in I}^{s\in\mathbb{N}}$. Due to Proposition 1.1(b), the PBWD monomials $\{\psi_{h_0}\}_{h_0\in\mathsf{H}_0}$ of (1.22) form a basis of the free $\mathbb{Z}[\boldsymbol{v},\boldsymbol{v}^{-1}]$-module $\mathfrak{U}_{\boldsymbol{v}}^0(L\mathfrak{sl}_n)$. Finally, we define the integral form $\mathfrak{U}_{\boldsymbol{v}}(L\mathfrak{sl}_n)$ as the $\mathbb{Z}[\boldsymbol{v},\boldsymbol{v}^{-1}]$-subalgebra of $U_{\boldsymbol{v}}(L\mathfrak{sl}_n)$ generated by $\{\widetilde{e}_\beta(r),\widetilde{f}_\beta(r),\psi_{i,\pm s}^{\pm}\}_{\beta\in\Delta^+,i\in I}^{r\in\mathbb{Z},s\in\mathbb{N}}$.

*Remark 1.2* Since $\mathfrak{U}_{\boldsymbol{v}}(L\mathfrak{sl}_n)$ is generated by $\mathfrak{U}_{\boldsymbol{v}}^0(L\mathfrak{sl}_n), \mathfrak{U}_{\boldsymbol{v}}^>(L\mathfrak{sl}_n), \mathfrak{U}_{\boldsymbol{v}}^<(L\mathfrak{sl}_n)$, it is independent of any choices (used to define $\widetilde{e}_\beta(r), \widetilde{f}_\beta(r)$), due to Theorem 1.3(a).

The following result is proved in [10, Theorem 3.24(b)] using the RTT realization of $U_{\boldsymbol{v}}(L\mathfrak{sl}_n)$ (though in [10] we worked with $\mathbb{C}[\boldsymbol{v}, \boldsymbol{v}^{-1}]$-subalgebras, the argument of *loc.cit.* applies verbatim to $\mathbb{Z}[\boldsymbol{v}, \boldsymbol{v}^{-1}]$-subalgebras instead):

**Theorem 1.4** *[10, Theorem 3.24] For the particular choice (1.20) of the PBWD basis elements $e_\beta(r), f_\beta(r)$ in (1.23), the elements*

$$\left\{ \widetilde{f}_{h_-} \cdot \psi_{h_0} \cdot \widetilde{e}_{h_+} \;\middle|\; h_- \in \mathsf{H}, h_0 \in \mathsf{H}_0, h_+ \in \mathsf{H} \right\}$$

*form a basis of the free $\mathbb{Z}[\boldsymbol{v}, \boldsymbol{v}^{-1}]$-module $\mathfrak{U}_{\boldsymbol{v}}(L\mathfrak{sl}_n)$.*

*Remark 1.3* Let us point out that it is often more convenient to work with the quantum loop algebra $U_{\boldsymbol{v}}(L\mathfrak{gl}_n)$ and its similarly defined integral form $\mathfrak{U}_{\boldsymbol{v}}(L\mathfrak{gl}_n)$. Following [10, Proposition 3.11], $\mathfrak{U}_{\boldsymbol{v}}(L\mathfrak{gl}_n)$ is identified with the RTT integral form $\mathfrak{U}_{\boldsymbol{v}}^{\mathrm{rtt}}(L\mathfrak{gl}_n)$ under the $\mathbb{C}(\boldsymbol{v})$-algebra isomorphism $U_{\boldsymbol{v}}(L\mathfrak{gl}_n) \simeq \mathfrak{U}_{\boldsymbol{v}}^{\mathrm{rtt}}(L\mathfrak{gl}_n) \otimes_{\mathbb{Z}[\boldsymbol{v}, \boldsymbol{v}^{-1}]} \mathbb{C}(\boldsymbol{v})$ of [4].

Combining Theorems 1.3–1.4, we get the triangular decomposition of $\mathfrak{U}_{\boldsymbol{v}}(L\mathfrak{sl}_n)$:

**Corollary 1.1** *The multiplication map*

$$m \colon \mathfrak{U}_{\boldsymbol{v}}^<(L\mathfrak{sl}_n) \otimes_{\mathbb{Z}[\boldsymbol{v}, \boldsymbol{v}^{-1}]} \mathfrak{U}_{\boldsymbol{v}}^0(L\mathfrak{sl}_n) \otimes_{\mathbb{Z}[\boldsymbol{v}, \boldsymbol{v}^{-1}]} \mathfrak{U}_{\boldsymbol{v}}^>(L\mathfrak{sl}_n) \longrightarrow \mathfrak{U}_{\boldsymbol{v}}(L\mathfrak{sl}_n)$$

*is an isomorphism of free $\mathbb{Z}[\boldsymbol{v}, \boldsymbol{v}^{-1}]$-modules.*

Combining this Corollary with Theorem 1.3 once again, we finally obtain:

**Theorem 1.5** *For any choices (1)–(3) made prior to (1.19), the elements*

$$\left\{ \widetilde{f}_{h_-} \cdot \psi_{h_0} \cdot \widetilde{e}_{h_+} \;\middle|\; h_- \in \mathsf{H}, h_0 \in \mathsf{H}_0, h_+ \in \mathsf{H} \right\}$$

*form a basis of the free $\mathbb{Z}[\boldsymbol{v}, \boldsymbol{v}^{-1}]$-module $\mathfrak{U}_{\boldsymbol{v}}(L\mathfrak{sl}_n)$.*

### 1.1.3 Integral form $\mathsf{U}_{\boldsymbol{v}}(L\mathfrak{sl}_n)$ and its PBWD bases

For $k \in \mathbb{N}$, recall $[k]_{\boldsymbol{v}}, [k]_{\boldsymbol{v}}! \in \mathbb{Z}[\boldsymbol{v}, \boldsymbol{v}^{-1}]$ defined via:

$$[k]_{\boldsymbol{v}} := (\boldsymbol{v}^k - \boldsymbol{v}^{-k})/(\boldsymbol{v} - \boldsymbol{v}^{-1}), \qquad [k]_{\boldsymbol{v}}! := [1]_{\boldsymbol{v}} \cdots [k]_{\boldsymbol{v}} . \tag{1.24}$$

Let us define the *divided powers*

$$\mathsf{e}_{i,r}^{(k)} := e_{i,r}^k / [k]_{\boldsymbol{v}}! \quad \text{and} \quad \mathsf{f}_{i,r}^{(k)} := f_{i,r}^k / [k]_{\boldsymbol{v}}! \quad \text{for} \quad i \in I, r \in \mathbb{Z}, k \in \mathbb{N} . \tag{1.25}$$

Following [11, §7.8], define integral forms $\mathsf{U}_v^>(L\mathfrak{sl}_n)$ and $\mathsf{U}_v^<(L\mathfrak{sl}_n)$ as the $\mathbb{Z}[v, v^{-1}]$-subalgebras of $U_v^>(L\mathfrak{sl}_n)$ and $U_v^<(L\mathfrak{sl}_n)$ generated by $\{e_{i,r}^{(k)}\}_{i \in I}^{r \in \mathbb{Z}, k \in \mathbb{N}}$ and $\{f_{i,r}^{(k)}\}_{i \in I}^{r \in \mathbb{Z}, k \in \mathbb{N}}$, respectively. Recall the PBWD basis elements $\{e_\beta(r), f_\beta(r)\}_{\beta \in \Delta^+}^{r \in \mathbb{Z}}$ of (1.19), depending on the choices (1)–(3) used in their definition. We define their divided powers

$$\mathsf{e}_\beta(r)^{(k)} := e_\beta(r)^k / [k]_v! \qquad \text{and} \qquad \mathsf{f}_\beta(r)^{(k)} := f_\beta(r)^k / [k]_v! \tag{1.26}$$

for $\beta \in \Delta^+, r \in \mathbb{Z}, k \in \mathbb{N}$.

**Proposition 1.2** *For any $\beta \in \Delta^+, r \in \mathbb{Z}, k \in \mathbb{N}$, we have:*

$$\mathsf{e}_\beta(r)^{(k)} \in \mathsf{U}_v^>(L\mathfrak{sl}_n) \qquad \text{and} \qquad \mathsf{f}_\beta(r)^{(k)} \in \mathsf{U}_v^<(L\mathfrak{sl}_n) \,. \tag{1.27}$$

*Proof* Let $U_v^{\mathrm{DJ}>}(\mathfrak{sl}_n)$ be the "positive" subalgebra of the Drinfeld-Jimbo quantum group of $\mathfrak{sl}_n$. Explicitly, $U_v^{\mathrm{DJ}>}(\mathfrak{sl}_n)$ is the associative $\mathbb{C}(v)$-algebra generated by $\{E_i\}_{i \in I}$ subject to the following quantum Serre relations:

$$[E_i, E_j] = 0 \text{ if } c_{ij} = 0, \qquad [E_i, [E_i, E_j]_{v^{-1}}]_v = 0 \text{ if } c_{ij} = -1 \,. \tag{1.28}$$

Recall the Lusztig integral form $\mathsf{U}_v^{\mathrm{DJ}>}(\mathfrak{sl}_n)$ defined as the $\mathbb{Z}[v, v^{-1}]$-subalgebra of $U_v^{\mathrm{DJ}>}(\mathfrak{sl}_n)$ generated by the divided powers $E_i^{(k)} = E_i^k / [k]_v!$ for all $i \in I, k \in \mathbb{N}$. Then, $E_{\alpha_j + \cdots + \alpha_i} = [\cdots [E_j, E_{j+1}]_v, \cdots, E_i]_v$ are specific Lusztig root generators of $U_v^{\mathrm{DJ}>}(\mathfrak{sl}_n)$ (for the convex order on $\Delta^+$ opposite to (1.17)) and their divided powers $E_{\alpha_j + \cdots + \alpha_i}^{(k)} := E_{\alpha_j + \cdots + \alpha_i}^k / [k]_v!$ belong to $\mathsf{U}_v^{\mathrm{DJ}>}(\mathfrak{sl}_n)$, due to [15, Theorem 6.6].

For any $r \in \mathbb{Z}$ and $j < i \in I$, the assignment $E_i \mapsto e_{i, r\delta_{ij}}$ gives rise to an algebra homomorphism $\sigma \colon U_v^{\mathrm{DJ}>}(\mathfrak{sl}_n) \to U_v^>(L\mathfrak{sl}_n)$ such that $\sigma(\mathsf{U}_v^{\mathrm{DJ}>}(\mathfrak{sl}_n)) \subset \mathsf{U}_v^>(L\mathfrak{sl}_n)$. This implies that $\mathsf{e}_\beta(r)^{(k)} \in \mathsf{U}_v^>(L\mathfrak{sl}_n)$ for the particular choice (1.20) of the PBWD basis elements $e_\beta(r)$ in (1.19), since $\mathsf{e}_{\alpha_j + \cdots + \alpha_i}(r)^{(k)} = \sigma(E_{\alpha_j + \cdots + \alpha_i}^{(k)}) \in \mathsf{U}_v^>(L\mathfrak{sl}_n)$. The proof of $\mathsf{f}_\beta(r)^{(k)} \in \mathsf{U}_v^<(L\mathfrak{sl}_n)$ for the particular choice (1.20) is similar.

The proof of the inclusions (1.27) for a more general construction (1.19) of $e_\beta(r)$ (and similarly also of $f_\beta(r)$) is derived in Section 1.2.4 below, see Lemma 1.13. □

The monomials

$$\mathsf{e}_h := \prod_{(\beta, r) \in \Delta^+ \times \mathbb{Z}}^{\rightarrow} \mathsf{e}_\beta(r)^{(h(\beta, r))} \quad \text{and} \quad \mathsf{f}_h := \prod_{(\beta, r) \in \Delta^+ \times \mathbb{Z}}^{\leftarrow} \mathsf{f}_\beta(r)^{(h(\beta, r))} \quad \text{for all} \ \ h \in \mathsf{H}$$

$$\tag{1.29}$$

will be called the *ordered PBWD monomials* of $\mathsf{U}_v^>(L\mathfrak{sl}_n)$ and $\mathsf{U}_v^<(L\mathfrak{sl}_n)$, respectively. Here, the arrows $\rightarrow$ and $\leftarrow$ over the product signs indicate that the products are ordered with respect to (1.18) or its opposite, respectively. The following result provides the PBWD bases for the above integral forms $\mathsf{U}_v^>(L\mathfrak{sl}_n)$ and $\mathsf{U}_v^<(L\mathfrak{sl}_n)$:

**Theorem 1.6** *The ordered PBWD monomials $\{\mathsf{e}_h\}_{h \in \mathsf{H}}$ and $\{\mathsf{f}_h\}_{h \in \mathsf{H}}$ form bases of the free $\mathbb{Z}[v, v^{-1}]$-modules $\mathsf{U}_v^>(L\mathfrak{sl}_n)$ and $\mathsf{U}_v^<(L\mathfrak{sl}_n)$, respectively.*

For $i \in I, r \in \mathbb{Z}, k > 0$, recall $\{h_{i, \pm k}\}$ of (1.9) and define $\begin{bmatrix} \psi_{i,0}^+; r \\ k \end{bmatrix} \in U_v^0(L\mathfrak{sl}_n)$ via:

$$\begin{bmatrix} \psi_{i,0}^+ ; r \\ k \end{bmatrix} := \prod_{l=1}^{k} \frac{\psi_{i,0}^+ \boldsymbol{\nu}^{r-l+1} - \psi_{i,0}^- \boldsymbol{\nu}^{-r+l-1}}{\boldsymbol{\nu}^l - \boldsymbol{\nu}^{-l}} \,. \tag{1.30}$$

Following [2, §3], we define $\mathsf{U}_{\boldsymbol{\nu}}^0(L\mathfrak{sl}_n)$ as the $\mathbb{Z}[\boldsymbol{\nu}, \boldsymbol{\nu}^{-1}]$-subalgebra of $U_{\boldsymbol{\nu}}^0(L\mathfrak{sl}_n)$ generated by $\left\{ \psi_{i,0}^{\pm}, \frac{h_{i,\pm k}}{[k]_{\boldsymbol{\nu}}}, \begin{bmatrix} \psi_{i,0}^+ ; r \\ k \end{bmatrix} \right\}_{i \in I}^{r \in \mathbb{Z}, k > 0}$. Finally, we define the integral form $\mathsf{U}_{\boldsymbol{\nu}}(L\mathfrak{sl}_n)$ as the $\mathbb{Z}[\boldsymbol{\nu}, \boldsymbol{\nu}^{-1}]$-subalgebra of $U_{\boldsymbol{\nu}}(L\mathfrak{sl}_n)$ generated by $\mathsf{U}_{\boldsymbol{\nu}}^<(L\mathfrak{sl}_n), \mathsf{U}_{\boldsymbol{\nu}}^0(L\mathfrak{sl}_n), \mathsf{U}_{\boldsymbol{\nu}}^>(L\mathfrak{sl}_n)$.

*Remark 1.4* Identifying $U_{\boldsymbol{\nu}}(L\mathfrak{sl}_n)$ with the Drinfeld-Jimbo quantum loop algebra $U_{\boldsymbol{\nu}}^{\mathrm{DJ}}(L\mathfrak{sl}_n)$ of [3] (a $\mathbb{C}(\boldsymbol{\nu})$-algebra generated by $\{E_i, F_i, K_i^{\pm 1}\}_{i \in \widetilde{I}}$ where $\widetilde{I} = I \cup \{0\}$ is the vertex set of the extended Dynkin diagram, cf. Remark 3.3(d)), the integral form $\mathsf{U}_{\boldsymbol{\nu}}(L\mathfrak{sl}_n)$ gets identified with the Lusztig form of $U_{\boldsymbol{\nu}}^{\mathrm{DJ}}(L\mathfrak{sl}_n)$ defined as the $\mathbb{Z}[\boldsymbol{\nu}, \boldsymbol{\nu}^{-1}]$-subalgebra generated by $\{K_i^{\pm 1}, E_i^{(k)}, F_i^{(k)}\}_{i \in \widetilde{I}}^{k \in \mathbb{N}}$, due to [2].

The following triangular decomposition of $\mathsf{U}_{\boldsymbol{\nu}}(L\mathfrak{sl}_n)$ is due to [2, Proposition 6.1]:

**Theorem 1.7** *[2, Proposition 6.1] The multiplication map*

$$m \colon \mathsf{U}_{\boldsymbol{\nu}}^<(L\mathfrak{sl}_n) \otimes_{\mathbb{Z}[\boldsymbol{\nu}, \boldsymbol{\nu}^{-1}]} \mathsf{U}_{\boldsymbol{\nu}}^0(L\mathfrak{sl}_n) \otimes_{\mathbb{Z}[\boldsymbol{\nu}, \boldsymbol{\nu}^{-1}]} \mathsf{U}_{\boldsymbol{\nu}}^>(L\mathfrak{sl}_n) \longrightarrow \mathsf{U}_{\boldsymbol{\nu}}(L\mathfrak{sl}_n)$$

*is an isomorphism of free $\mathbb{Z}[\boldsymbol{\nu}, \boldsymbol{\nu}^{-1}]$-modules.*

Combining Theorems 1.6–1.7, we thus obtain (similar to Theorems 1.2 and 1.5) a family of PBWD bases for $\mathsf{U}_{\boldsymbol{\nu}}(L\mathfrak{sl}_n)$.

## 1.2 Shuffle realizations and applications to the PBWD bases

In this section, we establish the shuffle algebra realization of $U_{\boldsymbol{\nu}}^>(L\mathfrak{sl}_n)$ and use it to prove Theorem 1.1. Likewise, we provide the shuffle realizations of both integral forms $\mathfrak{U}_{\boldsymbol{\nu}}^>(L\mathfrak{sl}_n)$ and $\mathsf{U}_{\boldsymbol{\nu}}^>(L\mathfrak{sl}_n)$ and use those to prove Theorem 1.3 and Theorem 1.6.

### 1.2.1 Shuffle algebra $S^{(n)}$

Let $\Sigma_k$ denote the symmetric group in $k$ elements, and set

$$\Sigma_{(k_1, \dots, k_{n-1})} := \Sigma_{k_1} \times \cdots \times \Sigma_{k_{n-1}}$$

for $k_1, \dots, k_{n-1} \in \mathbb{N}$. Consider an $\mathbb{N}^I$-graded $\mathbb{C}(\boldsymbol{\nu})$-vector space

$$\mathbb{S}^{(n)} = \bigoplus_{\underline{k} = (k_1, \dots, k_{n-1}) \in \mathbb{N}^I} \mathbb{S}_{\underline{k}}^{(n)},$$

where $\mathbb{S}_{\underline{k}}^{(n)}$ consists of $\Sigma_{\underline{k}}$-symmetric rational functions in the variables $\{x_{i,r}\}_{i \in I}^{1 \le r \le k_i}$.

We fix an $I \times I$ matrix of rational functions $(\zeta_{i,j}(z))_{i,j \in I} \in \mathrm{Mat}_{I \times I}(\mathbb{C}(\nu)(z))$ via:

$$\zeta_{i,j}(z) := \frac{z - \nu^{-c_{ij}}}{z - 1} \, . \tag{1.31}$$

Let us now introduce the bilinear *shuffle product* $\star$ on $\mathbb{S}^{(n)}$: given $F \in \mathbb{S}_{\underline{k}}^{(n)}$ and $G \in \mathbb{S}_{\underline{\ell}}^{(n)}$, we define $F \star G \in \mathbb{S}_{\underline{k}+\underline{\ell}}^{(n)}$ via

$$(F \star G)(x_{1,1}, \ldots, x_{1,k_1+\ell_1}; \ldots; x_{n-1,1}, \ldots, x_{n-1,k_{n-1}+\ell_{n-1}}) := \frac{1}{\underline{k}! \cdot \underline{\ell}!} \times$$

$$\mathrm{Sym}\left( F\left(\{x_{i,r}\}_{i \in I}^{1 \le r \le k_i}\right) G\left(\{x_{i',r'}\}_{i' \in I}^{k_{i'} < r' \le k_{i'}+\ell_{i'}}\right) \cdot \prod_{\substack{i \in I \\ r \le k_i}}^{i' \in I} \prod_{r' > k_{i'}}^{r' > k_{i'}} \zeta_{i,i'}\left(\frac{x_{i,r}}{x_{i',r'}}\right) \right). \tag{1.32}$$

Here, for $\underline{m} = (m_1, \ldots, m_{n-1}) \in \mathbb{N}^I$, we set

$$\underline{m}! = \prod_{i \in I} m_i!,$$

while the *symmetrization* of $f \in \mathbb{C}\left(\{x_{i,1}, \ldots, x_{i,m_i}\}_{i \in I}\right)$ is defined via:

$$\mathrm{Sym}(f)\Big(\{x_{i,1}, \ldots, x_{i,m_i}\}_{i \in I}\Big) := \sum_{(\sigma_1, \ldots, \sigma_{n-1}) \in \Sigma_{\underline{m}}} f\Big(\{x_{i,\sigma_i(1)}, \ldots, x_{i,\sigma_i(m_i)}\}_{i \in I}\Big).$$
$$\tag{1.33}$$

This endows $\mathbb{S}^{(n)}$ with a structure of an associative algebra with the unit $\mathbf{1} \in \mathbb{S}_{(0,\ldots,0)}^{(n)}$.

Let us consider the subspace of $\mathbb{S}^{(n)}$ defined by the *pole* and *wheel* conditions:

- We say that $F \in \mathbb{S}_{\underline{k}}^{(n)}$ satisfies the *pole conditions* if

$$F = \frac{f(x_{1,1}, \ldots, x_{n-1,k_{n-1}})}{\prod_{i=1}^{n-2} \prod_{r \le k_i}^{r' \le k_{i+1}} (x_{i,r} - x_{i+1,r'})}, \quad \text{where} \quad f \in \mathbb{C}(\nu)\left[\{x_{i,r}^{\pm 1}\}_{i \in I}^{1 \le r \le k_i}\right]^{\Sigma_{\underline{k}}}. \tag{1.34}$$

- We say that $F \in \mathbb{S}_{\underline{k}}^{(n)}$ satisfies the *wheel conditions*[1] if

$$F\left(\{x_{i,r}\}\right) = 0 \quad \text{once} \quad x_{i,r_1} = \nu x_{i+\epsilon,s} = \nu^2 x_{i,r_2} \quad \text{for some} \quad \epsilon, i, r_1, r_2, s, \tag{1.35}$$

where $\epsilon \in \{\pm 1\}$, $i, i + \epsilon \in I$, $1 \le r_1 \ne r_2 \le k_i$, $1 \le s \le k_{i+\epsilon}$.

Let $S_{\underline{k}}^{(n)} \subset \mathbb{S}_{\underline{k}}^{(n)}$ denote the subspace of all $F$ satisfying these two conditions and set

$$S^{(n)} := \bigoplus_{\underline{k} \in \mathbb{N}^I} S_{\underline{k}}^{(n)}.$$

It is straightforward to check that the subspace $S^{(n)} \subset \mathbb{S}^{(n)}$ is actually $\star$-closed:

---

[1] Following [6]–[8], the role of the wheel conditions is to replace intricate quantum Serre relations.

**Lemma 1.2** *For any $F \in S_{\underline{k}}^{(n)}$ and $G \in S_{\underline{\ell}}^{(n)}$, we have $F \star G \in S_{\underline{k}+\underline{\ell}}^{(n)}$.*

*Proof* If $F$ and $G$ satisfy (1.34), then $F \star G$ has at most first order poles at $x_{i,r} - x_{i+1,r'}$ or $x_{i,r} - x_{i,r'}$. But the latter cannot happen since $F \star G$ is symmetric in each family $\{x_{i,s}\}_{s=1}^{k_i + \ell_i}$. Thus, the shuffle product $F \star G$ also satisfies the pole conditions (1.34).

To verify the wheel conditions (1.35) for $F \star G$, we show that any summand from the symmetrization in the right-hand side of (1.32) vanishes under the specialization $x_{i,r_1} = \nu x_{i+\epsilon,s} = \nu^2 x_{i,r_2}$. To this end, either we shall encounter one of the vanishing $\zeta$-factors $\zeta_{i,i+\epsilon}(x_{i,r_1}/x_{i+\epsilon,s}) = \zeta_{i+\epsilon,i}(x_{i+\epsilon,s}/x_{i,r_2}) = \zeta_{i,i}(x_{i,r_2}/x_{i,r_1}) = 0$, or all three variables $x_{i,r_1}, x_{i+\epsilon,s}, x_{i,r_2}$ will get placed into the same function $F$ or $G$. In the latter case, the result will be zero again, due to the wheel conditions (1.35) for $F$ and $G$. □

The **shuffle algebra** $(S^{(n)}, \star)$ is linked to $U_\nu^>(L\mathfrak{sl}_n)$ via the following result:[2]

**Proposition 1.3** *There exists an injective $\mathbb{C}(\nu)$-algebra homomorphism*

$$\Psi \colon U_\nu^>(L\mathfrak{sl}_n) \longrightarrow S^{(n)} \tag{1.36}$$

*such that $e_{i,r} \mapsto x_{i,1}^r$ for any $i \in I, r \in \mathbb{Z}$.*

*Proof* Let us show that the assignment $\Psi \colon e_{i,r} \mapsto x_{i,1}^r$ ($i \in I, r \in \mathbb{Z}$) gives rise to an algebra homomorphism $\Psi \colon U_\nu^>(L\mathfrak{sl}_n) \to S^{(n)}$. To do so, it suffices to verify that $\Psi(e_i(z)) = \delta(\frac{x_{i,1}}{z})$ satisfy the defining relations (1.2, 1.7) of $U_\nu^>(L\mathfrak{sl}_n)$, see Proposition 1.1(b). To this end, let us first recall the key property of the $\delta$-function:

$$\delta(z/w) \cdot \gamma(z) = \delta(z/w) \cdot \gamma(w) \qquad \text{for any rational function} \quad \gamma(z). \tag{1.37}$$

- *Compatibility with (1.2) for $i \neq j$.*
  $(z - \nu^{c_{ij}} w)\delta(\frac{x_{i,1}}{z}) \star \delta(\frac{x_{j,1}}{w}) = \delta(\frac{x_{i,1}}{z})\delta(\frac{x_{j,1}}{w}) \cdot (x_{i,1} - \nu^{c_{ij}} x_{j,1}) \cdot \frac{x_{i,1} - \nu^{-c_{ij}} x_{j,1}}{x_{i,1} - x_{j,1}} = \delta(\frac{x_{j,1}}{w})\delta(\frac{x_{i,1}}{z}) \cdot (\nu^{c_{ij}} x_{i,1} - x_{j,1}) \cdot \frac{x_{j,1} - \nu^{-c_{ij}} x_{i,1}}{x_{j,1} - x_{i,1}} = (\nu^{c_{ij}} z - w)\delta(\frac{x_{j,1}}{w}) \star \delta(\frac{x_{i,1}}{z})$.

- *Compatibility with (1.2) for $i = j$.*
  $(z - \nu^{c_{ii}} w)\delta(\frac{x_{i,1}}{z})\star\delta(\frac{x_{i,1}}{w}) = \mathrm{Sym}\left(\delta(\frac{x_{i,1}}{z})\delta(\frac{x_{i,2}}{w}) \cdot (x_{i,1} - \nu^{c_{ii}} x_{i,2}) \frac{x_{i,1} - \nu^{-c_{ii}} x_{i,2}}{x_{i,1} - x_{i,2}}\right) = \mathrm{Sym}\left(\delta(\frac{x_{i,2}}{w})\delta(\frac{x_{i,1}}{z}) \cdot (\nu^{c_{ii}} x_{i,1} - x_{i,2}) \frac{x_{i,2} - \nu^{-c_{ii}} x_{i,1}}{x_{i,2} - x_{i,1}}\right) = (\nu^{c_{ii}} z - w)\delta(\frac{x_{i,1}}{w}) \star \delta(\frac{x_{i,1}}{z})$.

- *Compatibility with (1.7) for $c_{ij} = 0$.*
  $\delta(\frac{x_{i,1}}{z}) \star \delta(\frac{x_{j,1}}{w}) = \delta(\frac{x_{i,1}}{z})\delta(\frac{x_{j,1}}{w}) = \delta(\frac{x_{j,1}}{w}) \star \delta(\frac{x_{i,1}}{z})$.

- *Compatibility with (1.7) for $c_{ij} = -1$.*
  $[\delta(\frac{x_{i,1}}{z_1}), [\delta(\frac{x_{i,1}}{z_2}), \delta(\frac{x_{j,1}}{w})]_{\nu^{-1}}]_\nu = A(z_1, z_2, w) \cdot \left(\delta(\frac{x_{i,1}}{z_1})\delta(\frac{x_{i,2}}{z_2}) + \delta(\frac{x_{i,1}}{z_2})\delta(\frac{x_{i,2}}{z_1})\right) \cdot \delta(\frac{x_{j,1}}{w})$ with $A(z_1, z_2, w) = \zeta_{i,i}(\frac{z_1}{z_2})\zeta_{i,j}(\frac{z_1}{w})\zeta_{i,j}(\frac{z_2}{w}) - \nu^{-1}\zeta_{i,i}(\frac{z_2}{z_1})\zeta_{i,j}(\frac{z_1}{w})\zeta_{j,i}(\frac{w}{z_2}) - \nu\zeta_{i,i}(\frac{z_2}{z_1})\zeta_{i,j}(\frac{z_2}{w})\zeta_{j,i}(\frac{w}{z_1}) + \zeta_{i,i}(\frac{z_1}{z_2})\zeta_{j,i}(\frac{w}{z_1})\zeta_{j,i}(\frac{w}{z_2})$. Therefore, the compatibility with the cubic quantum Serre relations (1.7) follows from a straightforward equality $A(z_1, z_2, w) = -A(z_2, z_1, w)$.

---

[2] When working over $\mathbb{C}[[\hbar]]$ rather than over $\mathbb{C}(\nu)$, this construction goes back to [5, Corollary 1.4].

The quickest way to prove the injectivity of $\Psi$ is based on the analogue of (3.102):

$$\varphi(X, f_{i_1}(z_1) \ldots f_{i_k}(z_k)) = (\nu - \nu^{-1})^{-k} \frac{\Psi(X)(z_1, \ldots, z_k)}{\prod_{1 \leq a < b \leq k} \zeta_{i_a, i_b}(z_a/z_b)}, \qquad (1.38)$$

where $\varphi$ denotes the Hopf pairing $U_\nu^{\geq}(\widehat{\mathfrak{sl}}_n) \times U_\nu^{\leq}(\widehat{\mathfrak{sl}}_n) \to \mathbb{C}(\nu)$ as in Theorem 3.4 used to realize the quantum affine algebra $U_\nu(\widehat{\mathfrak{sl}}_n)$ as a Drinfeld double. Here, we note that the "positive" subalgebras $U_\nu^{>}(\widehat{\mathfrak{sl}}_n) \simeq U_\nu^{>}(L\mathfrak{sl}_n)$ are naturally isomorphic. Thus, the injectivity of $\Psi$ is an immediate consequence of the nondegeneracy of $\varphi$. $\square$

The next result follows from its much harder counterpart [16, Theorem 1.1], see Theorem 3.6 below, but we will provide its alternative simpler proof[3] in Section 1.2.2:

**Theorem 1.8** $\Psi: U_\nu^{>}(L\mathfrak{sl}_n) \xrightarrow{\sim} S^{(n)}$ of (1.36) is a $\mathbb{C}(\nu)$-algebra isomorphism.

Before proceeding to the simultaneous proofs of Theorem 1.1 and Theorem 1.8, let us introduce the key tool in our study of the shuffle algebra: the *specialization maps*. For any positive root $\beta = \alpha_j + \alpha_{j+1} + \cdots + \alpha_i$, define $j(\beta) := j, i(\beta) := i$, and let $[\beta]$ denote the integer interval $[j(\beta); i(\beta)]$. Consider a collection of the intervals $\{[\beta]\}_{\beta \in \Delta^+}$ each taken with a multiplicity $d_\beta \in \mathbb{N}$ and ordered with respect to the total order (1.17) on $\Delta^+$ (the order inside each group is irrelevant). Let $\underline{d} \in \mathbb{N}^{\Delta^+}$ denote the collection $\{d_\beta\}_{\beta \in \Delta^+}$. Define $\underline{\ell} = (\ell_1, \ldots, \ell_{n-1}) \in \mathbb{N}^I$, the $\mathbb{N}^I$-degree of $\underline{d}$, via:

$$\sum_{i \in I} \ell_i \alpha_i = \sum_{\beta \in \Delta^+} d_\beta \beta. \qquad (1.39)$$

Let us now define the **specialization map** $\phi_{\underline{d}}$ on the shuffle algebra $S^{(n)}$, which vanishes on degree $\neq \underline{\ell}$ components, and is thus determined by

$$\phi_{\underline{d}}: S_{\underline{\ell}}^{(n)} \longrightarrow \mathbb{C}(\nu) \left[ \{y_{\beta,s}^{\pm 1}\}_{\beta \in \Delta^+}^{1 \leq s \leq d_\beta} \right]^{\Sigma_{\underline{d}}}. \qquad (1.40)$$

To do so, we split the variables $\{x_{i,r}\}_{i \in I}^{1 \leq r \leq \ell_i}$ into $\sum_{\beta \in \Delta^+} d_\beta$ groups corresponding to the above intervals, and specialize those in the $s$-th copy of $[\beta]$ in the natural order to

$$\nu^{-j(\beta)} \cdot y_{\beta,s}, \ldots, \nu^{-i(\beta)} \cdot y_{\beta,s}$$

so that $x_{i,r}$ gets specialized to $\nu^{-i} y_{\beta,s}$. For $F = \frac{f(x_{1,1}, \ldots, x_{n-1, \ell_{n-1}})}{\prod_{i=1}^{n-2} \prod_{1 \leq r' \leq \ell_{i+1}}^{1 \leq r \leq \ell_i} (x_{i,r} - x_{i+1,r'})} \in S_{\underline{\ell}}^{(n)}$, we define $\phi_{\underline{d}}(F)$ as the corresponding specialization of $f$. Let us note right away that:

- $\phi_{\underline{d}}(F)$ is independent of our splitting of the variables $\{x_{i,r}\}_{i \in I}^{1 \leq r \leq \ell_i}$ into groups,
- $\phi_{\underline{d}}(F)$ is symmetric in $\{y_{\beta,s}\}_{s=1}^{d_\beta}$ for any $\beta \in \Delta^+$.

The key properties of the specialization maps $\phi_{\underline{d}}$ and their relevance to $\Psi(e_h)$, the images of the ordered PBWD monomials, are discussed in Lemmas 1.5–1.7 below.

---

[3] In particular, our proof will be directly generalized to establish the isomorphisms of Theorems 1.18 and 1.21 below, for which no analogue of Theorem 3.6 ([16, Theorem 1.1]) is presently known.

**Example 1.9** Let $F = \frac{x_{1,1}^a x_{2,1}^b x_{3,1}^c}{(x_{1,1}-x_{2,1})(x_{2,1}-x_{3,1})} \in S_{\underline{\ell}}^{(n)}$, where $\underline{\ell} = (1, 1, 1, 0, \ldots, 0) \in \mathbb{N}^I$ and $a, b, c \in \mathbb{Z}$. Let us compute the images of $F$ under all possible specialization maps:

(i) If $\underline{d}$ encodes a single positive root $\beta = \alpha_1 + \alpha_2 + \alpha_3$, then $\phi_{\underline{d}}(F)$ is a Laurent polynomial in a single variable $y_{\beta,1}$ and equals $(\nu^{-1} y_{\beta,1})^a (\nu^{-2} y_{\beta,1})^b (\nu^{-3} y_{\beta,1})^c$.

(ii) If $\underline{d}$ encodes two positive roots $\beta_1 = \alpha_1, \beta_2 = \alpha_2 + \alpha_3$, then $\phi_{\underline{d}}(F)$ is a Laurent polynomial in two variables $y_{\beta_1,1}, y_{\beta_2,1}$ and equals $(\nu^{-1} y_{\beta_1,1})^a (\nu^{-2} y_{\beta_2,1})^b (\nu^{-3} y_{\beta_2,1})^c$.

(iii) If $\underline{d}$ encodes two positive roots $\beta_1 = \alpha_1 + \alpha_2, \beta_2 = \alpha_3$, then $\phi_{\underline{d}}(F)$ is a Laurent polynomial in two variables $y_{\beta_1,1}, y_{\beta_2,1}$ and equals $(\nu^{-1} y_{\beta_1,1})^a (\nu^{-2} y_{\beta_1,1})^b (\nu^{-3} y_{\beta_2,1})^c$.

(iv) If $\underline{d}$ encodes roots $\beta_1 = \alpha_1, \beta_2 = \alpha_2, \beta_3 = \alpha_3$, then $\phi_{\underline{d}}(F)$ is a Laurent polynomial in three variables $y_{\beta_1,1}, y_{\beta_2,1}, y_{\beta_3,1}$ and equals $(\nu^{-1} y_{\beta_1,1})^a (\nu^{-2} y_{\beta_2,1})^b (\nu^{-3} y_{\beta_3,1})^c$.

### 1.2.2 Proofs of Theorem 1.1 and Theorem 1.8

Our proof of Theorem 1.1 will proceed in two steps: first, we shall establish the linear independence of the ordered PBWD monomials $\{e_h\}_{h \in \mathsf{H}}$ in Section 1.2.2.2, and then we will verify that they linearly span the entire algebra $U_\nu^>(L\mathfrak{sl}_n)$ in Section 1.2.2.3. The part of Theorem 1.1 for the "negative" subalgebra $U_\nu^<(L\mathfrak{sl}_n)$ can either be proved analogously, or can rather be deduced by utilizing a $\mathbb{C}(\nu)$-algebra anti-isomorphism

$$\varsigma \colon U_\nu^>(L\mathfrak{sl}_n) \longrightarrow U_\nu^<(L\mathfrak{sl}_n) \quad \text{defined by} \quad e_{i,r} \mapsto f_{i,r} \quad \forall\, i \in I, r \in \mathbb{Z}, \quad (1.41)$$

which maps the ordered PBWD monomials to the ordered PBWD monomials multiplied by constants in $\pm \nu^{\mathbb{Z}}$. Notably, we shall also prove Theorem 1.8 along the way.

First, we shall prove both theorems for $n = 2$ (the "rank 1" case) in Section 1.2.2.1.

#### 1.2.2.1 $n = 2$ case

For $n = 2$, we have $I = \{1\}$, hence, we shall denote the variables $x_{1,r}$ simply by $x_r$. We start from the following simple computation in the shuffle algebra $S^{(2)}$:

**Lemma 1.3** *For any $k \geq 1$ and $r \in \mathbb{Z}$, the $k$-th power of $x^r \in S_1^{(2)}$ equals*

$$\underbrace{x^r \star \cdots \star x^r}_{k \text{ times}} = \nu^{-\frac{k(k-1)}{2}} [k]_\nu! \cdot (x_1 \cdots x_k)^r. \quad (1.42)$$

*Proof* We prove this Lemma by induction on $k$, the base case $k = 1$ being trivial. Applying the induction assumption to the $(k-1)$-st power of $x^r$, the proof of (1.42) boils down to the verification of the following equality:

$$\sum_{p=1}^{k} \prod_{\substack{1 \le s \le k}}^{s \ne p} \frac{x_s - \nu^{-2} x_p}{x_s - x_p} = \nu^{1-k}[k]_{\nu}. \tag{1.43}$$

The left-hand side of (1.43) is a rational function in $\{x_p\}_{p=1}^{k}$ of degree 0 and has no poles (as it is symmetric and could only have simple poles along the diagonals), hence, must be a constant. To evaluate this constant, let $x_k \to \infty$: the last summand (corresponding to $p = k$) tends to $\nu^{-2(k-1)}$, while the sum of the first $k-1$ summands (corresponding to $1 \le p \le k-1$) tends to $1 + \nu^{-2} + \cdots + \nu^{-2(k-2)}$ by the induction assumption, thus, resulting in $1 + \nu^{-2} + \nu^{-4} + \cdots + \nu^{-2(k-1)} = \nu^{1-k}[k]_{\nu}$, as claimed.□

The following result simultaneously implies both Theorems 1.1 and 1.8 for $n = 2$:

**Proposition 1.4** *Fix any total order $\le$ on $\mathbb{Z}$. Then:*

*(a) the ordered monomials $\{x^{r_1} \star x^{r_2} \star \cdots \star x^{r_k}\}_{k \in \mathbb{N}}^{r_1 \le \cdots \le r_k}$ form a $\mathbb{C}(\nu)$-basis of $S^{(2)}$.*

*(b) the ordered monomials $\{e_{r_1} e_{r_2} \cdots e_{r_k}\}_{k \in \mathbb{N}}^{r_1 \le \cdots \le r_k}$ form a $\mathbb{C}(\nu)$-basis of $U_{\nu}^{>}(L\mathfrak{sl}_2)$.*

*Proof* For $r_1 = \cdots = r_{k_1} < r_{k_1+1} = \cdots = r_{k_1+k_2} < \cdots < r_{k_1+\cdots+k_{\ell-1}+1} = \cdots = r_{k_1+\cdots+k_\ell}$, set $k := k_1 + \cdots + k_\ell$ and choose $\sigma \in \Sigma_k$ so that $r_{\sigma(1)} \le \cdots \le r_{\sigma(k)}$. Then, the shuffle product $x^{r_1} \star \cdots \star x^{r_k}$ is a symmetric Laurent polynomial of the form

$$\nu_{\underline{r}} m_{(r_{\sigma(1)}, \ldots, r_{\sigma(k)})}(x_1, \ldots, x_k) + \sum \nu_{\underline{r}'} m_{\underline{r}'}(x_1, \ldots, x_k), \tag{1.44}$$

where

- $m_{(s_1, \ldots, s_k)}(x_1, \ldots, x_k)$ (with $s_1 \le \cdots \le s_k$) are the monomial symmetric polynomials (that is, the sums of all monomials $x_1^{r_1} \cdots x_k^{t_k}$ as $(t_1, \ldots, t_k)$ ranges over all distinct permutations of $(s_1, \ldots, s_k)$),

- the sum in (1.44) is over $\underline{r}' = (r'_1 \le \cdots \le r'_k) \in \mathbb{Z}^k$ distinct from $(r_{\sigma(1)}, \ldots, r_{\sigma(k)})$ and satisfying

  $$r_{\sigma(1)} \le r'_1 \le r'_k \le r_{\sigma(k)} \quad \text{as well as} \quad r'_1 + \cdots + r'_k = r_1 + \cdots + r_k,$$

- the coefficients $\nu_{\underline{r}'}$ in (1.44) are Laurent polynomials in $\nu$, that is, $\nu_{\underline{r}'} \in \mathbb{Z}[\nu, \nu^{-1}]$,

- due to Lemma 1.3, the coefficient $\nu_{\underline{r}}$ is explicitly given by

  $$\nu_{\underline{r}} = \prod_{p=1}^{\ell} \left( \nu^{-\frac{k_p(k_p-1)}{2}} [k_p]_{\nu}! \right). \tag{1.45}$$

Since $S_k^{(2)} \simeq \mathbb{C}(\nu)[\{x_p^{\pm 1}\}_{p=1}^{k}]^{\Sigma_k}$ as vector spaces and $\{m_{(s_1, \ldots, s_k)}(x_1, \ldots, x_k)\}_{s_1 \le \cdots \le s_k}$ form a $\mathbb{C}(\nu)$-basis of $\mathbb{C}(\nu)[\{x_p^{\pm 1}\}_{p=1}^{k}]^{\Sigma_k}$, we conclude that the ordered shuffle products $\{x^{r_1} \star x^{r_2} \star \cdots \star x^{r_k}\}_{k \in \mathbb{N}}^{r_1 \le \cdots \le r_k}$ form a $\mathbb{C}(\nu)$-basis of $S^{(2)}$. This proves part (a).

As the homomorphism $\Psi \colon U_{\nu}^{>}(L\mathfrak{sl}_2) \to S^{(2)}, e_r \mapsto x^r$, of (1.36) is injective by Proposition 1.3 and $\Psi(e_{r_1} e_{r_2} \cdots e_{r_k}) = x^{r_1} \star x^{r_2} \star \cdots \star x^{r_k}$, part (a) implies that $\{e_{r_1} e_{r_2} \cdots e_{r_k}\}_{k \in \mathbb{N}}^{r_1 \le \cdots \le r_k}$ form a $\mathbb{C}(\nu)$-basis of $U_{\nu}^{>}(L\mathfrak{sl}_2)$. This proves part (b). □

#### 1.2.2.2 Linear independence of $e_h$ and key properties of $\phi_{\underline{d}}$

For an ordered PBWD monomial $e_h$ ($h \in \mathsf{H}$) of (1.21), we define its *degree*

$$\deg(e_h) = \deg(h) \in \mathbb{N}^{\Delta^+} \tag{1.46}$$

as a collection of $d_\beta := \sum_{r \in \mathbb{Z}} h(\beta, r) \in \mathbb{N}$ ($\beta \in \Delta^+$) ordered with respect to the total order (1.17) on $\Delta^+$. Let us further consider the *anti-lexicographical order* on $\mathbb{N}^{\Delta^+}$:

$$\{d_\beta\}_{\beta \in \Delta^+} < \{d'_\beta\}_{\beta \in \Delta^+} \quad \text{iff} \quad \exists \gamma \in \Delta^+ \text{ s.t. } d_\gamma > d'_\gamma \text{ and } d_\beta = d'_\beta \ \forall \beta < \gamma . \tag{1.47}$$

In what follows, we shall need an explicit formula for the image $\Psi(e_\beta(r)) \in S^{(n)}$:

**Lemma 1.4** *For any $1 \le j < i < n$ and $r \in \mathbb{Z}$, we have:*

$$\Psi(e_{\alpha_j + \alpha_{j+1} + \cdots + \alpha_i}(r)) = (1 - \nu^2)^{i-j} \frac{p(x_{j,1}, \ldots, x_{i,1})}{(x_{j,1} - x_{j+1,1}) \cdots (x_{i-1,1} - x_{i,1})} , \tag{1.48}$$

*where $p(x_{j,1}, \ldots, x_{i,1})$ is a degree $r + i - j$ monomial, up to a sign and an integer power of $\nu$.*

*Proof* Straightforward computation. □

**Example 1.10** For the particular choice (1.20) of the elements $e_{\alpha_j + \alpha_{j+1} + \cdots + \alpha_i}(r)$, we have $p(x_{j,1}, \ldots, x_{i,1}) = x_{j,1}^{r+1} x_{j+1,1} \cdots x_{i-1,1}$ in the formula (1.48).

Our proof of the linear independence of $\{e_h\}_{h \in \mathsf{H}}$ is crucially based on the following key properties of the specialization maps $\phi_{\underline{d}}$ from (1.40):

**Lemma 1.5** *If $\deg(h) < \underline{d}$ while $\mathbb{N}^I$-degrees of $\underline{d}$ and $\deg(h)$ coincide, then*

$$\phi_{\underline{d}}(\Psi(e_h)) = 0 .$$

*Proof* The condition $\deg(h) < \underline{d}$ guarantees that $\phi_{\underline{d}}$-specialization of any summand from the symmetrization appearing in $\Psi(e_h)$ contains among all the $\zeta$-factors at least one factor of the form $\zeta_{i,i+1}(\nu) = 0$, hence, must vanish. The result follows. □

**Lemma 1.6** *For any $\underline{d} \in \mathbb{N}^{\Delta^+}$, the specializations $\left\{\phi_{\underline{d}}(\Psi(e_h))\right\}_{h \in \mathsf{H}}^{\deg(h) = \underline{d}}$ are linearly independent.*

*Proof* Consider the image $\Psi(e_h) \in S_{\underline{k}}^{(n)}$ (here, $\underline{k}$ is the $\mathbb{N}^I$-degree of $\underline{d} = \deg(h)$). It is a sum of $\underline{k}!$ terms, but most of them specialize to zero under $\phi_{\underline{d}}$, as in the proof of Lemma 1.5. The summands which do not specialize to zero are parametrized by $\Sigma_{\underline{d}} := \prod_{\beta \in \Delta^+} \Sigma_{d_\beta}$. More precisely, given $(\sigma_\beta)_{\beta \in \Delta^+} \in \Sigma_{\underline{d}}$, the associated summand corresponds to the case when for all $\beta \in \Delta^+$ and $1 \le s \le d_\beta$, the $(\sum_{\beta' < \beta} d_{\beta'} + s)$-th factor in the corresponding term of $\Psi(e_h)$ is evaluated at $\nu^{-j(\beta)} y_{\beta, \sigma_\beta(s)}, \ldots, \nu^{-i(\beta)} y_{\beta, \sigma_\beta(s)}$. The image of this summand under $\phi_{\underline{d}}$ equals

$$\prod_{\beta,\beta'\in\Delta^+}^{\beta<\beta'} G_{\beta,\beta'} \cdot \prod_{\beta\in\Delta^+} G_\beta \cdot \prod_{\beta\in\Delta^+} G_\beta^{(\sigma_\beta)}$$

(up to a common sign and a power of $\boldsymbol{v}$), with the factors given explicitly by:

$$
\begin{aligned}
G_{\beta,\beta'} &= \prod_{1\le s\le d_\beta}^{1\le s'\le d_{\beta'}} \left( \prod_{j\in[\beta],j'\in[\beta']}^{j=j'} (y_{\beta,s}-\boldsymbol{v}^{-2}y_{\beta',s'}) \cdot \prod_{j\in[\beta],j'\in[\beta']}^{j=j'+1} (y_{\beta,s}-\boldsymbol{v}^2 y_{\beta',s'}) \right) \times \\
&\qquad \prod_{1\le s\le d_\beta}^{1\le s'\le d_{\beta'}} (y_{\beta,s}-y_{\beta',s'})^{\delta_{j(\beta')>j(\beta)}\delta_{i(\beta)+1\in[\beta']}}, \\
G_\beta &= (1-\boldsymbol{v}^2)^{d_\beta(i(\beta)-j(\beta))} \cdot \prod_{1\le s\ne s'\le d_\beta} (y_{\beta,s}-\boldsymbol{v}^2 y_{\beta,s'})^{i(\beta)-j(\beta)} \cdot \prod_{1\le s\le d_\beta} y_{\beta,s}^{i(\beta)-j(\beta)}, \\
G_\beta^{(\sigma_\beta)} &= \prod_{s=1}^{d_\beta} y_{\beta,\sigma_\beta(s)}^{r_\beta(h,s)} \cdot \prod_{s<s'} \frac{y_{\beta,\sigma_\beta(s)}-\boldsymbol{v}^{-2}y_{\beta,\sigma_\beta(s')}}{y_{\beta,\sigma_\beta(s)}-y_{\beta,\sigma_\beta(s')}}.
\end{aligned}
$$

$$(1.49)$$

Here, for any $\beta\in\Delta^+$, the exponents $\{r_\beta(h,1),\dots,r_\beta(h,d_\beta)\}$ are obtained by listing every $r\in\mathbb{Z}$ with multiplicity $h(\beta,r)>0$ with respect to the total order $\le_\beta$ on $\mathbb{Z}$, see Section 1.1.1. We also use the standard *delta function* notation:

$$
\delta_{\mathrm{condition}} = \begin{cases} 1 & \text{if condition holds} \\ 0 & \text{if condition fails} \end{cases}.
$$

Note that the factors $\{G_{\beta,\beta'},G_\beta\}_{\beta<\beta'}$ in (1.49) are independent of $(\sigma_\beta)_{\beta\in\Delta^+}\in\Sigma_{\underline{d}}$. Therefore, the specialization $\phi_{\underline{d}}(\Psi(e_h))$ has the following form:

$$
\phi_{\underline{d}}(\Psi(e_h)) = c \cdot \prod_{\beta,\beta'\in\Delta^+}^{\beta<\beta'} G_{\beta,\beta'} \cdot \prod_{\beta\in\Delta^+} G_\beta \cdot \prod_{\beta\in\Delta^+} \left( \sum_{\sigma_\beta\in\Sigma_{d_\beta}} G_\beta^{(\sigma_\beta)} \right) \qquad (1.50)
$$

with $c\in\pm\boldsymbol{v}^{\mathbb{Z}}$. The key observation is that for any $\beta\in\Delta^+$, the sum $\sum_{\sigma_\beta\in\Sigma_{d_\beta}} G_\beta^{(\sigma_\beta)}$ coincides with the value of the shuffle product $x^{r_\beta(h,1)}\star\cdots\star x^{r_\beta(h,d_\beta)}\in S_{d_\beta}^{(2)}$ (in the shuffle algebra $S^{(2)}!$, we call this feature a "rank 1 reduction") evaluated at $\{y_{\beta,s}\}_{s=1}^{d_\beta}$. The latter elements are linearly independent, according to Proposition 1.4(a).

Combining this observation with (1.50) completes the proof of Lemma 1.6.  □

Now we can establish the linear independence of $\{e_h\}_{h\in\mathsf{H}}\subset U_{\boldsymbol{v}}^>(L\mathfrak{sl}_n)$:

**Proposition 1.5** *The ordered PBWD monomials $\{e_h\}_{h\in\mathsf{H}}$ are linearly independent.*

*Proof* Assuming the contradiction, pick a nontrivial linear combination $\sum_{h\in\mathsf{H}} c_h e_h$ that vanishes in $U_{\boldsymbol{v}}^>(L\mathfrak{sl}_n)$ (here, all but finitely many of $c_h$ are zero, but at least one of them is nonzero). Define $\underline{d} := \max\{\deg(h) | c_h \neq 0\}$. Applying the specialization map $\phi_{\underline{d}}$ to $\sum_{h \in \mathsf{H}} c_h \Psi(e_h) = 0$, we get $\sum_{h \in \mathsf{H}}^{\deg(h)=\underline{d}} c_h \phi_{\underline{d}}(\Psi(e_h)) = 0$ by Lemma 1.5. Furthermore, we get $c_h = 0$ for all $h \in \mathsf{H}$ of degree $\deg(h) = \underline{d}$, due to Lemma 1.6. This contradicts our choice of $\underline{d}$.                                                    □

### 1.2.2.3  Spanning property of $e_h$ and dominance property of $\phi_{\underline{d}}$

Let $M \subset S^{(n)}$ be the $\mathbb{C}(\nu)$-span of $\{\Psi(e_h)\}_{h \in \mathsf{H}}$, the images of the ordered PBWD monomials. For $\underline{k} \in \mathbb{N}^I$, let $T_{\underline{k}}$ denote a finite set consisting of all degree vectors $\underline{d} = \{d_\beta\}_{\beta \in \Delta^+} \in \mathbb{N}^{\Delta^+}$ such that $\sum_{\beta \in \Delta^+} d_\beta \beta = \sum_{i \in I} k_i \alpha_i$. We order $T_{\underline{k}}$ with respect to the anti-lexicographical order (1.47) on $\mathbb{N}^{\Delta^+}$. In particular, the minimal element $\underline{d}_{\min} = \{d_\beta\}_{\beta \in \Delta^+} \in T_{\underline{k}}$ is characterized by $d_{\alpha_i} = d_i$ and $d_\beta = 0$ when $\beta \in \Delta^+$ is not simple. The following "dominance property" of $\phi_{\underline{d}}$ is of crucial importance:

**Lemma 1.7** *Consider a shuffle element $F \in S^{(n)}_{\underline{k}}$ and a degree vector $\underline{d} \in T_{\underline{k}}$. If $\phi_{\underline{d}'}(F) = 0$ for all $\underline{d}' \in T_{\underline{k}}$ such that $\underline{d}' > \underline{d}$, then there exists $F_{\underline{d}} \in M \cap S^{(n)}_{\underline{k}}$ such that*

$$\phi_{\underline{d}}(F) = \phi_{\underline{d}}(F_{\underline{d}}) \qquad \text{and} \qquad \phi_{\underline{d}'}(F_{\underline{d}}) = 0 \quad \text{for all} \quad \underline{d}' > \underline{d} \, .$$

*Proof*  Consider the following total order on the set $\{(\beta, s) | \beta \in \Delta^+, 1 \leq s \leq d_\beta\}$:

$$(\beta, s) \leq (\beta', s') \quad \text{iff} \quad \beta < \beta' \ \text{ or } \ \beta = \beta' \text{ and } s \leq s' \, . \tag{1.51}$$

First, we note that the wheel conditions (1.35) for $F$ guarantee that $\phi_{\underline{d}}(F)$ (which is a Laurent polynomial in the variables $\{y_{\beta,s}\}_{\beta \in \Delta^+}^{s \leq d_\beta}$) vanishes, up to appropriate orders, under the following specializations:

(i) $y_{\beta,s} = \nu^{-2} y_{\beta',s'}$ for $(\beta, s) < (\beta', s')$,
(ii) $y_{\beta,s} = \nu^2 y_{\beta',s'}$ for $(\beta, s) < (\beta', s')$.

To evaluate the aforementioned orders of vanishing, let us view the specialization appearing in the definition of $\phi_{\underline{d}}$ as a step-by-step specialization in each interval $[\beta]$, ordered first in the non-increasing length order, while the intervals of the same length are ordered in the non-decreasing order of $j(\beta)$. As we specialize the variables in the $s$-th interval ($1 \leq s \leq \sum_{\beta \in \Delta^+} d_\beta$), we count only those wheel conditions that arise from the non-specialized yet variables. A straightforward case-by-case verification (treating each of the following cases separately): $j = j' = i = i'$, $j = j' = i < i'$, $j = j' < i = i'$, $j = j' < i < i'$, $j < j' \leq i < i'$, $j < j' < i' = i$, $j < j' < i < i'$, $j = i < j' < i'$, $j < i = j' < i'$, $j < j' = i < i'$, $j < i < j' < i'$, $j < j' < i < i'$, $j < i < i' \leq i'$, where we set $j := j(\beta)$, $j' := j(\beta')$, $i := i(\beta)$, $i' := i(\beta')$) shows that the corresponding orders of vanishing under the specializations (i) and (ii) are respectively equal to:

$$\#\{(j, j') \in [\beta] \times [\beta'] \,|\, j = j'\} - \delta_{\beta\beta'} \quad \text{and} \quad \#\{(j, j') \in [\beta] \times [\beta'] \,|\, j = j' + 1\} \, . \tag{1.52}$$

Second, we claim that $\phi_{\underline{d}}(F)$ vanishes under the following specializations:

(iii) $y_{\beta,s} = y_{\beta',s'}$ for $(\beta,s) < (\beta',s')$ such that $j(\beta) < j(\beta')$ and $i(\beta) + 1 \in [\beta']$.

Indeed, if $j(\beta) < j(\beta')$ and $i(\beta) + 1 \in [\beta']$, then there are positive roots $\gamma, \gamma' \in \Delta^+$ such that $j(\gamma) = j(\beta)$, $i(\gamma) = i(\beta')$, $j(\gamma') = j(\beta')$, $i(\gamma') = i(\beta)$. Consider the degree vector $\underline{d}' \in T_{\underline{k}}$ given by $d'_\alpha = d_\alpha + \delta_{\alpha\gamma} + \delta_{\alpha\gamma'} - \delta_{\alpha\beta} - \delta_{\alpha\beta'}$. As $\beta < \gamma < \gamma' < \beta'$, we have $\underline{d}' > \underline{d}$ and therefore $\phi_{\underline{d}'}(F) = 0$ by the assumption. The result follows.

Combining the above vanishing conditions (i)–(iii) and (1.52), we see that $\phi_{\underline{d}}(F)$ is divisible precisely by the product $\prod_{\beta<\beta'} G_{\beta,\beta'} \cdot \prod_\beta G_\beta$ of (1.49). Thus, we have:

$$\phi_{\underline{d}}(F) = \prod_{\beta,\beta'\in\Delta^+}^{\beta<\beta'} G_{\beta,\beta'} \cdot \prod_{\beta\in\Delta^+} G_\beta \cdot G \tag{1.53}$$

for some symmetric Laurent polynomial

$$G \in \mathbb{C}(\nu) \left[ \{y_{\beta,s}^{\pm 1}\}_{\beta\in\Delta^+}^{1\le s\le d_\beta} \right]^{\Sigma_{\underline{d}}} \simeq \bigotimes_{\beta\in\Delta^+} \mathbb{C}(\nu) \left[ \{y_{\beta,s}^{\pm 1}\}_{s=1}^{d_\beta} \right]^{\Sigma_{d_\beta}}. \tag{1.54}$$

Combining (1.53, 1.54) with Proposition 1.4(a), the formula (1.50) and the discussion of a "rank 1 reduction" after it, we see that there is a linear combination $F_{\underline{d}} = \sum_{h\in\mathsf{H}}^{\deg(h)=\underline{d}} c_h \Psi(e_h)$ such that $\phi_{\underline{d}}(F) = \phi_{\underline{d}}(F_{\underline{d}})$, with $c_h \in \mathbb{C}(\nu)$. Furthermore, we also have the equalities $\phi_{\underline{d}'}(F_{\underline{d}}) = 0$ for $\underline{d}' > \underline{d}$, due to Lemma 1.5.

This completes the proof of Lemma 1.7.    □

Now we are ready to prove the surjectivity of the homomorphism $\Psi$ of (1.36):

**Proposition 1.6** *We have $M = S^{(n)}$. In particular, the homomorphism $\Psi$ is surjective.*

Combining this with the injectivity of $\Psi$ (Proposition 1.3) implies Theorem 1.8.

*Proof* Let us show that any $F \in S_{\underline{k}}^{(n)}$ actually belongs to $M \cap S_{\underline{k}}^{(n)}$. Let $\underline{d}_{\max}$ and $\underline{d}_{\min}$ denote the maximal and the minimal elements of $T_{\underline{k}}$, respectively. The condition of Lemma 1.7 is vacuous for $\underline{d} = \underline{d}_{\max}$, hence, Lemma 1.7 applies. Replacing $F$ by $F - F_{\underline{d}_{\max}}$ and iterating this argument, we will eventually find $\widetilde{F} \in M$ such that $\phi_{\underline{d}}(F) = \phi_{\underline{d}}(\widetilde{F})$ for all $\underline{d} \in T_{\underline{k}}$. For $\underline{d} = \underline{d}_{\min}$, this yields $F = \widetilde{F}$, hence, $F \in M$.    □

Now we can prove that the ordered PBWD monomials $\{e_h\}_{h\in\mathsf{H}}$ span $U_\nu^>(L\mathfrak{sl}_n)$:

**Proposition 1.7** *The ordered PBWD monomials $\{e_h\}_{h\in\mathsf{H}}$ linearly span $U_\nu^>(L\mathfrak{sl}_n)$.*

*Proof* Invoking the injectivity of the homomorphism $\Psi\colon U_\nu^>(L\mathfrak{sl}_n) \to S^{(n)}$, Proposition 1.6 implies that the ordered PBWD monomials $\{e_h\}_{h\in\mathsf{H}}$ span $U_\nu^>(L\mathfrak{sl}_n)$.    □

Combining Propositions 1.5 and 1.7, we finally obtain:

**Corollary 1.2** $\{e_h\}_{h\in\mathsf{H}}$ *form a $\mathbb{C}(\nu)$-basis of $U_\nu^>(L\mathfrak{sl}_n)$.*

This completes our proof of Theorem 1.1 (as noted earlier, the claim for $U_\nu^<(L\mathfrak{sl}_n)$ can be easily deduced from Corollary 1.2 by using the anti-isomorphism (1.41)).

*Remark 1.5* We note that the above proof of Theorem 1.1 actually provides a much larger class of the PBWD bases for $U_\nu^>(L\mathfrak{sl}_n)$ and $U_\nu^<(L\mathfrak{sl}_n)$, with the PBWD basis elements given explicitly in the shuffle form (rather than long expressions in $e_{i,r}, f_{i,r}$).

### 1.2.3 Integral form $\mathfrak{S}^{(n)}$ and a proof of Theorem 1.3

Inspired by the above proof of Theorem 1.1 through the shuffle algebra realization of $U_v^>(L\mathfrak{sl}_n)$, let us now prove Theorem 1.3 by simultaneously providing the shuffle algebra realization of the integral form $\mathfrak{U}_v^>(L\mathfrak{sl}_n)$, which is of independent interest.

#### 1.2.3.1 Integral form $\mathfrak{S}^{(n)}$

For $\underline{k} \in \mathbb{N}^I$, consider a $\mathbb{Z}[\nu, \nu^{-1}]$-submodule $\widetilde{\mathfrak{S}}_{\underline{k}}^{(n)}$ of $S_{\underline{k}}^{(n)}$ consisting of those $F \in S_{\underline{k}}^{(n)}$ for which the numerator $f$ in (1.34) has the following form:

$$f \in (\nu - \nu^{-1})^{|\underline{k}|} \cdot \mathbb{Z}[\nu, \nu^{-1}]\left[\{x_{i,r}^{\pm 1}\}_{i \in I}^{1 \le r \le k_i}\right]^{\Sigma_{\underline{k}}}, \tag{1.55}$$

where $|\underline{k}| := \sum_{i \in I} k_i$. We define

$$\widetilde{\mathfrak{S}}^{(n)} := \bigoplus_{\underline{k} \in \mathbb{N}^I} \widetilde{\mathfrak{S}}_{\underline{k}}^{(n)},$$

which is easily seen to be a $\mathbb{Z}[\nu, \nu^{-1}]$-subalgebra of $S^{(n)}$. To this end, let us note that the shuffle product $F \star G$ of (1.32) still makes sense as $\frac{1}{\underline{k}! \cdot \underline{\ell}!} \cdot \text{Sym}(\cdots)$ may be simplified as a sum over so-called $(\underline{k}, \underline{\ell})$-shuffles, since $F$ and $G$ are symmetric.

Moreover, we note that $\Psi(\mathfrak{U}_v^>(L\mathfrak{sl}_n)) \subset \widetilde{\mathfrak{S}}^{(n)}$, due to the formula (1.48). The key goal of this subsection is to explicitly describe the image $\Psi(\mathfrak{U}_v^>(L\mathfrak{sl}_n))$ inside $S^{(n)}$.

**Example 1.11** Consider $F(x_1, x_2) = (\nu - \nu^{-1})^2(x_1 x_2)^r \in \widetilde{\mathfrak{S}}_2^{(2)}$ with $r \in \mathbb{Z}$. Then, due to Lemmas 1.8–1.9 below, we have $[2]_\nu \cdot F \in \Psi(\mathfrak{U}_v^>(L\mathfrak{sl}_2))$, but $F \notin \Psi(\mathfrak{U}_v^>(L\mathfrak{sl}_2))$.

We shall now introduce some more specialization maps. Pick $F \in \widetilde{\mathfrak{S}}_{\underline{k}}^{(n)}$. For any $\underline{d} = \{d_\beta\}_{\beta \in \Delta^+} \in \mathbb{N}^{\Delta^+}$ such that $\sum_{\beta \in \Delta^+} d_\beta \beta = \sum_{i \in I} k_i \alpha_i$, consider the specialization $\phi_{\underline{d}}(F) \in \mathbb{Z}[\nu, \nu^{-1}][\{y_{\beta,s}^{\pm 1}\}_{\beta \in \Delta^+}^{1 \le s \le d_\beta}]^{\Sigma_{\underline{d}}}$ of (1.40). First, we note that $\phi_{\underline{d}}(F)$ is divisible by

$$A_{\underline{d}} := (\nu - \nu^{-1})^{|\underline{k}|}. \tag{1.56}$$

Second, following the first part of the proof of Lemma 1.7 (see (i), (ii), and (1.52), which are consequences of the wheel conditions (1.35)), $\phi_{\underline{d}}(F)$ is also divisible by

$$B_{\underline{d}} := \prod_{(\beta,s)<(\beta',s')} (y_{\beta,s} - \nu^{-2} y_{\beta',s'})^{\#\{(j,j') \in [\beta] \times [\beta'] | j=j'\} - \delta_{\beta\beta'}} \times$$
$$\prod_{(\beta,s)<(\beta',s')} (y_{\beta,s} - \nu^2 y_{\beta',s'})^{\#\{(j,j') \in [\beta] \times [\beta'] | j=j'+1\}}, \tag{1.57}$$

where we use the order (1.51) on the set $\{(\beta,s)\}_{\beta \in \Delta^+}^{1 \le s \le d_\beta}$. Using these observations, we define the *reduced specialization map* (vanishing on degree $\ne \underline{k}$ components):

$$\varphi_{\underline{d}} \colon \widetilde{\mathfrak{S}}^{(n)} \longrightarrow \mathbb{Z}[\boldsymbol{\nu}, \boldsymbol{\nu}^{-1}][\{y_{\beta,s}^{\pm 1}\}_{\beta \in \Delta^+}^{1 \leq s \leq d_\beta}]^{\Sigma_{\underline{d}}} \qquad \text{via} \qquad \varphi_{\underline{d}}(F) := \frac{\phi_{\underline{d}}(F)}{A_{\underline{d}} B_{\underline{d}}} \,. \tag{1.58}$$

Let us introduce another type of specialization maps. Given a collection of positive integers $\underline{t} = \{t_{\beta,r}\}_{\beta \in \Delta^+}^{1 \leq r \leq \ell_\beta}$ ($\ell_\beta \in \mathbb{N}$), define a degree vector $\underline{d} = \{d_\beta\}_{\beta \in \Delta^+} \in \mathbb{N}^{\Delta^+}$ via $d_\beta := \sum_{r=1}^{\ell_\beta} t_{\beta,r}$. Let us now define the *vertical specialization map*

$$\varpi_{\underline{t}} \colon \mathbb{Z}[\boldsymbol{\nu}, \boldsymbol{\nu}^{-1}]\Big[\{y_{\beta,s}^{\pm 1}\}_{\beta \in \Delta^+}^{1 \leq s \leq d_\beta}\Big]^{\Sigma_{\underline{d}}} \longrightarrow \mathbb{Z}[\boldsymbol{\nu}, \boldsymbol{\nu}^{-1}][\{z_{\beta,r}^{\pm 1}\}_{\beta \in \Delta^+}^{1 \leq r \leq \ell_\beta}] \,. \tag{1.59}$$

To do so, for each $\beta \in \Delta^+$ we split the variables $\{y_{\beta,s}\}_{s=1}^{d_\beta}$ into $\ell_\beta$ groups of size $t_{\beta,r}$ each ($1 \leq r \leq \ell_\beta$) and specialize the variables in the $r$-th group to

$$\boldsymbol{\nu}^{-2} z_{\beta,r} \,, \ \boldsymbol{\nu}^{-4} z_{\beta,r} \,, \ \boldsymbol{\nu}^{-6} z_{\beta,r} \,, \ \dots \,, \ \boldsymbol{\nu}^{-2t_{\beta,r}} z_{\beta,r} \,.$$

For $K \in \mathbb{Z}[\boldsymbol{\nu}, \boldsymbol{\nu}^{-1}][\{y_{\beta,s}^{\pm 1}\}_{\beta \in \Delta^+}^{1 \leq s \leq d_\beta}]^{\Sigma_{\underline{d}}}$, we define $\varpi_{\underline{t}}(K)$ as the corresponding specialization of $K$ (it is independent of our splitting of the variables $\{y_{\beta,s}\}_{s=1}^{d_\beta}$ into groups).

Finally, let us combine (1.58) and (1.59) to produce the key tool in the study of the integral form $\Psi(\mathfrak{U}_{\boldsymbol{\nu}}^{>}(L\mathfrak{sl}_n))$ of $S^{(n)}$. To this end, given $\underline{k} \in \mathbb{N}^I$, $\underline{d} = \{d_\beta\}_{\beta \in \Delta^+} \in \mathbb{N}^{\Delta^+}$, and a collection of positive integers $\underline{t} = \{t_{\beta,r}\}_{\beta \in \Delta^+}^{1 \leq r \leq \ell_\beta}$ satisfying the compatibilities:

$$\sum_{\beta \in \Delta^+} d_\beta \beta = \sum_{i \in I} k_i \alpha_i \qquad \text{and} \qquad \sum_{r=1}^{\ell_\beta} t_{\beta,r} = d_\beta \,, \tag{1.60}$$

we define the **cross specialization map** (vanishing on degree $\neq \underline{k}$ components):

$$\Upsilon_{\underline{d},\underline{t}} \colon \widetilde{\mathfrak{S}}^{(n)} \longrightarrow \mathbb{Z}[\boldsymbol{\nu}, \boldsymbol{\nu}^{-1}]\Big[\{z_{\beta,r}^{\pm 1}\}_{\beta \in \Delta^+}^{1 \leq r \leq \ell_\beta}\Big] \qquad \text{via} \qquad \Upsilon_{\underline{d},\underline{t}}(F) := \varpi_{\underline{t}}(\varphi_{\underline{d}}(F)) \,. \tag{1.61}$$

**Definition 1.1** We call $F \in S_{\underline{k}}^{(n)}$ **integral** if $F \in \widetilde{\mathfrak{S}}_{\underline{k}}^{(n)}$ and $\Upsilon_{\underline{d},\underline{t}}(F)$ is divisible by $\prod_{\beta \in \Delta^+}^{1 \leq r \leq \ell_\beta} [t_{\beta,r}]_{\boldsymbol{\nu}}!$ (the product of $\boldsymbol{\nu}$-factorials) for all possible $\underline{d}$ and $\underline{t}$ satisfying (1.60).

**Example 1.12** In the simplest "rank 1" case of $n = 2$, a symmetric Laurent polynomial $F \in S_k^{(2)}$ is integral if and only if it has the form $F = (\boldsymbol{\nu} - \boldsymbol{\nu}^{-1})^k \cdot \tilde{F}$ with $\tilde{F} \in \mathbb{Z}[\boldsymbol{\nu}, \boldsymbol{\nu}^{-1}][\{x_p^{\pm 1}\}_{p=1}^k]^{\Sigma_k}$ satisfying the following divisibility condition:

$$\tilde{F}(\boldsymbol{\nu}^{-2} z_1, \boldsymbol{\nu}^{-4} z_1, \dots, \boldsymbol{\nu}^{-2k_1} z_1, \dots, \boldsymbol{\nu}^{-2} z_\ell, \dots, \boldsymbol{\nu}^{-2k_\ell} z_\ell) \text{ is divisible by } [k_1]_{\boldsymbol{\nu}}! \cdots [k_\ell]_{\boldsymbol{\nu}}! \tag{1.62}$$

for any decomposition $k = k_1 + \cdots + k_\ell$ into a sum of positive integers.

Let $\mathfrak{S}_{\underline{k}}^{(n)} \subset \widetilde{\mathfrak{S}}_{\underline{k}}^{(n)}$ denote the $\mathbb{Z}[\boldsymbol{\nu}, \boldsymbol{\nu}^{-1}]$-submodule of all integral elements and set

$$\mathfrak{S}^{(n)} := \bigoplus_{\underline{k} \in \mathbb{N}^I} \mathfrak{S}_{\underline{k}}^{(n)} \,.$$

The following is the key result of this section:

**Theorem 1.13** *The $\mathbb{C}(\nu)$-algebra isomorphism $\Psi\colon U_\nu^>(L\mathfrak{sl}_n) \overset{\sim}{\longrightarrow} S^{(n)}$ of Theorem 1.8 gives rise to a $\mathbb{Z}[\nu,\nu^{-1}]$-algebra isomorphism $\Psi\colon \mathfrak{U}_\nu^>(L\mathfrak{sl}_n) \overset{\sim}{\longrightarrow} \mathfrak{S}^{(n)}$.*

Before proceeding to the proof of this result, let us record the following corollary:

**Corollary 1.3** *(a) $\mathfrak{S}^{(n)}$ is a $\mathbb{Z}[\nu,\nu^{-1}]$-subalgebra of $S^{(n)}$.*

*(b) The subalgebra $\mathfrak{U}_\nu^>(L\mathfrak{sl}_n)$ is independent of the choices (1)–(3) preceding (1.19).*

### 1.2.3.2  Proofs of Theorem 1.3 and Theorem 1.13

Both Theorem 1.3(b) and Theorem 1.13 follow from the following two statements:

(I)  For any $k \geq 1$, $\{\beta_p\}_{p=1}^k \subset \Delta^+$, $\{r_p\}_{p=1}^k \subset \mathbb{Z}$, we have $\Psi(\widetilde{e}_{\beta_1}(r_1)\cdots\widetilde{e}_{\beta_k}(r_k)) \in \mathfrak{S}^{(n)}$.

(II)  Any $F \in \mathfrak{S}^{(n)}$ may be written as a $\mathbb{Z}[\nu,\nu^{-1}]$-linear combination of $\{\Psi(\widetilde{e}_h)\}_{h\in\mathsf{H}}$.

The proof of (I) shall be easily deduced from our definition of $\mathfrak{S}^{(n)}$, while the proof of (II) will follow from Lemma 1.7 and the validity of (II) for $n=2$, see Lemma 1.9.

We start by establishing both (I), (II) in the simplest "rank 1" case of $n=2$, for which the integral form $\mathfrak{S}^{(2)}$ has been already described in the Example above. Set $\widetilde{e}_r := (\nu - \nu^{-1})e_r \in U_\nu^>(L\mathfrak{sl}_2)$ for $r \in \mathbb{Z}$. The following lemma implies (I) for $n=2$:

**Lemma 1.8** *For any $k \geq 1$ and $r_1,\ldots,r_k \in \mathbb{Z}$, we have $\Psi(\widetilde{e}_{r_1}\cdots\widetilde{e}_{r_k}) \in \mathfrak{S}_k^{(2)}$.*

*Proof* Pick any decomposition $k = k_1 + \cdots + k_\ell$ with all $k_p \geq 1$. According to the above Example, it suffices to prove that as we specialize the variables $x_1,\ldots,x_k$ to $\{\nu^{-2r}z_p\}_{1\leq p\leq \ell}^{1\leq r\leq k_p}$, the image of any summand from the symmetrization appearing in $\Psi(e_{r_1}\cdots e_{r_k}) \in S_k^{(2)}$ will be divisible by $\prod_{p=1}^\ell [k_p]_\nu!$, as required in (1.62).

To this end, let us fix $1 \leq p \leq \ell$ and consider the relative position of the variables $\nu^{-2}z_p, \nu^{-4}z_p, \ldots, \nu^{-2k_p}z_p$. If there is an index $1 \leq r < k_p$ such that $\nu^{-2(r+1)}z_p$ is placed to the left of $\nu^{-2r}z_p$ (that is, $\Psi(e_{r_a})$ is evaluated at $\nu^{-2(r+1)}z_p$ and $\Psi(e_{r_b})$ is evaluated at $\nu^{-2r}z_p$ for $a < b$), then the above specialization of the corresponding $\zeta$-factor equals $\frac{\nu^{-2(r+1)}z_p - \nu^{-2r}\cdot\nu^{-2r}z_p}{\nu^{-2(r+1)}z_p - \nu^{-2r}z_p} = 0$. However, if $\nu^{-2r}z_p$ is placed to the left of $\nu^{-2(r+1)}z_p$ for all $1 \leq r < k_p$, then the total contribution of the specializations of the corresponding $\zeta$-factors equals

$$\prod_{1\leq r<s\leq k_p} \frac{\nu^{-2r}z_p - \nu^{-2}\cdot\nu^{-2s}z_p}{\nu^{-2r}z_p - \nu^{-2s}z_p} = \nu^{-\frac{k_p(k_p-1)}{2}}[k_p]_\nu! \tag{1.63}$$

Combining this for all $p \in \{1,2,\ldots,\ell\}$, we see that $\prod_{p=1}^\ell [k_p]_\nu!$ indeed divides the above specialization of $\Psi(e_{r_1}\cdots e_{r_k})$. This completes our proof of Lemma 1.8.   $\square$

To simplify our exposition, let us assume that the total order $\leq$ on $\mathbb{Z}$ is the usual $\leq$. The following result establishes (II) for $n=2$:

**Lemma 1.9** *Any symmetric Laurent polynomial $\bar{F} \in \mathbb{Z}[\boldsymbol{\nu}, \boldsymbol{\nu}^{-1}][\{x_p^{\pm 1}\}_{p=1}^k]^{\Sigma_k}$ satisfying the divisibility condition (1.62) is a $\mathbb{Z}[\boldsymbol{\nu}, \boldsymbol{\nu}^{-1}]$-linear combination of $\{\Psi(e_h)\}_{h \in \mathsf{H}}$.*

*Proof* We may assume that $\bar{F}$ is homogeneous of the total degree $N$. Let $V_N$ denote the set of all ordered $k$-tuples of integers $\underline{r} = (r_1, r_2, \ldots, r_k)$, $r_1 \leq \cdots \leq r_k$, such that $r_1 + \cdots + r_k = N$. This set is totally ordered with respect to the lexicographical order:

$$\underline{r} < \underline{r}' \quad \text{iff} \quad \exists\, 1 \leq p \leq k \ \text{s.t.} \ r_p < r_p' \ \text{and} \ r_s = r_s' \ \text{for all} \ s > p\,.$$

Let us present $\bar{F}$ as a linear combination of the monomial symmetric polynomials:

$$\bar{F}(x_1, \ldots, x_k) = \sum_{\underline{r} \in V_N} \mu_{\underline{r}} m_{\underline{r}}(x_1, \ldots, x_k) \quad \text{with} \quad \mu_{\underline{r}} \in \mathbb{Z}[\boldsymbol{\nu}, \boldsymbol{\nu}^{-1}]\,.$$

Pick the maximal $\underline{r}_{\max} = (r_1, \ldots, r_k)$ from the finite set $V_N(\bar{F}) := \{\underline{r} \in V_N | \mu_{\underline{r}} \neq 0\}$ and consider a decomposition $k = k_1 + \cdots + k_\ell$ such that

$$r_1 = \cdots = r_{k_1} < r_{k_1+1} = \cdots = r_{k_1+k_2} < \cdots < r_{k_1+\cdots+k_{\ell-1}+1} = \cdots = r_k\,.$$

Evaluating $\bar{F}$ at the resulting specialization $\{\boldsymbol{\nu}^{-2s} z_p\}_{1 \leq p \leq \ell}^{1 \leq s \leq k_p}$, we see that the coefficient of the lexicographically largest monomial in the variables $\{z_p\}_{p=1}^\ell$ equals $\mu_{\underline{r}_{\max}}$, up to an integer power of $\boldsymbol{\nu}$. Therefore, the divisibility condition (1.62) implies that

$$\frac{\mu_{\underline{r}_{\max}}}{\prod_{p=1}^\ell [k_p]_{\boldsymbol{\nu}}!} \in \mathbb{Z}[\boldsymbol{\nu}, \boldsymbol{\nu}^{-1}]\,. \tag{1.64}$$

Set $\bar{F}^{(0)} := \bar{F}$ and define $\bar{F}^{(1)} \in \mathbb{Z}[\boldsymbol{\nu}, \boldsymbol{\nu}^{-1}][\{x_p^{\pm 1}\}_{p=1}^k]^{\Sigma_k}$ via

$$\bar{F}^{(1)} := \bar{F}^{(0)} - \boldsymbol{\nu}^{\sum_{p=1}^\ell k_p (k_p-1)/2} \frac{\mu_{\underline{r}_{\max}}}{\prod_{p=1}^\ell [k_p]_{\boldsymbol{\nu}}!} \Psi(e_{r_1} \cdots e_{r_k})\,. \tag{1.65}$$

We note that $\bar{F}^{(1)}$ satisfies the divisibility condition (1.62), due to (1.64) and Lemma 1.8. Applying the same argument to $\bar{F}^{(1)}$ in place of $F = \bar{F}^{(0)}$, we obtain $\bar{F}^{(2)}$ also satisfying (1.62). Iterating this process, we thus construct a sequence of symmetric Laurent polynomials $\{\bar{F}^{(s)}\}_{s \in \mathbb{N}}$ satisfying the divisibility condition (1.62).

However, invoking our proof of Proposition 1.4, specifically the formula (1.45), we note that the sequence $\underline{r}_{\max}^{(s)} \in V_N$ of the maximal elements of $V_N(\bar{F}^{(s)})$ strictly decreases. Meanwhile, the sequence of the minimal powers of any variable in $\bar{F}^{(s)}$ is a non-decreasing sequence. Hence, $\bar{F}^{(s)} = 0$ for some $s \in \mathbb{N}$.

This completes our proof of Lemma 1.9. □

Let us now generalize the above arguments to prove (I) and (II) for any $n > 2$. The proof of the former is quite similar (though is more tedious) to that of Lemma 1.8:

**Lemma 1.10** *For any $m \geq 1$, $\{\beta_p\}_{p=1}^m \subset \Delta^+$, and $\{r_p\}_{p=1}^m \subset \mathbb{Z}$, we have:*

$$\Psi\left(\widetilde{e}_{\beta_1}(r_1) \cdots \widetilde{e}_{\beta_m}(r_m)\right) \in \mathfrak{S}^{(n)}\,.$$

*Proof* Define $\underline{k} \in \mathbb{N}^I$ via $\sum_{i \in I} k_i \alpha_i = \sum_{p=1}^{m} \beta_p$, so that $F := \Psi\left(\widetilde{e}_{\beta_1}(r_1) \cdots \widetilde{e}_{\beta_m}(r_m)\right)$ is in $S_{\underline{k}}^{(n)}$. First, we note that $F \in \widetilde{\mathfrak{S}}_{\underline{k}}^{(n)}$ as $(\boldsymbol{v} - \boldsymbol{v}^{-1})^{|\underline{k}|}$ divides $F$, due to Lemma 1.4.

It remains to show that $\Upsilon_{\underline{d},\underline{t}}(F)$ satisfies the divisibility condition of Definition 1.1 for any collections $\underline{d} = \{d_\beta\}_{\beta \in \Delta^+} \in \mathbb{N}^{\Delta^+}$ and $\underline{t} = \{t_{\beta,s}\}_{\beta \in \Delta^+}^{1 \le s \le \ell_\beta}$ satisfying (1.60). To this end, we recall that the cross specialization $\Upsilon_{\underline{d},\underline{t}}(F)$ is computed in three steps:

- first, we specialize $x_{*,*}$'s in $f$ of (1.34) to $\boldsymbol{v}^?$-multiples of $y_{*,*}$'s as in (1.40);
- second, we divide that specialization by the product of appropriate powers of $\boldsymbol{v} - \boldsymbol{v}^{-1}$ and linear terms in $y_{*,*}$'s due to the wheel conditions, see (1.56)–(1.58);
- finally, we specialize $y_{*,*}$-variables to $\boldsymbol{v}^?$-multiples of $z_{*,*}$-variables as in (1.59).

Fix $\beta \in \Delta^+$ and $1 \le s \le \ell_\beta$, and consider those $x_{*,*}$'s that eventually got specialized to $\boldsymbol{v}^? z_{\beta,s}$. Without loss of generality, we may assume those are $\{x_{i,r}\}_{j(\beta) \le i \le i(\beta)}^{1 \le r \le t_{\beta,s}}$. We may also assume that $x_{i,r}$ got specialized to $\boldsymbol{v}^{-i} y_{\beta,r}$ under the first specialization (1.40), while $y_{\beta,r}$ got further specialized to $\boldsymbol{v}^{-2r} z_{\beta,s}$ under the second specialization (1.59), for any $j(\beta) \le i \le i(\beta)$ and $1 \le r \le t_{\beta,s}$.

For $j(\beta) \le i < i(\beta)$ and $1 \le r \ne r' \le t_{\beta,s}$, consider the relative position of the variables $x_{i,r}, x_{i,r'}, x_{i+1,r'}$ in any summand from the symmetrization appearing in $F$. As $x_{i,r}, x_{i,r'}$ cannot enter the same function $\Psi(\widetilde{e}_*(*))$, $x_{i,r}$ is placed either to the left of $x_{i,r'}$ or to the right. In the former case, we gain the factor $\zeta_{i,i}(x_{i,r}/x_{i,r'})$, which upon the specialization $\phi_{\underline{d}}$ contributes the factor $(y_{\beta,r} - \boldsymbol{v}^{-2} y_{\beta,r'})$. Likewise, if $x_{i+1,r'}$ is placed to the left of $x_{i,r}$, we gain the factor $\zeta_{i+1,i}(x_{i+1,r'}/x_{i,r})$, which upon the specialization $\phi_{\underline{d}}$ contributes the factor $(y_{\beta,r} - \boldsymbol{v}^{-2} y_{\beta,r'})$ as well. In the remaining case, when $x_{i,r'}$ is placed to the left of $x_{i,r}$ while $x_{i+1,r'}$ is not, we gain the factor $\zeta_{i,i+1}(x_{i,r'}/x_{i+1,r'})$, which specializes to 0 under the specialization $\phi_{\underline{d}}$. As $i$ ranges from $j(\beta)$ up to $i(\beta) - 1$, we thus gain the $(i(\beta) - j(\beta))$-th power of $(y_{\beta,r} - \boldsymbol{v}^{-2} y_{\beta,r'})$. The latter precisely coincides with the power of $(y_{\beta,r} - \boldsymbol{v}^{-2} y_{\beta,r'})$ in $B_{\underline{d}}$ of (1.57), by which we divide $\phi_{\underline{d}}(F)$ to obtain $\varphi_{\underline{d}}(F)$ of (1.58).

However, we have not used yet $\zeta$-factors $\zeta_{i(\beta),i(\beta)}(x_{i(\beta),r}/x_{i(\beta),r'})$ for $x_{i(\beta),r}$ placed to the left of $x_{i(\beta),r'}$. If $x_{i(\beta),r+1}$ is placed to the left of $x_{i(\beta),r}$ for some $1 \le r < t_{\beta,s}$, then $\zeta_{i(\beta),i(\beta)}(x_{i(\beta),r+1}/x_{i(\beta),r})$ specializes to 0 under (1.59). In the remaining case, when each $x_{i(\beta),r}$ is placed to the left of $x_{i(\beta),r+1}$, the total contribution of the specializations of the corresponding $\zeta$-factors equals $\boldsymbol{v}^{-t_{\beta,s}(t_{\beta,s}-1)/2}[t_{\beta,s}]_{\boldsymbol{v}}!$, cf. (1.63).

This completes the proof of Lemma 1.10. □

Let $\widetilde{M} \subset S^{(n)}$ denote the $\mathbb{Z}[\boldsymbol{v},\boldsymbol{v}^{-1}]$-span of $\{\Psi(\widetilde{e}_h)\}_{h \in \mathbb{H}}$, the images of the ordered PBWD monomials, and define $T_{\underline{k}}$ as before. Generalizing Lemma 1.7, we have:

**Lemma 1.11** *Consider $F \in \widetilde{\mathfrak{S}}_{\underline{k}}^{(n)}$ and a degree vector $\underline{d} \in T_{\underline{k}}$. If $\phi_{\underline{d}'}(F) = 0$ for all $\underline{d}' \in T_{\underline{k}}$ such that $\underline{d}' > \underline{d}$, then there exists an element $F_{\underline{d}} \in \widetilde{M} \cap S_{\underline{k}}^{(n)}$ such that*

$$\phi_{\underline{d}}(F) = \phi_{\underline{d}}(F_{\underline{d}}) \qquad \text{and} \qquad \phi_{\underline{d}'}(F_{\underline{d}}) = 0 \quad \text{for all} \quad \underline{d}' > \underline{d}.$$

*Proof* The proof is completely analogous to that of Lemma 1.7. More precisely, combining the formulas (1.53, 1.54), the condition $F \in \widetilde{\mathfrak{S}}_{\underline{k}}^{(n)}$ together with the formula (1.50) and the discussion of a "rank 1 reduction" after it, the result follows from its $n = 2$ counterpart. The latter has been already established in Lemma 1.9. □

Combining Lemmas 1.10–1.11, we obtain (II) similarly to Proposition 1.6:

**Proposition 1.8** *We have $\widetilde{M} = \mathfrak{S}^{(n)}$ and $\Psi \colon \mathfrak{U}_v^>(L\mathfrak{sl}_n) \to \mathfrak{S}^{(n)}$ is surjective.*

Combining this with the injectivity of $\Psi$ (Proposition 1.3) implies Theorem 1.13.

**Corollary 1.4** *$\{\widetilde{e}_h\}_{h \in \mathsf{H}}$ is a basis of the free $\mathbb{Z}[v, v^{-1}]$-module $\mathfrak{U}_v^>(L\mathfrak{sl}_n)$.*

This completes our proof of Theorem 1.3(b) (as before, the claim for $\mathfrak{U}_v^<(L\mathfrak{sl}_n)$ follows from Corollary 1.4 by using the algebra anti-isomorphism (1.41)). Finally, the result of Theorems 1.3(a) has been already established in Corollary 1.3(b).

### 1.2.4 Integral form $\mathsf{S}^{(n)}$ and a proof of Theorem 1.6

Likewise, let us now prove Theorem 1.6 by simultaneously providing the shuffle algebra realization of the integral form $\mathsf{U}_v^>(L\mathfrak{sl}_n)$, which is of independent interest.

**Definition 1.2** We call $F \in S_{\underline{k}}^{(n)}$ **good** if it satisfies the following two properties:

(i) $F = \frac{f(x_{1,1}, \ldots, x_{n-1,k_{n-1}})}{\prod_{i=1}^{n-2} \prod_{r \le k_i}^{r' \le k_{i+1}} (x_{i,r} - x_{i+1,r'})}$ with $f \in \mathbb{Z}[v, v^{-1}] \left[ \{x_{i,r}^{\pm 1}\}_{i \in I}^{1 \le r \le k_i} \right]^{\Sigma_{\underline{k}}}$;

(ii) the specialization $\phi_{\underline{d}}(F)$ of (1.40) is divisible by $(v - v^{-1})^{\sum_{\beta \in \Delta^+} d_\beta (i(\beta) - j(\beta))}$ for any degree vector $\underline{d} = \{d_\beta\}_{\beta \in \Delta^+} \in \mathbb{N}^{\Delta^+}$ satisfying $\sum_{\beta \in \Delta^+} d_\beta \beta = \sum_{i \in I} k_i \alpha_i$.

**Example 1.14** In the simplest "rank 1" case of $n = 2$, we have:

$$F \in S_k^{(2)} \text{ is good} \iff F \in \mathbb{Z}[v, v^{-1}][\{x_p^{\pm 1}\}_{p=1}^k]^{\Sigma_k}.$$

Let $\mathsf{S}_{\underline{k}}^{(n)} \subset S_{\underline{k}}^{(n)}$ denote the $\mathbb{Z}[v, v^{-1}]$-submodule of all good elements and set

$$\mathsf{S}^{(n)} := \bigoplus_{\underline{k} \in \mathbb{N}^I} \mathsf{S}_{\underline{k}}^{(n)}.$$

The following is the key result of this section:

**Theorem 1.15** *The $\mathbb{C}(v)$-algebra isomorphism $\Psi \colon U_v^>(L\mathfrak{sl}_n) \xrightarrow{\sim} S^{(n)}$ of Theorem 1.8 gives rise to a $\mathbb{Z}[v, v^{-1}]$-algebra isomorphism $\Psi \colon \mathsf{U}_v^>(L\mathfrak{sl}_n) \xrightarrow{\sim} \mathsf{S}^{(n)}$.*

**Corollary 1.5** *$\mathsf{S}^{(n)}$ is a $\mathbb{Z}[v, v^{-1}]$-subalgebra of $S^{(n)}$.*

Both Theorem 1.6 and Theorem 1.15 follow from the following three results:

(I) For any $N \in \mathbb{N}$, $\{(i_p, r_p, k_p)\}_{p=1}^N \subset I \times \mathbb{Z} \times \mathbb{N}$, we have $\Psi(\mathsf{e}_{i_1, r_1}^{(k_1)} \cdots \mathsf{e}_{i_N, r_N}^{(k_N)}) \in \mathsf{S}^{(n)}$.

(II) For any $h \in \mathsf{H}$ and $\mathsf{e}_h$ defined in (1.29), we have $\Psi(\mathsf{e}_h) \in \mathsf{S}^{(n)}$.

(III) Any $F \in \mathsf{S}^{(n)}$ may be written as a $\mathbb{Z}[\nu, \nu^{-1}]$-linear combination of $\{\Psi(\mathsf{e}_h)\}_{h \in \mathsf{H}}$.

The proof of the first statement (I) is straightforward:

**Lemma 1.12** *For any $N \in \mathbb{N}$ and $\{(i_p, r_p, k_p)\}_{p=1}^{N} \subset I \times \mathbb{Z} \times \mathbb{N}$, we have:*

$$\Psi\left(\mathsf{e}_{i_1,r_1}^{(k_1)} \cdots \mathsf{e}_{i_N,r_N}^{(k_N)}\right) \in \mathsf{S}^{(n)} \,.$$

*Proof* The validity of (i) for $F = \Psi(\mathsf{e}_{i_1,r_1}^{(k_1)} \cdots \mathsf{e}_{i_N,r_N}^{(k_N)})$ follows immediately from the equality $\Psi(\mathsf{e}_{i,r}^{(k)}) = \nu^{-\frac{k(k-1)}{2}} (x_1 \cdots x_k)^r$, due to Lemma 1.3. To verify (ii), pick any degree vector $\underline{d} = \{d_\beta\}_{\beta \in \Delta^+} \in \mathbb{N}^{\Delta^+}$ such that $\sum_{\beta \in \Delta^+} d_\beta \beta = \sum_{p=1}^{N} k_p \alpha_{i_p}$. For any $\beta \in \Delta^+$ and $1 \leq s \leq d_\beta$, consider $\zeta$-factors between those pairs of $x_{*,*}$-variables that are specialized to $\nu^{-\ell} y_{\beta,s}$ and $\nu^{-\ell-1} y_{\beta,s}$ with $j(\beta) \leq \ell < i(\beta)$ under the specialization map (1.40). Each of them contributes a multiple of $(\nu - \nu^{-1})$ into $\phi_{\underline{d}}(F)$, and there are exactly $\sum_{\beta \in \Delta^+} d_\beta (i(\beta) - j(\beta))$ of such pairs. Hence, (ii) indeed holds. $\qquad \square$

The above lemma can be rephrased as:

$$\Psi\left(\mathsf{U}_\nu^>(L\mathfrak{sl}_n)\right) \subseteq \mathsf{S}^{(n)} \,. \tag{1.66}$$

Thus, the second result (II) follows from (previously stated Proposition 1.2):

**Lemma 1.13** *For any $\beta \in \Delta^+$, $r \in \mathbb{Z}$, $k \in \mathbb{N}$, we have $\mathsf{e}_\beta(r)^{(k)} \in \mathsf{U}_\nu^>(L\mathfrak{sl}_n)$.*

*Proof* For $\beta = \alpha_j + \alpha_{j+1} + \cdots + \alpha_i$ and $r \in \mathbb{Z}$, consider the PBWD basis element $e_\beta(r)$ of (1.19), and let $\bar{e}_\beta(r)$ denote the particular choice of (1.20). Combining Lemma 1.4 with the Example following after it, we get $\Psi(e_\beta(r)) = \Psi(\bar{e}_\beta(r)) \cdot c \cdot x_{j,1}^{m_j} x_{j+1,1}^{m_{j+1}} \cdots x_{i,1}^{m_i}$ for some $\underline{m} = (m_j, m_{j+1}, \ldots, m_i) \in \mathbb{Z}^{i-j+1}$ satisfying $m_j + \cdots + m_i = 0$ and $c \in \pm \nu^\mathbb{Z}$. Thus, we have:

$$\Psi\left(\mathsf{e}_\beta(r)^{(k)}\right) = \Psi\left(\bar{\mathsf{e}}_\beta(r)^{(k)}\right) \cdot c^k \cdot \prod_{\iota=j}^{i} (x_{\iota,1} \cdots x_{\iota,k})^{m_\iota} \,. \tag{1.67}$$

Consider an algebra automorphism $\sigma_{\underline{m}}$ of $U_\nu^>(L\mathfrak{sl}_n)$ that maps $e_{\iota,k} \mapsto e_{\iota,k+m_\iota}$ if $j \leq \iota \leq i$ and $e_{\iota,k} \mapsto e_{\iota,k}$ otherwise. Since $\sigma_{\underline{m}}$ maps $\mathsf{e}_{\iota,r}^{(k)}$ to $\mathsf{e}_{\iota,r+m_\iota}^{(k)}$ or $\mathsf{e}_{\iota,r}^{(k)}$, it does induce the same-named automorphism of $\mathsf{U}_\nu^>(L\mathfrak{sl}_n)$. Then, combining (1.67) and Theorem 1.8, we see that $\mathsf{e}_\beta(r)^{(k)} = c^k \sigma_{\underline{m}} \bar{\mathsf{e}}_\beta(r)^{(k)}$. However, the inclusion $\bar{\mathsf{e}}_\beta(r)^{(k)} \in \mathsf{U}_\nu^>(L\mathfrak{sl}_n)$ has been already established in our proof of Proposition 1.2. Therefore, we get $\mathsf{e}_\beta(r)^{(k)} \in \mathsf{U}_\nu^>(L\mathfrak{sl}_n)$, which completes our proof of Lemma 1.13 (and hence also of Proposition 1.2). $\qquad \square$

Let $\mathsf{M} \subset \mathsf{S}^{(n)}$ denote the $\mathbb{Z}[\nu, \nu^{-1}]$-span of $\{\Psi(\mathsf{e}_h)\}_{h \in \mathsf{H}}$ and define $T_{\underline{k}}$ as before. Generalizing Lemma 1.7, we have:

**Lemma 1.14** *Consider $F \in \mathsf{S}^{(n)}_{\underline{k}}$ and a degree vector $\underline{d} \in T_{\underline{k}}$. If $\phi_{\underline{d}'}(F) = 0$ for all $\underline{d}' \in T_{\underline{k}}$ such that $\underline{d}' > \underline{d}$, then there exists an element $F_{\underline{d}} \in \mathsf{M} \cap S^{(n)}_{\underline{k}}$ such that*

$$\phi_{\underline{d}}(F) = \phi_{\underline{d}}(F_{\underline{d}}) \qquad \text{and} \qquad \phi_{\underline{d}'}(F_{\underline{d}}) = 0 \quad \text{for all} \quad \underline{d}' > \underline{d} \,.$$

*Proof* The proof is completely analogous to that of Lemma 1.7. More precisely, combining the formula (1.49) with the condition $F \in \mathsf{S}^{(n)}_{\underline{k}}$, we obtain the formula (1.53) with the factor $G$ necessarily of the form $G \in \mathbb{Z}[\nu, \nu^{-1}] \left[ \{y^{\pm 1}_{\beta, s}\}^{1 \leq s \leq d_\beta}_{\beta \in \Delta^+} \right]^{\Sigma_{\underline{d}}}$. To this end, it suffices to note that the power of $(\nu - \nu^{-1})$ in the condition (ii) of Definition 1.2 precisely coincides with the overall power of $(\nu - \nu^{-1})$ arising from the factors $G_\beta$ in (1.49). Therefore, the result follows from its "rank 1" $n = 2$ counterpart, which states that $\{\frac{x^{r_1}}{[r_1]_\nu!} \star \cdots \star \frac{x^{r_N}}{[r_N]_\nu!}\}^{r_1 < \cdots < r_N}_{N \in \mathbb{N}}$ is a basis of the free $\mathbb{Z}[\nu, \nu^{-1}]$-module $\mathsf{S}^{(2)}$. The proof of this claim is completely analogous to that of Proposition 1.4(a): the key point is that the formula (1.45) will get replaced by $\nu_{\underline{r}} = \nu^{-\sum^N_{p=1} \frac{r_p(r_p - 1)}{2}}$. $\qquad \square$

Combining Lemmas 1.12–1.14, we obtain (III) similarly to Proposition 1.6:

**Proposition 1.9** *We have $\mathsf{M} = \mathsf{S}^{(n)}$ and $\Psi \colon \mathsf{U}^>_\nu(L\mathfrak{sl}_n) \to \mathsf{S}^{(n)}$ is surjective.*

This completes our proof of Theorem 1.15 and Theorem 1.6 (as before, the claim for $\mathsf{U}^<_\nu(L\mathfrak{sl}_n)$ follows by using the algebra anti-isomorphism (1.41)).

*Remark 1.6* We note that an interesting duality between the integral forms $\mathsf{U}^>_\nu(L\mathfrak{sl}_n)$ and $\mathfrak{U}^>_\nu(L\mathfrak{sl}_n)$ (similarly, between $\mathsf{U}^<_\nu(L\mathfrak{sl}_n)$ and $\mathfrak{U}^>_\nu(L\mathfrak{sl}_n)$) was established in [18] using the shuffle realizations of these integral forms from Theorems 1.13 and 1.15.

## 1.3 Generalizations to two-parameter quantum loop algebras

The two-parameter quantum loop algebra $U_{\nu_1, \nu_2}(L\mathfrak{sl}_n)$ was first introduced in [13][4] as a generalization of $U_\nu(L\mathfrak{sl}_n)$, the latter recovered by setting $\nu_1 = \nu_2^{-1} = \nu$ and identifying some Cartan elements, see [13, Remark 3.3(4)]. The key results of [13] are:

- [13, Theorem 3.12]: the Drinfeld-Jimbo type realization of $U_{\nu_1, \nu_2}(L\mathfrak{sl}_n)$,
- [13, Theorem 3.11]: the PBW basis of its subalgebras $U^>_{\nu_1, \nu_2}(L\mathfrak{sl}_n), U^<_{\nu_1, \nu_2}(L\mathfrak{sl}_n)$.

However, we note that [13, Theorem 3.11] was stated without any glimpse of a proof.

The primary goal of this section is to generalize Theorem 1.1 to the case of $U^>_{\nu_1, \nu_2}(L\mathfrak{sl}_n)$, thus finally proving [13, Theorem 3.11]. Along the way, we also generalize Theorem 1.8 by providing the shuffle algebra realization of $U^>_{\nu_1, \nu_2}(L\mathfrak{sl}_n)$, which is of independent interest. The latter is also used to construct the PBWD bases for the integral form $\mathfrak{U}^>_{\nu_1, \nu_2}(L\mathfrak{sl}_n)$ of $U^>_{\nu_1, \nu_2}(L\mathfrak{sl}_n)$, thus generalizing Theorem 1.3.

---

[4] To be more precise, this recovers the algebra of *loc.cit.* with the trivial central charges. However, we can ignore this since their "positive" subalgebras are naturally isomorphic.

### 1.3.1 Quantum loop algebra $U^>_{\nu_1,\nu_2}(L\mathfrak{sl}_n)$ and its PBWD bases

To simplify our exposition, we shall consider only the "positive" subalgebra $U^>_{\nu_1,\nu_2}(L\mathfrak{sl}_n)$ of $U_{\nu_1,\nu_2}(L\mathfrak{sl}_n)$. Let $\nu_1, \nu_2$ be two independent formal variables and set $\mathbb{K} := \mathbb{C}(\nu_1^{1/2}, \nu_2^{1/2})$. Following [13, Definition 3.1], we define $U^>_{\nu_1,\nu_2}(L\mathfrak{sl}_n)$ as the associative $\mathbb{K}$-algebra generated by $\{e_{i,r}\}^{r\in\mathbb{Z}}_{i\in I}$ with the following defining relations:

$$\Big(z - (\langle j,i\rangle\langle i,j\rangle)^{1/2}w\Big) e_i(z)e_j(w) = \Big(\langle j,i\rangle z - (\langle j,i\rangle\langle i,j\rangle^{-1})^{1/2}w\Big) e_j(w)e_i(z) \tag{1.68}$$

as well as quantum Serre relations:

$$\begin{aligned}
&e_i(z)e_j(w) = e_j(w)e_i(z) \quad \text{if} \quad c_{ij} = 0, \\
&[e_i(z_1), [e_i(z_2), e_{i+1}(w)]_{\nu_2}]_{\nu_1} + [e_i(z_2), [e_i(z_1), e_{i+1}(w)]_{\nu_2}]_{\nu_1} = 0, \\
&[e_i(z_1), [e_i(z_2), e_{i-1}(w)]_{\nu_2^{-1}}]_{\nu_1^{-1}} + [e_i(z_2), [e_i(z_1), e_{i-1}(w)]_{\nu_2^{-1}}]_{\nu_1^{-1}} = 0,
\end{aligned} \tag{1.69}$$

where $e_i(z) = \sum_{r\in\mathbb{Z}} e_{i,r}z^{-r}$ and $\langle i,j\rangle \in \mathbb{K}$ is defined via

$$\langle i,j\rangle := \nu_1^{\delta_{ij}-\delta_{i+1,j}}\nu_2^{\delta_{i-1,j}-\delta_{ij}}.$$

Let us follow our previous notations, except that now $(\lambda_1,\dots,\lambda_{p-1}) \in \{\nu_1,\nu_2\}^{p-1}$. Similarly to (1.19), we define the *PBWD basis elements* $e_\beta(r) \in U^>_{\nu_1,\nu_2}(L\mathfrak{sl}_n)$ via:

$$e_\beta(r) := [\cdots[[e_{i_1,r_1}, e_{i_2,r_2}]_{\lambda_1}, e_{i_3,r_3}]_{\lambda_2}, \cdots, e_{i_p,r_p}]_{\lambda_{p-1}}. \tag{1.70}$$

Following (1.21), the monomials

$$e_h := \prod^{\rightarrow}_{(\beta,r)\in\Delta^+\times\mathbb{Z}} e_\beta(r)^{h(\beta,r)} \qquad \text{with} \qquad h \in \mathsf{H} \tag{1.71}$$

will be called the *ordered PBWD monomials* of $U^>_{\nu_1,\nu_2}(L\mathfrak{sl}_n)$, where the arrow $\rightarrow$ over the product sign indicates that the product is ordered with respect to (1.18).

Our first main result establishes the PBWD property of $U^>_{\nu_1,\nu_2}(L\mathfrak{sl}_n)$ (the proof is outlined in Section 1.3.4 and is based on the shuffle approach):

**Theorem 1.16** $\{e_h\}_{h\in\mathsf{H}}$ *of* (1.71) *form a $\mathbb{K}$-basis of* $U^>_{\nu_1,\nu_2}(L\mathfrak{sl}_n)$.

*Remark 1.7* We note that the PBWD basis elements introduced in [13, (3.14)] are

$$e_{\alpha_j+\alpha_{j+1}+\cdots+\alpha_i}(r) := [\cdots[[e_{j,r}, e_{j+1,0}]_{\nu_1}, e_{j+2,0}]_{\nu_1}, \cdots, e_{i,0}]_{\nu_1}. \tag{1.72}$$

In this setup, Theorem 1.16 recovers Theorem 3.11 of [13] presented without a proof.

*Remark 1.8* We note that the entire two-parameter quantum loop algebra $U_{\nu_1,\nu_2}(L\mathfrak{sl}_n)$ admits a triangular decomposition similar to that of $U_\nu(L\mathfrak{sl}_n)$ from Proposition 1.1(a), cf. [13, §3.2]. Hence, a natural analogue of Theorem 1.2 holds for $U_{\nu_1,\nu_2}(L\mathfrak{sl}_n)$ as well, thus providing a family of PBWD $\mathbb{K}$-bases for $U_{\nu_1,\nu_2}(L\mathfrak{sl}_n)$.

### 1.3.2 Integral form $\mathfrak{U}^{>}_{v_1,v_2}(L\mathfrak{sl}_n)$ and its PBWD bases

Following (1.23), for any $(\beta, r) \in \Delta^+ \times \mathbb{Z}$, we define $\widetilde{e}_\beta(r) \in U^{>}_{v_1,v_2}(L\mathfrak{sl}_n)$ via:

$$\widetilde{e}_\beta(r) := \left( v_1^{1/2} v_2^{-1/2} - v_1^{-1/2} v_2^{1/2} \right) e_\beta(r). \tag{1.73}$$

We also define $\widetilde{e}_h \in U^{>}_{v_1,v_2}(L\mathfrak{sl}_n)$ via (1.71) but using $\widetilde{e}_\beta(r)$ of (1.73) instead of $e_\beta(r)$. Finally, let us define an integral form $\mathfrak{U}^{>}_{v_1,v_2}(L\mathfrak{sl}_n)$ as the $\mathbb{Z}[v_1^{1/2}, v_2^{1/2}, v_1^{-1/2}, v_2^{-1/2}]$-subalgebra of $U^{>}_{v_1,v_2}(L\mathfrak{sl}_n)$ generated by $\{\widetilde{e}_\beta(r)\}_{\beta \in \Delta^+}^{r \in \mathbb{Z}}$.

The following natural counterpart of Theorem 1.3 strengthens our Theorem 1.16:

**Theorem 1.17** *(a) The subalgebra $\mathfrak{U}^{>}_{v_1,v_2}(L\mathfrak{sl}_n)$ is independent of all our choices.*

*(b) $\{\widetilde{e}_h\}_{h \in \mathsf{H}}$ form a basis of the free $\mathbb{C}[v_1^{1/2}, v_2^{1/2}, v_1^{-1/2}, v_2^{-1/2}]$-module $\mathfrak{U}^{>}_{v_1,v_2}(L\mathfrak{sl}_n)$.*

The proof of Theorem 1.17 follows easily from that of Theorem 1.16 sketched below in the same way as we deduced the proof of Theorem 1.3 in Section 1.2.3.2 from that of Theorem 1.1; we leave details to the interested reader.

*Remark 1.9* Similarly to Remark 1.3, let us point out that it is often more convenient to work with the two-parameter quantum loop algebra $U_{v_1,v_2}(L\mathfrak{gl}_n)$ and its similarly defined integral form $\mathfrak{U}_{v_1,v_2}(L\mathfrak{gl}_n)$. Following the arguments of [10, Proposition 3.11], $\mathfrak{U}_{v_1,v_2}(L\mathfrak{gl}_n)$ is identified with the RTT integral form $\mathfrak{U}^{\mathrm{rtt}}_{v_1,v_2}(L\mathfrak{gl}_n)$ under the $\mathbb{K}$-algebra isomorphism $U_{v_1,v_2}(L\mathfrak{gl}_n) \simeq \mathfrak{U}^{\mathrm{rtt}}_{v_1,v_2}(L\mathfrak{gl}_n) \otimes_{\mathbb{Z}[v_1^{1/2}, v_2^{1/2}, v_1^{-1/2}, v_2^{-1/2}]} \mathbb{K}$ of [14]. Therefore, the analogue of [10, Theorem 3.24] (cf. our Theorem 1.4) provides a family of Poincaré-Birkhoff-Witt-Drinfeld bases for $\mathfrak{U}_{v_1,v_2}(L\mathfrak{gl}_n)$, cf. Theorem 1.5.

### 1.3.3 Shuffle algebra $\widetilde{S}^{(n)}$

We define the **shuffle algebra** $(\widetilde{S}^{(n)}, \star)$ similarly to $(S^{(n)}, \star)$ with three modifications:

(1) all vector spaces are now defined over $\mathbb{K}$ (rather than over $\mathbb{C}(v)$);

(2) the shuffle factors (1.31) used in the shuffle product (1.32) are now replaced with:

$$\zeta_{i,j}(z) = \left( \frac{z - v_1^{1/2} v_2^{-1/2}}{z-1} \right)^{\delta_{j,i-1}} \left( \frac{z - v_1^{-1} v_2}{z-1} \right)^{\delta_{ji}} \left( v_1^{1/2} v_2^{1/2} \cdot \frac{z - v_1^{1/2} v_2^{-1/2}}{z-1} \right)^{\delta_{j,i+1}};$$

(3) the wheel conditions (1.35) for $F$ are replaced with

$$F\left(\{x_{i,r}\}\right) = 0 \quad \text{once} \quad x_{i,r_1} = v_1^{1/2} v_2^{-1/2} x_{i+\epsilon,s} = v_1 v_2^{-1} x_{i,r_2}$$

for some $\epsilon \in \{\pm 1\}$, $i$, $r_1 \neq r_2$, $s$.

The following result relates $(\widetilde{S}^{(n)}, \star)$ to $U^>_{\nu_1, \nu_2}(L\mathfrak{sl}_n)$, generalizing Proposition 1.3:

**Proposition 1.10** *There exists an injective $\mathbb{K}$-algebra homomorphism*

$$\Psi \colon U^>_{\nu_1, \nu_2}(L\mathfrak{sl}_n) \longrightarrow \widetilde{S}^{(n)} \tag{1.74}$$

*such that $e_{i,r} \mapsto x^r_{i,1}$ for any $i \in I, r \in \mathbb{Z}$.*

Our proof of Theorem 1.16 also implies the following counterpart of Theorem 1.8:

**Theorem 1.18** $\Psi \colon U^>_{\nu_1, \nu_2}(L\mathfrak{sl}_n) \overset{\sim}{\longrightarrow} \widetilde{S}^{(n)}$ *of (1.74) is a $\mathbb{K}$-algebra isomorphism.*

### 1.3.4  Proofs of Theorem 1.16 and Theorem 1.18

The proofs of Theorems 1.16, 1.18 are completely analogous to our proofs of Theorems 1.1, 1.8, and utilize the embedding $\Psi \colon U^>_{\nu_1, \nu_2}(L\mathfrak{sl}_n) \hookrightarrow \widetilde{S}^{(n)}$ of Proposition 1.10. To this end, we note that:

- The linear independence of $\{e_h\}_{h \in \mathsf{H}}$ is deduced exactly as in Proposition 1.5 with the only modification of the *specialization maps* $\phi_{\underline{d}} \colon \widetilde{S}^{(n)}_{\underline{\ell}} \to \mathbb{K}[\{y^{\pm 1}_{\beta,s}\}^{1 \le s \le d_\beta}_{\beta \in \Delta^+}]^{\Sigma_{\underline{d}}}$ of (1.40) by replacing $\nu^{-i}$ with $(\nu_1^{1/2}\nu_2^{-1/2})^{-i}$ everywhere. Then, the results of Lemmas 1.5 and 1.6 still hold, thus implying the linear independence of $\{e_h\}_{h \in \mathsf{H}}$.

- The fact that $\{e_h\}_{h \in \mathsf{H}}$ span $U^>_{\nu_1, \nu_2}(L\mathfrak{sl}_n)$ is deduced exactly as in Proposition 1.7. To be more precise, Lemma 1.7 still holds and its iterative application immediately implies that any shuffle element $F \in \widetilde{S}^{(n)}$ belongs to the $\mathbb{K}$-span of $\{\Psi(e_h)\}_{h \in \mathsf{H}}$.

- The latter observation also implies the surjectivity of the algebra embedding $\Psi$ from Proposition 1.10, thus proving Theorem 1.18.

## 1.4  Generalizations to quantum loop superalgebras

The quantum loop superalgebra $U_\nu(L\mathfrak{sl}(m|n))$[5] was introduced in [20], both in the Drinfeld-Jimbo and new Drinfeld realizations, see [20, Theorem 8.5.1] for an identification of those. The representation theory of these algebras was partially studied in [21] by crucially utilizing a weak version of the PBW theorem for $U^>_\nu(L\mathfrak{sl}(m|n))$, see [21, Theorem 3.12]. Inspired by [13, Theorem 3.11] discussed above, the author also conjectured the PBW theorem for $U^>_\nu(L\mathfrak{sl}(m|n))$, see [21, Remark 3.13(2)].

The primary goal of this section is to generalize Theorem 1.1 to the case of $U^>_\nu(L\mathfrak{sl}(m|n))$, thus proving the aforementioned conjecture of [21]. Along the way, we also generalize Theorem 1.8 by providing the shuffle realization of $U_v^>(L\mathfrak{sl}(m|n))$, which is of independent interest. The latter is also used to construct the PBWD bases for the integral form $\mathfrak{U}_v^>(L\mathfrak{sl}(m|n))$ of $U_v^>(L\mathfrak{sl}(m|n))$, thus generalizing Theorem 1.3.

---

[5] Here, one actually needs to use the classical Lie superalgebra $A(m-1, n-1)$ in place of $\mathfrak{sl}(m|n)$, which do coincide only for $m \ne n$. However, we shall ignore this difference, since we will consider only the "positive" subalgebras and those are isomorphic: $U^>_\nu(L\mathfrak{sl}(m|n)) \simeq U^>_\nu(LA(m-1, n-1))$.

*Remark 1.10* (a) Let us emphasize right away that in our exposition below we consider the distinguished Dynkin diagram with a single simple positive odd root. The generalization of all our results to an arbitrary Dynkin diagram is carried out in [17].

(b) The new Drinfeld realization of quantum loop superalgebras surprisingly has not been properly developed yet. Thus, the shuffle approach, involving both symmetric and skew-symmetric functions, provides a promising new toolkit for such a treatment.

## 1.4.1 Quantum loop superalgebra $U_v^>(L\mathfrak{sl}(m|n))$

To simplify our exposition, we shall only consider the "positive" subalgebra $U_v^>(L\mathfrak{sl}(m|n))$ of $U_v(L\mathfrak{sl}(m|n))$. From now on, we redefine $I = \{1, \ldots, m+n-1\}$. Let us consider a free $\mathbb{Z}$-module $\bigoplus_{i=1}^{m+n} \mathbb{Z}\epsilon_i$ with the bilinear form $(\cdot, \cdot)$ determined by $(\epsilon_i, \epsilon_j) = (-1)^{\delta_{i>m}}\delta_{ij}$. Let $v$ be a formal variable and define $\{v_i\}_{i \in I} \subset \{v, v^{-1}\}$ via $v_i := v^{(\epsilon_i, \epsilon_i)}$. For $i, j \in I$, we set $\bar{c}_{ij} := (\alpha_i, \alpha_j)$ with $\alpha_i := \epsilon_i - \epsilon_{i+1}$.

Following [20] (cf. [21, Theorem 3.3]), we define $U_v^>(L\mathfrak{sl}(m|n))$ as the associative $\mathbb{C}(v)$-superalgebra generated by $\{e_{i,r}\}_{i \in I}^{r \in \mathbb{Z}}$, with the $\mathbb{Z}_2$-grading defined by:

$$|e_{m,r}| = \bar{1}, \quad |e_{i,r}| = \bar{0} \quad \text{for all} \quad i \neq m, \, r \in \mathbb{Z}, \tag{1.75}$$

and subject to the following defining relations:

$$(z - v^{\bar{c}_{ij}}w)e_i(z)e_j(w) = (v^{\bar{c}_{ij}}z - w)e_j(w)e_i(z) \quad \text{if} \quad \bar{c}_{ij} \neq 0, \tag{1.76}$$

quadratic and cubic quantum Serre relations:

$$[e_i(z), e_j(w)] = 0 \quad \text{if} \quad \bar{c}_{ij} = 0,$$
$$[e_i(z_1), [e_i(z_2), e_j(w)]_{v^{-1}}]_v + [e_i(z_2), [e_i(z_1), e_j(w)]_{v^{-1}}]_v = 0 \text{ if } \bar{c}_{ij} = \pm 1 \text{ and } i \neq m, \tag{1.77}$$

as well as quartic quantum Serre relations:

$$[[[e_{m-1}(w), e_m(z_1)]_{v^{-1}}, e_{m+1}(u)]_v, e_m(z_2)] +$$
$$[[[e_{m-1}(w), e_m(z_2)]_{v^{-1}}, e_{m+1}(u)]_v, e_m(z_1)] = 0. \tag{1.78}$$

Here, $e_i(z) = \sum_{r \in \mathbb{Z}} e_{i,r} z^{-r}$ and we use the super-bracket notations:

$$[a, b]_x := ab - (-1)^{|a||b|} x \cdot ba,$$
$$[a, b] := [a, b]_1 \tag{1.79}$$

for $\mathbb{Z}_2$-homogeneous elements $a$ and $b$ (with conventions $(-1)^{\bar{0}} = 1$, $(-1)^{\bar{1}} = -1$).

### 1.4.2 PBWD bases of $U_v^>(L\mathfrak{sl}(m|n))$

Let $\Delta^+ = \{\alpha_j + \alpha_{j+1} + \cdots + \alpha_i\}_{1 \le j \le i < m+n}$. For $\beta \in \Delta^+$, define its parity $|\beta| \in \mathbb{Z}_2$ by:

$$|\beta| = \begin{cases} \bar{1} & \text{if } m \in [\beta] \\ \bar{0} & \text{if } m \notin [\beta] \end{cases} . \tag{1.80}$$

Then, following the notations of Section 1.1.1, we define the *PBWD basis elements* $e_\beta(r) \in U_v^>(L\mathfrak{sl}(m|n))$ via (1.19), but with $[\cdot, \cdot]_{\lambda_k}$ denoting the super-bracket (1.79).

Let $\bar{\mathsf{H}}$ be the set of all functions $h : \Delta^+ \times \mathbb{Z} \to \mathbb{N}$ with finite support and such that

$$h(\beta, r) \le 1 \qquad \text{if} \qquad |\beta| = \bar{1} .$$

Following (1.21), the monomials

$$e_h := \overset{\rightarrow}{\prod_{(\beta,r) \in \Delta^+ \times \mathbb{Z}}} e_\beta(r)^{h(\beta,r)} \qquad \text{with} \qquad h \in \bar{\mathsf{H}} \tag{1.81}$$

will be called the *ordered PBWD monomials* of $U_v^>(L\mathfrak{sl}(m|n))$, where the arrow $\rightarrow$ over the product sign indicates that the product is ordered with respect to (1.18).

Our first main result establishes the PBWD property of $U_v^>(L\mathfrak{sl}(m|n))$ (the proof is outlined in Section 1.4.5 and is based on the shuffle approach):

**Theorem 1.19** $\{e_h\}_{h \in \bar{\mathsf{H}}}$ *of* (1.81) *form a* $\mathbb{C}(v)$-*basis of* $U_v^>(L\mathfrak{sl}(m|n))$.

*Remark 1.11* We note that the PBWD basis elements introduced in [21, (3.12)] are

$$e_{\alpha_j + \alpha_{j+1} + \cdots + \alpha_i}(r) := [\cdots [[e_{j,r}, e_{j+1,0}]_{v_{j+1}}, e_{j+2,0}]_{v_{j+2}}, \cdots, e_{i,0}]_{v_i} . \tag{1.82}$$

In this setup, Theorem 1.19 recovers the conjecture from Remark 3.13(2) of [21].

*Remark 1.12* We note that the entire quantum loop superalgebra $U_v(L\mathfrak{sl}(m|n))$ admits a triangular decomposition similar to that of $U_v(L\mathfrak{sl}_n)$ from Proposition 1.1(a), see [21, Theorem 3.3]. Hence, a natural analogue of Theorem 1.2 holds for $U_v(L\mathfrak{sl}(m|n))$ as well, thus providing a family of PBWD $\mathbb{C}(v)$-bases for $U_v(L\mathfrak{sl}(m|n))$.

### 1.4.3 Integral form $\mathfrak{U}_v^>(L\mathfrak{sl}(m|n))$ and its PBWD bases

Following (1.23), for any $(\beta, r) \in \Delta^+ \times \mathbb{Z}$, we define $\widetilde{e}_\beta(r) \in U_v^>(L\mathfrak{sl}(m|n))$ via:

$$\widetilde{e}_\beta(r) := (v - v^{-1}) e_\beta(r) . \tag{1.83}$$

For $h \in \bar{\mathsf{H}}$, we also define $\widetilde{e}_h \in U_v^>(L\mathfrak{sl}(m|n))$ via (1.81) but using $\widetilde{e}_\beta(r)$ of (1.83) instead of $e_\beta(r)$. Finally, let us define an integral form $\mathfrak{U}_v^>(L\mathfrak{sl}(m|n))$ as the $\mathbb{Z}[v, v^{-1}]$-subalgebra of $U_v^>(L\mathfrak{sl}(m|n))$ generated by $\{\widetilde{e}_\beta(r)\}_{\beta \in \Delta^+}^{r \in \mathbb{Z}}$.

The following natural counterpart of Theorem 1.3 strengthens our Theorem 1.19:

**Theorem 1.20** *(a) The subalgebra* $\mathfrak{U}_{\boldsymbol{v}}^{>}(L\mathfrak{sl}(m|n))$ *is independent of all our choices.*

*(b) The elements* $\{\widetilde{e}_h\}_{h\in\bar{\mathsf{H}}}$ *form a basis of the free* $\mathbb{Z}[\boldsymbol{v},\boldsymbol{v}^{-1}]$-*module* $\mathfrak{U}_{\boldsymbol{v}}^{>}(L\mathfrak{sl}(m|n))$.

The proof of Theorem 1.20 follows easily from that of Theorem 1.19 outlined below in the same way as we deduced the proof of Theorem 1.3 in Section 1.2.3.2 from that of Theorem 1.1; we leave details to the interested reader.

*Remark 1.13* Similarly to Remarks 1.3 and 1.9, let us note that it is often more convenient to work with the quantum loop superalgebra $U_{\boldsymbol{v}}(L\mathfrak{gl}(m|n))$ and its similarly defined integral form $\mathfrak{U}_{\boldsymbol{v}}(L\mathfrak{gl}(m|n))$. Following the arguments of [10, Proposition 3.11], $\mathfrak{U}_{\boldsymbol{v}}(L\mathfrak{gl}(m|n))$ is identified with the RTT integral form $\mathfrak{U}_{\boldsymbol{v}}^{\mathrm{rtt}}(L\mathfrak{gl}(m|n))$, under the $\mathbb{C}(\boldsymbol{v})$-superalgebra isomorphism $U_{\boldsymbol{v}}(L\mathfrak{gl}(m|n)) \simeq \mathfrak{U}_{\boldsymbol{v}}^{\mathrm{rtt}}(L\mathfrak{gl}(m|n)) \otimes_{\mathbb{Z}[\boldsymbol{v},\boldsymbol{v}^{-1}]} \mathbb{C}(\boldsymbol{v})$, see [22, Definition 3.1]. Hence, the analogue of [10, Theorem 3.24] (cf. our Theorem 1.4) provides a family of PBWD bases for $\mathfrak{U}_{\boldsymbol{v}}(L\mathfrak{gl}(m|n))$, cf. Theorem 1.5.

## 1.4.4 Shuffle algebra $S^{(m|n)}$

Consider an $\mathbb{N}^I$-graded $\mathbb{C}(\boldsymbol{v})$-vector superspace

$$\mathbb{S}^{(m|n)} = \bigoplus_{\underline{k}=(k_1,\ldots,k_{m+n-1})\in\mathbb{N}^I} \mathbb{S}_{\underline{k}}^{(m|n)},$$

$\mathbb{Z}_2$-graded by the parity of $k_m$, where $\mathbb{S}_{\underline{k}}^{(m|n)}$ consists of **supersymmetric** rational functions in $\{x_{i,r}\}_{i\in I}^{1\le r\le k_i}$, that is, symmetric in $\{x_{i,r}\}_{r=1}^{k_i}$ for every $i\ne m$ and skew-symmetric in $\{x_{m,r}\}_{r=1}^{k_m}$. We also fix an $I\times I$ matrix of rational functions $\zeta_{i,j}(z)$ via:

$$\zeta_{i,j}(z) = \frac{z - \boldsymbol{v}^{-\bar{c}_{ij}}}{z - 1}. \tag{1.84}$$

We equip $\mathbb{S}^{(m|n)}$ with a structure of an associative algebra with the shuffle product defined similar to (1.32), but Sym of (1.33) being replaced with a supersymmetrization SSym. Here, the *supersymmetrization* of $f\in\mathbb{C}(\{x_{i,1},\ldots,x_{i,s_i}\}_{i\in I})$ is defined via:

$$\mathrm{SSym}(f)\left(\{x_{i,1},\ldots,x_{i,s_i}\}_{i\in I}\right) :=$$
$$\sum_{(\sigma_1,\ldots,\sigma_{m+n-1})\in\Sigma_{\underline{s}}} \mathrm{sgn}(\sigma_m) f\left(\{x_{i,\sigma_i(1)},\ldots,x_{i,\sigma_i(s_i)}\}_{i\in I}\right). \tag{1.85}$$

As before, let us consider the subspace of $\mathbb{S}^{(m|n)}$ defined by the *pole* and *wheel* conditions (though now there will be two kinds of the wheel conditions):

• We say that $F\in\mathbb{S}_{\underline{k}}^{(m|n)}$ satisfies the *pole conditions* if

$$F = \frac{f(x_{1,1}, \ldots, x_{m+n-1,k_{m+n-1}})}{\prod_{i=1}^{m+n-2} \prod_{r \leq k_i}^{r' \leq k_{i+1}} (x_{i,r} - x_{i+1,r'})}, \tag{1.86}$$

where $f \in \mathbb{C}(\nu)[\{x_{i,r}^{\pm 1}\}_{i \in I}^{1 \leq r \leq k_i}]$ is a supersymmetric Laurent polynomial, that is, symmetric in $\{x_{i,r}\}_{r=1}^{k_i}$ for every $i \neq m$ and skew-symmetric in $\{x_{m,r}\}_{r=1}^{k_m}$.

- We say that $F \in \mathbb{S}_{\underline{k}}^{(m|n)}$ satisfies the *first kind wheel conditions* if

$$F\left(\{x_{i,r}\}\right) = 0 \quad \text{once} \quad x_{i,r_1} = \nu_i x_{i+\epsilon,s} = \nu_i^2 x_{i,r_2} \tag{1.87}$$

  for some $\epsilon \in \{\pm 1\}$, $i \in I \setminus \{m\}$, $1 \leq r_1 \neq r_2 \leq k_i$, $1 \leq s \leq k_{i+\epsilon}$.

- We say that $F \in \mathbb{S}_{\underline{k}}^{(m|n)}$ satisfies the *second kind wheel conditions* if

$$F\left(\{x_{i,r}\}\right) = 0 \quad \text{once} \quad x_{m-1,s} = \nu x_{m,r_1} = x_{m+1,s'} = \nu^{-1} x_{m,r_2} \tag{1.88}$$

  for some $1 \leq r_1 \neq r_2 \leq k_m$, $1 \leq s \leq k_{m-1}$, $1 \leq s' \leq k_{m+1}$.

Let $S_{\underline{k}}^{(m|n)} \subset \mathbb{S}_{\underline{k}}^{(m|n)}$ denote the subspace of all $F$ satisfying these conditions and set

$$S^{(m|n)} := \bigoplus_{\underline{k} \in \mathbb{N}^I} S_{\underline{k}}^{(m|n)}.$$

It is straightforward to check that $S^{(m|n)} \subset \mathbb{S}^{(m|n)}$ is $\star$-closed, cf. Lemma 1.2. Similar to Proposition 1.3, the **shuffle algebra** $\left(S^{(m|n)}, \star\right)$ is related to $U_\nu^>(L\mathfrak{sl}(m|n))$ via:

**Proposition 1.11** *There exists an injective $\mathbb{C}(\nu)$-superalgebra homomorphism*

$$\Psi \colon U_\nu^>(L\mathfrak{sl}(m|n)) \longrightarrow S^{(m|n)} \tag{1.89}$$

*such that $e_{i,r} \mapsto x_{i,1}^r$ for any $i \in I, r \in \mathbb{Z}$.*

In fact, we have the following generalization of Theorem 1.8:

**Theorem 1.21** *$\Psi \colon U_\nu^>(L\mathfrak{sl}(m|n)) \overset{\sim}{\longrightarrow} S^{(m|n)}$ of (1.89) is a $\mathbb{C}(\nu)$-superalgebra isomorphism.*

### 1.4.5 Proofs of Theorem 1.19 and Theorem 1.21

The proofs of Theorems 1.19, 1.21 are completely analogous to those of Theorems 1.1, 1.8, and utilize the embedding $\Psi \colon U_\nu^>(L\mathfrak{sl}(m|n)) \hookrightarrow S^{(m|n)}$ of Proposition 1.11. Let us first establish both theorems in the simplest "rank 1" $m = n = 1$ case:

**Proposition 1.12** *Fix any total order $\leq$ on $\mathbb{Z}$. Then:*

*(a) the ordered monomials $\{x^{r_1} \star x^{r_2} \star \cdots \star x^{r_k}\}_{\underline{k} \in \mathbb{N}}^{r_1 < \cdots < r_k}$ form a $\mathbb{C}(\nu)$-basis of $S^{(1|1)}$.*

*(b) the ordered monomials $\{e_{r_1} e_{r_2} \cdots e_{r_k}\}_{\underline{k} \in \mathbb{N}}^{r_1 < \cdots < r_k}$ form a $\mathbb{C}(\nu)$-basis of $U_\nu^>(L\mathfrak{sl}(1|1))$.*

*Proof* To prove (a), it suffices to note the superalgebra isomorphism $S^{(1|1)} \simeq \oplus_k \Lambda_k$, where $\Lambda_k$ denotes the vector space of skew-symmetric Laurent polynomials in $k$ variables, while the algebra structure on $\oplus_k \Lambda_k$ arises via the skew-symmetrization maps $\Lambda_k \otimes \Lambda_\ell \to \Lambda_{k+\ell}$. Combining part (a) with Proposition 1.11 implies part (b). □

Let us now treat the general $m, n \in \mathbb{N}$. For a degree vector $\underline{d} = \{d_\beta\}_{\beta \in \Delta^+} \in \mathbb{N}^{\Delta^+}$, define its $\mathbb{N}^I$-degree $\underline{\ell} \in \mathbb{N}^I$ via $\sum_{\beta \in \Delta^+} d_\beta \beta = \sum_{i \in I} \ell_i \alpha_i$. The **specialization map** $\phi_{\underline{d}}$, vanishing on degree $\neq \underline{\ell}$ components, and thus determined by

$$\phi_{\underline{d}} \colon S^{(m|n)}_{\underline{\ell}} \longrightarrow \mathbb{C}(v) \left[ \{y^{\pm 1}_{\beta, s}\}^{1 \le s \le d_\beta}_{\beta \in \Delta^+} \right] \tag{1.90}$$

is defined alike (1.40), with the only modification that the variable $x_{i,r}$ in the $s$-th copy of $[\beta]$ is specialized to $v^{-i} y_{\beta, s}$ if $i \le m$ and to $v^{i-2m} y_{\beta, s}$ if $i > m$. Let us note right away that $\phi_{\underline{d}}(F)$ is a supersymmetric Laurent polynomial in $\{y_{\beta, s}\}^{d_\beta}_{s=1}$ for any $\beta \in \Delta^+$, that is, symmetric in $\{y_{\beta, s}\}^{d_\beta}_{s=1}$ if $|\beta| = \bar{0}$ and skew-symmetric if $|\beta| = \bar{1}$.

Then, the results of Lemmas 1.5 and 1.6 still hold, with $\bar{\mathsf{H}}$ being used instead of $\mathsf{H}$, thus implying the linear independence of $\{e_h\}_{h \in \bar{\mathsf{H}}}$. Furthermore, for any $h \in \bar{\mathsf{H}}$ with $\deg(h) = \underline{d}$, the formula (1.50) still holds but with the factors generalizing (1.49):

$$G_{\beta, \beta'} = \prod^{1 \le s' \le d_{\beta'}}_{1 \le s \le d_\beta} (y_{\beta, s} - v^{-2} y_{\beta', s'})^{\nu^-(\beta, \beta')} \cdot (y_{\beta, s} - v^2 y_{\beta', s'})^{\nu^+(\beta, \beta')} \times$$

$$\prod^{1 \le s' \le d_{\beta'}}_{1 \le s \le d_\beta} (y_{\beta, s} - y_{\beta', s'})^{\delta_{j(\beta') > j(\beta)} \delta_{i(\beta)+1 \in [\beta']} + \delta_{m \in [\beta]} \delta_{m \in [\beta']}},$$

$$G_\beta = (1 - v^2)^{d_\beta (i(\beta) - j(\beta))} \cdot \prod_{1 \le s \ne s' \le d_\beta} (y_{\beta, s} - v^2 y_{\beta, s'})^{i(\beta) - j(\beta)} \cdot \prod_{1 \le s \le d_\beta} y^{i(\beta) - j(\beta)}_{\beta, s},$$

$$G^{(\sigma_\beta)}_\beta = \prod^{d_\beta}_{s=1} y^{r_\beta(h,s)}_{\beta, \sigma_\beta(s)} \cdot \begin{cases} \prod_{s < s'} \frac{y_{\beta, \sigma_\beta(s)} - v^{-2} y_{\beta, \sigma_\beta(s')}}{y_{\beta, \sigma_\beta(s)} - y_{\beta, \sigma_\beta(s')}} & \text{if } m > i(\beta) \\ \prod_{s < s'} \frac{y_{\beta, \sigma_\beta(s)} - v^2 y_{\beta, \sigma_\beta(s')}}{y_{\beta, \sigma_\beta(s)} - y_{\beta, \sigma_\beta(s')}} & \text{if } m < j(\beta) \\ \operatorname{sgn}(\sigma_\beta) & \text{if } m \in [\beta] \end{cases} \cdot$$

Here, the exponents $\nu^-(\beta, \beta')$ and $\nu^+(\beta, \beta')$ are given respectively by:

$$\#\{(j, j') \in [\beta] \times [\beta'] \mid j = j' < m\} + \#\{(j, j') \in [\beta] \times [\beta'] \mid j = j' + 1 > m\},$$

$$\#\{(j, j') \in [\beta] \times [\beta'] \mid j = j' > m\} + \#\{(j, j') \in [\beta] \times [\beta'] \mid j = j' + 1 \le m\}.$$

The key observation is that for any $\beta \in \Delta^+$, the sum $\sum_{\sigma_\beta \in \Sigma_{d_\beta}} G^{(\sigma_\beta)}_\beta$ coincides with the value of the shuffle product $x^{r_\beta(h,1)} \star \cdots \star x^{r_\beta(h, d_\beta)}$, viewed as an element of $S^{(2|0)}$ if $m > i(\beta)$, $S^{(0|2)}$ if $m < j(\beta)$, or $S^{(1|1)}$ if $m \in [\beta]$ (a "rank 1 reduction" feature), evaluated at $\{y_{\beta, s}\}^{d_\beta}_{s=1}$. The latter elements are linearly independent, due to Propositions 1.4(a) and 1.12(a), hence so are $\{e_h\}_{h \in \bar{\mathsf{H}}}$ of (1.81), cf. Proposition 1.5.

The fact that the ordered PBWD monomials $\{e_h\}_{h \in \check{H}}$ span $U_v^>(L\mathfrak{sl}(m|n))$ follows from the validity of Lemma 1.7 in the present setup. The latter result also implies the surjectivity of the homomorphism $\Psi$ from (1.89), thus proving Theorem 1.21.

# References

1. J. Beck: Convex bases of PBW type for quantum affine algebras. Commun. Math. Phys. **165**, no. 1, 193–199 (1994)

2. V. Chari, A. Pressley: Quantum affine algebras at roots of unity. Represent. Theory **1**, 280–328 (1997)

3. V. Drinfeld: A New realization of Yangians and quantized affine algebras. Sov. Math. Dokl. **36**, no. 2, 212–216 (1988)

4. J. Ding, I. Frenkel: Isomorphism of two realizations of quantum affine algebra $U_q(\widehat{\mathfrak{gl}(n)})$. Commun. Math. Phys. **156**, no. 2, 277–300 (1993)

5. B. Enriquez: PBW and duality theorems for quantum groups and quantum current algebras. J. Lie Theory **13**, no. 1, 21–64 (2003)

6. B. Feigin, A. Odesskii: Sklyanin's elliptic algebras. (Russian) Funktsional. Anal. i Prilozhen. **23**, no. 3, 45–54 (1989); translation in Funct. Anal. Appl. **23**, no. 3, 207–214 (1990)

7. B. Feigin, A. Odesskii: Elliptic deformations of current algebras and their representations by difference operators. (Russian) Funktsional. Anal. i Prilozhen. **31**, no. 3, 57–70 (1997); translation in Funct. Anal. Appl. **31**, no. 3, 193–203 (1998)

8. B. Feigin, A. Odesskii: Quantized moduli spaces of the bundles on the elliptic curve and their applications. Integrable structures of exactly solvable two-dimensional models of quantum field theory (Kiev, 2000), 123–137; NATO Sci. Ser. II Math. Phys. Chem., 35, Kluwer Acad. Publ., Dordrecht (2001)

9. L. Faddeev, N. Reshetikhin, L. Takhtadzhyan: Quantization of Lie groups and Lie algebras. (Russian) Algebra i Analiz **1**, no. 1, 178–206 (1989); translation in Leningrad Math. J. **1**, no. 1, 193–225 (1990)

10. M. Finkelberg, A. Tsymbaliuk: Shifted quantum affine algebras: integral forms in type $A$ (with appendices by A. Tsymbaliuk, A. Weekes). Arnold Math. J. **5**, no. 2-3, 197–283 (2019)

11. I. Grojnowski: Affinizing quantum algebras: from $D$-modules to $K$-theory, unpublished manuscript available at https://www.dpmms.cam.ac.uk/~groj/char.ps (1994)

12. D. Hernandez: Representations of quantum affinizations and fusion product. Transform. Groups **10**, no. 2, 163–200 (2005)

13. N. Hu, M. Rosso, H. Zhang: Two-parameter quantum affine algebra $U_{r,s}(\widehat{\mathfrak{sl}}_n)$, Drinfel'd realization and quantum affine Lyndon basis. Commun. Math. Phys. **278**, no. 2, 453–486 (2008)

14. N. Jing, M. Liu: $R$-matrix realization of two-parameter quantum affine algebra $U_{r,s}(\widehat{\mathfrak{gl}}_n)$. J. Algebra **488**, 1–28 (2017)

15. G. Lusztig: Quantum groups at roots of 1. Geom. Dedicata **35**, no. 1-3, 89–113 (1990)

16. A. Neguţ: Quantum toroidal and shuffle algebras. Adv. Math. **372**, Paper No. 107288 (2020)

17. A. Tsymbaliuk: Shuffle algebra realizations of type $A$ super Yangians and quantum affine superalgebras for all Cartan data. Lett. Math. Physics **110**, no. 8, 2083–2111 (2020)

18. A. Tsymbaliuk: Duality of Lusztig and RTT integral forms of $U_v(L\mathfrak{sl}_n)$. J. Pure Appl. Algebra **225**, no. 1, Paper No. 106469 (2021)

19. A. Tsymbaliuk: PBWD bases and shuffle algebra realizations for $U_v(L\mathfrak{sl}_n)$, $U_{v_1,v_2}(L\mathfrak{sl}_n)$, $U_v(L\mathfrak{sl}(m|n))$ and their integral forms. Selecta Math. (N.S.) **27**, no. 3, Paper No. 35 (2021)

20. H. Yamane: On defining relations of affine Lie superalgebras and affine quantized universal enveloping superalgebras. Publ. Res. Inst. Math. Sci. **35**, no. 3, 321–390 (1999)

21. H. Zhang: Representations of quantum affine superalgebras. Math. Z. **278**, no. 3-4, 663–703 (2014)

22. H. Zhang: RTT realization of quantum affine superalgebras and tensor products. Int. Math. Res. Not. IMRN, no. 4, 1126–1157 (2016)

# Chapter 2
# Quantum toroidal $\mathfrak{gl}_1$, its representations, and geometric realization

**Abstract** In this chapter, we recall the basic results on the quantum toroidal algebra of $\mathfrak{gl}_1$. This algebra has an important interpretation [18, 19] as an extended elliptic Hall algebra of [1], which provides its "90 degree rotation" automorphism $\varpi$ (first discovered in [13]). We also establish the shuffle realization of its "positive" subalgebra and its particular commutative subalgebra, due to [17] and [6], respectively. Following [4, 5, 7], we discuss a combinatorial approach to the construction of several important families of representations. We recall the shuffle realization [8, 22] of these combinatorial modules, as well as explain how vertex-type representations of [6] can be naturally related to them via a $\varpi$-twist. Finally, following [9, 19, 21], we provide some flavor of the applications to the geometry by realizing Fock modules and their tensor products via equivariant $K$-theory of the Gieseker moduli spaces, as well as evoking the $K$-theoretic version of the Nakajima's construction from [15].

## 2.1 Quantum toroidal algebras of $\mathfrak{gl}_1$

In this section, we recall the definition of the quantum toroidal algebra $\ddot{U}_{q_1,q_2,q_3}(\mathfrak{gl}_1)$, first studied in [13] and rediscovered later in [9, 4]. We provide the shuffle algebra realization of $\ddot{U}^>_{q_1,q_2,q_3}(\mathfrak{gl}_1)$ established in [17]. Following [18, 19], we evoke the identification of $\ddot{U}_{q_1,q_2,q_3}(\mathfrak{gl}_1)$ with a (centrally extended) elliptic Hall algebra of [1], and the resulting "90 degree rotation" automorphism of [13] interchanging the horizontal and vertical Heisenberg subalgebras. We also provide the shuffle realization [6] of the commutative "positive half" $\mathcal{A}$ of the horizontal Heisenberg.

### 2.1.1 Quantum toroidal $\mathfrak{gl}_1$

We fix $q_1, q_2, q_3 \in \mathbb{C}^\times$ satisfying the relation $q_1 q_2 q_3 = 1$. Choosing a square root $q$ of $q_2$, we can equivalently use $q, d \in \mathbb{C}^\times$ related to the above via:

$$q_1 := q^{-1}d, \quad q_2 := q^2, \quad q_3 := q^{-1}d^{-1}. \tag{2.1}$$

Throughout this section, we assume that none of $q_1, q_2, q_3$ is a root of 1. We define

$$\beta_r := (1 - q_1^r)(1 - q_2^r)(1 - q_3^r) \in \mathbb{C}^\times \quad \text{for all} \quad r > 0, \tag{2.2}$$

and consider the rational function

$$g(z) := \frac{(1 - q_1 z)(1 - q_2 z)(1 - q_3 z)}{(1 - q_1^{-1} z)(1 - q_2^{-1} z)(1 - q_3^{-1} z)} \tag{2.3}$$

satisfying the equality $g(1/z) = 1/g(z)$. Following [4, 7], we define the quantum toroidal algebra of $\mathfrak{gl}_1$, denoted by $\ddot{U}_{q_1,q_2,q_3}(\mathfrak{gl}_1)$, to be the associative $\mathbb{C}$-algebra generated by $\{e_r, f_r, \psi_r, \psi_0^{-1}, \gamma^{\pm 1/2}, q^{\pm d_1}, q^{\mp d_2}\}_{r \in \mathbb{Z}}$ with the following defining relations:

$$[\psi^\pm(z), \psi^\pm(w)] = 0, \quad \gamma^{\pm 1/2} - \text{central}, \tag{t0.1}$$

$$\psi_0^{\pm 1} \cdot \psi_0^{\mp 1} = \gamma^{\pm 1/2} \cdot \gamma^{\mp 1/2} = q^{\pm d_1} \cdot q^{\mp d_1} = q^{\pm d_2} \cdot q^{\mp d_2} = 1, \tag{t0.2}$$

$$q^{d_1} e(z) q^{-d_1} = e(qz), \quad q^{d_1} f(z) q^{-d_1} = f(qz), \quad q^{d_1} \psi^\pm(z) q^{-d_1} = \psi^\pm(qz), \tag{t0.3}$$

$$q^{d_2} e(z) q^{-d_2} = qe(z), \quad q^{d_2} f(z) q^{-d_2} = q^{-1} f(z), \quad q^{d_2} \psi^\pm(z) q^{-d_2} = \psi^\pm(z), \tag{t0.4}$$

$$g(\gamma^{-1} z/w) \psi^+(z) \psi^-(w) = g(\gamma z/w) \psi^-(w) \psi^+(z), \tag{t1}$$

$$e(z) e(w) = g(z/w) e(w) e(z), \tag{t2}$$

$$f(z) f(w) = g(z/w)^{-1} f(w) f(z), \tag{t3}$$

$$\beta_1 \cdot [e(z), f(w)] = \delta(\gamma w/z) \psi^+(\gamma^{1/2} w) - \delta(\gamma z/w) \psi^-(\gamma^{1/2} z), \tag{t4}$$

$$\psi^\pm(z) e(w) = g(\gamma^{\pm 1/2} z/w) e(w) \psi^\pm(z), \tag{t5}$$

$$\psi^\pm(z) f(w) = g(\gamma^{\mp 1/2} z/w)^{-1} f(w) \psi^\pm(z), \tag{t6}$$

$$\text{Sym} \frac{z_2}{z_3} [e(z_1), [e(z_2), e(z_3)]] = 0, \tag{t7}$$

$$\text{Sym} \frac{z_2}{z_3} [f(z_1), [f(z_2), f(z_3)]] = 0, \tag{t8}$$

where the generating series are defined via:

$$e(z) := \sum_{r \in \mathbb{Z}} e_r z^{-r}, \quad f(z) := \sum_{r \in \mathbb{Z}} f_r z^{-r}, \quad \psi^\pm(z) := \psi_0^{\pm 1} + \sum_{r > 0} \psi_{\pm r} z^{\mp r},$$

and Sym denotes the symmetrization with respect to all permutations of $z_1, z_2, z_3$. This construction is manifestly symmetric with respect to the parameters $q_1, q_2, q_3$.

*Remark 2.1* In [4, 5], this algebra without $q^{\pm d_1}, q^{\pm d_2}$ and with $\gamma^{1/2} = 1$ was referred to as the "quantum continuous $\mathfrak{gl}_\infty$". A version of the latter algebra without the Serre relations (t7, t8) was referred to as the "Ding-Iohara algebra" in [9].

It is convenient to use the generators $\{h_{\pm r}\}_{r>0}$ instead of $\{\psi_{\pm r}\}_{r>0}$, defined via:

$$\exp\left(\mp\sum_{r>0}\frac{\beta_r}{r}h_{\pm r}z^{\mp r}\right)=\bar\psi^{\pm}(z):=\psi_0^{\mp 1}\psi^{\pm}(z) \qquad (2.4)$$

with $\beta_r\in\mathbb{C}^{\times}$ introduced in (2.2). Then, the relations (t5, t6) are equivalent to:

$$\psi_0 e_r=e_r\psi_0,\qquad [h_k,e_r]=\gamma^{-|k|/2}e_{r+k}\quad(k\neq 0),\qquad\text{(t5$'$)}$$

$$\psi_0 f_r=f_r\psi_0,\qquad [h_k,f_r]=-\gamma^{|k|/2}f_{r+k}\quad(k\neq 0),\qquad\text{(t6$'$)}$$

as follows from the equality

$$\log\left(\frac{(z-q_1^{-1}\gamma^{-1/2}w)(z-q_2^{-1}\gamma^{-1/2}w)(z-q_3^{-1}\gamma^{-1/2}w)}{(z-q_1\gamma^{-1/2}w)(z-q_2\gamma^{-1/2}w)(z-q_3\gamma^{-1/2}w)}\right)=\sum_{r>0}-\frac{\beta_r}{r}\cdot\gamma^{-r/2}\frac{w^r}{z^r},$$

with the left-hand side expanded in $w/z$. Likewise, the relation (t1) is equivalent to:

$$\psi_0 h_k=h_k\psi_0,\qquad [h_k,h_{-l}]=-\delta_{kl}\frac{k(\gamma^k-\gamma^{-k})}{\beta_k}\quad(k,l\neq 0),\qquad\text{(t1$'$)}$$

which we refer to as the Heisenberg algebra relations.

Finally, we note that in view of the relations (t0.1)–(t6), the Serre relations (t7, t8) can be replaced by equivalent simpler ones:

**Proposition 2.1** *If the relations (t5$'$, t6$'$) hold, then* (t7) *and* (t8) *are equivalent to:*

$$[e_0,[e_1,e_{-1}]]=0,\qquad\text{(t7$'$)}$$

$$[f_0,[f_1,f_{-1}]]=0\,.\qquad\text{(t8$'$)}$$

*Proof* The relation (t7$'$) arises by comparing the free terms in (t7). Let us now prove that (t5$'$, t7$'$) imply (t7). To this end, we consider endomorphisms $B_r$ of $\ddot{U}_{q_1,q_2,q_3}(\mathfrak{gl}_1)$ given by $B_0:=3\mathrm{Id}$ and $B_r\colon x\mapsto[\gamma^{|r|/2}h_r,x]$ for $r\neq 0$. Then, combining (t5$'$) with the algebraic equality $[t,[a,[b,c]]]=[[t,a],[b,c]]+[a,[[t,b],c]]+[a,[b,[t,c]]]$, we get that $B_{r_1,r_2,r_3}:=B_{r_1}B_{r_2}B_{r_3}-B_{r_1+r_2}B_{r_3}-B_{r_1+r_3}B_{r_2}-B_{r_2+r_3}B_{r_1}+2B_{r_1+r_2+r_3}$ maps

$$B_{r_1,r_2,r_3}\colon[e_0,[e_1,e_{-1}]]\mapsto\mathrm{Sym}\,[e_{r_1},[e_{r_2+1},e_{r_3-1}]],$$

where Sym denotes the symmetrization with respect to all permutations of $r_1,r_2,r_3$. Thus, (t5$'$) and (t7$'$) imply:

$$\mathrm{Sym}\,[e_{r_1},[e_{r_2+1},e_{r_3-1}]]=0\,.\qquad(2.5)$$

Multiplying (2.5) by $z_1^{-r_1}z_2^{-r_2}z_3^{-r_3}$ and summing over all $r_1,r_2,r_3\in\mathbb{Z}$, we obtain (t7). The proof of the equivalence of (t8) and (t8$'$) modulo (t6$'$) is analogous. $\qquad\square$

Let $\ddot{U}_{q_1,q_2,q_3}^<$, $\ddot{U}_{q_1,q_2,q_3}^>$, and $\ddot{U}_{q_1,q_2,q_3}^0$ be the $\mathbb{C}$-subalgebras of $\ddot{U}_{q_1,q_2,q_3}(\mathfrak{gl}_1)$ generated by $\{f_r\}_{r\in\mathbb{Z}}$, $\{e_r\}_{r\in\mathbb{Z}}$, and $\{\psi_r,\psi_0^{-1},\gamma^{\pm 1/2},q^{\pm d_1},q^{\pm d_2}\}_{r\in\mathbb{Z}}$, respectively. The following result is completely analogous to Proposition 1.1 (though its proof is slightly more tedious, due to the presence of the generators $\gamma^{\pm 1/2},q^{\pm d_1},q^{\pm d_2}$):

**Proposition 2.2** *(a) (Triangular decomposition of $\ddot{U}_{q_1,q_2,q_3}(\mathfrak{gl}_1)$) The multiplication map*

$$m\colon \ddot{U}^<_{q_1,q_2,q_3} \otimes \ddot{U}^0_{q_1,q_2,q_3} \otimes \ddot{U}^>_{q_1,q_2,q_3} \longrightarrow \ddot{U}_{q_1,q_2,q_3}(\mathfrak{gl}_1)$$

*is an isomorphism of $\mathbb{C}$-vector spaces.*

*(b) The algebra $\ddot{U}^>_{q_1,q_2,q_3}$ (resp. $\ddot{U}^0_{q_1,q_2,q_3}$ and $\ddot{U}^<_{q_1,q_2,q_3}$) is isomorphic to the associative $\mathbb{C}$-algebra generated by $\{e_r\}_{r\in\mathbb{Z}}$ (resp. $\{\psi_r, \psi_0^{-1}, \gamma^{\pm 1/2}, q^{\pm d_1}, q^{\pm d_2}\}_{r\in\mathbb{Z}}$ and $\{f_r\}_{r\in\mathbb{Z}}$) with the defining relations (t2, t7) (resp. (t0.1, t0.2, t1) and (t3, t8)).*

*Remark 2.2* We warn the reader that in Proposition 2.2(b), the relations (t7) and (t8) cannot be replaced with (t7′) and (t8′), respectively, in contrast to Proposition 2.1.

We equip $\ddot{U}_{q_1,q_2,q_3}(\mathfrak{gl}_1)$ with a $\mathbb{Z}^2$-grading via the following assignment:

$$\begin{aligned}
&\deg(e_r) := (1,r), \quad \deg(f_r) := (-1,r), \quad \deg(\psi_r) := (0,r) \quad \text{for} \quad r\in\mathbb{Z},\\
&\deg(x) := (0,0) \quad \text{for} \quad x = \psi_0^{-1}, \gamma^{\pm 1/2}, q^{\pm d_1}, q^{\pm d_2}\,.
\end{aligned} \tag{2.6}$$

We note that both $\gamma^{1/2}$ and $\psi_0$ are central elements of $\ddot{U}_{q_1,q_2,q_3}(\mathfrak{gl}_1)$. Furthermore, we shall also need the following modifications of the above algebra $\ddot{U}_{q_1,q_2,q_3}(\mathfrak{gl}_1)$:

- Let $\ddot{U}'_{q_1,q_2,q_3}(\mathfrak{gl}_1)$ be obtained from $\ddot{U}_{q_1,q_2,q_3}(\mathfrak{gl}_1)$ by "ignoring" the generators $q^{\pm d_2}$ and taking a quotient by the ideal $(\psi_0 - 1)$.

- Let $'\ddot{U}_{q_1,q_2,q_3}(\mathfrak{gl}_1)$ be obtained from $\ddot{U}_{q_1,q_2,q_3}(\mathfrak{gl}_1)$ by "ignoring" the generators $q^{\pm d_1}$ and taking a quotient by the ideal $(\gamma^{1/2} - 1)$.

- Let $U^{\mathrm{tor}}_{q_1,q_2,q_3}(\mathfrak{gl}_1)$ be obtained from $\ddot{U}_{q_1,q_2,q_3}(\mathfrak{gl}_1)$ by "ignoring" $q^{\pm d_1}$ and $q^{\pm d_2}$.

The above algebras $\ddot{U}'_{q_1,q_2,q_3}(\mathfrak{gl}_1)$, $'\ddot{U}_{q_1,q_2,q_3}(\mathfrak{gl}_1)$, $U^{\mathrm{tor}}_{q_1,q_2,q_3}(\mathfrak{gl}_1)$ satisfy the analogues of Proposition 2.2 with the corresponding subalgebras denoted by $\ddot{U}'^<_{q_1,q_2,q_3}$, $\ddot{U}'^>_{q_1,q_2,q_3}$, $\ddot{U}'^0_{q_1,q_2,q_3}$, $'\ddot{U}^<_{q_1,q_2,q_3}$, $'\ddot{U}^>_{q_1,q_2,q_3}$, $'\ddot{U}^0_{q_1,q_2,q_3}$, and $U^{\mathrm{tor}<}_{q_1,q_2,q_3}$, $U^{\mathrm{tor}>}_{q_1,q_2,q_3}$, $U^{\mathrm{tor}0}_{q_1,q_2,q_3}$, respectively. In particular, their "positive" subalgebras are naturally isomorphic:

$$\ddot{U}^>_{q_1,q_2,q_3} \simeq \ddot{U}'^>_{q_1,q_2,q_3} \simeq {}'\ddot{U}^>_{q_1,q_2,q_3} \simeq U^{\mathrm{tor}>}_{q_1,q_2,q_3}\,. \tag{2.7}$$

Finally, let us note that the algebra $\ddot{U}_{q_1,q_2,q_3}(\mathfrak{gl}_1)$ is endowed with a topological Hopf algebra structure given by the same formulas (H1)–(H3) as in Theorem 3.3, see [3], and can be realized as a Drinfeld double similarly to Theorem 3.4(b). However, we do not need these features in the present chapter besides for the discussion of $'\ddot{U}_{q_1,q_2,q_3}(\mathfrak{gl}_1)$-representations for which a simplified coproduct is recalled in (2.41).

### 2.1.2  Elliptic Hall algebra

We recall the (centrally extended) elliptic Hall algebra of [1]. First, let us introduce the following notation:

- Define $(\mathbb{Z}^2)^* := \mathbb{Z}^2\setminus\{(0,0)\}$ and $(\mathbb{Z}^2)^\pm := \{(a,b)|\pm a > 0 \text{ or } a = 0, \pm b > 0\}$.
- For $\mathbf{x} = (a,b) \in (\mathbb{Z}^2)^*$, we define $\deg(\mathbf{x}) := \gcd(a,b) \in \mathbb{Z}_{\geq 1}$.
- For $\mathbf{x} \in (\mathbb{Z}^2)^*$, we define $\epsilon_{\mathbf{x}} := 1$ if $\mathbf{x} \in (\mathbb{Z}^2)^+$ and $\epsilon_{\mathbf{x}} := -1$ if $\mathbf{x} \in (\mathbb{Z}^2)^-$.
- For noncollinear $\mathbf{x}, \mathbf{y} \in (\mathbb{Z}^2)^*$, we define $\epsilon_{\mathbf{x},\mathbf{y}} := \operatorname{sign}(\det(\mathbf{x},\mathbf{y})) \in \{\pm 1\}$.
- For noncollinear $\mathbf{x}, \mathbf{y} \in (\mathbb{Z}^2)^*$, we use $\triangle_{\mathbf{x},\mathbf{y}}$ to denote the triangle with vertices $\{(0,0), \mathbf{x}, \mathbf{x}+\mathbf{y}\}$.
- We call the triangle $\triangle_{\mathbf{x},\mathbf{y}}$ *empty* if it has no interior lattice points.
- For $r \geq 1$, we define $\alpha_r := -\frac{\beta_r}{r} = \frac{(1-q_1^{-r})(1-q_2^{-r})(1-q_3^{-r})}{r}$.

**Definition 2.1** [1, Definition 6.4] The (central extension of) elliptic Hall algebra $\widetilde{\mathcal{E}}$ is the associative algebra generated by $\{u_{\mathbf{x}}, \kappa_{\mathbf{y}} \,|\, \mathbf{x} \in (\mathbb{Z}^2)^*, \mathbf{y} \in \mathbb{Z}^2\}$ with the following defining relations:

$$\kappa_{\mathbf{x}} - \text{central}, \quad \kappa_{\mathbf{x}}\kappa_{\mathbf{y}} = \kappa_{\mathbf{x}+\mathbf{y}}, \quad \kappa_{0,0} = 1, \tag{2.8}$$

$$[u_{\mathbf{y}}, u_{\mathbf{x}}] = \delta_{\mathbf{x},-\mathbf{y}} \cdot \frac{\kappa_{\mathbf{x}} - \kappa_{\mathbf{x}}^{-1}}{\alpha_{\deg(\mathbf{x})}} \quad \text{if} \quad \mathbf{x}, \mathbf{y} \text{ are collinear}, \tag{2.9}$$

$$[u_{\mathbf{y}}, u_{\mathbf{x}}] = \epsilon_{\mathbf{x},\mathbf{y}}\kappa_{\alpha(\mathbf{x},\mathbf{y})}\frac{\theta_{\mathbf{x}+\mathbf{y}}}{\alpha_1} \quad \text{if} \quad \triangle_{\mathbf{x},\mathbf{y}} \text{ is empty and} \quad \deg(\mathbf{x}) = 1, \tag{2.10}$$

where

$$\alpha(\mathbf{x},\mathbf{y}) = \begin{cases} \epsilon_{\mathbf{x}}(\epsilon_{\mathbf{x}}\mathbf{x} + \epsilon_{\mathbf{y}}\mathbf{y} - \epsilon_{\mathbf{x}+\mathbf{y}}(\mathbf{x}+\mathbf{y}))/2 & \text{if } \epsilon_{\mathbf{x},\mathbf{y}} = 1 \\ \epsilon_{\mathbf{y}}(\epsilon_{\mathbf{x}}\mathbf{x} + \epsilon_{\mathbf{y}}\mathbf{y} - \epsilon_{\mathbf{x}+\mathbf{y}}(\mathbf{x}+\mathbf{y}))/2 & \text{if } \epsilon_{\mathbf{x},\mathbf{y}} = -1 \end{cases}$$

and the elements $\{\theta_{\mathbf{x}} \,|\, \mathbf{x} \in \mathbb{Z}^2\}$ are defined via

$$\sum_{r\geq 0} \theta_{r\mathbf{x}_0} z^r = \exp\left(\sum_{r>0} \alpha_r u_{r\mathbf{x}_0} z^r\right) \quad \text{for} \quad \mathbf{x}_0 \in (\mathbb{Z}^2)^* \quad \text{with} \quad \deg(\mathbf{x}_0) = 1. \tag{2.11}$$

This algebra is related to the quantum toroidal of $\mathfrak{gl}_1$ via the following result:

**Theorem 2.1** *[18, 19] The assignment*

$$u_{1,r} \mapsto e_r, \quad u_{-1,r} \mapsto f_r, \quad \theta_{0,\pm k} \mapsto \gamma^{\pm k/2}\psi_{\pm k}\psi_0^{\mp 1}, \quad \kappa_{1,0} \mapsto \psi_0, \quad \kappa_{0,1} \mapsto \gamma,$$

*for all $r \in \mathbb{Z}$ and $k \geq 1$, gives rise to an algebra isomorphism*

$$\Xi : \widetilde{\mathcal{E}}[\kappa_{0,1}^{\pm 1/2}] \xrightarrow{\sim} U_{q_1,q_2,q_3}^{\text{tor}}(\mathfrak{gl}_1). \tag{2.12}$$

*Remark 2.3* Since the algebra $\widetilde{\mathcal{E}}$ is naturally $\mathbb{Z}^2$-graded via $\deg(u_{\mathbf{x}}) = \mathbf{x}$, $\deg(\kappa_{\mathbf{y}}) = \mathbf{y}$, one can enlarge it to $\ddot{\mathcal{E}}$ by adjoining degree generators $q^{\pm d_1}$ and $q^{\pm d_2}$ satisfying:

$$\begin{aligned} q^{d_1} u_{a,b} q^{-d_1} &= q^{-b} u_{a,b}, \quad q^{d_2} u_{a,b} q^{-d_2} = q^a u_{a,b}, \\ q^{\pm d_r} \cdot q^{\mp d_r} &= 1, \quad q^{d_r}\kappa_{a,b}q^{-d_r} = \kappa_{a,b} \quad \text{for} \quad r \in \{1,2\}, \, a,b \in \mathbb{Z}. \end{aligned} \tag{2.13}$$

Then, (2.12) naturally extends to an algebra isomorphism $\ddot{\mathcal{E}}[\kappa_{0,1}^{\pm 1/2}] \xrightarrow{\sim} \ddot{U}_{q_1,q_2,q_3}(\mathfrak{gl}_1)$.

### 2.1.3  "90 degree rotation" automorphism

In this section, we shall enlarge the algebra $U^{\mathrm{tor}}_{q_1,q_2,q_3}(\mathfrak{gl}_1)$ by formally adjoining $\psi_0^{\pm 1/2}$, a square root of the central element $\psi_0$. After this upgrade, the algebra $U^{\mathrm{tor}}_{q_1,q_2,q_3}(\mathfrak{gl}_1)$ admits an important automorphism constructed first in [13]:

**Theorem 2.2** *[13, Proposition 3.2] (a) The assignment*

$$
\begin{aligned}
\psi_0^{\pm 1/2} &\mapsto \gamma^{\pm 1/2}, \qquad \gamma^{\pm 1/2} \mapsto \psi_0^{\mp 1/2}, \\
e_0 &\mapsto -h_1 \gamma^{-1/2}, \quad f_0 \mapsto h_{-1}\gamma^{1/2}, \quad h_1 \mapsto -f_0 \psi_0^{-1/2} \quad h_{-1} \mapsto e_0 \psi_0^{1/2},
\end{aligned}
\tag{2.14}
$$

*gives rise to an order four automorphism of $U^{\mathrm{tor}}_{q_1,q_2,q_3}(\mathfrak{gl}_1)$:*

$$
\varpi\colon U^{\mathrm{tor}}_{q_1,q_2,q_3}(\mathfrak{gl}_1) \xrightarrow{\ \sim\ } U^{\mathrm{tor}}_{q_1,q_2,q_3}(\mathfrak{gl}_1), \qquad \varpi^4 = \mathrm{Id}\,.
\tag{2.15}
$$

*(b) If $\deg(x) = (a,b) \in \mathbb{Z}^2$, see (2.6), then $\deg(\varpi(x)) = (-b,a) \in \mathbb{Z}^2$.*

*(c) For $k \geq 1$, we have:*

$$
\varpi\colon Y_k^+ \mapsto \beta_1^{-1}\gamma^{-k/2}\cdot\psi_k\psi_0^{-1}, \quad Y_k^- \mapsto \beta_1^{-1}\gamma^{k/2}\cdot\psi_{-k}\psi_0,
\tag{2.16}
$$

*where we set*

$$
Y_k^+ := \underbrace{[e_{-1},[e_0,\cdots,[e_0,e_1]\cdots]]}_{k\ \mathrm{factors}}, \qquad Y_k^- := \underbrace{[f_1,[f_0,\cdots,[f_0,f_{-1}]\cdots]]}_{k\ \mathrm{factors}}
\tag{2.17}
$$

*for $k \geq 2$, while $Y_1^+ := e_0$ and $Y_1^- := f_0$.*

Evoking the identification of $U^{\mathrm{tor}}_{q_1,q_2,q_3}(\mathfrak{gl}_1)$ and $\widetilde{\mathcal{E}}$ from Theorem 2.1, we shall refer to $\varpi$ of (2.15) as (counterclockwise) "90 degree rotation" automorphism. We note that $Y_k^\pm$ from (2.17) get identified with nonzero multiples of $\theta_{(\pm k,0)}$ from (2.11).

*Remark 2.4* (a) In [13], the parameters $q,\gamma$ are related to our $q_1,q_2$ via $q \leftrightarrow q_1^{1/2}$, $\gamma \leftrightarrow q_2^{1/2}$, while the generators $\{a_{\pm r}, X_k^\pm, C^{\pm 1}, C^{\pm 1}\}_{k \in \mathbb{Z}}^{r \geq 1}$ of $\mathcal{U}_{q,\gamma}$ are related to ours via:

$$
\begin{aligned}
C^{\pm 1} &\leftrightarrow \gamma^{\pm 1}, \quad C^{\pm 1} \leftrightarrow \psi_0^{\pm 1}, \quad X_k^+ \leftrightarrow (q_3^{-1/2}-q_3^{1/2})e_k, \quad X_k^- \leftrightarrow (q_3^{-1/2}-q_3^{1/2})f_k, \\
a_{\pm r} &\leftrightarrow \mp q_1^{1/2} q_2^{r/2}(1-q_1^r)(1-q_3^r)/r(1-q_1)\cdot h_{\pm r}\gamma^{\mp r/2}\,.
\end{aligned}
$$

(b) The automorphism $\psi$ from [13, Proposition 3.2] is the inverse of $\varpi$ from (2.15).

In the next section we will need a version of $\varpi$, the algebra isomorphism:

$$
\varpi\colon {}'\ddot{U}_{q_1,q_2,q_3}(\mathfrak{gl}_1) \xrightarrow{\ \sim\ } \ddot{U}'_{q_1,q_2,q_3}(\mathfrak{gl}_1)
\tag{2.18}
$$

determined by (2.14) and the assignment $q^{d_2} \mapsto q^{-d_1}$, which is compatible with Theorem 2.2(b). According to Theorem 2.2(c), the $\varpi$-preimage of $\ddot{U}'^0_{q_1,q_2,q_3}$ provides a *horizontal* Heisenberg subalgebra of ${}'\ddot{U}^0_{q_1,q_2,q_3}$ generated by $\{Y_k^\pm, \psi_0^{\pm 1/2}\}_{k \geq 1}$.

### 2.1.4 Shuffle algebra realization

Consider an $\mathbb{N}$-graded $\mathbb{C}$-vector space $\mathbb{S} = \bigoplus_{k \in \mathbb{N}} \mathbb{S}_k$, where $\mathbb{S}_k$ consists of $\Sigma_k$-symmetric rational functions in the variables $\{x_r\}_{r=1}^k$. We also fix a rational function

$$\zeta(z) := \frac{(q_1 z - 1)(q_2 z - 1)(q_3 z - 1)}{(z-1)^3} \,. \tag{2.19}$$

Similarly to (1.32), let us now introduce the bilinear shuffle product $\star$ on $\mathbb{S}$: given $F \in \mathbb{S}_k$ and $G \in \mathbb{S}_\ell$, we define $F \star G \in \mathbb{S}_{k+\ell}$ via

$$(F \star G)(x_1, \ldots, x_{k+\ell}) := \frac{1}{k!\,\ell!} \operatorname{Sym}\left(F(x_1, \ldots, x_k)G(x_{k+1}, \ldots, x_{k+\ell}) \prod_{\substack{r \le k}}^{\substack{r' > k}} \zeta\left(\frac{x_r}{x_{r'}}\right)\right).$$

This endows $\mathbb{S}$ with a structure of an associative $\mathbb{C}$-algebra with the unit $\mathbf{1} \in \mathbb{S}_0$.

As before, we consider the subspace of $\mathbb{S}$ defined by the *pole* and *wheel* conditions:

- We say that $F \in \mathbb{S}_k$ satisfies the *pole conditions* if

$$F(x_1, \ldots, x_k) = \frac{f(x_1, \ldots, x_k)}{\prod_{r \ne r'}(x_r - x_{r'})}, \quad \text{where} \quad f \in \mathbb{C}[x_1^{\pm 1}, \ldots, x_k^{\pm 1}]^{\Sigma_k} \,. \tag{2.20}$$

- We say that $F \in \mathbb{S}_k$ satisfies the *wheel conditions* if

$$F(x_1, \ldots, x_k) = 0 \ \text{ once } \ \left\{\frac{x_{r_1}}{x_{r_2}}, \frac{x_{r_2}}{x_{r_3}}, \frac{x_{r_3}}{x_{r_1}}\right\} = \{q_1, q_2, q_3\} \ \text{ for some } \ r_1, r_2, r_3 \,. \tag{2.21}$$

Let $S_k \subset \mathbb{S}_k$ denote the subspace of all $F$ satisfying these two conditions and set

$$S := \bigoplus_{k \in \mathbb{N}} S_k \,.$$

Similarly to Lemma 1.2, it is straightforward to check that $S \subset \mathbb{S}$ is actually $\star$-closed:

**Lemma 2.1** *For any $F \in S_k$ and $G \in S_\ell$, we have $F \star G \in S_{k+\ell}$.*

The algebra $(S, \star)$ is called the **shuffle algebra** (of type $\widehat{\mathfrak{gl}}_1$). It is $\mathbb{Z}^2$-graded via:

$$S = \bigoplus_{(k,d) \in \mathbb{N} \times \mathbb{Z}} S_{k,d} \qquad \text{with} \qquad S_{k,d} := \left\{F \in S_k \,\middle|\, \text{tot.deg}(F) = d\right\} \,. \tag{2.22}$$

Similarly to Proposition 1.3, the shuffle algebra $(S, \star)$ is linked to the "positive" subalgebra of the quantum toroidal $\mathfrak{gl}_1$, cf. (2.7), via:

**Proposition 2.3** *There exists an injective $\mathbb{C}$-algebra homomorphism*

$$\Psi \colon \ddot{U}^{>}_{q_1, q_2, q_3} \longrightarrow S \tag{2.23}$$

*such that $e_r \mapsto x_1^r$ for any $r \in \mathbb{Z}$.*

Let $\widetilde{\mathcal{E}}^>$ be the "positive" subalgebra of $\widetilde{\mathcal{E}}$ generated by $\{u_{k,d}\}_{k \geq 1}^{d \in \mathbb{Z}}$. Then, composing (2.23) with the restriction of (2.12) to the "positive" subalgebras, we obtain:

$$\Upsilon \colon \widetilde{\mathcal{E}}^> \hookrightarrow S. \tag{2.24}$$

The following result was conjectured in [9, 19] and proved in [17, Proposition 3.2]:

**Theorem 2.3** *[17]* $\Upsilon \colon \widetilde{\mathcal{E}}^> \xrightarrow{\sim} S$ *is an isomorphism of* $\mathbb{N} \times \mathbb{Z}$*-graded* $\mathbb{C}$*-algebras.*

Thus, the shuffle algebra $S$ is generated by its first component $S_1$, and we have:

$$\widetilde{\mathcal{E}}^> \simeq \ddot{U}_{q_1,q_2,q_3}^> \simeq S. \tag{2.25}$$

### 2.1.5 Neguț's proof of Theorem 2.3

In this section, we provide a proof of Theorem 2.3 following the original argument of [17] that crucially relies on the notion of slope $\leq \mu$ subalgebras:

**Definition 2.2** For $\mu \in \mathbb{R}$, define $S^\mu = \bigoplus_{(k,d) \in \mathbb{N} \times \mathbb{Z}} S_{k,d}^\mu$ with

$$S_{k,d}^\mu = \left\{ F \in S_{k,d} \,\Big|\, \lim_{\xi \to \infty} \frac{F(\xi \cdot x_1, \ldots, \xi \cdot x_\ell, x_{\ell+1}, \ldots, x_k)}{\xi^{\ell\mu}} \text{ exists } \forall\, 0 \leq \ell \leq k \right\}. \tag{2.26}$$

Using the equalities $\lim\limits_{\xi \to \infty} \zeta(\xi) = \lim\limits_{\xi \to \infty} \zeta(1/\xi) = 1$, we immediately get:

**Lemma 2.2** *[17, Proposition 2.3] The subspace* $S^\mu$ *is a subalgebra of* $S$.

The subspaces $S_{k,d}^\mu$ give an increasing filtration of (an infinite dimensional) $S_{k,d}$: $S_{k,d} = \bigcup_{\mu \in \mathbb{R}} S_{k,d}^\mu$ and $S_{k,d}^\mu \subset S_{k,d}^{\mu'}$ for $\mu < \mu'$. The key point is that $S_{k,d}^\mu$ are actually finite dimensional, with an explicit upper bound on their dimensions:

**Proposition 2.4** *[17, Proposition 2.4]* $\dim(S_{k,d}^\mu) \leq$ *number of unordered collections*

$$(k_1, d_1), \ldots, (k_t, d_t) \text{ with } k_1 + \cdots + k_t = k, \, d_1 + \cdots + d_t = d, \, d_r \leq \mu k_r, \tag{2.27}$$

*where* $t \geq 1$, $k_r \geq 1$, *and* $d_r \in \mathbb{Z}$ *for* $1 \leq r \leq t$.

*Proof* Similarly to (1.40), for any partition $\pi = \{k_1 \geq k_2 \geq \cdots \geq k_t > 0\}$ of $k$ let us define the **specialization map**

$$\phi_\pi \colon S_{k,d}^\mu \longrightarrow \mathbb{C}[y_1^{\pm 1}, \ldots, y_t^{\pm 1}]. \tag{2.28}$$

Explicitly, for $F$ of the form (2.20), we define $\phi_\pi(F)$ as a specialization of $f$ under:

$$x_{k_1 + \cdots + k_{r-1} + a} \mapsto q_2^a y_r \quad \text{for} \quad 1 \leq r \leq t, \, 1 \leq a \leq k_r.$$

This gives rise to the filtration (called the *Gordon filtration* in [6]) on $S_{k,d}^{\mu}$ via:

$$S_{k,d}^{\mu,\pi} := \bigcap_{\widetilde{\pi} > \pi} \mathrm{Ker}(\phi_{\widetilde{\pi}}), \tag{2.29}$$

where $>$ denotes the dominance partial order on the set of all partitions of $k$. To prove the desired upper bound on $\dim(S_{k,d}^{\mu})$, it thus suffices to show:

$$\dim \phi_{\pi}\left(S_{k,d}^{\mu,\pi}\right) \leq \#\left\{(d_1, \ldots, d_t) \in \mathbb{Z}^t \mid d_1 + \cdots + d_t = d, \, d_r \leq \mu k_r \,\forall\, r\right\}, \tag{2.30}$$

where the order of $d_r$ and $d_{r'}$ is disregarded whenever $k_r = k_{r'}$.

First, we note that the wheel conditions (2.21) for $F$ guarantee that $\phi_{\pi}(F)$, which is a Laurent polynomial in $\{y_1, \ldots, y_t\}$, vanishes under the following specializations:

(i) $q_2^x y_r = q_1 \cdot q_2^{x'} y_{r'}$ for any $1 \leq r < r' \leq t$, $1 \leq x < k_r$, $1 \leq x' \leq k_{r'}$,

(ii) $q_2^x y_r = q_3 \cdot q_2^{x'} y_{r'}$ for any $1 \leq r < r' \leq t$, $1 \leq x < k_r$, $1 \leq x' \leq k_{r'}$.

Second, the condition $\phi_{\widetilde{\pi}}(F) = 0$ for any $\widetilde{\pi} > \pi$ also implies that $\phi_{\pi}(F)$ vanishes under the following specializations:

(iii) $y_r = q_2^{x'} y_{r'}$ for any $1 \leq r < r' \leq t$, $1 \leq x' \leq k_{r'}$,

(iv) $q_2^{k_r+1} y_r = q_2^{x'} y_{r'}$ for any $1 \leq r < r' \leq t$, $1 \leq x' \leq k_{r'}$.

Combining the above vanishing conditions (i)–(iv) for $\phi_{\pi}(F)$ with $F \in S_{k,d}^{\mu,\pi}$, we conclude that $\phi_{\pi}(F)$ is divisible by $Q_{\pi} \in \mathbb{C}[y_1, \ldots, y_t]$, defined as a product of the linear terms in $y_r$'s arising from (i)–(iv), counted with the correct multiplicities.

The total degree of $Q_{\pi}$ as well as the degree with respect to each variable $y_r$ are:

$$\mathrm{tot.deg}(Q_{\pi}) = k^2 - \sum_{1 \leq r \leq t} k_r^2 \quad \text{and} \quad \deg_{y_r}(Q_{\pi}) = 2k k_r - 2k_r^2. \tag{2.31}$$

On the other hand, the specialization $\phi_{\pi}(F)$ satisfies:

$$\mathrm{tot.deg}(\phi_{\pi}(F)) = d + k^2 - k \quad \text{and} \quad \deg_{y_r}(\phi_{\pi}(F)) \leq \mu k_r + 2k k_r - k_r^2 - k_r. \tag{2.32}$$

To see the above inequality one should multiply by $\xi$ only those $x_*$-variables that got specialized to $q_2^2$-multiples of $y_r$ and use the condition that the slope of $F$ is $\leq \mu$.

Define $r_{\pi} := \phi_{\pi}(F)/Q_{\pi} \in \mathbb{C}[y_1^{\pm 1}, \ldots, y_t^{\pm 1}]$. Combining (2.31) and (2.32), we get:

$$\mathrm{tot.deg}(r_{\pi}) = d + \sum_{r=1}^{t} k_r(k_r - 1) \quad \text{and} \quad \deg_{y_r}(r_{\pi}) \leq k_r^2 - k_r + \mu k_r. \tag{2.33}$$

The space of Laurent polynomials satisfying (2.33) is spanned by the monomials

$$y_1^{k_1^2 - k_1 + d_1} \cdots y_t^{k_t^2 - k_t + d_t} \quad \text{with} \quad d_1 + \cdots + d_t = d \text{ and } d_r \leq \mu k_r \,\forall\, r.$$

Finally, we note that both $\phi_{\pi}(F)$ and $Q_{\pi}$ are symmetric in $y_r$ and $y_{r'}$ whenever $k_r = k_{r'}$, hence, so is $r_{\pi}$. This implies (2.30) and thus completes the proof. □

Now we are ready to complete the proof of Theorem 2.3. To this end, for any $k \geq 1$, $d \in \mathbb{Z}$, $\mu \in \mathbb{R}$, we consider the following subspace $\widetilde{\mathcal{E}}_{k,d}^{\mu}$ of $\widetilde{\mathcal{E}}$:

$$\widetilde{\mathcal{E}}_{k,d}^{\mu} = \text{Span}\left\{ u_{k_1,d_1} \cdots u_{k_t,d_t} | k_1 + \cdots + k_t = k, d_1 + \cdots + d_t = d, d_r \leq \mu k_r, k_r \geq 1 \right\}$$

According to [1, Lemma 5.6], the "positive" subalgebra $\widetilde{\mathcal{E}}^{>}$ generated by $\{u_{a,b}\}_{a>0}^{b \in \mathbb{Z}}$ has a basis consisting of convex paths in $\mathbf{Conv}^{+}$ (cf. [1]), while $\widetilde{\mathcal{E}}_{k,d}^{\mu}$ has a basis consisting of such convex paths of slope $\leq \mu$ starting from $(0,0)$ and ending at $(k,d)$. The latter are in bijection with the collections (2.27). Combining this observation with the injectivity of $\Upsilon$ from (2.24) and Proposition 2.4, it thus remains to prove the following result:

**Proposition 2.5** *[17]* $\Upsilon(\widetilde{\mathcal{E}}_{k,d}^{\mu}) \subseteq S_{k,d}^{\mu}$.

*Proof* Due to Lemma 2.2, it suffices to verify that elements $F_{k,d} := \Upsilon(u_{k,d})$ have slope $\leq \frac{d}{k}$. This is proved by an induction on $k$, the base case $k = 1$ being trivial. To prove the induction step, for a fixed $(k,d)$ pick an empty triangle $\triangle$ with vertices $\{(0,0),(k_2,d_2),(k,d)\}$ such that $\frac{d_2}{k_2} < \frac{d}{k} < \frac{d_1}{k_1}$ and $\gcd(k_1,d_1) = 1 = \gcd(k_2,d_2)$ with $(k_1,d_1) = (k - k_2, d - d_2)$ and $k_1, k_2 > 0$ (for example, pick a triangle with the minimal area). By the induction hypothesis and the relation (2.10), it remains to show that the commutator $[F_{k_1,d_1}, F_{k_2,d_2}]$ has a slope $\leq \frac{d}{k}$.

By the induction hypothesis, a slope of $F_{k_1,d_1}$ is $\leq \frac{d_1}{k_1}$, that is, as we multiply $\ell_1$ variables by $\xi$ the result grows at the speed $\xi^{\leq \lfloor \frac{d_1}{k_1} \ell_1 \rfloor}$, and similarly for $F_{k_2,d_2}$. Let us now consider any summand from the symmetrization appearing in $F_{k_1,d_1} \star F_{k_2,d_2}$ and $F_{k_2,d_2} \star F_{k_1,d_1}$. As we multiply $\ell$ variables by $\xi$, let $\ell_1$ and $\ell_2$ denote the number of those $x_*$-variables that got placed into $F_{k_1,d_1}$ and $F_{k_2,d_2}$, respectively, so that $\ell = \ell_1 + \ell_2$. Then, the corresponding summand grows at the speed $\xi^{\leq \lfloor \frac{d_1}{k_1} \ell_1 \rfloor + \lfloor \frac{d_2}{k_2} \ell_2 \rfloor}$. However, the choice of the triangle $\triangle$ guarantees that $\lfloor \frac{d_1}{k_1} \ell_1 \rfloor + \lfloor \frac{d_2}{k_2} \ell_2 \rfloor \leq \frac{d}{k}(\ell_1 + \ell_2)$ for any $0 \leq \ell_1 \leq k_1$ and $0 \leq \ell_2 \leq k_2$, with the only exception when $\ell_1 = k_1, \ell_2 = 0$. But in the latter case, the corresponding summands in both $F_{k_1,d_1} \star F_{k_2,d_2}$ and $F_{k_2,d_2} \star F_{k_1,d_1}$ with $\ell$ variables multiplied by $\xi$ have the form

$$\xi^{d_1} \cdot F_{k_1,d_1}(x_1,\ldots,x_{k_1}) \cdot F_{k_2,d_2}(x_{k_1+1},\ldots,x_k) + O(\xi^{\frac{k_1 d}{k}}),$$

and hence the commutator $[F_{k_1,d_1}, F_{k_2,d_2}]$ indeed has a slope $\leq \frac{d}{k}$, as claimed.  $\square$

*Remark 2.5* (a) The result above establishes a vector space isomorphism $\widetilde{\mathcal{E}}_{k,d}^{\mu} \simeq S_{k,d}^{\mu}$ which implies that the inequality of Proposition 2.4 is actually an equality.

(b) The proof of Proposition 2.4 implies that the subspace of all $F \in S_{k,d}$ satisfying

$$\lim_{\xi \to \infty} \frac{F(\xi \cdot x_1,\ldots,\xi \cdot x_\ell, x_{\ell+1},\ldots,x_k)}{\xi^{\frac{d}{k}\mu}} = 0 \qquad \forall\, 0 < \ell < k \qquad (2.34)$$

is at most one-dimensional. It is proved in [17, §5] that $\Upsilon(u_{k,d})$ satisfies this property, while an explicit formula for $\Upsilon(u_{k,d}) \in S_{k,d}$ is derived in [17, Theorem 1.1, §6].

### 2.1.6 Commutative subalgebra $\mathcal{A}$

Following [6], let us consider an important $\mathbb{N}$-graded subspace $\mathcal{A} = \bigoplus_{k \geq 0} \mathcal{A}_k$ of $S$:

$$\mathcal{A}_k = \left\{ F \in S_{k,0} \, \Big| \, \partial^{(\infty,\ell)} F \text{ exists } \forall \, 1 \leq \ell \leq k \right\}, \tag{2.35}$$

where

$$\partial^{(\infty,\ell)} F := \lim_{\xi \to \infty} F(\xi \cdot x_1, \ldots, \xi \cdot x_\ell, x_{\ell+1}, \ldots, x_k). \tag{2.36}$$

Equivalently, we have $\mathcal{A}_k = S_{k,0}^0$, see (2.26). By Lemma 2.2, $\mathcal{A}$ is a subalgebra of $S$. This subalgebra obviously contains the following family of elements:

$$K_1(x_1) = x_1^0, \quad K_k(x_1, \ldots, x_k) = \prod_{1 \leq r \neq r' \leq k} \frac{x_r - q_2 x_{r'}}{x_r - x_{r'}} \quad \forall \, k \geq 2. \tag{2.37}$$

The following description of $\mathcal{A}$ was established in [6, Proposition 2.20]:

**Theorem 2.4** *[6] The subalgebra $\mathcal{A} \subset S$ is shuffle-generated by $\{K_k\}_{k \geq 1}$ and is a polynomial algebra in them. In particular, $\mathcal{A}$ is a commutative subalgebra of $S$.*

*Remark 2.6* The same result obviously holds if $q_2$ in (2.37) is replaced with $q_1$ or $q_3$.

*Proof* According to Remark 2.5(a), we have:

$$\dim \mathcal{A}_k = p(k) = \# \{\text{partitions of } k\}. \tag{2.38}$$

We shall identify partitions $\pi$ with Young diagrams $\lambda = \{\lambda_1 \geq \lambda_2 \geq \cdots \geq \lambda_t > 0\}$ and use $\lambda'$ for the transposed Young diagram. We define $K_\lambda := K_{\lambda_1} \star K_{\lambda_2} \star \cdots \star K_{\lambda_t}$, so that $K_\lambda \in \mathcal{A}_{|\lambda|}$ by Lemma 2.2. Analogously to Lemmas 1.5–1.6, we have:

**Lemma 2.3** $\phi_\lambda(K_{\lambda'}) \neq 0$ and $\phi_\mu(K_{\lambda'}) = 0$ for any $\mu > \lambda$ with $|\mu| = |\lambda|$. $\qquad \square$

Combining this lemma with the dimension formula (2.38), we immediately see that $\{K_\lambda\}_{|\lambda|=k}$ is a basis of $\mathcal{A}_k$. Thus, it remains only to prove the commutativity:

**Lemma 2.4** $[K_{m_1}, K_{m_2}] = 0$ for any $m_1, m_2 \geq 1$. $\qquad \square$

*Proof* It is sufficient to prove that $\phi_{(2,1^{m_1+m_2-2})}([K_{m_1}, K_{m_2}]) = 0$. Indeed, this would imply that the numerator $f$ from (2.20) in $[K_{m_1}, K_{m_2}]$ is divisible by $\prod_{r \neq r'} (x_r - q_2 x_{r'})$, and hence for degree reasons $[K_{m_1}, K_{m_2}] = v_{m_1, m_2} K_{m_1+m_2}$ with $v_{m_1, m_2} \in \mathbb{C}$. Let us multiply both sides of this equality by $\prod_{r \neq r'} (x_r - x_{r'})$ and consider a specialization $x_r \mapsto y$ for all $1 \leq r \leq m_1 + m_2$. The left-hand side will clearly specialize to 0, while the right-hand side will specialize to $((1 - q_2)y)^{(m_1+m_2)(m_1+m_2-1)} \cdot v_{m_1, m_2}$. Thus, $v_{m_1, m_2} = 0$ and so $[K_{m_1}, K_{m_2}] = 0$.

Let us now prove the stated equality $\phi_{(2,1^{m_1+m_2-2})}([K_{m_1}, K_{m_2}]) = 0$ for $m_1, m_2 \geq 0$, where we set $K_0 := 1 \in S_0$. The proof is by induction on $\min\{m_1, m_2\}$, with the base of induction ($m_1 = 0$ or $m_2 = 0$) being obvious. By a direct computation, we get

$$\phi_{(2,1^{m_1+m_2-2})}(K_{m_1} \star K_{m_2}) = \frac{1}{(m_1-1)!(m_2-1)!} \operatorname{Sym}(A_1 \cdot B_1)$$

with the symmetrization Sym taken over all permutations of $\{y_r\}_{r=2}^{m_1+m_2-1}$, where $A_1$ is an $(m_1, m_2)$-independent rational function in $\{y_r\}_{r=1}^{m_1+m_2-1}$ which is symmetric in $\{y_r\}_{r=2}^{m_1+m_2-1}$ and $B_1$ is explicitly given by:

$$B_1 = \prod_{2 \le r \ne r' \le m_1} \frac{y_r - q_2 y_{r'}}{y_r - y_{r'}} \prod_{m_1 < r \ne r' < m_1+m_2} \frac{y_r - q_2 y_{r'}}{y_r - y_{r'}} \prod_{2 \le r \le m_1}^{m_1 < r' < m_1+m_2} \zeta\left(\frac{y_r}{y_{r'}}\right).$$

Therefore, we have:

$$\phi_{(2,1^{m_1+m_2-2})}(K_{m_1} \star K_{m_2}) = A_1 \cdot \left(K_{m_1-1} \star K_{m_2-1}\right)(y_2,\ldots,y_{m_1+m_2-1}),$$

while permuting $m_1 \leftrightarrow m_2$ we also obtain:

$$\phi_{(2,1^{m_1+m_2-2})}(K_{m_2} \star K_{m_1}) = A_1 \cdot \left(K_{m_2-1} \star K_{m_1-1}\right)(y_2,\ldots,y_{m_1+m_2-1}).$$

Thus, we get:

$$\phi_{(2,1^{m_1+m_2-2})}([K_{m_1}, K_{m_2}]) = A_1 \cdot [K_{m_1-1}, K_{m_2-1}](y_2,\ldots,y_{m_1+m_2-1}) = 0$$

by the induction hypothesis. This proves the induction step, hence, also the claim. □

This completes the proof of Theorem 2.4.                                      □

On the other hand, we also have another important set of generators for $\mathcal{A}$:

**Proposition 2.6** *The algebra $\mathcal{A}$ is a polynomial algebra in $\{Y_k^+\}_{k \ge 1}$ of (2.17).*

*Proof* According to Theorem 2.1 and the formula (2.10), $Y_k^+$ are nonzero multiples of $\Xi(\theta_{k,0})$. Due to (2.9, 2.11), it thus suffices to prove that unordered monomials in pairwise commuting $\{u_{k,0}\}_{k \ge 1}$ form a basis of $\bigoplus_{k \in \mathbb{N}} \widetilde{\mathcal{E}}_{k,0}^0$. But the former are in bijection with the convex paths in $\mathbf{Conv}^+$ of slope $0$ staring from $(0,0)$ and ending in $(k,0)$ with $k > 0$, which do form a basis of $\bigoplus_{k \in \mathbb{N}} \widetilde{\mathcal{E}}_{k,0}^0$ by [1], as mentioned above. □

## 2.2  Representations of quantum toroidal $\mathfrak{gl}_1$

In this section, we introduce the appropriate categories $O^\pm$ of $'\ddot{U}_{q_1,q_2,q_3}(\mathfrak{gl}_1)$-modules and recall several families of $'\ddot{U}_{q_1,q_2,q_3}(\mathfrak{gl}_1)$-modules from [4, 5, 7] that admit nice combinatorial bases. Twisting one of them by the "90 degree rotation" isomorphism $\varpi^{-1}$ of (2.18), we recover "level 1 vertex-type" $\ddot{U}_{q_1,q_2,q_3}'(\mathfrak{gl}_1)$-modules of [6]. We also recall the shuffle realization of these combinatorial modules following [8, 22].

In this section, we shall assume that $q_1$ and $q_3$ are *generic*:

$$q_1^a q_3^b = 1 \ (a, b \in \mathbb{Z}) \iff a = b = 0 \,. \tag{2.39}$$

### 2.2.1 Categories $O^{\pm}$

In this section, we introduce the categories $O^{\pm}$ for ${}'\ddot{U}_{q_1,q_2,q_3}(\mathfrak{gl}_1)$ and classify their lowest and highest weight modules in a standard way. Let us recall the $\mathbb{Z}$-grading on ${}'\ddot{U}_{q_1,q_2,q_3}(\mathfrak{gl}_1)$ of (2.6) arising via the adjoint action $x \mapsto q^{d_2} x q^{-d_2}$, see (t0.4).

**Definition 2.3** The category $O^-$ of ${}'\ddot{U}_{q_1,q_2,q_3}(\mathfrak{gl}_1)$-modules is the full subcategory of the category of $\mathbb{Z}$-graded modules $L$ satisfying the following two conditions:

(i) for any $v \in L$ there exists $N \in \mathbb{Z}$ such that ${}'\ddot{U}_{q_1,q_2,q_3}(\mathfrak{gl}_1)_{\leq N}(v) = 0$,
(ii) $L$ is of *finite type*, i.e., the graded components $\{L_k\}_{k \in \mathbb{Z}}$ are all finite dimensional.

An important family of modules satisfying (i) are the lowest weight modules:

**Definition 2.4** A ${}'\ddot{U}_{q_1,q_2,q_3}(\mathfrak{gl}_1)$-module $L$ is called a *lowest weight module* if it is generated by $v_0 \in L_k$ satisfying $f_r(v_0) = 0$, $\psi_r(v_0) = \phi_r \cdot v_0$, $q^{d_2}(v_0) = \nu \cdot v_0$ for all $r \in \mathbb{Z}$ with $\nu \in \mathbb{C}^{\times}$, $\phi_r \in \mathbb{C}$, $\phi_0 \neq 0$. The vector $v_0$ is called the *lowest weight vector*.

Without loss of generality, we can assume that $v_0 \in L_0$. Let us encode $\{\phi_r\}_{r \in \mathbb{Z}}$ by two series $\phi^{\pm}(z) := \phi_0^{\pm 1} + \sum_{r>0} \phi_{\pm r} z^{\mp r} \in \mathbb{C}[[z^{\mp 1}]]$, so that $\psi^{\pm}(z) v_0 = \phi^{\pm}(z) \cdot v_0$. For any such series $\phi^{\pm}(z) \in \mathbb{C}[[z^{\mp 1}]]$ and $\nu \in \mathbb{C}^{\times}$, there is a universal lowest weight module $M_{\phi^+,\phi^-,\nu}$, defined as the quotient of ${}'\ddot{U}_{q_1,q_2,q_3}(\mathfrak{gl}_1)$ by the left ideal generated by $\{f_r, \psi_r - \phi_r, q^{d_2} - \nu\}_{r \in \mathbb{Z}}$ with $v_0 = 1$. By standard arguments, it has a unique irreducible quotient $V_{\phi^+,\phi^-,\nu}$, which clearly satisfies the condition (i) of Definition 2.3. The next result (which first appeared in [13, Theorem 6.1]) provides a standard criterion for $V_{\phi^+,\phi^-,\nu}$ to satisfy (ii):

**Proposition 2.7** $V_{\phi^+,\phi^-,\nu}$ *is in $O^-$ if and only if there is a rational function $\phi(z) \in \mathbb{C}(z)$ such that $\phi(0)\phi(\infty) = 1$ and $\phi^{\pm}(z) = \phi(z)^{\pm}$–the expansions of $\phi(z)$ in $z^{\mp 1}$.*

*Proof* The proof is similar to that for quantum loop algebras, see [2], so we shall only sketch it. Define constants $\{\bar{\phi}_r\}_{r \in \mathbb{Z}}$ as $\phi_r$ (for $r > 0$), $-\phi_r$ (for $r < 0$), and $\phi_0 - \phi_0^{-1}$ (for $r = 0$). To prove the "only if" part, choose indices $k \in \mathbb{Z}, l \in \mathbb{N}$ such that $\{e_k(v_0), \ldots, e_{k+l}(v_0)\}$ span the component $(V_{\phi^+,\phi^-,\nu})_1$ and $a_0 e_k(v_0) + a_1 e_{k+1}(v_0) + \cdots + a_l e_{k+l}(v_0) = 0$ for some $a_0, \ldots, a_l \in \mathbb{C}$ with $a_l \neq 0$. Applying $f_{r-k}$ to the above equality and using (t4) in the form $f_r e_s(v_0) = -\beta_1^{-1} \bar{\phi}_{r+s} \cdot v_0$, we get $a_0 \bar{\phi}_r + a_1 \bar{\phi}_{r+1} + \cdots + a_l \bar{\phi}_{r+l} = 0$ for all $r \in \mathbb{Z}$. Therefore, the collection $\{\bar{\phi}_r\}_{r \in \mathbb{Z}}$ satisfies a simple recurrence relation. Solving this recurrence and using the condition $\bar{\phi}_0 = \phi_0 - \phi_0^{-1}$, we immediately see that $\phi^{\pm}(z)$ are expansions in $z^{\mp 1}$ of the same rational function.

To prove the "if" part, let us assume that $\phi^{\pm}(z) = \phi(z)^{\pm}$ for a rational function $\phi(z)$. Reversing the above arguments, we get $\dim(V_{\phi^+,\phi^-,\nu})_1 < \infty$. Combining this with the relation (t2), a simple induction implies that $\dim(V_{\phi^+,\phi^-,\nu})_l < \infty$. □

All the above admits natural opposite counterparts. We define the category $O^+$ as in Definition 2.3, but now ${}'\ddot{U}_{q_1,q_2,q_3}(\mathfrak{gl}_1)_{\geq N}(v) = 0$ in (i). We also define the highest weight ${}'\ddot{U}_{q_1,q_2,q_3}(\mathfrak{gl}_1)$-modules as in Definition 2.4, but now $e_r(v_0) = 0$ for all $r \in \mathbb{Z}$. Then, the analogue of Proposition 2.7 still holds. Moreover, we have natural equivalences $O^{\pm} \to O^{\mp}$ given by twists with the automorphism $\vartheta$ of ${}'\ddot{U}_{q_1,q_2,q_3}(\mathfrak{gl}_1)$:

$$\vartheta \colon e(z) \mapsto f(1/z), \ f(z) \mapsto e(1/z), \ \psi^{\pm}(z) \mapsto \psi^{\mp}(1/z), \ q^{d_2} \mapsto q^{-d_2}. \quad (2.40)$$

### 2.2.2  Vector, Fock, Macmahon modules, and their tensor products

In this section, we recall some important families of $'\ddot{U}_{q_1,q_2,q_3}(\mathfrak{gl}_1)$-modules in $O^-$ with nice combinatorial bases and explicit action. They arise via a semiinfinite construction from the simplest *vector representations* which however are not in $O^-$.

- *Tensor product construction*

    In what follows, we shall use the formal coproduct $\Delta$ on $'\ddot{U}_{q_1,q_2,q_3}(\mathfrak{gl}_1)$, cf. (H1):

$$\Delta\colon e(z) \mapsto e(z) \otimes 1 + \psi^-(z) \otimes e(z), \quad f(z) \mapsto 1 \otimes f(z) + f(z) \otimes \psi^+(z),$$
$$\psi^\pm(z) \mapsto \psi^\pm(z) \otimes \psi^\pm(z), \quad q^{d_2} \mapsto q^{d_2} \otimes q^{d_2}. \tag{2.41}$$

We note that one cannot naively use (2.41) to define the tensor product of two modules, since $\Delta(e_r)$ and $\Delta(f_r)$ contain infinite sums. However, for the modules in the present section "with general spectral parameters" this issue can be fixed as we explain now. To this end, let us recall the key property of the $\delta$-functions, cf. (1.37):

$$\gamma(z)\delta(u/z) = \gamma(u)\delta(u/z) \quad \text{for any} \quad \gamma(z) \in \mathbb{C}(z) \quad \text{regular at} \quad u \in \mathbb{C}^\times. \tag{2.42}$$

All the $'\ddot{U}_{q_1,q_2,q_3}(\mathfrak{gl}_1)$-modules considered below admit bases $\{w_\iota\}_{\iota \in \mathcal{I}}$ such that

$$e(z)(w_\iota) = \sum_{J \in \mathcal{I}} c_{\iota,J}\delta(\chi_{\iota,J}/z)w_J, \quad f(z)(w_\iota) = \sum_{J \in \mathcal{I}} d_{J,\iota}\delta(\chi_{J,\iota}/z)w_J \tag{2.43}$$

with some $\chi_{\iota,J}$ and $c_{\iota,J}, d_{J,\iota}$ (all but finitely many of which vanish for any $\iota \in \mathcal{I}$), and

$$\psi^\pm(z)(w_\iota) = \gamma_\iota(z)^\pm \cdot w_\iota \tag{2.44}$$

for some rational functions $\gamma_\iota(z)$. Let $V^{(1)}, V^{(2)}$ be two such modules, and we shall use the corresponding superscripts for $c_{\bullet,\bullet}$ etc. Then, applying (2.41), we get:

$$\Delta(e(z))(w_{\iota_1}^{(1)} \otimes w_{\iota_2}^{(2)}) = e(z)(w_{\iota_1}^{(1)}) \otimes w_{\iota_2}^{(2)} + \psi^-(z)(w_{\iota_1}^{(1)}) \otimes e(z)(w_{\iota_2}^{(2)}). \tag{2.45}$$

The first summand above is well-defined, and we express the second one using (2.42):

$$\psi^-(z)(w_{\iota_1}^{(1)}) \otimes e(z)(w_{\iota_2}^{(2)}) = \sum_{J_2 \in \mathcal{I}^{(2)}} c_{\iota_2,J_2}^{(2)}\gamma_{\iota_1}^{(1)}(\chi_{\iota_2,J_2}^{(2)})\delta(\chi_{\iota_2,J_2}^{(2)}/z) \cdot w_{\iota_1}^{(1)} \otimes w_{J_2}^{(2)}, \tag{2.46}$$

assuming that the rational function $\gamma_\bullet^{(1)}$ is regular at $\{\chi_{\iota_2,J_2}^{(2)}|c_{\iota_2,J_2}^{(2)} \neq 0\}$. The action of $f_r$'s on $V^{(1)} \otimes V^{(2)}$ is defined analogously via:

$$f(z)(w_{\iota_1}^{(1)} \otimes w_{\iota_2}^{(2)}) = w_{\iota_1}^{(1)} \otimes f(z)w_{\iota_2}^{(2)} + \sum_{J_1 \in \mathcal{I}^{(1)}} d_{J_1,\iota_1}^{(1)}\gamma_{\iota_2}^{(2)}(\chi_{J_1,\iota_1}^{(1)})\delta(\chi_{J_1,\iota_1}^{(1)}/z) \cdot w_{J_1}^{(1)} \otimes w_{\iota_2}^{(2)},$$
$$\tag{2.47}$$

assuming that the rational function $\gamma_\bullet^{(2)}$ is regular at $\{\chi_{J_1,\iota_1}^{(1)}|d_{J_1,\iota_1}^{(1)} \neq 0\}$. Finally, the action of $\psi_r$'s and $q^{d_2}$ on $V^{(1)} \otimes V^{(2)}$ is well-defined via (2.41), in particular:

$$\psi^{\pm}(z)(w_{t_1}^{(1)} \otimes w_{t_2}^{(2)}) = \gamma_{t_1}^{(1)}(z)^{\pm} \cdot \gamma_{t_2}^{(2)}(z)^{\pm} \cdot w_{t_1}^{(1)} \otimes w_{t_2}^{(2)} \,. \tag{2.48}$$

Thus, we shall always equip $V^{(1)} \otimes V^{(2)}$ with a $´\ddot{U}_{q_1, q_2, q_3}(\mathfrak{gl}_1)$-module structure through the formulas (2.45)–(2.48) whenever they are well-defined, and it will always be straightforward to check that this indeed defines an action.

● *Vector representations $V(u)$ and their tensor products*

The "building block" of all the constructions below is the family of *vector representations* $\{V(u)\}_{u \in \mathbb{C}^{\times}}$ with a basis parametrized by $\mathbb{Z}$, see [4, Proposition 3.1]:

**Proposition 2.8 (Vector representations)** *For $u \in \mathbb{C}^{\times}$, let $V(u)$ be a $\mathbb{C}$-vector space with a basis $\{[u]_r\}_{r \in \mathbb{Z}}$. The following formulas define a $´\ddot{U}_{q_1, q_2, q_3}(\mathfrak{gl}_1)$-action on it:*

$$\begin{aligned}
e(z)[u]_r &= (1 - q_1)^{-1} \delta(q_1^r u/z) \cdot [u]_{r+1} \,, \\
f(z)[u]_r &= (q_1^{-1} - 1)^{-1} \delta(q_1^{r-1} u/z) \cdot [u]_{r-1} \,, \\
q^{\pm d_2}[u]_r &= q^{\pm r} \cdot [u]_r \,, \\
\psi^{\pm}(z)[u]_r &= \left( \frac{(z - q_1^r q_2 u)(z - q_1^r q_3 u)}{(z - q_1^r u)(z - q_1^{r-1} u)} \right)^{\pm} \cdot [u]_r \,,
\end{aligned} \tag{2.49}$$

*where $\gamma(z)^{\pm}$ denotes the expansion of a rational function $\gamma(z)$ in $z^{\mp 1}$, respectively.*

*Proof* The proof is straightforward and is based on (2.42) and the following result:

**Lemma 2.5** *For any rational function $\gamma(z)$ with simple poles $\{x_t\} \subset \mathbb{C}^{\times}$ and possibly poles of higher order at $0, \infty$, the following equality holds:*

$$\gamma(z)^{+} - \gamma(z)^{-} = \sum_t \delta\left(\frac{z}{x_t}\right) \mathrm{Res}_{z=x_t} \gamma(z) \frac{dz}{z} \,. \tag{2.50}$$

*Proof* Consider the partial fraction decomposition of $\gamma(z)$:

$$\gamma(z) = P(z) + \sum_t \frac{\nu_t}{z - x_t} \quad \text{with} \quad P(z) \in \mathbb{C}[z, z^{-1}], \ \nu_t \in \mathbb{C} \,.$$

Then, $P(z)^{\pm} = P(z)$ and so $P(z)^{+} - P(z)^{-} = 0$. Meanwhile, we also have:

$$\left( \frac{\nu_t}{z - x_t} \right)^{+} = \frac{\nu_t}{z} + \frac{\nu_t x_t}{z^2} + \frac{\nu_t x_t^2}{z^3} + \dots \text{ and } \left( \frac{\nu_t}{z - x_t} \right)^{-} = -\frac{\nu_t}{x_t} - \frac{\nu_t z}{x_t^2} - \frac{\nu_t z^2}{x_t^3} - \dots$$

so that

$$\left( \frac{\nu_t}{z - x_t} \right)^{+} - \left( \frac{\nu_t}{z - x_t} \right)^{-} = \frac{\nu_t}{x_t} \delta\left( \frac{z}{x_t} \right) = \delta\left( \frac{z}{x_t} \right) \cdot \mathrm{Res}_{z=x_t} \frac{\nu_t}{z - x_t} \frac{dz}{z} \,.$$

This implies (2.50). □

Let $g^{\pm}(z,w) = (z - q_1^{\pm 1}w)(z - q_2^{\pm 1}w)(z - q_3^{\pm 1}w)$ so that $\frac{g^-(z,w)}{g^+(z,w)} = g(z/w)$ of (2.3). Then:

$$g^+(z,w)e(z)e(w)[u]_r = (1-q_1)^{-2}g^+(z,w)\delta(q_1^{r+1}u/z)\delta(q_1^r u/w)[u]_{r+2} = 0,$$

$$g^-(z,w)e(w)e(z)[u]_r = (1-q_1)^{-2}g^-(z,w)\delta(q_1^{r+1}u/w)\delta(q_1^r u/z)[u]_{r+2} = 0,$$

due to (2.42) and $g^{\pm}(q_1^{r\pm 1}u, q_1^r u) = 0$. This implies a compatibility of (2.49) with (t2). The compatibility of (2.49) with the defining relations (t3, t5, t6) is verified similarly.

To verify the compatibility with (t4), it suffices to note that

$$[e(z), f(w)][u]_r = \frac{-\left(\delta\left(\frac{q_1^{r-1}u}{z}\right)\delta\left(\frac{q_1^{r-1}u}{w}\right) - \delta\left(\frac{q_1^r u}{z}\right)\delta\left(\frac{q_1^r u}{w}\right)\right)}{(1-q_1)(1-q_1^{-1})} \cdot [u]_r =$$

$$\frac{1}{(1-q_1)(1-q_1^{-1})}\delta\left(\frac{z}{w}\right)\left(\delta\left(\frac{q_1^r u}{z}\right) - \delta\left(\frac{q_1^{r-1}u}{z}\right)\right) \cdot [u]_r,$$

while applying (2.50) we also get

$$(\psi^+(z) - \psi^-(z))[u]_r = \frac{(1-q_2)(1-q_3)}{1-q_1^{-1}}\left(\delta\left(\frac{z}{q_1^r u}\right) - \delta\left(\frac{z}{q_1^{r-1}u}\right)\right) \cdot [u]_r.$$

Finally, it is straightforward to check that $[e_0, [e_1, e_{-1}]][u]_r = 0 = [f_0, [f_1, f_{-1}]][u]_r$ for all $r$, hence, (t7′, t8′) hold. This completes our proof of the Proposition. □

*Remark 2.7* Replacing $q_1$ with $q_2$ or $q_3$ in the formula (2.49) also gives rise to a $\ddot{U}_{q_1,q_2,q_3}(\mathfrak{gl}_1)$-action on $V(u)$.

Since the formulas (2.49) are of the form (2.43, 2.44), we can now form tensor products of such $V(u_i)$'s with "$u_i$'s in general position" by applying (2.45)–(2.48):

**Lemma 2.6** *[4, Lemma 3.6] If $u_1,\ldots,u_N \in \mathbb{C}^{\times}$ satisfy $\frac{u_i}{u_j} \notin q_1^{\mathbb{Z}}$ for any $i < j$, then the formulas (2.45)–(2.48) define a $\ddot{U}_{q_1,q_2,q_3}(\mathfrak{gl}_1)$-action on $V(u_1) \otimes \cdots \otimes V(u_N)$.*

*Proof* The action of $e(z)$ in the basis $[u_1]_{r_1} \otimes \cdots \otimes [u_N]_{r_N}$ is given by:

$$(1-q_1)e(z)([u_1]_{r_1} \otimes [u_2]_{r_2} \otimes \cdots \otimes [u_N]_{r_N}) =$$

$$\sum_{l=1}^N \prod_{k=1}^{l-1} \gamma_{r_k,u_k}(q_1^{r_l}u_l)\delta\left(\frac{q_1^{r_l}u_l}{z}\right)[u_1]_{r_1} \otimes \cdots \otimes [u_l]_{r_l+1} \otimes \cdots \otimes [u_N]_{r_N} \quad (2.51)$$

in accordance with (2.45, 2.46, 2.48), where

$$\gamma_{r,u}(z) := \frac{(z - q_1^r q_2 u)(z - q_1^r q_3 u)}{(z - q_1^r u)(z - q_1^{r-1}u)}.$$

Here, we note that the above assumption $u_i/u_j \notin q_1^{\mathbb{Z}}$ for $i < j$ guarantees that all $\gamma_{r_k,u_k}(q_1^{r_l}u_l)$ with $k < l$ are well-defined. Likewise, the action of $f(z)$ is well-defined and is explicitly given by:

$$(q_1^{-1} - 1)f(z)([u_1]_{r_1} \otimes [u_2]_{r_2} \otimes \cdots \otimes [u_N]_{r_N}) =$$

$$\sum_{l=1}^{N} \prod_{k=l+1}^{N} \gamma_{r_k, u_k}(q_1^{r_l-1}u_l)\delta\left(\frac{q_1^{r_l-1}u_l}{z}\right)[u_1]_{r_1} \otimes \cdots \otimes [u_l]_{r_l-1} \otimes \cdots \otimes [u_N]_{r_N}.$$

The rest of the proof is straightforward. □

In particular, for generic parameters $q_1, q_3$, we obtain a $\ddot{U}_{q_1, q_2, q_3}(\mathfrak{gl}_1)$-action on

$$V^N(u) := V(u) \otimes V(q_2^{-1}u) \otimes V(q_2^{-2}u) \otimes \cdots \otimes V(q_2^{-N+1}u). \tag{2.52}$$

Let $W^N(u)$ be the subspace of $V^N(u)$ spanned by

$$|\lambda\rangle_u := [u]_{\lambda_1} \otimes [q_2^{-1}u]_{\lambda_2-1} \otimes [q_2^{-2}u]_{\lambda_3-2} \otimes \cdots \otimes [q_2^{-N+1}u]_{\lambda_N-N+1}, \qquad \lambda \in \mathcal{P}^N,$$

where $\mathcal{P}^N = \{\lambda = (\lambda_1, \ldots, \lambda_N) \in \mathbb{Z}^N | \lambda_1 \geq \cdots \geq \lambda_N\}$. We introduce the notation:

$$\lambda \pm 1_l := (\lambda_1, \ldots, \lambda_{l-1}, \lambda_l \pm 1, \lambda_{l+1}, \ldots, \lambda_N).$$

Then, we have:

**Lemma 2.7** $W^N(u)$ is a $\ddot{U}_{q_1, q_2, q_3}(\mathfrak{gl}_1)$-submodule of $V^N(u)$ for generic $q_1, q_3$ (2.39).

*Proof* Evoking the explicit formula (2.51), we obtain:

$$e(z)|\lambda\rangle = \sum_{l=1}^{N} \prod_{k=1}^{l-1} \frac{(1 - q_1^{\lambda_l - \lambda_k} q_3^{l-k-1})(1 - q_1^{\lambda_l - \lambda_k + 1} q_3^{l-k+1})}{(1 - q_1^{\lambda_l - \lambda_k} q_3^{l-k})(1 - q_1^{\lambda_l - \lambda_k + 1} q_3^{l-k})} \cdot \frac{\delta(q_1^{\lambda_l} q_3^{l-1} u/z)}{1 - q_1} \cdot |\lambda + 1_l\rangle. \tag{2.53}$$

We note that $|\lambda + 1_l\rangle \in \mathcal{P}^N$ unless $\lambda_{l-1} = \lambda_l$. But if $\lambda_{l-1} = \lambda_l$, then (2.53) does not contain $|\lambda + 1_l\rangle$ as $1 - q_1^{\lambda_l - \lambda_{l-1}} q_3^{l-(l-1)-1} = 0$. Thus, $e_r(W^N(u)) \subseteq W^N(u)$ for all $r$. The proof of $f_r(W^N(u)) \subseteq W^N(u)$ is analogous and is based on the explicit formula:

$$f(z)|\lambda\rangle = \sum_{l=1}^{N} \prod_{k=l+1}^{N} \frac{(1 - q_1^{\lambda_k - \lambda_l + 1} q_3^{k-l+1})(1 - q_1^{\lambda_k - \lambda_l} q_3^{k-l-1})}{(1 - q_1^{\lambda_k - \lambda_l + 1} q_3^{k-l})(1 - q_1^{\lambda_k - \lambda_l} q_3^{k-l})} \times$$

$$\frac{\delta(q_1^{\lambda_l - 1} q_3^{l-1} u/z)}{q_1^{-1} - 1} \cdot |\lambda - 1_l\rangle. \tag{2.54}$$

Let us finally record the action of $\psi^{\pm}(z)$ and $q^{d_2}$ that are diagonal in the basis $|\lambda\rangle$:

$$\psi^{\pm}(z)|\lambda\rangle = \left(\prod_{l=1}^{N} \frac{(z - q_1^{\lambda_l} q_3^l u)(z - q_1^{\lambda_l - 1} q_3^{l-2} u)}{(z - q_1^{\lambda_l} q_3^{l-1} u)(z - q_1^{\lambda_l - 1} q_3^{l-1} u)}\right)^{\pm} |\lambda\rangle, \quad q^{d_2}|\lambda\rangle = q^{\lambda_1 + \cdots + \lambda_N}|\lambda\rangle. \tag{2.55}$$

This completes our proof of the lemma. □

• *Fock representations $\mathcal{F}(u)$ and their tensor products*

While $V(u), W^N(u)$ are not in $O^-$, we shall now construct $\overset{'}{\ddot{U}}_{q_1,q_2,q_3}(\mathfrak{gl}_1)$-modules $\mathcal{F}(u) \in O^-$ as inductive limits of the subspaces (not submodules!) of $W^N(u)$.

*Remark 2.8* Following [4], we note that the formulas for the action of $\overset{'}{\ddot{U}}_{q_1,q_2,q_3}(\mathfrak{gl}_1)$ on $V(u)$ and $W^N(u)$ are structurally similar to the action of $\mathfrak{gl}_\infty$ on its vector representation $V$ and exterior powers $\Lambda^N V$, but with an additional "continuous" parameter $u$. Meanwhile, the modules $\mathcal{F}(u)$ will be reminiscent of the $\mathfrak{gl}_\infty$-representation $\Lambda^{\infty/2} V$.

Let $\mathcal{P}^{N,+} = \{\lambda = (\lambda_1, \ldots, \lambda_N) \in \mathcal{P}^N | \lambda_N \geq 0\}$ and define $W^{N,+}(u)$ as a subspace of $W^N(u)$ spanned by $|\lambda\rangle$ with $\lambda \in \mathcal{P}^{N,+}$. Then, the natural maps

$$\tau_N \colon \mathcal{P}^{N,+} \longrightarrow \mathcal{P}^{N+1,+} \quad \text{given by} \quad (\lambda_1, \ldots, \lambda_N) \mapsto (\lambda_1, \ldots, \lambda_N, 0), \qquad (2.56)$$

induce embeddings $W^{N,+}(u) \overset{\tau_N}{\hookrightarrow} W^{N+1,+}(u)$, and we consider the inductive limit

$$\mathcal{F}(u) := \lim_{N \to \infty} W^{N,+}(u). \qquad (2.57)$$

Thus, $\mathcal{F}(u)$ has a basis $|\lambda\rangle$ parametrized by all possible Young diagrams:

$$\mathcal{P}^+ = \left\{\lambda = (\lambda_1, \lambda_2, \ldots) \,|\, \lambda_1 \geq \lambda_2 \geq \cdots \quad \text{and} \quad \lambda_l = 0 \ \text{for} \ l \gg 0\right\}. \qquad (2.58)$$

To endow it with a $\overset{'}{\ddot{U}}_{q_1,q_2,q_3}(\mathfrak{gl}_1)$-action, we first modify the operators (2.53)–(2.55). Set

$$\phi(t) = \frac{q^{-1}t - q}{t - 1} \quad \text{where} \quad q_2 = q^2. \qquad (2.59)$$

Then, we consider the operators acting on $W^N(u)$:

$$e^{[N]}(z) = e(z), \ f^{[N]}(z) = \phi(z/q_3^N u) \cdot f(z), \ \psi^{\pm[N]}(z) = \phi(z/q_3^N u)^{\pm} \cdot \psi^{\pm}(z). \ (2.60)$$

This renormalization is chosen to satisfy the following "stabilization" property:

$$\tau_N(x^{[N]}(z)|\lambda\rangle) = x^{[N+1]}(z)\tau_N(|\lambda\rangle) \quad \text{for} \quad \lambda \in \mathcal{P}^{N,+} \ \text{s.t.} \ \lambda_N = 0, \qquad (2.61)$$

with $x = e, f, \psi^{\pm}$, see [4, Lemma 4.2]. Thus, we can finally define operators on $\mathcal{F}(u)$:

$$\tau_u(x(z))|\lambda\rangle = \lim_{N \to \infty} x^{[N]}(z)|\lambda_1, \ldots, \lambda_N\rangle \in \mathcal{F}(u) \quad \text{with} \quad x = e, f, \psi^{\pm}, \qquad (2.62)$$

for any $\lambda = (\lambda_1, \lambda_2, \ldots) \in \mathcal{P}^+$. The following is [4, Theorem 4.3, Corollary 4.4]:

**Proposition 2.9 (Fock representations)** *Assume that $q_1$ and $q_3$ are generic* (2.39).

*(a) The formula* (2.62) *does endow $\mathcal{F}(u)$ with a $\overset{'}{\ddot{U}}_{q_1,q_2,q_3}(\mathfrak{gl}_1)$-module structure.*

*(b) $\mathcal{F}(u)$ is an irreducible lowest weight $\overset{'}{\ddot{U}}_{q_1,q_2,q_3}(\mathfrak{gl}_1)$-module in the category $O^-$ with the lowest weight $\phi(z/u)$, and $|\emptyset\rangle = (0, 0, \ldots)$ is its lowest weight vector.*

*(c) Explicitly, $\overset{'}{\ddot{U}}_{q_1,q_2,q_3}(\mathfrak{gl}_1)$-action $\tau_u$ is given by the following matrix coefficients:*

$$\langle\lambda+1_l|e(z)|\lambda\rangle = \frac{1}{1-q_1}\prod_{k=1}^{l-1}\phi(q_1^{\lambda_k-\lambda_l-1}q_3^{k-l})\phi(q_1^{\lambda_l-\lambda_k}q_3^{l-k})\cdot\delta(q_1^{\lambda_l}q_3^{l-1}u/z),$$

$$\langle\lambda|f(z)|\lambda+1_l\rangle = \frac{1}{q_1^{-1}-1}\prod_{k=l+1}^{\infty}\phi(q_1^{\lambda_k-\lambda_l-1}q_3^{k-l})\phi(q_1^{\lambda_l-\lambda_k}q_3^{l-k})\cdot\delta(q_1^{\lambda_l}q_3^{l-1}u/z),$$

$$\langle\lambda|\psi^{\pm}(z)|\lambda\rangle = \left(\prod_{l=1}^{\infty}\phi(q_1^{\lambda_l-1}q_3^{l-1}u/z)\phi(q_1^{\lambda_l-1}q_3^{l-2}u/z)^{-1}\right)^{\pm}, \quad \langle\lambda|q^{\pm d_2}|\lambda\rangle = q^{\pm|\lambda|},$$

(2.63)

*where $|\lambda| = \sum_{l\geq 1}\lambda_l$ and all other matrix coefficients are zero.*

*Proof* Due to Lemma 2.7 and the renormalization (2.60), the only nontrivial check is the verification of (t4). The latter follows from Lemma 2.5 as the pole of $\phi(z/q_3^N u)$ is cancelled by zeros of the eigenvalues of $\psi(z)$ on $|\lambda\rangle$ when $\lambda_N=0$ (viewed as rational functions). The formulas in part (c) are direct consequences of the formulas (2.53)–(2.55). Finally, part (b) follows from the fact that $\{\psi_r\}_{r\in\mathbb{Z}}$ act diagonally in the basis $|\lambda\rangle$ with a simple joint spectrum as follows from (2.63), see [5, Lemma 2.3]. □

*Remark 2.9* (a) We note that: the infinite products in part (c) stabilize, the matrix coefficient $\langle\lambda+1_l|e(z)|\lambda\rangle$ vanishes if $\lambda$ is a Young diagram and $\lambda+1_l$ is not, the matrix coefficient $\langle\lambda|f(z)|\lambda+1_l\rangle$ vanishes if $\lambda+1_l$ is a Young diagram and $\lambda$ is not.

(b) If $q_1, q_3$ are not generic, then the formulas (2.63) are ill-defined. The particular $(k,r)$-*resonance* case "$q_1^a q_3^b = 1$ $(a, b\in\mathbb{Z})\Leftrightarrow a=(1-r)c, b=(k+1)c$ for some $c\in\mathbb{Z}$" (with $k\geq 1, r\geq 2$) was considered in [4, §6].

Since the formulas (2.63) are of the form (2.43), (2.44), we can now form tensor products of such $\mathcal{F}(u_i)$'s with "$u_i$'s in general position" by applying (2.45)–(2.48):

**Lemma 2.8** *[5, Lemma 3.1] If $q_1, q_3, u_1,\ldots,u_N\in\mathbb{C}^{\times}$ are generic in the sense*

$$q_1^a q_3^b u_1^{c_1}\cdots u_N^{c_N} = 1\ (a, b, c_i\in\mathbb{Z}) \Longleftrightarrow a=b=c_1=\ldots=c_N=0, \quad (2.64)$$

*then the formulas (2.45)–(2.48) define a $'\ddot{U}_{q_1,q_2,q_3}(\mathfrak{gl}_1)$-action on $\mathcal{F}(u_1)\otimes\cdots\otimes\mathcal{F}(u_N)$. Moreover, it is an irreducible lowest weight $'\ddot{U}_{q_1,q_2,q_3}(\mathfrak{gl}_1)$-module in the category $O^-$ with the lowest weight $\prod_{i=1}^N\phi(z/u_i)$, and $|\emptyset\rangle\otimes\cdots\otimes|\emptyset\rangle$ is its lowest weight vector.*

*Proof* The well-definedness of (2.45)–(2.47) is guaranteed by (2.64). The second part again follows from the observation that $\{\psi_r\}_{r\in\mathbb{Z}}$ act diagonally in the basis $|\lambda^1\rangle\otimes\cdots\otimes|\lambda^N\rangle$ with a simple joint spectrum, due to (2.48, 2.63, 2.64). □

*Remark 2.10* For non-generic $q_1, q_3, u_1,\ldots,u_N$ such an action of $'\ddot{U}_{q_1,q_2,q_3}(\mathfrak{gl}_1)$ on $\mathcal{F}(u_1)\otimes\cdots\otimes\mathcal{F}(u_N)$ may be ill-defined. The particular resonance case of $q_1, q_3, u_1$ being generic while $u_k = u_{k+1}q_1^{a_k+1}q_3^{b_k+1}$ for $1\leq k<N$ and $a_k, b_k\in\mathbb{N}$ was considered in [5, §3.2]. As shown in *loc.cit.* there is still an action of $'\ddot{U}_{q_1,q_2,q_3}(\mathfrak{gl}_1)$ on

$$\text{Span}\left\{|\lambda^1\rangle\otimes\cdots\otimes|\lambda^N\rangle\ \Big|\ \lambda_l^k\geq\lambda_{l+b_k}^{k+1}-a_k\ \text{ for }\ 1\leq k<N\text{ and }l\geq 1\right\}. \quad (2.65)$$

If moreover $q_1, q_3$ are not generic, then this action is ill-defined, but one gets an action on a subspace similar to (2.65) but with $1 \leq k \leq N$, see [5, §3.3] for more details.

● *Macmahon representations $\mathcal{M}(u, K)$*

Following Remark 2.10 with all $a_k = b_k = 0$, one obtains a ${}'\ddot{U}_{q_1, q_2, q_3}(\mathfrak{gl}_1)$-action on the subspace of $\mathcal{F}(u) \otimes \mathcal{F}(q_2 u) \otimes \cdots \otimes \mathcal{F}(q_2^{N-1} u)$ spanned by $|\lambda^1\rangle \otimes \cdots \otimes |\lambda^N\rangle$ satisfying $\lambda_l^k \geq \lambda_l^{k+1}$ for all $1 \leq k < N$ and $l \geq 1$. Considering the $N \to \infty$ limit once again, one can construct ${}'\ddot{U}_{q_1, q_2, q_3}(\mathfrak{gl}_1)$-actions on the vector space with a basis labeled by plane partitions $\bar{\lambda}$:

$$\bar{\lambda} = (\lambda^1, \lambda^2, \ldots): \quad \lambda^k \in \mathcal{P}^+, \quad \lambda_l^k \geq \lambda_l^{k+1} \; \forall k, l > 0, \quad \lambda^k = \emptyset \text{ for } k \gg 0, \quad (2.66)$$

that depend on an extra parameter $K \in \mathbb{C}^\times$. The resulting modules $\mathcal{M}(u, K)$ are called *(vacuum) Macmahon modules*. As we will not presently need explicit formulas for the action, we rather refer the reader to [7, §3.1] for more details. We shall only need:

**Proposition 2.10** *[7, §3.1] Assuming that $q_1, q_3$ are generic (2.39) and $K^2 \notin q_1^{\mathbb{Z}} q_3^{\mathbb{Z}}$, the module $\mathcal{M}(u, K)$ can be invariantly described as the lowest weight ${}'\ddot{U}_{q_1, q_2, q_3}(\mathfrak{gl}_1)$-module in the category $\mathcal{O}^-$ with the lowest weight $\phi^K(z/u) := \frac{K^{-1}z/u - K}{z/u - 1}$. In the above realization, the empty plane partition $\bar{\emptyset} = (\emptyset, \emptyset, \ldots)$ is its lowest weight vector.*

### 2.2.3 Vertex representations and their relation to Fock modules

In this section, we recall an important family of vertex-type ${}'\ddot{U}_{q_1, q_2, q_3}(\mathfrak{gl}_1)$-modules. To this end, we consider a Heisenberg algebra $\mathfrak{A}$ generated by $\{H_k\}_{k \in \mathbb{Z}}$ subject to:

$$[H_0, H_k] = 0, \qquad [H_k, H_l] = k \frac{1 - q_1^{|k|}}{1 - q_3^{-|k|}} \delta_{k, -l} H_0 \quad (k, l \neq 0). \qquad (2.67)$$

Let $\mathfrak{A}^{\geq}$ be the subalgebra generated by $\{H_k\}_{k \geq 0}$, and $\mathbb{C}v_0$ be the $\mathfrak{A}^{\geq}$-representation with $H_k v_0 = \delta_{k0} v_0$. The induced representation $F := \text{Ind}_{\mathfrak{A}^{\geq}}^{\mathfrak{A}}(\mathbb{C}v_0)$ is called the *Fock representation* of $\mathfrak{A}$. We also define the degree operator $\mathsf{d}$ acting on $F$ via

$$\mathsf{d}(H_{-k_1} H_{-k_2} \cdots H_{-k_N} v_0) = (k_1 + k_2 + \cdots + k_N) \cdot H_{-k_1} H_{-k_2} \cdots H_{-k_N} v_0.$$

The following result provides a natural structure of a ${}'\ddot{U}_{q_1, q_2, q_3}(\mathfrak{gl}_1)$-module on $F$:

**Proposition 2.11** *[6, Proposition A.6] For $c, u \in \mathbb{C}^\times$, the following formulas define an action of ${}'\ddot{U}_{q_1, q_2, q_3}(\mathfrak{gl}_1)$ on $F$:*

$$\rho_{u,c}(e(z)) =$$

$$\frac{-q_1 q_3 c}{(1 - q_1)^2 (1 - q_3)^2} \cdot \exp\left(\sum_{k > 0} \frac{1 - q_3^k}{k} H_{-k} \left(\frac{z}{u}\right)^k\right) \exp\left(-\sum_{k > 0} \frac{1 - q_3^{-k}}{k} H_k \left(\frac{z}{u}\right)^{-k}\right),$$

$$\rho_{u,c}(f(z)) =$$

$$c^{-1} \cdot \exp\left(-\sum_{k>0} \frac{1-q_3^k}{k} q_2^{k/2} H_{-k} \left(\frac{z}{u}\right)^k\right) \exp\left(\sum_{k>0} \frac{1-q_3^{-k}}{k} q_2^{k/2} H_k \left(\frac{z}{u}\right)^{-k}\right),$$

$$\rho_{u,c}(\psi^\pm(z)) = \exp\left(\mp \sum_{k>0} \frac{1-q_3^{\mp k}}{k} (1-q_2^k) q_2^{-k/4} H_{\pm k} \left(\frac{z}{u}\right)^{\mp k}\right),$$

$$\rho_{u,c}(\gamma^{\pm 1/2}) = q_2^{\pm 1/4}, \qquad \rho_{u,c}(q^{\pm d_1}) = q^{\pm d}.$$

*Remark 2.11* This module is irreducible since so is its restriction to $\ddot{U}'^{,0}_{q_1,q_2,q_3}$. To this end, we note the following compatibility between the relations (t1′) and (2.67):

$$\frac{-\beta_k^2}{k^2} \cdot \frac{-k(\gamma^k - \gamma^{-k})_{|\gamma \mapsto q_2^{1/2}}}{\beta_k} = \frac{-(1-q_3^{-k})(1-q_3^k)(1-q_2^k)^2 q_2^{-k/2}}{k^2} \cdot \frac{k(1-q_1^k)}{1-q_3^{-k}}.$$

Let us now relate the above vertex-type $\ddot{U}'_{q_1,q_2,q_3}(\mathfrak{gl}_1)$-modules $\rho_{v,c}$ to the Fock $'\ddot{U}_{q_1,q_2,q_3}(\mathfrak{gl}_1)$-modules $\tau_u$. To this end, we recall that $\tau_u$ is the lowest weight module in $O^-$ with the lowest weight $\phi(z/u)$. Let $\tau^\vartheta_{1/q_2 u}$ denote the $\vartheta$-twisted action of $'\ddot{U}_{q_1,q_2,q_3}(\mathfrak{gl}_1)$ on $\mathcal{F}(1/q_2 u)$ given by $\tau^\vartheta_{1/q_2 u}(x) = \tau_{1/q_2 u}(\vartheta(x))$, with the automorphism $\vartheta \colon '\ddot{U}_{q_1,q_2,q_3}(\mathfrak{gl}_1) \overset{\sim}{\longrightarrow} '\ddot{U}_{q_1,q_2,q_3}(\mathfrak{gl}_1)$ of (2.40). As pointed out in the end of Section 2.2.1, $\tau^\vartheta_{1/q_2 u}$ is a highest weight module in the category $O^+$ and its highest weight is $\phi(q^2 u/z) = \phi(z/u)^{-1}$. We also consider $\rho^\varpi_{v,c}$, the $\varpi$-twisted action of $\ddot{U}_{q_1,q_2,q_3}(\mathfrak{gl}_1)$ on $F$ given by $\rho^\varpi_{v,c}(x) = \rho_{v,c}(\varpi(x))$, with the isomorphism $\varpi \colon '\ddot{U}_{q_1,q_2,q_3}(\mathfrak{gl}_1) \overset{\sim}{\longrightarrow} \ddot{U}_{q_1,q_2,q_3}(\mathfrak{gl}_1)$ of (2.18).

Then, we have the following result:

**Theorem 2.5** *For any $v, c \in \mathbb{C}^\times$, we have a $'\ddot{U}_{q_1,q_2,q_3}(\mathfrak{gl}_1)$-module isomorphism:*

$$\rho^\varpi_{v,c} \simeq \tau^\vartheta_{1/q_2 u}, \tag{2.68}$$

*where $u = -(1-q_1)(1-q_3)/c$.*

*Proof* Since both $\rho_{v,c}$ and $\tau_{1/q_2 u}$ are irreducible modules, so are $'\ddot{U}_{q_1,q_2,q_3}(\mathfrak{gl}_1)$-modules $\rho^\varpi_{v,c}$ and $\tau^\vartheta_{1/q_2 u}$. Thus, it suffices to verify that $v_0 \in \rho^\varpi_{v,c}$ and $|\emptyset\rangle \in \tau^\vartheta_{1/q_2 u}$ are the highest weight vectors of the same highest weight.

As noted right before the theorem, $\tau^\vartheta_{1/q_2 u}$ is a highest weight module with the highest weight $\phi(z/u)^{-1}$ and the highest weight vector $|\emptyset\rangle$. Therefore:

$$\tau^\vartheta_{1/q_2 u}(e_k)|\emptyset\rangle = 0, \quad \tau^\vartheta_{1/q_2 u}(\psi_{\pm r})|\emptyset\rangle = \pm(q-q^{-1})(q^2 u)^{\pm r}|\emptyset\rangle \quad \forall k \in \mathbb{Z}, r > 0.$$

On the other hand, we also have $\rho^\varpi_{v,c}(e_k)v_0 = 0$ for any $k \in \mathbb{Z}$, which follows from Theorem 2.2(b) as all eigenvalues of the degree operator $d$ on $F$ are nonnegative. Thus, it just remains to prove:

$$\rho^\varpi_{v,c}(\psi_{\pm r})(v_0) = \pm(q-q^{-1})(q^2 u)^{\pm r} v_0 \quad \forall r > 0. \tag{2.69}$$

According to (t4), we have:

$$\rho^{\varpi}_{v,c}(\psi_{\pm r})(v_0) = \pm\beta_1 \cdot [\rho^{\varpi}_{v,c}(e_{\pm r}), \rho^{\varpi}_{v,c}(f_0)](v_0) = \pm\beta_1 \cdot \rho^{\varpi}_{v,c}(e_{\pm r})(\rho^{\varpi}_{v,c}(f_0)(v_0)) \tag{2.70}$$

for all $r > 0$. Since $\varpi(f_0) = h_{-1}\gamma^{1/2} = \beta_1^{-1}\psi_{-1}\psi_0\gamma^{1/2}$ by (2.14), we get:

$$\rho^{\varpi}_{v,c}(f_0)(v_0) = \frac{1}{v(1-q_1)}H_{-1}v_0 \,. \tag{2.71}$$

On the other hand, we have $\varpi(e_0) = -h_1\gamma^{-1/2} = \beta_1^{-1}\psi_1\psi_0^{-1}\gamma^{-1/2}$ by (2.14), hence:

$$\rho^{\varpi}_{v,c}(e_0)\,(H_{-1}v_0) = -\frac{v}{q(1-q_3)}v_0 \,. \tag{2.72}$$

We can now evaluate all $\rho^{\varpi}_{v,c}(e_{\pm r})(H_{-1}v_0)$ for $r > 0$ by applying (t5′):

$$\rho^{\varpi}_{v,c}(e_{\pm(r+1)})(H_{-1}v_0) = [\rho^{\varpi}_{v,c}(h_{\pm 1}), \rho^{\varpi}_{v,c}(e_{\pm r})](H_{-1}v_0) \,. \tag{2.73}$$

To this end, we note that $\varpi(h_1) = -f_0$ and $\varpi(h_{-1}) = e_0$ by (2.14), hence we have:

$$\begin{aligned}
\rho^{\varpi}_{v,c}(h_1)(v_0) &= -c^{-1}v_0 \,, \\
\rho^{\varpi}_{v,c}(h_{-1})(v_0) &= -\frac{q_1 q_3 c}{(1-q_1)^2(1-q_3)^2}v_0 \,, \\
\rho^{\varpi}_{v,c}(h_1)(H_{-1}v_0) &= -c^{-1}\left(1 - q_2(1-q_1)(1-q_3)\right)H_{-1}v_0 \,, \\
\rho^{\varpi}_{v,c}(h_{-1})(H_{-1}v_0) &= -\frac{q_1 q_3 c}{(1-q_1)^2(1-q_3)^2}\left(1 - (1-q_1)(1-q_3)\right)H_{-1}v_0 \,.
\end{aligned} \tag{2.74}$$

Combining (2.72)–(2.74), we thus get by induction on $r$:

$$\rho^{\varpi}_{v,c}(e_{\pm r})(H_{-1}v_0) = -\frac{v}{q(1-q_3)} \cdot \left(-q_2(1-q_1)(1-q_3)/c\right)^{\pm r}v_0 \quad \forall\, r > 0 \,. \tag{2.75}$$

Finally, combining (2.70, 2.71, 2.75), we obtain (2.69) with $u = -(1-q_1)(1-q_3)/c$. □

Following [10], we note that the vertex representations $\rho_{v,c}$ allow to interpret the elements $K_k$ of (2.37) generating the commutative subalgebra $\mathcal{A}$ and their "elliptic" generalizations via the transfer matrices as explained in Section 3.5, see Remark 3.17.

*Remark 2.12* The above observation together with the shuffle realization of the Fock modules (as discussed in the next subsection) played the crucial role in [8].

*Remark 2.13* Analogously to Theorem 2.5, one can show that $\tau_{-(1-q_1)(1-q_3)/c} \simeq \rho^{*,\varpi}_{v,c}$. Here, $\rho^{*,\varpi}_{v,c}$ is the $\varpi$-twist of the $\ddot{U}'_{q_1,q_2,q_3}(\mathfrak{gl}_1)$-module $\rho^*_{v,c}$, dual to $\rho_{v,c}$, which can be naturally realized on $F$ by similar vertex-type formulas: $\rho^*_{v,c}(e(z)) = \frac{q_1 q_3 c}{(1-q_1)^2(1-q_3)^2} \cdot$ $\exp(-\sum_{k>0} \frac{(1-q_3^{-k})q_2^{k/2}}{k}H_{-k}\left(\frac{z}{v}\right)^{-k})\exp(\sum_{k>0} \frac{(1-q_3^k)q_2^{k/2}}{k}H_k\left(\frac{z}{v}\right)^k)$, $\rho^*_{v,c}(f(z)) = -c^{-1} \cdot$ $\exp(\sum_{k>0} \frac{1-q_3^{-k}}{k}H_{-k}\left(\frac{z}{v}\right)^{-k})\exp(-\sum_{k>0} \frac{1-q_3^k}{k}H_k\left(\frac{z}{v}\right)^k)$, $\rho^*_{v,c}(\gamma^{\pm 1/2}) = q_2^{\mp 1/4}$, as well as $\rho^*_{v,c}(q^{\pm d_1}) = q^{\mp d}$, $\rho^*_{v,c}(\psi^{\pm}(z)) = \exp(\pm\sum_{k>0} \frac{(1-q_3^{\mp k})(1-q_2^k)q_2^{-k/4}}{k}H_{\mp k}\left(\frac{z}{v}\right)^{\mp k})$.

### 2.2.4 Shuffle bimodules and their relation to Fock modules

Following [8], we recall the shuffle realization of the Fock modules and their tensor products by constructing natural $S$-bimodules that become ${}'\ddot{U}_{q_1,q_2,q_3}(\mathfrak{gl}_1)$-modules.

For $u \in \mathbb{C}^\times$, consider an $\mathbb{N}$-graded $\mathbb{C}$-vector space $S_1(u) = \bigoplus_{k \in \mathbb{N}} S_1(u)_k$, where the degree $k$ component $S_1(u)_k$ consists of $\Sigma_k$-symmetric rational functions $F$ in the variables $\{x_r\}_{r=1}^k$ satisfying the following three conditions:

(i) *Pole conditions* on $F$:

$$F = \frac{f(x_1,\dots,x_k)}{\prod_{r \neq r'}(x_r - x_{r'}) \cdot \prod_{r=1}^k (x_r - u)}, \quad f \in \mathbb{C}\left[x_1^{\pm 1},\dots,x_k^{\pm 1}\right]^{\Sigma_k}. \quad (2.76)$$

(ii) *First kind wheel conditions* on $F$:

$$F(x_1,\dots,x_k) = 0 \quad \text{once} \quad \left\{ \frac{x_{r_1}}{x_{r_2}}, \frac{x_{r_2}}{x_{r_3}}, \frac{x_{r_3}}{x_{r_1}} \right\} = \{q_1, q_2, q_3\} \quad \text{for some} \quad r_1, r_2, r_3.$$

(iii) *Second kind wheel conditions* on $F$ (with $f$ from (2.76)):

$$f(x_1,\dots,x_k) = 0 \quad \text{once} \quad x_{r_1} = u \quad \text{and} \quad x_{r_2} = q_2 u \quad \text{for some} \quad 1 \leq r_1 \neq r_2 \leq k.$$

Given $F \in S_k$ and $G \in S_1(u)_\ell$, we define $F \star G, G \star F \in S_1(u)_{k+\ell}$ via:

$$(F \star G)(x_1,\dots,x_{k+\ell}) := \frac{1}{k! \cdot \ell!} \times$$
$$\text{Sym}\left( F(x_1,\dots,x_k)G(x_{k+1},\dots,x_{k+\ell}) \prod_{r \leq k}^{r' > k} \zeta\left(\frac{x_r}{x_{r'}}\right) \prod_{r=1}^k \phi\left(\frac{x_r}{u}\right) \right), \quad (2.77)$$

$$(G \star F)(x_1,\dots,x_{k+\ell}) := \frac{1}{k! \cdot \ell!} \times$$
$$\text{Sym}\left( G(x_1,\dots,x_\ell)F(x_{\ell+1},\dots,x_{k+\ell}) \prod_{r \leq \ell}^{r' > \ell} \zeta\left(\frac{x_r}{x_{r'}}\right) \right). \quad (2.78)$$

These formulas clearly endow $S_1(u)$ with an $S$-bimodule structure. Identifying $S$ with ${}'\ddot{U}_{q_1,q_2,q_3}^>$ via (2.25), we thus get two commuting ${}'\ddot{U}_{q_1,q_2,q_3}^>$-actions on $S_1(u)$.

In fact, the left action can be extended to an action of the entire ${}'\ddot{U}_{q_1,q_2,q_3}(\mathfrak{gl}_1)$:

**Proposition 2.12** *[8, Proposition 3.2] The following formulas define a* ${}'\ddot{U}_{q_1,q_2,q_3}(\mathfrak{gl}_1)$-*action on $S_1(u)$:*

$$\pi_u(\psi_0)G = q^{-1} \cdot G, \qquad \pi_u(q^{d_2})G = q^n \cdot G, \qquad \pi_u(e_k)G = x^k \star G,$$

$$\pi_u(h_l)G = \left( \sum_{r=1}^n x_r^l - \frac{u^l}{(1-q_1^l)(1-q_3^l)} \right) \cdot G,$$

$$\pi_u(f_k)G = \frac{n}{\beta_1} \left( \operatorname*{Res}_{z=0} + \operatorname*{Res}_{z=\infty} \right) \left\{ \frac{z^k G(x_1, \ldots, x_{n-1}, z)}{\prod_{r=1}^{n-1} \zeta(\frac{x_r}{z})} \frac{dz}{z} \right\},$$

*for $k \in \mathbb{Z}$, $l \neq 0$, and $G \in S_1(u)_n$.*

To relate this to the Fock modules $\mathcal{F}(u)$, let us consider the following subspace:

$$J_0 = \operatorname{span}_{\mathbb{C}} \{ G \star F \mid G \in S_1(u), F \in S_k \text{ with } k \geq 1 \} \subset S_1(u).$$

**Proposition 2.13** *[8, Proposition 3.3] (a) $J_0$ is a $'\ddot{U}_{q_1,q_2,q_3}(\mathfrak{gl}_1)$-submodule of $\pi_u$.*

*(b) We have a $'\ddot{U}_{q_1,q_2,q_3}(\mathfrak{gl}_1)$-module isomorphism $S_1(u)/J_0 \simeq \mathcal{F}(u)$.*

The above admits a *higher rank* counterpart. For any $m \geq 1$ and generic $\underline{u} = (u_1, \ldots, u_m) \in (\mathbb{C}^\times)^m$ satisfying (2.64), so that $\mathcal{F}(u_1) \otimes \cdots \otimes \mathcal{F}(u_m)$ is irreducible by Lemma 2.8, define $S(\underline{u}) = \bigoplus_{k \in \mathbb{N}} S(\underline{u})_k$ similarly to $S_1(u)$ with the following changes:

(i′) *Pole conditions* for $F$ from a degree $k$ component should read as follows:

$$F = \frac{f(x_1, \ldots, x_k)}{\prod_{r \neq r'}(x_r - x_{r'}) \cdot \prod_{s=1}^m \prod_{r=1}^k (x_r - u_s)}, \quad f \in \mathbb{C}\left[ x_1^{\pm 1}, \ldots, x_k^{\pm 1} \right]^{\Sigma_k}.$$

(iii′) *Second kind wheel conditions* for such $F$ should read as follows:

$$f(x_1, \ldots, x_k) = 0 \text{ once } x_{r_1} = u_s \text{ and } x_{r_2} = q_2 u_s \text{ for some } s, 1 \leq r_1 \neq r_2 \leq k.$$

We endow $S(\underline{u})$ with an $S$-bimodule structure by the formulas (2.77, 2.78) with

$$\prod_{1 \leq r \leq k} \phi(x_r/u) \rightsquigarrow \prod_{1 \leq s \leq m} \prod_{1 \leq r \leq k} \phi(x_r/u_s).$$

The resulting left action of $'\ddot{U}_{q_1,q_2,q_3}^{>} \simeq S$ on $S(\underline{u})$ can be extended to the action of the entire $'\ddot{U}_{q_1,q_2,q_3}(\mathfrak{gl}_1)$, denoted by $\pi_{\underline{u}}$, cf. Section 3.3.3. Consider the subspace

$$J_0 = \operatorname{span}_{\mathbb{C}} \{ G \star F \mid G \in S(\underline{u}), F \in S_k \text{ with } k \geq 1 \} \subset S(\underline{u}).$$

**Proposition 2.14** *[8, Proposition 6.1] $J_0$ is a $'\ddot{U}_{q_1,q_2,q_3}(\mathfrak{gl}_1)$-submodule of $\pi_{\underline{u}}$. Moreover, we have a $'\ddot{U}_{q_1,q_2,q_3}(\mathfrak{gl}_1)$-module isomorphism $S(\underline{u})/J_0 \simeq \mathcal{F}(u_1) \otimes \cdots \otimes \mathcal{F}(u_m)$.*

Yet another natural generalization of $S_1(u)$ is provided by the $S$-bimodules $S_1^K(u)$. As a vector space, $S_1^K(u)$ is defined similarly to $S_1(u)$ but without imposing the second kind wheel conditions. The $S$-bimodule structure on $S_1^K(u)$ is defined by the formulas (2.77, 2.78) with the following modification in (2.77):

$$\phi(t) \rightsquigarrow \phi^K(t) := (K^{-1} \cdot t - K)/(t - 1). \tag{2.79}$$

The resulting left action of $'\ddot{U}_{q_1,q_2,q_3}^{>} \simeq S$ on $S_1^K(u)$ can be extended to the action of the entire $'\ddot{U}_{q_1,q_2,q_3}(\mathfrak{gl}_1)$, denoted by $\pi_u^K$. Consider the subspace

$$J_0 = \operatorname{span}_{\mathbb{C}} \{ G \star F \mid G \in S_1^K(u), F \in S_k \text{ with } k \geq 1 \} \subset S_1^K(u).$$

**Proposition 2.15** *[22] $J_0$ is a $'\ddot{U}_{q_1,q_2,q_3}(\mathfrak{gl}_1)$-submodule of $\pi_u^K$. Moreover, we have a $'\ddot{U}_{q_1,q_2,q_3}(\mathfrak{gl}_1)$-module isomorphism $S_1^K(u)/J_0 \simeq \mathcal{M}(u,K)$ if $K^2 \notin q_1^{\mathbb{Z}} q_3^{\mathbb{Z}}$.*

One can naturally define the *higher rank* counterparts $S^{\underline{K}}(\underline{u})$ of $S_1^K(u)$ similarly to $S(\underline{u})$. Then, the analogue of the Proposition above provides a shuffle realization of the tensor product $\mathcal{M}(u_1,K_1) \otimes \cdots \otimes \mathcal{M}(u_N,K_N)$ of vacuum Macmahon modules (for general parameters $\underline{u} = (u_1,\ldots,u_N)$ and $\underline{K} = (K_1,\ldots,K_N)$).

## 2.3 Geometric realizations

In this section, we provide some flavour of the geometric applications of the quantum toroidal and shuffle algebras by recalling the realization of the Fock modules $\mathcal{F}(u)$ and their tensor products through the geometry of the Hilbert schemes of points on the plane $(\mathbb{A}^2)^{[n]}$ and more generally the Gieseker varieties $M(r,n)$. We also discuss how the commutative subalgebra $\mathcal{A} \subset S$ from Section 2.1.6 and its "negative" version $\mathcal{A}^{\mathrm{op}} \subset S^{\mathrm{op}}$ give rise to an action of the Heisenberg algebra on the localized equivariant $K$-theory of $(\mathbb{A}^2)^{[n]}$, providing a $K$-theoretic counterpart of the famous Nakajima's construction from [15]. Our exposition closely follows that of [9, 21].

### 2.3.1 Correspondences and fixed points for $(\mathbb{A}^2)^{[n]}$

Throughout this section $X = \mathbb{A}^2$. Let $X^{[n]}$ be the Hilbert scheme of $n$ points on $X$. Its $\mathbb{C}$-points are the codimension $n$ ideals $J \subset \mathbb{C}[x,y]$. Let $P[k] \subset \bigsqcup_n X^{[n]} \times X^{[n+k]}$ be the Nakajima-Grojnowski correspondence. For $k > 0$, $P[k] = \bigsqcup_n P[k]_n$ with

$$P[k]_n = \left\{ (J_1,J_2) \in X^{[n]} \times X^{[n+k]} \,\middle|\, J_2 \subset J_1 \,,\, \mathrm{supp}(J_1/J_2) \text{ is a single point} \right\}. \quad (2.80)$$

It is known that $P[1]$ is a smooth variety. It is equipped with two natural projections:

$$X^{[n]} \xleftarrow{\mathbf{p}} P[1]_n \xrightarrow{\mathbf{q}} X^{[n+1]},$$
$$J_1 \leftarrowtail (J_1,J_2) \mapsto J_2. \quad (2.81)$$

Let $L$ be the *tautological* line bundle on $P[1]$ whose fiber at $(J_1,J_2)$ equals $J_1/J_2$.

Consider a natural action of a torus $\mathbb{T} = \mathbb{C}^{\times} \times \mathbb{C}^{\times}$ on each $X^{[n]}$, induced from the action on $X$ given by the formula $(t_1,t_2) \cdot (x,y) = (t_1 x, t_2 y)$. The set $(X^{[n]})^{\mathbb{T}}$ of $\mathbb{T}$-fixed points in $X^{[n]}$ is finite and is in bijection with size $n$ Young diagrams. For such a Young diagram $\lambda = (\lambda_1,\ldots,\lambda_l)$, the corresponding ideal $J_\lambda \in (X^{[n]})^{\mathbb{T}}$ is given by:

$$J_\lambda = \mathbb{C}[x,y] \cdot \left( \mathbb{C}x^{\lambda_1} y^0 \oplus \cdots \oplus \mathbb{C}x^{\lambda_l} y^{l-1} \oplus \mathbb{C}y^l \right).$$

For a Young diagram $\lambda$, let $\lambda'$ denote its transpose and $|\lambda| = \lambda_1 + \cdots + \lambda_l$, as before.

### 2.3.2  Geometric action I

Let $'M$ be a direct sum of equivariant (complexified) $K$-groups: $'M = \bigoplus_n K^{\mathbb{T}}(X^{[n]})$. It is a module over $K^{\mathbb{T}}(\mathrm{pt}) = \mathbb{C}[t_1^{\pm 1}, t_2^{\pm 1}]$, the equivariant (complexified) $K$-group of a point. Set $t_3 := t_1^{-1} t_2^{-1} \in \mathbb{C}(t_1, t_2)$ and consider $\mathbb{F} = \mathbb{C}(t_1, t_2)[t]/(t^2 - t_3)$. We define

$$M := {}'M \otimes_{K^{\mathbb{T}}(\mathrm{pt})} \mathbb{F}. \tag{2.82}$$

It has a natural $\mathbb{N}$-grading:

$$M = \bigoplus_{n \in \mathbb{N}} M_n, \quad M_n = K^{\mathbb{T}}(X^{[n]}) \otimes_{K^{\mathbb{T}}(\mathrm{pt})} \mathbb{F}.$$

According to the Thomason localization theorem, restriction to the $\mathbb{T}$-fixed point set induces an isomorphism

$$K^{\mathbb{T}}(X^{[n]}) \otimes_{K^{\mathbb{T}}(\mathrm{pt})} \mathbb{F} \xrightarrow{\sim} K^{\mathbb{T}}((X^{[n]})^{\mathbb{T}}) \otimes_{K^{\mathbb{T}}(\mathrm{pt})} \mathbb{F}.$$

The structure sheaves of $\mathbb{T}$-fixed points $J_\lambda$ form a basis in $\bigoplus_n K^{\mathbb{T}}((X^{[n]})^{\mathbb{T}}) \otimes_{K^{\mathbb{T}}(\mathrm{pt})} \mathbb{F}$, denoted by $\{\lambda\}$. Since embedding of a point $J_\lambda$ into $X^{[|\lambda|]}$ is a proper morphism, the direct image in the equivariant $K$-theory is well-defined, and we denote by $[\lambda] \in M_{|\lambda|}$ the direct image of the structure sheaf $\{\lambda\}$. Thus, the set $\{[\lambda]\}$ is an $\mathbb{F}$-basis of $M$.

Let $\mathfrak{F}$ be the *tautological* vector bundle on $X^{[n]}$ whose fiber at a point $J \in X^{[n]}$ equals the quotient $\mathbb{C}[x, y]/J$. Define the generating series $\mathbf{a}(z), \mathbf{c}(z) \in M(z)$ via:

$$\mathbf{a}(z) := \Lambda^\bullet_{-1/z}(\mathfrak{F}) = \sum_{k \geq 0} [\Lambda^k(\mathfrak{F})](-1/z)^k, \tag{2.83}$$

$$\mathbf{c}(z) := \frac{\mathbf{a}(zt_1)\mathbf{a}(zt_2)\mathbf{a}(zt_3)}{\mathbf{a}(zt_1^{-1})\mathbf{a}(zt_2^{-1})\mathbf{a}(zt_3^{-1})}. \tag{2.84}$$

We also define the linear operators $e_k, f_k, \psi_k, \psi_0^{-1}, d_2 \ (k \in \mathbb{Z})$ acting on $M$ via:

$$e_k = \mathbf{q}_*(L^{\otimes k} \otimes \mathbf{p}^*) \colon M_n \longrightarrow M_{n+1}, \tag{2.85}$$

$$f_k = -t^{-1} \cdot \mathbf{p}_*(L^{\otimes(k-1)} \otimes \mathbf{q}^*) \colon M_n \longrightarrow M_{n-1}, \tag{2.86}$$

$$\psi^\pm(z) = \psi_0^{\pm 1} + \sum_{k=1}^\infty \psi_{\pm k} z^{\mp k} := \left( \frac{t^{-1}z - t}{z - 1} \mathbf{c}(z) \right)^\pm \in \prod_n \mathrm{End}(M_n)[[z^{\mp 1}]], \tag{2.87}$$

$$d_2 = n \cdot \mathrm{Id} \colon M_n \longrightarrow M_n, \tag{2.88}$$

where $\gamma(z)^\pm$ denotes the expansion of a rational function $\gamma(z)$ in $z^{\mp 1}$, respectively. The following is the key result of this section:

**Theorem 2.6** *The operators (2.85)–(2.88) give rise to an action of $'\ddot{U}_{q_1, q_2, q_3}(\mathfrak{gl}_1)$ on $M$ from (2.82) with the parameters $q_1 = t_1, q_2 = t_3, q_3 = t_2$.*

This theorem was proved independently and simultaneously in [9] and [19]. We refer the interested reader to [9, §4] for a straightforward verification of the defining relations that crucially uses the properties (2.42, 2.50) of the $\delta$-functions.

In fact, this representation recovers the previously constructed Fock module:

**Proposition 2.16** *We have a* $'\ddot{U}_{q_1,q_2,q_3}(\mathfrak{gl}_1)$ *-module isomorphism* $M \simeq \mathcal{F}(1)$.

*Proof* Such an isomorphism is provided by the assignment $[\lambda] \mapsto c_\lambda |\lambda\rangle$ for a unique choice of nonzero scalars $c_\lambda$ with $c_\emptyset = 1$, see [4, Corollary 4.5]. □

### 2.3.3 Heisenberg algebra action on the equivariant $K$-theory

Following [9], let us now construct an action of a Heisenberg algebra on $M$ of (2.82), whose cohomological version [11, 15] utilizes higher order correspondences $P[k]$.

To this end, we first recall the notion of Macdonald polynomials following [12]. The algebra $\Lambda$ of symmetric polynomials in the variables $\{x_k\}_{k=1}^{\infty}$ over $F = \mathbb{C}(q,t)$ is freely generated by the power-sum symmetric polynomials $p_r = \sum_k x_k^r$, that is, $\Lambda = F[p_1, p_2, \ldots]$. For a Young diagram $\lambda = (\lambda_1, \ldots, \lambda_l) = (1^{m_1} 2^{m_2} \cdots)$, we define:

$$p_\lambda := p_{\lambda_1} \cdots p_{\lambda_l} \in \Lambda, \qquad z_\lambda := \prod_{r \geq 1} r^{m_r} m_r! \in \mathbb{N}.$$

Consider the Macdonald inner product $(\cdot, \cdot)_{q,t}$ on $\Lambda$ such that:

$$(p_\lambda, p_\mu)_{q,t} = \delta_{\lambda\mu} z_\lambda \prod_{1 \leq r \leq l} \frac{1 - q^{\lambda_r}}{1 - t^{\lambda_r}}.$$

**Definition 2.5** The Macdonald polynomials $\{P_\lambda\}_{\lambda \in \mathcal{P}^+} \subset \Lambda$ are defined as the unique family of polynomials such that:

- $P_\lambda = m_\lambda + \sum_{\mu < \lambda} \eta_{\lambda\mu} m_\mu$ with some $\eta_{\lambda\mu} \in F$ and $<$ being the dominance order,
- $(P_\lambda, P_\mu)_{q,t} = 0$ if $\lambda \neq \mu$.

Identifying $q \leftrightarrow t_1$ and $t \leftrightarrow t_2^{-1}$, we define $\Lambda_{\mathbb{F}} := \Lambda \otimes_F \mathbb{F}$. As both $\Lambda_{\mathbb{F}}$ and $M$ have bases parametrized by Young diagrams, we consider a vector space isomorphism

$$\Theta \colon \Lambda_{\mathbb{F}} \xrightarrow{\sim} M \quad \text{given by} \quad P_\lambda \mapsto \langle \lambda \rangle := c_\lambda \cdot [\lambda] \tag{2.89}$$

with the constants $c_\lambda$ defined via:

$$c_\lambda := (-(1-t_2)^{-1} t_2)^{-|\lambda|} t_1^{\sum_r (\lambda_r (\lambda_r - 1))/2} \prod_{\square \in \lambda} \left(1 - t_1^{l(\square)} t_2^{-a(\square)-1}\right)^{-1}. \tag{2.90}$$

Here, for a box $\square = (a, b)$ located in the $a$-th row and $b$-th column of $\lambda$, we define $l(\square) := \lambda_a - b$, $a(\square) := \max\{k | \lambda_k \geq b\} - a$ as the standard *leg* and *arm* functions. The specific choice (2.90) is determined by the $r = 1$ case of the following result:

**Theorem 2.7** *[9, Theorem 6.1] The isomorphism* (2.89) *intertwines the multiplication by the $r$-th elementary symmetric polynomial $e_r$ on $\Lambda_{\overline{\mathbb{F}}}$ with the action of $\widetilde{K}_r$ on $M$, where $\widetilde{K}_r(x_1, \ldots, x_r) = t_1^{-r(r-1)/2}(1-t_1)^r(1-t_2)^r \cdot \prod_{1 \leq k \neq k' \leq r} \frac{x_k - t_1 x_{k'}}{x_k - x_{k'}} \in \mathcal{A} \cap S_r$.*

*Proof* This is a straightforward computation based on matching the matrix coefficients of $\widetilde{K}_r$ in the basis $\langle\lambda\rangle$ with the Pieri formulas of [12, §VI.6] for the multiplication by $e_r$ in the basis of Macdonald polynomials $P_\lambda$. □

Evoking the well-known equality of generating functions

$$1 + \sum_{r \geq 1} e_r z^r = \exp\left(\sum_{r > 0} \frac{(-1)^{r-1}}{r} p_r z^r\right),$$

we introduce $\{\mathfrak{h}_r\}_{r>0} \subset S \simeq {}'\ddot{U}^{>}_{t_1, t_3, t_2}$ via $1 + \sum_{r \geq 1} \widetilde{K}_r z^r = \exp\left(\sum_{r>0} \frac{(-1)^{r-1}}{r} \mathfrak{h}_r z^r\right)$. We also define $\{\mathfrak{h}_{-r}\}_{r>0} \subset S^{\mathrm{op}} \simeq {}'\ddot{U}^{<}_{t_1, t_3, t_2}$ by the same formula but perceiving $\widetilde{K}_r$ as elements of $S^{\mathrm{op}} \simeq {}'\ddot{U}^{<}_{t_1, t_3, t_2}$. Then, the elements $\{\mathfrak{h}_r\}_{r \neq 0}$ give rise to the desired action of a Heisenberg algebra on $M$, which is naturally identified with the Fock module over this Heisenberg algebra.

### 2.3.4 Correspondences and fixed points for $M(r, n)$

The purpose of this section is to provide *higher rank* counterparts of the previous results by replacing $(\mathbb{A}^2)^{[n]}$ with the Gieseker moduli spaces $M(r, n)$.

First, let us recall some basics on the Gieseker framed moduli spaces $M(r, n)$ of torsion free sheaves on $\mathbb{P}^2$ of rank $r$ (with the second Chern class) $c_2 = n$. Its $\mathbb{C}$-points are the isomorphism classes of pairs $\{(E, \Phi)\}$, where $E$ is a torsion free sheaf on $\mathbb{P}^2$ of rank $r$ with $c_2(E) = n$, and which is locally free in a neighborhood of the line $l_\infty = \{(0 : z_1 : z_2)\} \subset \mathbb{P}^2$, while $\Phi$ is a trivialization at this line $\Phi: E|_{l_\infty} \xrightarrow{\sim} O^{\oplus r}_{l_\infty}$.

This space has an alternative quiver description (see [15, Ch. 2] for details):

$$M(r, n) = \mathcal{M}(r, n)/GL_n(\mathbb{C}), \quad \mathcal{M}(r, n) = \left\{(B_1, B_2, i, j)\,\middle|\,[B_1, B_2] + ij = 0\right\}^s,$$

where $B_1, B_2 \in \mathrm{End}(\mathbb{C}^n), i \in \mathrm{Hom}(\mathbb{C}^r, \mathbb{C}^n), j \in \mathrm{Hom}(\mathbb{C}^n, \mathbb{C}^r)$, the group $GL_n(\mathbb{C})$ acts via $g \cdot (B_1, B_2, i, j) = (gB_1 g^{-1}, gB_2 g^{-1}, gi, jg^{-1})$, and the superscript $s$ denotes the stability condition "there is no subspace $S \subsetneq \mathbb{C}^n$ such that $B_1(S), B_2(S), \mathrm{Im}(i) \subseteq S$". Let $\mathfrak{F}_r$ be the *tautological* rank $n$ vector bundle on $M(r, n)$.

Consider a natural action of a torus $\mathbb{T}_r = (\mathbb{C}^\times)^2 \times (\mathbb{C}^\times)^r$ on $M(r, n)$, where $(\mathbb{C}^\times)^2$ acts on $\mathbb{P}^2$ via $(t_1, t_2) \cdot ([z_0 : z_1 : z_2]) = [z_0 : t_1 z_1 : t_2 z_2]$, while $(\mathbb{C}^\times)^r$ acts by rescaling the trivialization $\Phi$. The set $M(r, n)^{\mathbb{T}_r}$ of $\mathbb{T}_r$-fixed points in $M(r, n)$ is finite and is in bijection with $r$-tuples of Young diagrams $\underline{\lambda} = (\lambda^1, \ldots, \lambda^r)$ satisfying $|\underline{\lambda}| := |\lambda^1| + \cdots + |\lambda^r| = n$, denoted by $\underline{\lambda} \vdash n$ (see [16, Proposition 2.9]). For such $\underline{\lambda}$, the corresponding point $J_{\underline{\lambda}} \in M(r, n)^{\mathbb{T}_r}$ is represented by $(E_{\underline{\lambda}}, \Phi_{\underline{\lambda}})$, where $E_{\underline{\lambda}} = J_{\lambda^1} \oplus \cdots \oplus J_{\lambda^r}$ and $\Phi_{\underline{\lambda}}$ is a sum of natural inclusions $J_{\lambda^k}|_{l_\infty} \hookrightarrow O_{l_\infty}$ $(1 \leq k \leq r)$.

Following [14, Section 5], we recall the *Hecke correspondences* generalizing $P[1]$ for $r = 1$. Consider $\mathcal{M}(r; n, n + 1) \subset \mathcal{M}(r, n) \times \mathcal{M}(r, n + 1)$ consisting of pairs of tuples $\{(B_1^{(k)}, B_2^{(k)}, i^{(k)}, j^{(k)})\}_{k=n,n+1}$ such that there exists $\varrho \colon \mathbb{C}^{n+1} \to \mathbb{C}^n$ satisfying

$$\varrho B_1^{(n+1)} = B_1^{(n)} \varrho, \quad \varrho B_2^{(n+1)} = B_2^{(n)} \varrho, \quad \varrho i^{(n+1)} = i^{(n)}, \quad j^{(n+1)} = j^{(n)} \varrho.$$

The stability condition implies the surjectivity of $\varrho$, hence, $S := \mathrm{Ker}(\varrho) \subset \mathbb{C}^{n+1}$ is a one-dimensional subspace of $\mathrm{Ker}(j^{(n+1)})$ invariant with respect to $\{B_a^{(n+1)}\}_{a=1,2}$. This identifies $\mathcal{M}(r; n, n+1)$ with pairs of $(B_1^{(n+1)}, B_2^{(n+1)}, i^{(n+1)}, j^{(n+1)}) \in \mathcal{M}(r, n+1)$ and a one-dimensional subspace $S \subset \mathbb{C}^{n+1}$ satisfying the above conditions. Taking the latter viewpoint, the Hecke correspondence $M(r; n, n+1)$ is defined as the quotient:

$$M(r; n, n + 1) = \mathcal{M}(r; n, n + 1)/GL_{n+1}(\mathbb{C}).$$

It is equipped with two natural projections

$$M(r, n) \xleftarrow{\mathbf{p}_r} M(r; n, n + 1) \xrightarrow{\mathbf{q}_r} M(r, n + 1) \tag{2.91}$$

and the *tautological* line bundle $L_r$. The set $M(r; n, n+1)^{\mathbb{T}_r}$ of $\mathbb{T}_r$-fixed points in $M(r; n, n+1)$ is in bijection with pairs of $r$-tuples of Young diagrams $\underline{\lambda} \vdash n, \underline{\mu} \vdash n + 1$ such that $\lambda^k \subseteq \mu^k$ for $1 \leq k \leq r$.

### 2.3.5 Geometric action II

Let $'M^r$ be the direct sum of equivariant (complexified) $K$-groups: $'M^r = \bigoplus_n K^{\mathbb{T}_r}(M(r, n))$. It is a module over $K^{\mathbb{T}_r}(\mathrm{pt}) = \mathbb{C}[\mathbb{T}_r] = \mathbb{C}[t_1^{\pm 1}, t_2^{\pm 1}, \chi_1^{\pm 1}, \dots, \chi_r^{\pm 1}]$. Consider $\mathbb{F}_r = \mathrm{Frac}(K^{\mathbb{T}_r}(\mathrm{pt}))[t]/(t^2 - t_3)$ with $t_3 := t_1^{-1} t_2^{-1}$, and define

$$M^r := {}'M^r \otimes_{K^{\mathbb{T}_r}(\mathrm{pt})} \mathbb{F}_r.$$

It has a natural $\mathbb{N}$-grading:

$$M^r = \bigoplus_{n \in \mathbb{N}} M_n^r, \quad M_n^r = K^{\mathbb{T}_r}(M(r, n)) \otimes_{K^{\mathbb{T}_r}(\mathrm{pt})} \mathbb{F}_r.$$

According to the Thomason localization theorem, restriction to the $\mathbb{T}_r$-fixed point set induces an isomorphism

$$K^{\mathbb{T}_r}(M(r, n)) \otimes_{K^{\mathbb{T}_r}(\mathrm{pt})} \mathbb{F}_r \xrightarrow{\sim} K^{\mathbb{T}_r}(M(r, n)^{\mathbb{T}_r}) \otimes_{K^{\mathbb{T}_r}(\mathrm{pt})} \mathbb{F}_r.$$

Again, the structure sheaves $\{\underline{\lambda}\}$ of the $\mathbb{T}_r$-fixed points $J_{\underline{\lambda}}$ form a basis in $\bigoplus_n K^{\mathbb{T}_r}(M(r, n)^{\mathbb{T}_r}) \otimes_{K^{\mathbb{T}_r}(\mathrm{pt})} \mathbb{F}_r$. Since embedding of a point $J_{\underline{\lambda}}$ into $M(r, |\underline{\lambda}|)$ is a proper morphism, the direct image in the equivariant $K$-theory is well-defined, and we denote by $[\underline{\lambda}] \in M_{|\underline{\lambda}|}^r$ the direct image of $\{\underline{\lambda}\}$. Thus, the set $\{[\underline{\lambda}]\}$ is a basis of $M^r$.

Generalizing (2.83, 2.84), define the generating series $\mathbf{a}_r(z), \mathbf{c}_r(z) \in M^r(z)$ via:

$$\mathbf{a}_r(z) := \Lambda^{\bullet}_{-1/z}(\mathfrak{F}_r) = \sum_{k \geq 0} [\Lambda^k(\mathfrak{F}_r)](-1/z)^k \,, \tag{2.92}$$

$$\mathbf{c}_r(z) := \frac{\mathbf{a}_r(zt_1)\mathbf{a}_r(zt_2)\mathbf{a}_r(zt_3)}{\mathbf{a}_r(zt_1^{-1})\mathbf{a}_r(zt_2^{-1})\mathbf{a}_r(zt_3^{-1})} \,. \tag{2.93}$$

We also define the linear operators $e_k, f_k, \psi_k, \psi_0^{\pm 1}, d_2$ ($k \in \mathbb{Z}$) acting on $M^r$ via:

$$e_k = \mathbf{q}_{r*}(L_r^{\otimes k} \otimes \mathbf{p}_r^*) \colon M_n^r \longrightarrow M_{n+1}^r, \tag{2.94}$$

$$f_k = (-t)^{r-2} \chi_1^{-1} \cdots \chi_r^{-1} \cdot \mathbf{p}_{r*}(L_r^{\otimes (k-r)} \otimes \mathbf{q}_r^*) \colon M_n^r \longrightarrow M_{n-1}^r, \tag{2.95}$$

$$\psi^{\pm}(z) = \psi_0^{\pm 1} + \sum_{k=1}^{\infty} \psi_{\pm k} z^{\mp k} := \left( \prod_{a=1}^r \frac{t^{-1}z - t\chi_a^{-1}}{z - \chi_a^{-1}} \cdot \mathbf{c}_r(z) \right)^{\pm} \in \prod_n \mathrm{End}(M_n^r)[[z^{\mp 1}]], \tag{2.96}$$

$$d_2 = n \cdot \mathrm{Id} \colon M_n^r \longrightarrow M_n^r \,. \tag{2.97}$$

The following is the key result of this section:

**Theorem 2.8** *The operators (2.94)–(2.97) give rise to an action of $\ddot{U}_{q_1,q_2,q_3}(\mathfrak{gl}_1)$ on $M^r$ with the parameters $q_1 = t_1, q_2 = t_3, q_3 = t_2$.*

This higher rank generalization of Theorem 2.6 first appeared in [19] and can be proved similarly to Theorem 2.6 as presented in [9]. We refer the interested reader to [21, Appendix B] for such a straightforward proof, which only uses the properties (2.42, 2.50) of the $\delta$-functions as well as the following explicit formulas for the matrix coefficients of $e_p, f_p, \psi^{\pm}(z)$ in the fixed point basis $\{[\underline{\lambda}]\}$ of $M^r$:

**Lemma 2.9** *Consider the fixed point basis $\{[\underline{\lambda}]\}$ of $M^r$. Define $\chi_k^{(a)} := t_1^{\lambda_k^a - 1} t_2^{k-1} \chi_a^{-1}$.*

*(a) The only nonzero matrix coefficients of the operators $e_p, f_p$ are as follows:*

$$\langle \underline{\lambda} | e_p | \underline{\lambda} - 1_l^b \rangle = \frac{(\chi_l^{(b)})^p}{1 - t_1 t_2} \cdot \prod_{a=1}^r \prod_{k=1}^{\infty} \frac{1 - t_1 t_2 \chi_k^{(a)}/\chi_l^{(b)}}{1 - t_1 \chi_k^{(a)}/\chi_l^{(b)}} \,,$$

$$\langle \underline{\lambda} | f_p | \underline{\lambda} + 1_l^b \rangle = (-t)^{r-2} \chi_1^{-1} \cdots \chi_r^{-1} \cdot \frac{(t_1 \chi_l^{(b)})^{p-r}}{1 - t_1 t_2} \cdot \prod_{a=1}^r \prod_{k=1}^{\infty} \frac{1 - t_1 t_2 \chi_l^{(b)}/\chi_k^{(a)}}{1 - t_1 \chi_l^{(b)}/\chi_k^{(a)}} \,,$$

*where $\underline{\lambda} \pm 1_l^b$ denotes the $r$-tuple of diagrams $(\lambda^1, \ldots, \lambda^{b-1}, \lambda^b \pm 1_l, \lambda^{b+1}, \ldots, \lambda^r)$.*

*(b) $\psi^{\pm}(z)$ is diagonal in the fixed point basis with the eigenvalue on $[\underline{\lambda}]$ given by:*

$$\langle \underline{\lambda} | \psi^{\pm}(z) | \underline{\lambda} \rangle =$$

$$\left( \prod_{a=1}^r \frac{t^{-1} - t\chi_a^{-1}/z}{1 - \chi_a^{-1}/z} \prod_{a=1}^r \prod_{\square \in \lambda^a} \frac{(1 - t_1^{-1}\chi(\square)/z)(1 - t_2^{-1}\chi(\square)/z)(1 - t_3^{-1}\chi(\square)/z)}{(1 - t_1\chi(\square)/z)(1 - t_2\chi(\square)/z)(1 - t_3\chi(\square)/z)} \right)^{\pm} ,$$

*where* $\chi(\square^a_{k,l}) := t_1^{k-1} t_2^{l-1} \chi_a^{-1}$ *for a box* $\square^a_{k,l}$ *in the l-th row and k-th column of* $\lambda^a$.

In fact, this representation recovers one of the previously constructed ones:

**Proposition 2.17** *We have a* $'\ddot{U}_{q_1,q_2,q_3}(\mathfrak{gl}_1)$*-module isomorphism:*

$$M^r \simeq \mathcal{F}(\chi_1^{-1}) \otimes \cdots \otimes \mathcal{F}(\chi_r^{-1}). \tag{2.98}$$

*Proof* Such an isomorphism is provided by the assignment

$$[\underline{\lambda}] = [(\lambda^1, \ldots, \lambda^r)] \mapsto c_{\underline{\lambda}} \cdot |\lambda^1\rangle \otimes \cdots \otimes |\lambda^r\rangle$$

for a unique choice of nonzero scalars $c_{\underline{\lambda}}$ with $c_{(\emptyset,\ldots,\emptyset)} = 1$, see [21, Theorem 4.6]. □

### 2.3.6 Whittaker vector

We conclude this section with a brief discussion of an important property of the sum

$$\mathfrak{w}_r := \sum_{n \geq 0} \left[ O_{M(r,n)} \right] \in \widehat{M}^r \tag{2.99}$$

called the *Whittaker vector*, which is an element of the completion $\widehat{M}^r := \prod_{n=1}^{\infty} M_n^r$.

To this end, we consider the following elements $\{K_{l,n}\}_{n \in \mathbb{N}}^{l \in \mathbb{Z}}$ of $S$, cf. (2.37):

$$K_{l,n}(x_1, \ldots, x_n) := \prod_{1 \leq a < b \leq n} \frac{(x_a - q_1 x_b)(x_b - q_1 x_a)}{(x_a - x_b)^2} \prod_{1 \leq a \leq n} x_a^l. \tag{2.100}$$

Let $\{K_{l,-n}\}_{n \in \mathbb{N}}^{l \in \mathbb{Z}}$ be the analogous elements in the opposite algebra $S^{\mathrm{op}} \simeq {}'\ddot{U}_{q_1,q_2,q_3}^<$. The above terminology is motivated by the following result of [21, Theorem 8.5] (cf. [20, §9.2] for a cohomological counterpart):

**Theorem 2.9** $\mathfrak{w}_r$ *of (2.99) is an eigenvector with respect to* $\{K_{l,-n} | 0 \leq l \leq r, n > 0\}$:

$$K_{l,-n}(\mathfrak{w}_r) = C_{l,-n} \cdot \mathfrak{w}_r,$$

*where*

$$\begin{aligned}
C_{0,-n} &= \frac{(-1)^n q^{(r-2)n}(-q_1)^{n(n-1)/2}}{(1-q_1)^n(1-q_3)(1-q_3^2)\cdots(1-q_3^n)}, \\
C_{1,-n} &= \ldots = C_{r-1,-n} = 0, \\
C_{r,-n} &= \frac{(-q)^{(r-2)n}(-q_1 q_3)^{n(n-1)/2}}{(1-q_1)^n(1-q_3)(1-q_3^2)\cdots(1-q_3^n)\cdot(\chi_1\cdots\chi_r)^n}.
\end{aligned} \tag{2.101}$$

*Remark 2.14* Evoking Theorems 2.1 and 2.3, this result states that the vector $\mathfrak{w}_r \in \widehat{M}^r$ is a common eigenvector of the subalgebra of $\widetilde{\mathcal{E}}$ generated by $\{u_{-k,d}\}_{k \geq 1}^{0 \leq d \leq rk}$.

# References

1. I. Burban, O. Schiffmann: On the Hall algebra of an elliptic curve, I. Duke Math. J. **161**, no. 7, 1171–1231 (2012)

2. V. Chari, A. Pressley: Quantum affine algebras. Commun. Math. Phys. **142**, no. 2, 261–283 (1991)

3. J. Ding, K. Iohara: Generalization of Drinfeld quantum affine algebras. Lett. Math. Phys. **41**, no. 2, 181–193 (1997)

4. B. Feigin, E. Feigin, M. Jimbo, T. Miwa, E. Mukhin: Quantum continuous $\mathfrak{gl}_\infty$: semiinfinite construction of representations. Kyoto J. Math. **51**, no. 2, 337–364 (2011)

5. B. Feigin, E. Feigin, M. Jimbo, T. Miwa, E. Mukhin: Quantum continuous $\mathfrak{gl}_\infty$: tensor products of Fock modules and $\mathcal{W}_n$-characters. Kyoto J. Math. **51**, no. 2, 365–392 (2011)

6. B. Feigin, K. Hashizume, A. Hoshino, J. Shiraishi, S. Yanagida: A commutative algebra on degenerate $\mathbb{CP}^1$ and Macdonald polynomials. J. Math. Phys. **50**, no. 9, Paper No. 095215 (2009)

7. B. Feigin, M. Jimbo, T. Miwa, E. Mukhin: Quantum toroidal $\mathfrak{gl}_1$-algebra: Plane partitions. Kyoto J. Math. **52**, no. 3, 621–659 (2012)

8. B. Feigin, M. Jimbo, T. Miwa, E. Mukhin: Quantum toroidal $\mathfrak{gl}_1$ and Bethe ansatz. J. Phys. A **48**, no. 24, Paper No. 244001 (2015)

9. B. Feigin, A. Tsymbaliuk: Equivariant $K$-theory of Hilbert schemes via shuffle algebra. Kyoto J. Math. **51**, no. 4, 831–854 (2011)

10. B. Feigin, A. Tsymbaliuk: Bethe subalgebras of $U_q(\widehat{\mathfrak{gl}_n})$ via shuffle algebras. Selecta Math. (N.S.) **22**, no. 2, 979–1011 (2016)

11. W. Li, Z. Qin, W. Wang: The cohomology rings of Hilbert schemes via Jack polynomials. In Algebraic Structures and Moduli Spaces, CRM Proc. Lecture Notes **38**, Amer. Math. Soc., Providence, 249–258 (2004)

12. I. Macdonald: Symmetric Functions and Hall Polynomials, 2nd ed. Oxford Math. Monogr., Oxford Univ. Press, New York (1995)

13. K. Miki: A $(q, \gamma)$ analog of the $W_{1+\infty}$ algebra. J. Math. Phys. **48**, no. 12, Paper No. 123520 (2007)

14. H. Nakajima: Quiver varieties and Kac-Moody algebras. Duke Math. J. **91**, no. 3, 515–560 (1998)

15. H. Nakajima: Lectures on Hilbert schemes of points on surfaces. Univ. Lecture Ser. **18**, Amer. Math. Soc., Providence (1999)

16. H. Nakajima, K. Yoshioka: Instanton counting on blowup. I. 4-dimensional pure gauge theory. Invent. Math. **162**, no. 2, 313–355 (2005)

17. A. Neguţ: The shuffle algebra revisited. Int. Math. Res. Not. IMRN, no. 22, 6242–6275 (2014)

18. O. Schiffmann: Drinfeld realization of the elliptic Hall algebra. J. Algebraic Combin. **35**, no. 2, 237–262 (2012)

19. O. Schiffmann, E. Vasserot: The elliptic Hall algebra and the $K$-theory of the Hilbert scheme of $\mathbb{A}^2$. Duke Math. J. **162**, no. 2, 279–366 (2013)

20. O. Schiffmann, E. Vasserot: Cherednik algebras, $W$-algebras and the equivariant cohomology of the moduli space of instantons on $\mathbb{A}^2$. Publ. Math. Inst. Hautes Études Sci. **118**, 213– 342 (2013)

21. A. Tsymbaliuk: The affine Yangian of $\mathfrak{gl}_1$ revisited. Adv. Math. **304**, 583–645 (2017)

22. A. Tsymbaliuk: Several realizations of Fock modules for toroidal $\ddot{U}_{q,d}(\mathfrak{sl}_n)$. Algebr. Represent. Theory **22**, no. 1, 177–209 (2019)

# Chapter 3
# Quantum toroidal $\mathfrak{sl}_n$, its representations, and Bethe subalgebras

**Abstract** In this chapter, we generalize the results of Chapter 2 in the context of the quantum toroidal algebras of $\mathfrak{sl}_n$ (for $n \geq 3$). Besides for presenting the natural counterparts of the main constructions, we also recall the shuffle realizations of:

(1) the "positive" subalgebra of the quantum toroidal algebra $\ddot{U}_{q,d}(\mathfrak{sl}_n)$, due to [20];
(2) the various combinatorial modules, such as Fock and Macmahon, due to [23];
(3) the commutative (Bethe) subalgebras of $\ddot{U}_{q,d}(\mathfrak{sl}_n)$, due to [12].

## 3.1 Quantum toroidal algebras of $\mathfrak{sl}_n$

In this section, we recall the notion of the quantum toroidal algebras of $\mathfrak{sl}_n$ ($n \geq 3$) as well as their basic algebraic properties that will be relevant for the upcoming results. In particular, we evoke Miki's isomorphism that interchanges vertical and horizontal quantum affine subalgebras and is analogous to (2.18) for the quantum toroidal $\mathfrak{gl}_1$.

### 3.1.1 Quantum toroidal $\mathfrak{sl}_n$ ($n \geq 3$)

For an integer $n \geq 3$, consider index sets $[n] := \{0, 1, \ldots, n-1\}$ (also viewed as a set of residues modulo $n$) and $[n]^{\times} := [n] \backslash \{0\}$ (naturally identified with $I$ of Section 1.1). We define two matrices $(a_{ij})_{i,j \in [n]}$ (the Cartan matrix of $\widehat{\mathfrak{sl}}_n$) and $(m_{ij})_{i,j \in [n]}$ via:

$$a_{ii} = 2, \quad a_{i,i\pm 1} = -1, \quad m_{i,i\pm 1} = \mp 1, \quad \text{and} \quad a_{ij} = 0 = m_{ij} \quad \text{otherwise}. \quad (3.1)$$

We fix $q, d \in \mathbb{C}^{\times}$ which, similarly to (2.1), can be encoded by three parameters:

$$q_1 := q^{-1}d, \quad q_2 := q^2, \quad q_3 := q^{-1}d^{-1} \qquad \text{satisfying} \qquad q_1 q_2 q_3 = 1. \quad (3.2)$$

For $m \in \mathbb{Z}$, we also consider the rational function $g_m(z) \in \mathbb{C}(z)$ defined by:

$$g_m(z) := \frac{q^m z - 1}{z - q^m} \, . \tag{3.3}$$

Following [19], we define the quantum toroidal algebra of $\mathfrak{sl}_n$, denoted by $\ddot{U}_{q,d}(\mathfrak{sl}_n)$, to be the associative $\mathbb{C}$-algebra generated by $\{e_{i,r}, f_{i,r}, \psi_{i,r}, \psi_{i,0}^{-1}, \gamma^{\pm 1/2}, q^{\pm d_1}, q^{\pm d_2}\}_{i \in [n]}^{r \in \mathbb{Z}}$ with the following defining relations:

$$[\psi_i^\pm(z), \psi_j^\pm(w)] = 0, \quad \gamma^{\pm 1/2} - \text{central}, \tag{T0.1}$$

$$\psi_{i,0}^{\pm 1} \cdot \psi_{i,0}^{\mp 1} = \gamma^{\pm 1/2} \cdot \gamma^{\mp 1/2} = q^{\pm d_1} \cdot q^{\mp d_1} = q^{\pm d_2} \cdot q^{\mp d_2} = 1, \tag{T0.2}$$

$$q^{d_1} e_i(z) q^{-d_1} = e_i(qz), \quad q^{d_1} f_i(z) q^{-d_1} = f_i(qz), \quad q^{d_1} \psi_i^\pm(z) q^{-d_1} = \psi_i^\pm(qz), \tag{T0.3}$$

$$q^{d_2} e_i(z) q^{-d_2} = q e_i(z), \quad q^{d_2} f_i(z) q^{-d_2} = q^{-1} f_i(z), \quad q^{d_2} \psi_i^\pm(z) q^{-d_2} = \psi_i^\pm(z), \tag{T0.4}$$

$$g_{a_{ij}}(\gamma^{-1} d^{m_{ij}} z/w) \psi_i^+(z) \psi_j^-(w) = g_{a_{ij}}(\gamma d^{m_{ij}} z/w) \psi_j^-(w) \psi_i^+(z), \tag{T1}$$

$$e_i(z) e_j(w) = g_{a_{ij}}(d^{m_{ij}} z/w) e_j(w) e_i(z), \tag{T2}$$

$$f_i(z) f_j(w) = g_{a_{ij}}(d^{m_{ij}} z/w)^{-1} f_j(w) f_i(z), \tag{T3}$$

$$(q - q^{-1})[e_i(z), f_j(w)] = \delta_{ij}\left(\delta(\gamma w/z)\psi_i^+(\gamma^{1/2}w) - \delta(\gamma z/w)\psi_i^-(\gamma^{1/2}z)\right), \tag{T4}$$

$$\psi_i^\pm(z) e_j(w) = g_{a_{ij}}(\gamma^{\pm 1/2} d^{m_{ij}} z/w) e_j(w) \psi_i^\pm(z), \tag{T5}$$

$$\psi_i^\pm(z) f_j(w) = g_{a_{ij}}(\gamma^{\mp 1/2} d^{m_{ij}} z/w)^{-1} f_j(w) \psi_i^\pm(z), \tag{T6}$$

$$[e_i(z), e_j(w)] = 0 \text{ if } a_{ij} = 0,$$
$$[e_i(z_1), [e_i(z_2), e_j(w)]_{q^{-1}}]_q + [e_i(z_2), [e_i(z_1), e_j(w)]_{q^{-1}}]_q = 0 \text{ if } a_{ij} = -1, \tag{T7}$$

$$[f_i(z), f_j(w)] = 0 \text{ if } a_{ij} = 0,$$
$$[f_i(z_1), [f_i(z_2), f_j(w)]_{q^{-1}}]_q + [f_i(z_2), [f_i(z_1), f_j(w)]_{q^{-1}}]_q = 0 \text{ if } a_{ij} = -1, \tag{T8}$$

where, as before, $[a, b]_x = ab - x \cdot ba$ and the generating series are defined via:

$$e_i(z) := \sum_{r \in \mathbb{Z}} e_{i,r} z^{-r}, \quad f_i(z) := \sum_{r \in \mathbb{Z}} f_{i,r} z^{-r}, \quad \psi_i^\pm(z) := \psi_{i,0}^{\pm 1} + \sum_{r > 0} \psi_{i,\pm r} z^{\mp r} \, .$$

It is convenient to use the generators $\{h_{i,r}\}_{r \neq 0}$ instead of $\{\psi_{i,r}\}_{r \neq 0}$, defined via:

$$\exp\left(\pm(q - q^{-1})\sum_{r > 0} h_{i,\pm r} z^{\mp r}\right) = \bar{\psi}_i^\pm(z) := \psi_{i,0}^{\mp 1}\psi_i^\pm(z) \tag{3.4}$$

with $h_{i,\pm r} \in \mathbb{C}[\psi_{i,0}^{\mp 1}, \psi_{i,\pm 1}, \psi_{i,\pm 2}, \ldots]$. Then, (T1, T5, T6) are equivalent to:

$$\psi_{i,0} h_{j,l} = h_{j,l}\psi_{i,0}, \quad [h_{i,k}, h_{j,-l}] = \delta_{kl} d^{-km_{ij}} \frac{[ka_{ij}]_q}{k}\frac{\gamma^k - \gamma^{-k}}{q - q^{-1}} \ (k, l \neq 0), \tag{T1$'$}$$

$$\psi_{i,0}e_{j,l} = q^{a_{ij}}e_{j,l}\psi_{i,0}, \quad [h_{i,k},e_{j,l}] = d^{-km_{ij}}\gamma^{-|k|/2}\frac{[ka_{ij}]_q}{k}e_{j,l+k} \quad (k \neq 0), \quad \text{(T5')}$$

$$\psi_{i,0}f_{j,l} = q^{-a_{ij}}f_{j,l}\psi_{i,0}, \quad [h_{i,k},f_{j,l}] = -d^{-km_{ij}}\gamma^{|k|/2}\frac{[ka_{ij}]_q}{k}f_{j,l+k} \quad (k \neq 0), \quad \text{(T6')}$$

with $[m]_q = \frac{q^m-q^{-m}}{q-q^{-1}}$, cf. (1.24). We also introduce $h_{i,0}$ and central elements $c,c'$ via:

$$\psi_{i,0} = q^{h_{i,0}}, \qquad \gamma^{1/2} = q^c, \qquad c' = \sum_{i \in [n]} h_{i,0}.$$

When needed below, we will also formally add elements $q^{\frac{h_{i,0}}{N}}, q^{\frac{c}{N}}, q^{\frac{d_1}{N}}, q^{\frac{d_2}{N}}$ $(N > 0)$.

Let $\ddot{U}_{q,d}^<$, $\ddot{U}_{q,d}^>$, and $\ddot{U}_{q,d}^0$ be the $\mathbb{C}$-subalgebras of $\ddot{U}_{q,d}(\mathfrak{sl}_n)$ generated by $\{f_{i,r}\}_{i \in [n]}^{r \in \mathbb{Z}}$, $\{e_{i,r}\}_{i \in [n]}^{r \in \mathbb{Z}}$, and $\{\psi_{i,r}, \psi_{i,0}^{-1}, \gamma^{\pm 1/2}, q^{\pm d_1}, q^{\pm d_2}\}_{i \in [n]}^{r \in \mathbb{Z}}$, respectively. The following result is completely analogous to Propositions 1.1 and 2.2:

**Proposition 3.1** *(a) (Triangular decomposition of $\ddot{U}_{q,d}(\mathfrak{sl}_n)$) The multiplication map*

$$m\colon \ddot{U}_{q,d}^< \otimes \ddot{U}_{q,d}^0 \otimes \ddot{U}_{q,d}^> \longrightarrow \ddot{U}_{q,d}(\mathfrak{sl}_n)$$

*is an isomorphism of $\mathbb{C}$-vector spaces.*

*(b) The algebra $\ddot{U}_{q,d}^>$ (resp. $\ddot{U}_{q,d}^0$ and $\ddot{U}_{q,d}^<$) is isomorphic to the associative $\mathbb{C}$-algebra generated by $\{e_{i,r}\}_{i \in [n]}^{r \in \mathbb{Z}}$ (resp. $\{\psi_{i,r}, \psi_{i,0}^{-1}, \gamma^{\pm 1/2}, q^{\pm d_1}, q^{\pm d_2}\}_{i \in [n]}^{r \in \mathbb{Z}}$ and $\{f_{i,r}\}_{i \in [n]}^{r \in \mathbb{Z}}$) with the defining relations (T2, T7) (resp. (T0.1, T0.2, T1) and (T3, T8)).*

We equip $\ddot{U}_{q,d}(\mathfrak{sl}_n)$ with a $\mathbb{Z}^{[n]} \times \mathbb{Z}$-grading via the following assignment:

$$\deg(e_{i,r}) := (1_i, r), \quad \deg(f_{i,r}) := (-1_i, r), \quad \deg(\psi_{i,r}) := (0, r),$$
$$\deg(x) := (0,0) \text{ for } x = \psi_{i,0}^{-1}, \gamma^{\pm 1/2}, q^{\pm d_1}, q^{\pm d_2} \text{ and any } i \in [n], r \in \mathbb{Z}, \tag{3.5}$$

where $1_j \in \mathbb{Z}^{[n]}$ has the $j$-th coordinate equal to 1 and all other coordinates being zero.

We shall also need the following modifications of the above algebra $\ddot{U}_{q,d}(\mathfrak{sl}_n)$.

- Let $\ddot{U}'_{q,d}(\mathfrak{sl}_n)$ be obtained from $\ddot{U}_{q,d}(\mathfrak{sl}_n)$ by "ignoring" the generators $q^{\pm d_2}$ and taking a quotient by the ideal $(c')$, that is, setting $\prod_{i \in [n]} \psi_{i,0} = 1$.

- Let ${}'\ddot{U}_{q,d}(\mathfrak{sl}_n)$ be obtained from $\ddot{U}_{q,d}(\mathfrak{sl}_n)$ by "ignoring" the generators $q^{\pm d_1}$ and taking a quotient by the ideal $(c)$, that is, setting $\gamma^{1/2} = 1$.

- Let $U_{q,d}^{\text{tor}}(\mathfrak{sl}_n)$ be obtained from $\ddot{U}_{q,d}(\mathfrak{sl}_n)$ by "ignoring" $q^{\pm d_1}$ and $q^{\pm d_2}$.

The above algebras $\ddot{U}'_{q,d}(\mathfrak{sl}_n)$, ${}'\ddot{U}_{q,d}(\mathfrak{sl}_n)$, $U_{q,d}^{\text{tor}}(\mathfrak{sl}_n)$ are also $\mathbb{Z}^{[n]} \times \mathbb{Z}$-graded via (3.5). Furthermore, they also satisfy the analogues of Proposition 3.1 with the corresponding subalgebras denoted by $\ddot{U}_{q,d}'^<, \ddot{U}_{q,d}'^>, \ddot{U}_{q,d}'^0, {}'\ddot{U}_{q,d}^<, {}'\ddot{U}_{q,d}^>, {}'\ddot{U}_{q,d}^0, U_{q,d}^{\text{tor}<}, U_{q,d}^{\text{tor}>}, U_{q,d}^{\text{tor}0}$, respectively. In particular, their "positive" subalgebras are naturally isomorphic:

$$\ddot{U}_{q,d}^> \simeq \ddot{U}_{q,d}'^> \simeq {}'\ddot{U}_{q,d}^> \simeq U_{q,d}^{\text{tor}>}. \tag{3.6}$$

### 3.1.2  Hopf pairing, Drinfeld double, and a universal $R$-matrix

Let us recall the general notion of a Hopf pairing, following [17, Chapter 3]. Given two Hopf $\mathbb{C}$-algebras $A$ and $B$ with invertible antipodes $S_A$ and $S_B$, a bilinear map

$$\varphi \colon A \times B \longrightarrow \mathbb{C}$$

is called a *Hopf pairing* if it satisfies the following properties:

$$\begin{aligned}
\varphi(a, bb') &= \varphi(a_1, b)\varphi(a_2, b') &&\forall\, a \in A,\ b, b' \in B\,, \\
\varphi(aa', b) &= \varphi(a, b_2)\varphi(a', b_1) &&\forall\, a, a' \in A,\ b \in B\,,
\end{aligned} \tag{3.7}$$

as well as

$$\varphi(a, 1_B) = \epsilon_A(a),\ \varphi(1_A, b) = \epsilon_B(b),\ \varphi(S_A(a), b) = \varphi(a, S_B^{-1}(b)) \quad \forall\, a \in A,\ b \in B\,,$$

where we use the Swedler notation (suppressing the summation) for the coproduct:

$$\Delta(x) = x_1 \otimes x_2\,. \tag{3.8}$$

For such a data, one defines the *generalized Drinfeld double* $D_\varphi(A, B)$ as follows ($D_\varphi(A, B)$ is called the *Drinfeld double* when $B = A^{\star, \mathrm{cop}}$ and $\varphi$ is the natural pairing):

**Theorem 3.1** *[17, Theorem 3.2] There is a unique Hopf algebra $D_\varphi(A, B)$ such that:*

*(i) $D_\varphi(A, B) \simeq A \otimes B$ as coalgebras.*

*(ii) Under the natural inclusions*

$$\begin{aligned}
A &\hookrightarrow D_\varphi(A, B) &&\text{given by} \quad a \mapsto a \otimes 1_B\,, \\
B &\hookrightarrow D_\varphi(A, B) &&\text{given by} \quad b \mapsto 1_A \otimes b\,,
\end{aligned}$$

*$A$ and $B$ are Hopf subalgebras of $D_\varphi(A, B)$.*

*(iii) For any $a \in A$ and $b \in B$, we have:*

$$(a \otimes 1_B) \cdot (1_A \otimes b) = a \otimes b\,,$$

$$(1_A \otimes b) \cdot (a \otimes 1_B) = \varphi(S_A^{-1}(a_1), b_1)\varphi(a_3, b_3)\, a_2 \otimes b_2\,,$$

*where we use the Swedler notation again (suppressing the summation):*

$$(\mathrm{Id} \otimes \Delta)(\Delta(x)) = (\Delta \otimes \mathrm{Id})(\Delta(x)) = x_1 \otimes x_2 \otimes x_3\,.$$

A Hopf algebra $A$ is called *quasitriangular* if there is an invertible element $R \in A \otimes A$ (or in a proper completion $R \in A \widehat{\otimes} A$) satisfying the following properties:

$$\begin{aligned}
(\Delta \otimes \mathrm{Id})(R) &= R^{13}R^{23}\,, \quad (\mathrm{Id} \otimes \Delta)(R) = R^{13}R^{12}\,, \\
R\Delta(x) &= \Delta^{\mathrm{op}}(x)R \quad \text{for all} \quad x \in A\,,
\end{aligned} \tag{3.9}$$

where for decomposable tensors (and extended to all tensors by linearity) we define

$$(a \otimes b)^{12} = a \otimes b \otimes 1, \quad (a \otimes b)^{13} = a \otimes 1 \otimes b, \quad (a \otimes b)^{23} = 1 \otimes a \otimes b \in A \otimes A \otimes A$$

as well as $(a \otimes b)^{\mathrm{op}} = b \otimes a$. Such an element $R$ is called a *universal R-matrix* of $A$.

The basic property of the generalized Drinfeld doubles is their quasitriangularity:

**Theorem 3.2** *[17, Theorem 3.2] For a nondegenerate Hopf pairing $\varphi \colon A \times B \to \mathbb{C}$, the generalized Drinfeld double $D_\varphi(A, B)$ of Theorem 3.1 is quasitriangular with a universal R-matrix given by:*

$$R = \sum_{i \in \mathfrak{I}} e_i \otimes e_i^*,$$

*where $\{e_i\}_{i \in \mathfrak{I}}$ and $\{e_i^*\}_{i \in \mathfrak{I}}$ are bases of $A$ and $B$, respectively, dual with respect to $\varphi$.*

To apply the above constructions to the quantum toroidal algebra $\ddot{U}_{q,d}(\mathfrak{sl}_n)$ and its subalgebras, we first need to endow the former with a Hopf algebra structure. This was carried out in [5, Theorem 2.1] in a more general setup:

**Theorem 3.3** $\ddot{U}_{q,d}(\mathfrak{sl}_n)$ *is endowed with a (topological) Hopf algebra structure via:*

coproduct $\Delta \colon e_i(z) \mapsto e_i(z) \otimes 1 + \psi_i^-(\gamma_{(1)}^{1/2} z) \otimes e_i(\gamma_{(1)} z),$

$\qquad\qquad f_i(z) \mapsto 1 \otimes f_i(z) + f_i(\gamma_{(2)} z) \otimes \psi_i^+(\gamma_{(2)}^{1/2} z),$

$\qquad\qquad \psi_i^\pm(z) \mapsto \psi_i^\pm(\gamma_{(2)}^{\pm 1/2} z) \otimes \psi_i^\pm(\gamma_{(1)}^{\mp 1/2} z),$ $\qquad$ (H1)

$\qquad\qquad x \mapsto x \otimes x \quad \text{for} \quad x = \gamma^{\pm 1/2}, q^{\pm d_1}, q^{\pm d_2},$

counit $\epsilon \colon e_i(z) \mapsto 0, \; f_i(z) \mapsto 0, \; \psi_i^\pm(z) \mapsto 1,$

$\qquad\qquad x \mapsto 1 \quad \text{for} \quad x = \gamma^{\pm 1/2}, q^{\pm d_1}, q^{\pm d_2},$ $\qquad$ (H2)

antipode $S \colon e_i(z) \mapsto -\psi_i^-(\gamma^{-1/2} z)^{-1} e_i(\gamma^{-1} z),$

$\qquad\qquad f_i(z) \mapsto -f_i(\gamma^{-1} z) \psi_i^+(\gamma^{-1/2} z)^{-1},$ $\qquad$ (H3)

$\qquad\qquad x \mapsto x^{-1} \quad \text{for} \quad x = \gamma^{\pm 1/2}, q^{\pm d_1}, q^{\pm d_2}, \psi_i^\pm(z),$

*where we set $\gamma_{(1)}^{1/2} := \gamma^{1/2} \otimes 1$ and $\gamma_{(2)}^{1/2} := 1 \otimes \gamma^{1/2}$.*

Let $\ddot{U}_{q,d}^{\geq}, \ddot{U}_{q,d}^{\leq}$ (resp. $\ddot{U}_{q,d}'^{\geq}, \ddot{U}_{q,d}'^{\leq}$ or $'\ddot{U}_{q,d}^{\geq}, '\ddot{U}_{q,d}^{\leq}$) be the subalgebras of $\ddot{U}_{q,d}(\mathfrak{sl}_n)$ (resp. $\ddot{U}_{q,d}'(\mathfrak{sl}_n)$ or $'\ddot{U}_{q,d}(\mathfrak{sl}_n)$) generated by $\{e_{i,r}, \psi_{i,l}, \psi_{i,0}^{\pm 1}, \gamma^{\pm 1/2}, q^{\pm d_1}, q^{\pm d_2}\}_{i \in [n]}^{r \in \mathbb{Z}, l \in -\mathbb{N}}$ and $\{f_{i,r}, \psi_{i,l}, \psi_{i,0}^{\pm 1}, \gamma^{\pm 1/2}, q^{\pm d_1}, q^{\pm d_2}\}_{i \in [n]}^{r \in \mathbb{Z}, l \in \mathbb{N}}$. Then, we have the following result:

**Theorem 3.4** *(a) There exists a unique Hopf algebra pairing*

$$\varphi \colon \ddot{U}_{q,d}^{\geq} \times \ddot{U}_{q,d}^{\leq} \longrightarrow \mathbb{C} \qquad\qquad (3.10)$$

*satisfying*

$$\varphi(e_i(z), f_j(w)) = \frac{\delta_{ij}}{q - q^{-1}} \cdot \delta\left(\frac{z}{w}\right), \quad \varphi(\psi_i^-(z), \psi_j^+(w)) = g_{a_{ij}}(d^{m_{ij}}z/w), \quad \text{(P1)}$$

$$\varphi(e_i(z), x^-) = \varphi(x^+, f_i(z)) = 0 \quad \text{for} \quad x^\pm = \psi_j^\mp(w), \gamma^{1/2}, q^{d_1}, q^{d_2}, \quad \text{(P2)}$$

$$\varphi(\gamma^{1/2}, q^{d_1}) = \varphi(q^{d_1}, \gamma^{1/2}) = q^{-1/2}, \; \varphi(\psi_i^-(z), q^{d_2}) = q^{-1}, \; \varphi(q^{d_2}, \psi_i^+(z)) = q, \quad \text{(P3)}$$

$$\varphi(\psi_i^-(z), x) = \varphi(x, \psi_i^+(z)) = 1 \quad \text{for} \quad x = \gamma^{1/2}, q^{d_1}, \quad \text{(P4)}$$

$$\varphi(\gamma^{1/2}, q^{d_2}) = \varphi(q^{d_2}, \gamma^{1/2}) = \varphi(\gamma^{1/2}, \gamma^{1/2}) = 1, \quad \text{(P5)}$$

$$\varphi(q^{d_1}, q^{d_1}) = \varphi(q^{d_1}, q^{d_2}) = \varphi(q^{d_2}, q^{d_1}) = 1, \quad \varphi(q^{d_2}, q^{d_2}) = q^{\frac{n^3 - n}{12}}. \quad \text{(P6)}$$

(b) *The natural Hopf algebra homomorphism* $D_\varphi(\ddot{U}_{q,d}^\geq, \ddot{U}_{q,d}^\leq) \to \ddot{U}_{q,d}(\mathfrak{sl}_n)$ *gives rise to the isomorphism*

$$D_\varphi(\ddot{U}_{q,d}^\geq, \ddot{U}_{q,d}^\leq)/\mathcal{J} \xrightarrow{\sim} \ddot{U}_{q,d}(\mathfrak{sl}_n), \quad \mathcal{J} := (x \otimes 1 - 1 \otimes x \mid x = \psi_{i,0}^{\pm 1}, \gamma^{\pm 1/2}, q^{\pm d_1}, q^{\pm d_2}).$$

(c) *Likewise, the algebras* $'\ddot{U}_{q,d}(\mathfrak{sl}_n)$, $\ddot{U}_{q,d}'(\mathfrak{sl}_n)$ *admit the Drinfeld double realizations via* $D_{'\varphi}('\ddot{U}_{q,d}^\geq, '\ddot{U}_{q,d}^\leq)$, $D_{\varphi'}(\ddot{U}_{q,d}'^\geq, \ddot{U}_{q,d}'^\leq)$ *with similarly defined pairings* $'\varphi, \varphi'$.

(d) *The pairings* $\varphi, \varphi', '\varphi$ *are nondegenerate iff* $q, qd, qd^{-1}$ *are not roots of unity.*

(e) *If* $q, qd, qd^{-1}$ *are not roots of unity, then* $\ddot{U}_{q,d}(\mathfrak{sl}_n), \ddot{U}_{q,d}'(\mathfrak{sl}_n), '\ddot{U}_{q,d}(\mathfrak{sl}_n)$ *admit universal R-matrices* $R, R', 'R$ *associated to the Hopf pairings* $\varphi, \varphi', '\varphi$, *respectively.*

*Remark 3.1* The existence of such Hopf pairings and the identification of the algebras given by generators and relations with the Drinfeld doubles are always verified by direct computations. The hardest part (often used as "folklore") is establishing the nondegeneracy of such pairings. For general quantum affine algebras this was proved in [NT, Proposition 6.7], while in the present toroidal setup it follows from [N].

*Remark 3.2* (a) We note that $\varphi(\psi_i^-(z), \psi_j^+(w)) = g_{a_{ij}}(d^{m_{ij}}z/w)$ is equivalent to:

$$\varphi(\psi_{i,0}^-, \psi_{j,0}^+) = q^{-a_{ij}}, \quad \varphi(h_{i,-k}, h_{j,l}) = \delta_{kl} \frac{d^{km_{ij}}[ka_{ij}]_q}{k(q - q^{-1})} \quad (k, l > 0). \quad (3.11)$$

(b) The specific choice of values in (P6) is only to simplify the formulas of [23, §4].

### 3.1.3 Horizontal and vertical copies of quantum affine $\mathfrak{sl}_n$

Following [7], we define the quantum affine algebra of $\mathfrak{sl}_n$ (in the new Drinfeld realization), denoted by $U_q(\widehat{\mathfrak{sl}}_n)$, to be the associative $\mathbb{C}$-algebra generated by $\{e_{i,r}, f_{i,r}, \psi_{i,r}, \psi_{i,0}^{-1}, C^{\pm 1}, \widetilde{D}^{\pm 1}\}_{i \in [n]^\times}^{r \in \mathbb{Z}}$ with the defining relations similar to (T0.1)–(T8):

$$[\psi_i^\pm(z), \psi_j^\pm(w)] = 0, \quad C^{\pm 1} - \text{central},$$

$$\psi_{i,0}^{\pm 1} \cdot \psi_{i,0}^{\mp 1} = C^{\pm 1} \cdot C^{\mp 1} = \widetilde{D}^{\pm 1} \cdot \widetilde{D}^{\mp 1} = 1,$$

$$\widetilde{D}e_i(z)\widetilde{D}^{-1} = qe_i(q^{-n}z), \quad \widetilde{D}f_i(z)\widetilde{D}^{-1} = q^{-1}f_i(q^{-n}z), \quad \widetilde{D}\psi_i^{\pm}(z)\widetilde{D}^{-1} = \psi_i^{\pm}(q^{-n}z),$$

$$g_{a_{ij}}(C^{-1}z/w)\psi_i^+(z)\psi_j^-(w) = g_{a_{ij}}(Cz/w)\psi_j^-(w)\psi_i^+(z),$$

$$e_i(z)e_j(w) = g_{a_{ij}}(z/w)e_j(w)e_i(z),$$

$$f_i(z)f_j(w) = g_{a_{ij}}(z/w)^{-1}f_j(w)f_i(z),$$

$$(q - q^{-1})[e_i(z), f_j(w)] = \delta_{ij}\left(\delta(Cw/z)\psi_i^+(Cw) - \delta(Cz/w)\psi_i^-(Cz)\right),$$

$$\psi_i^+(z)e_j(w) = g_{a_{ij}}(z/w)e_j(w)\psi_i^+(z), \quad \psi_i^-(z)e_j(w) = g_{a_{ij}}(C^{-1}z/w)e_j(w)\psi_i^-(z),$$

$$\psi_i^+(z)f_j(w) = g_{a_{ij}}(C^{-1}z/w)^{-1}f_j(w)\psi_i^+(z), \quad \psi_i^-(z)f_j(w) = g_{a_{ij}}(z/w)^{-1}f_j(w)\psi_i^-(z),$$

$$e_i(z)e_j(w) = e_j(w)e_i(z) \text{ if } a_{ij} = 0,$$
$$[e_i(z_1), [e_i(z_2), e_j(w)]_q]_q + [e_i(z_2), [e_i(z_1), e_j(w)]_{q^{-1}}]_q = 0 \text{ if } a_{ij} = -1,$$

$$f_i(z)f_j(w) = f_j(w)f_i(z) \text{ if } a_{ij} = 0,$$
$$[f_i(z_1), [f_i(z_2), f_j(w)]_q]_q + [f_i(z_2), [f_i(z_1), f_j(w)]_{q^{-1}}]_q = 0 \text{ if } a_{ij} = -1,$$

where the generating series $\{e_i(z), f_i(z), \psi_i^{\pm}(z)\}_{i \in [n]^\times}$ are defined as before. Let $U_q^{\text{aff}}(\mathfrak{sl}_n)$ be the subalgebra of $U_q(\widehat{\mathfrak{sl}}_n)$ obtained by "ignoring" the generators $\widetilde{D}^{\pm 1}$.

The above algebra $U_q(\widehat{\mathfrak{sl}}_n)$ is known to admit the original *Drinfeld-Jimbo realization* of [6, 15]. To present this identification explicitly, let $U_q^{\text{DJ}}(\widehat{\mathfrak{sl}}_n)$ be the associative $\mathbb{C}$-algebra generated by $\{E_i, F_i, K_i^{\pm 1}, D^{\pm 1}\}_{i \in [n]}$ with the following defining relations:

$$D^{\pm 1}D^{\mp 1} = 1, \quad DK_iD^{-1} = K_i, \quad DE_iD^{-1} = qE_i, \quad DF_iD^{-1} = q^{-1}F_i,$$

$$K_i^{\pm 1}K_i^{\mp 1} = 1, \quad K_iK_j = K_jK_i, \quad K_iE_jK_i^{-1} = q^{a_{ij}}E_j, \quad K_iF_jK_i^{-1} = q^{-a_{ij}}F_j,$$

$$(q - q^{-1})[E_i, F_j] = \delta_{ij}(K_i - K_i^{-1}),$$

$$E_iE_j = E_jE_i \text{ if } a_{ij} = 0, \quad E_i^2E_j - (q + q^{-1})E_iE_jE_i + E_jE_i^2 = 0 \text{ if } a_{ij} = -1,$$

$$F_iF_j = F_jF_i \text{ if } a_{ij} = 0, \quad F_i^2F_j - (q + q^{-1})F_iF_jF_i + F_jF_i^2 = 0 \text{ if } a_{ij} = -1.$$

According to [7], there is a $\mathbb{C}$-algebra isomorphism

$$\Phi: U_q^{\text{DJ}}(\widehat{\mathfrak{sl}}_n) \xrightarrow{\sim} U_q(\widehat{\mathfrak{sl}}_n) \tag{3.12}$$

given by:

$$
\begin{aligned}
&E_i \mapsto e_{i,0}, \quad F_i \mapsto f_{i,0}, \quad K_i^{\pm 1} \mapsto \psi_{i,0}^{\pm 1} \quad \text{for} \quad 1 \le i \le n-1, \\
&K_0 \mapsto C \cdot (\psi_{1,0} \cdots \psi_{n-1,0})^{-1}, \quad D \mapsto \widetilde{D}, \\
&E_0 \mapsto C(\psi_{1,0} \cdots \psi_{n-1,0})^{-1} \cdot [\cdots [f_{1,1}, f_{2,0}]_q, \cdots, f_{n-1,0}]_q, \\
&F_0 \mapsto [e_{n-1,0}, \cdots, [e_{2,0}, e_{1,-1}]_{q^{-1}} \cdots]_{q^{-1}} \cdot (\psi_{1,0} \cdots \psi_{n-1,0})C^{-1}.
\end{aligned} \tag{3.13}
$$

*Remark 3.3* (a) The isomorphism (3.12, 3.13) was presented in [7] without a proof. The inverse $\Phi^{-1}$ was constructed in [1] by using Lusztig's braid group action, while the direct verification of (3.13) giving rise to a $\mathbb{C}$-algebra homomorphism $\Phi \colon U_q^{\mathrm{aff}}(\widehat{\mathfrak{sl}}_n) \to U_q(\widehat{\mathfrak{sl}}_n)$ appeared in [16]. However, proofs of injectivity of $\Phi^{-1}, \Phi$ in [1, 16] had gaps that were fixed more recently in [3].

(b) In some standard literature, the *grading* elements $D$ and $\widetilde{D}$ satisfy slightly different relations. However, the present conventions are better adapted to fit into the toroidal story (in particular, our discussion in the next section) and follow that of [19].

(c) The quantum loop algebra $U_q(L\mathfrak{sl}_n)$ from Chapter 1 is the quotient of $U_q^{\mathrm{aff}}(\widehat{\mathfrak{sl}}_n)$ by the ideal $(C-1)$, thus corresponding to a trivial central charge.

(d) Likewise, the algebra $U_q^{\mathrm{DI}}(L\mathfrak{sl}_n)$ from Remark 1.4 is obtained from $U_q^{\mathrm{DI}}(\widehat{\mathfrak{sl}}_n)$ by "ignoring" the generators $D^{\pm 1}$ and taking a quotient by the ideal $(K_0 K_1 \cdots K_{n-1} - 1)$.

Following [24], let us introduce the *vertical* and *horizontal* copies of the quantum affine $U_q(\widehat{\mathfrak{sl}}_n)$ inside $\ddot{U}_{q,d}(\mathfrak{sl}_n)$. To this end, we consider two algebra homomorphisms

$$U_q(\widehat{\mathfrak{sl}}_n) \overset{\mathrm{v}}{\hookrightarrow} \ddot{U}_{q,d}(\mathfrak{sl}_n) \overset{\mathrm{h}}{\hookleftarrow} U_q(\widehat{\mathfrak{sl}}_n), \tag{3.14}$$

defined by:

$$\Phi \circ \mathrm{h} \colon E_i \mapsto e_{i,0}, \quad F_i \mapsto f_{i,0}, \quad K_i \mapsto \psi_{i,0}, \quad D \mapsto q^{d_2} \quad \text{for} \quad i \in [n] \tag{3.15}$$

with $\Phi$ of (3.12, 3.13), and

$$\begin{aligned}
&\mathrm{v} \colon e_{i,r} \mapsto d^{ir} e_{i,r}, \quad f_{i,r} \mapsto d^{ir} f_{i,r}, \quad \psi_{i,r} \mapsto d^{ir} \gamma^{r/2} \psi_{i,r}, \\
&C \mapsto \gamma, \quad \widetilde{D} \mapsto q^{-nd_1} \cdot q^{\sum_{j=1}^{n-1} \frac{j(n-j)}{2} h_{j,0}} \quad \text{for} \quad 0 \le i < n, r \in \mathbb{Z},
\end{aligned} \tag{3.16}$$

where we follow the conventions of Section 3.1.1 and add elements $q^{\frac{h_{j,0}}{2}}$ to $\ddot{U}_{q,d}(\mathfrak{sl}_n)$. According to [24], these algebra homomorphisms are injective. The corresponding images of h and v are called the *horizontal* and the *vertical* copies of $U_q(\widehat{\mathfrak{sl}}_n)$.

Following Remark 3.3(c, d), we also note that $\ddot{U}'_{q,d}(\mathfrak{sl}_n)$, $'\ddot{U}_{q,d}(\mathfrak{sl}_n)$, $U_{q,d}^{\mathrm{tor}}(\mathfrak{sl}_n)$ admit similar constructions, that is, we have the following algebra embeddings:

$$U_q(\widehat{\mathfrak{sl}}_n) \overset{\mathrm{v}'}{\hookrightarrow} \ddot{U}'_{q,d}(\mathfrak{sl}_n) \overset{\mathrm{h}'}{\hookleftarrow} U_q(L\mathfrak{sl}_n), \tag{3.17}$$

$$U_q(L\mathfrak{sl}_n) \overset{'\mathrm{v}}{\hookrightarrow} {}'\ddot{U}_{q,d}(\mathfrak{sl}_n) \overset{'\mathrm{h}}{\hookleftarrow} U_q(\widehat{\mathfrak{sl}}_n), \tag{3.18}$$

$$U_q^{\mathrm{aff}}(\mathfrak{sl}_n) \overset{\mathrm{v}}{\hookrightarrow} U_{q,d}^{\mathrm{tor}}(\mathfrak{sl}_n) \overset{\mathrm{h}}{\hookleftarrow} U_q^{\mathrm{aff}}(\mathfrak{sl}_n). \tag{3.19}$$

By abuse of notations, we shall use $\dot{U}_q^{\mathrm{h}}(\mathfrak{sl}_n)$ and $\dot{U}_q^{\mathrm{v}}(\mathfrak{sl}_n)$ to denote the corresponding horizontal and vertical copies of $U_q(\mathfrak{sl}_n)$, $U_q^{\mathrm{aff}}(\mathfrak{sl}_n)$, $U_q(L\mathfrak{sl}_n)$ in one of the algebras $\ddot{U}_{q,d}(\mathfrak{sl}_n)$, $\ddot{U}'_{q,d}(\mathfrak{sl}_n)$, $'\ddot{U}_{q,d}(\mathfrak{sl}_n)$, $U_{q,d}^{\mathrm{tor}}(\mathfrak{sl}_n)$ from the embeddings (3.14, 3.17–3.19).

### 3.1.4 Miki's isomorphism

In this section, we recall the beautiful result of K. Miki that provides an algebra isomorphism $'\ddot{U}_{q,d}(\mathfrak{sl}_n) \xrightarrow{\sim} \ddot{U}'_{q,d}(\mathfrak{sl}_n)$ intertwining the *vertical* and *horizontal* embeddings of the quantum affine/loop algebras of $\mathfrak{sl}_n$ in (3.17, 3.18), similar to (2.18).

To state this result explicitly, we shall need some more notations:

- Let $\sigma$ be the algebra anti-automorphism of $U_q(\widehat{\mathfrak{sl}}_n) \simeq U_q^{\mathrm{DJ}}(\widehat{\mathfrak{sl}}_n)$ determined by:

$$\sigma \colon E_i \mapsto E_i, \quad F_i \mapsto F_i, \quad K_i \mapsto K_i^{-1}, \quad D \mapsto D^{-1} \quad \text{for} \quad i \in [n].$$

- Let $\eta$ be the algebra anti-automorphism of $U_q(\widehat{\mathfrak{sl}}_n)$ determined by:

$$\eta \colon e_{i,r} \mapsto e_{i,-r}, \quad f_{i,r} \mapsto f_{i,-r}, \quad h_{i,l} \mapsto -C^l h_{i,-l}, \quad \psi_{i,0} \mapsto \psi_{i,0}^{-1},$$
$$C \mapsto C, \quad \widetilde{D} \mapsto \widetilde{D} \cdot \prod_{1 \le i \le n-1} \psi_{i,0}^{-i(n-i)} \quad \text{for} \quad 1 \le i < n, \, r \in \mathbb{Z}, \, l \ne 0.$$

- We shall also use the corresponding anti-automorphisms $\sigma, \eta$ of $U_q(L\mathfrak{sl}_n)$.

- Let $'Q$ be the algebra automorphism of $'\ddot{U}_{q,d}(\mathfrak{sl}_n)$ determined by:

$$'Q \colon e_{i,r} \mapsto (-d)^r e_{i+1,r}, \quad f_{i,r} \mapsto (-d)^r f_{i+1,r}, \quad h_{i,l} \mapsto (-d)^l h_{i+1,l},$$
$$\psi_{i,0} \mapsto \psi_{i+1,0}, \quad q^{d_2} \mapsto q^{d_2} \quad \text{for} \quad i \in [n], \, r \in \mathbb{Z}, \, l \ne 0.$$

- Let $Q'$ be the algebra automorphism of $\ddot{U}'_{q,d}(\mathfrak{sl}_n)$ such that it maps the generators other than $\gamma^{\pm 1/2}, q^{\pm d_1}$ as $'Q$, while

$$Q' \colon \gamma^{1/2} \mapsto \gamma^{1/2}, \quad q^{d_1} \mapsto q^{d_1} \cdot \gamma^{-1}.$$

- Let $'\mathcal{Y}_j \ (1 \le j \le n)$ be the algebra automorphism of $'\ddot{U}_{q,d}(\mathfrak{sl}_n)$ determined by:

$$'\mathcal{Y}_j \colon h_{i,l} \mapsto h_{i,l}, \quad \psi_{i,0} \mapsto \psi_{i,0}, \quad q^{d_2} \mapsto q^{d_2},$$
$$e_{i,r} \mapsto (-d)^{-n\delta_{i0}\delta_{jn}-i\tilde{\delta}_{ij}+i\delta_{i,j-1}} e_{i,r-\tilde{\delta}_{ij}+\delta_{i,j-1}} \quad \text{for} \quad i \in [n], \, r \in \mathbb{Z}, \, l \ne 0,$$
$$f_{i,r} \mapsto (-d)^{n\delta_{i0}\delta_{jn}+i\tilde{\delta}_{ij}-i\delta_{i,j-1}} f_{i,r+\tilde{\delta}_{ij}-\delta_{i,j-1}}, \quad \tilde{\delta}_{ij} := \begin{cases} 1 & \text{if } j \equiv i \pmod{n} \\ 0 & \text{otherwise} \end{cases}$$

- Let $\mathcal{Y}'_j \ (1 \le j \le n)$ be the algebra automorphism of $\ddot{U}'_{q,d}(\mathfrak{sl}_n)$ such that it maps the generators other than $\psi_{i,0}^{\pm 1}, \gamma^{\pm 1/2}, q^{\pm d_1}$ as $'\mathcal{Y}_j$, while

$$\mathcal{Y}'_j \colon \gamma^{1/2} \mapsto \gamma^{1/2}, \quad \psi_{i,0} \mapsto \gamma^{-\tilde{\delta}_{ij}+\delta_{i,j-1}} \psi_{i,0}, \quad q^{d_1} \mapsto q^{d_1} \cdot \gamma^{-\frac{n+1}{2n}} \cdot \widetilde{K}_j,$$

where $\widetilde{K}_j = \prod_{l=1}^{j-1} q^{\frac{l}{n} h_{l,0}} \cdot \prod_{l=j}^{n-1} q^{\frac{l-n}{n} h_{l,0}}$ and we follow the conventions of Section 3.1.1 in that we add elements $\gamma^{\frac{1}{2n}}, q^{\frac{h_{j,0}}{2n}}$ to $\ddot{U}'_{q,d}(\mathfrak{sl}_n)$.

**Theorem 3.5** *[19, Proposition 1] There exists an algebra isomorphism*

$$\varpi : \, 'U''_{q,d}(\mathfrak{sl}_n) \stackrel{\sim}{\longrightarrow} \ddot{U}'_{q,d}(\mathfrak{sl}_n) \tag{3.20}$$

*satisfying the following properties:*

$$\varpi \circ 'h = v', \quad \varpi \circ 'v \circ \eta \circ \sigma = h', \quad Q' \circ \mathcal{Y}'_n \circ \varpi = \varpi \circ {}'\mathcal{Y}_1^{-1} \circ {}'Q \, . \tag{3.21}$$

*Remark 3.4* (a) The construction of $\varpi$ in [19] was based on the previous work [18], where an analogous algebra automorphism $\overline{\varpi}$ of $U^{\mathrm{tor}}_{q,d}(\mathfrak{sl}_n)$ was constructed.

(b) In [19], the parameters $q, \xi$ are related to our $q, d$ via $q \leftrightarrow q$ and $\xi \leftrightarrow d^{-n}$, while the generators $\{x^{\pm}_{i,r}, h_{i,l}, k^{\pm 1}_i, C^{\pm 1}, D^{\pm 1}, \widetilde{D}^{\pm 1}\}$ of [19, §B.1] are related to ours via:

$$x^+_{i,r} \leftrightarrow d^{ir} e_{i,r}, \quad x^-_{i,r} \leftrightarrow d^{ir} f_{i,r}, \quad h_{i,l} \leftrightarrow d^{il} \gamma^{l/2} h_{i,l},$$

$$k^{\pm 1}_i \leftrightarrow \psi^{\pm 1}_{i,0}, \quad C^{\pm 1} \leftrightarrow \gamma^{\pm 1}, \quad D^{\pm 1} \leftrightarrow q^{\pm d_2}, \quad \widetilde{D}^{\pm 1} \leftrightarrow q^{\mp n d_1} \cdot q^{\pm \sum_{j=1}^{n-1} \frac{j(n-j)}{2} h_{j,0}} \, .$$

(c) While the generators of [19] do not require one to add $q^{\frac{h_{j,0}}{N}}, q^{\frac{c}{N}}, q^{\frac{d_1}{N}}, q^{\frac{d_2}{N}}$, the present choice is more symmetric and is better suited for the rest of our exposition.

We conclude this section by computing explicitly $\varpi$-images of some generators:

**Proposition 3.2** *(a) We have*

$$\varpi : e_{i,0} \mapsto e_{i,0}, \quad f_{i,0} \mapsto f_{i,0}, \quad \psi^{\pm 1}_{i,0} \mapsto \psi^{\pm 1}_{i,0} \quad \text{for} \quad i \in [n]^{\times},$$

$$\psi^{\pm 1}_{0,0} \mapsto \gamma^{\pm 1} \cdot \psi^{\pm 1}_{0,0}, \quad q^{\pm d_2} \mapsto q^{\mp n d_1} \cdot q^{\pm \sum_{j=1}^{n-1} \frac{j(n-j)}{2} h_{j,0}},$$

$$e_{0,0} \mapsto d \cdot \gamma \psi_{0,0} \cdot [\cdots [f_{1,1}, f_{2,0}]_q, \cdots, f_{n-1,0}]_q \, ,$$

$$f_{0,0} \mapsto d^{-1} \cdot [e_{n-1,0}, \cdots, [e_{2,0}, e_{1,-1}]_{q^{-1}} \cdots]_{q^{-1}} \cdot \psi^{-1}_{0,0} \gamma^{-1} \, .$$

*(b) For $1 \leq i < n$, we have*

$$\varpi : h_{i,1} \mapsto (-1)^{i+1} d^{-i} \times$$
$$[[\cdots [[\cdots [f_{0,0}, f_{n-1,0}]_q, \cdots, f_{i+1,0}]_q, f_{1,0}]_q, \cdots, f_{i-1,0}]_q, f_{i,0}]_{q^2} \, ,$$
$$h_{i,-1} \mapsto (-1)^{i+1} d^i \times$$
$$[e_{i,0}, [e_{i-1,0}, \cdots, [e_{1,0}, [e_{i+1,0}, \cdots, [e_{n-1,0}, e_{0,0}]_{q^{-1}} \cdots]_{q^{-1}}]_{q^{-1}} \cdots]_{q^{-1}}]_{q^{-2}} \, .$$

*(c) For $i = 0$, we have*

$$\varpi : h_{0,1} \mapsto (-1)^n d^{1-n} \cdot [[\cdots [f_{1,1}, f_{2,0}]_q, \cdots, f_{n-1,0}]_q, f_{0,-1}]_{q^2} \, ,$$
$$h_{0,-1} \mapsto (-1)^n d^{n-1} \cdot [e_{0,1}, [e_{n-1,0}, \cdots, [e_{2,0}, e_{1,-1}]_{q^{-1}} \cdots]_{q^{-1}}]_{q^{-2}} \, .$$

*(d) We have*

$$\varpi : e_{0,-1} \mapsto (-d)^n e_{0,1} \, , \quad f_{0,1} \mapsto (-d)^{-n} f_{0,-1} \, .$$

*Proof* Part (a) follows directly by applying the equality $\varpi \circ {}'\mathsf{h} = \mathsf{v}'$ to the explicit formulas for $\Phi(E_i), \Phi(F_i), \Phi(K_i), \Phi(D)$ with $\Phi$ being the isomorphism of (3.12, 3.13).

(b) First, let us express $h_{i,\pm 1} \in U_q(\widehat{\mathfrak{sl}}_n)$ in the Drinfeld-Jimbo presentation:

$$\Phi^{-1}(h_{i,1}) = \tag{3.22}$$
$$(-1)^i[[\cdots[[\cdots[E_0, E_{n-1}]_{q^{-1}}, \cdots, E_{i+1}]_{q^{-1}}, E_1]_{q^{-1}}, \cdots, E_{i-1}]_{q^{-1}}, E_i]_{q^{-2}},$$

$$\Phi^{-1}(h_{i,-1}) = \tag{3.23}$$
$$(-1)^i[F_i, [F_{i-1}, \cdots, [F_1, [F_{i+1}, \cdots, [F_{n-1}, F_0]_q \cdots]_q]_q \cdots]_q]_{q^2}.$$

These formulas are proved by iteratively applying two useful identities on $q$-brackets (we leave details to the interested reader):

$$[a, [b, c]_u]_v = [[a, b]_x, c]_{uv/x} + x \cdot [b, [a, c]_{v/x}]_{u/x},$$

$$[[a, b]_u, c]_v = [a, [b, c]_x]_{uv/x} + x \cdot [[a, c]_{v/x}, b]_{u/x}.$$

Applying the equality $\varpi \circ {}'\mathsf{v} \circ \eta = \mathsf{h}' \circ \sigma^{-1}$ to the formulas (3.22) and (3.23), we obtain the stated formulas for the images $\varpi(h_{i,1})$ and $\varpi(h_{i,-1})$, respectively.

(c) Applying the equality $Q' \circ \mathcal{Y}'_n \circ \varpi = \varpi \circ {}'\mathcal{Y}_1^{-1} \circ {}'Q$ to $h_{n-1,\pm 1}$ and using the formulas for $\varpi(h_{n-1,\pm 1})$ from part (b), we obtain the stated formulas for $\varpi(h_{0,\pm 1})$.

(d) The explicit formulas for $\varpi(e_{0,-1})$ and $\varpi(f_{0,1})$ follow by applying the equality $Q' \circ \mathcal{Y}'_n \circ \varpi = \varpi \circ {}'\mathcal{Y}_1^{-1} \circ {}'Q$ to $e_{n-1,0}$ and $f_{n-1,0}$, respectively. □

### 3.1.5 Horizontal and vertical extra Heisenberg subalgebras

Following [9], let us now extend the vertical and horizontal $\dot{U}_q^{\mathsf{v}}(\mathfrak{sl}_n)$, $\dot{U}_q^{\mathsf{h}}(\mathfrak{sl}_n)$ to $\dot{U}_q^{\mathsf{v}}(\mathfrak{gl}_n)$, $\dot{U}_q^{\mathsf{h}}(\mathfrak{gl}_n)$ by adjoining vertical and horizontal Heisenberg algebras $\mathfrak{h}_1^{\mathsf{v}}$ and $\mathfrak{h}_1^{\mathsf{h}}$.

For $r \neq 0$, pick a nontrivial solution $\{p_{i,r}\}_{i \in [n]}$ of the following linear system:

$$\sum_{i \in [n]} p_{i,r} d^{-rm_{ij}} [ra_{ij}]_q = 0 \quad \forall j \in [n]^\times, \tag{3.24}$$

and let $\mathfrak{h}_1^{\mathsf{v}}$ be the subspace of $U_{q,d}^{\mathrm{tor}}(\mathfrak{sl}_n)$ spanned by

$$h_r^{\mathsf{v}} = \begin{cases} \sum_{i \in [n]} p_{i,r} h_{i,r} & \text{if } r \neq 0 \\ \gamma^{1/2} & \text{if } r = 0 \end{cases}. \tag{3.25}$$

We note that $\mathfrak{h}_1^{\mathsf{v}}$ is well-defined (as (3.24) has a one-dimensional space of solutions) and commutes with $\dot{U}_q^{\mathsf{v}}(\mathfrak{sl}_n)$ due to (T1′, T5′, T6′). Since $[h_k^{\mathsf{v}}, h_l^{\mathsf{v}}] = 0$ for $l \neq -k$, $\mathfrak{h}_1^{\mathsf{v}}$ is isomorphic to a Heisenberg algebra. Let $\dot{U}_q^{\mathsf{v}}(\mathfrak{gl}_n)$ be the subalgebra of $U_{q,d}^{\mathrm{tor}}(\mathfrak{sl}_n)$ generated by $\dot{U}_q^{\mathsf{v}}(\mathfrak{sl}_n)$ and $\mathfrak{h}_1^{\mathsf{v}}$, so that $\dot{U}_q^{\mathsf{v}}(\mathfrak{gl}_n) \simeq U_q^{\mathrm{aff}}(\mathfrak{gl}_n)$–the quantum affine $\mathfrak{gl}_n$. Let $\dot{U}_q^{\mathsf{v}}(\mathfrak{h}_1)$, $\dot{U}_q^{\mathsf{v}}(\mathfrak{h}_n)$ be the subalgebras generated by $\{h_r^{\mathsf{v}}, \gamma^{\pm 1/2}\}_{r \neq 0}$ and $\{h_{i,r}, \gamma^{\pm 1/2}\}_{i \in [n]}^{r \neq 0}$.

Let us recall the algebra automorphism $\overline{\varpi}$ of $U_{q,d}^{\mathrm{tor}}(\mathfrak{sl}_n)$ from Remark 3.4(a) satisfying

$$\overline{\varpi} \colon \dot{U}_q^{\mathrm{v}}(\mathfrak{sl}_n) \mapsto \dot{U}_q^{\mathrm{h}}(\mathfrak{sl}_n), \quad \dot{U}_q^{\mathrm{h}}(\mathfrak{sl}_n) \mapsto \dot{U}_q^{\mathrm{v}}(\mathfrak{sl}_n), \quad q^{c'} \mapsto q^{2c}, \quad q^c \mapsto q^{-c'/2}.$$

Then, we define $\dot{U}_q^{\mathrm{h}}(\mathfrak{gl}_n)$ as the subalgebra of $U_{q,d}^{\mathrm{tor}}(\mathfrak{sl}_n)$ generated by $\dot{U}_q^{\mathrm{h}}(\mathfrak{sl}_n)$ and the *horizontal* extra Heisenberg $\mathfrak{h}_1^{\mathrm{h}} := \overline{\varpi}(\mathfrak{h}_1^{\mathrm{v}})$. Thus $\dot{U}_q^{\mathrm{h}}(\mathfrak{gl}_n) = \overline{\varpi}(\dot{U}_q^{\mathrm{v}}(\mathfrak{gl}_n)) \simeq U_q^{\mathrm{aff}}(\mathfrak{gl}_n)$. We define $\dot{U}_q^{\mathrm{h}}(\mathfrak{b}_1) := \overline{\varpi}(\dot{U}_q^{\mathrm{v}}(\mathfrak{b}_1))$, $\dot{U}_q^{\mathrm{h}}(\mathfrak{b}_n) := \overline{\varpi}(\dot{U}_q^{\mathrm{v}}(\mathfrak{b}_n))$, and consider their "positive" subalgebras $\dot{U}_q^{\mathrm{h}>}(\mathfrak{b}_1)$ and $\dot{U}_q^{\mathrm{h}>}(\mathfrak{b}_n)$ generated by $\{\overline{\varpi}(h_{-r}^{\mathrm{v}})\}_{r>0}$ and $\{\overline{\varpi}(h_{i,-r})\}_{i\in[n]}^{r>0}$.

*Remark 3.5* (a) The other three versions of the quantum toroidal algebra of $\mathfrak{sl}_n$ admit natural counterparts of these constructions, but with $U_q^{\mathrm{aff}}(\mathfrak{gl}_n)$ being replaced either by $U_q(L\mathfrak{gl}_n)$ (the *trivial central charge* version) or $U_q(\widehat{\mathfrak{gl}}_n)$ (obtained by adjoining the *degree* generator). To be more precise, we have: $U_q^{\mathrm{h}}(\widehat{\mathfrak{gl}}_n), U_q^{\mathrm{v}}(\widehat{\mathfrak{gl}}_n) \subset \ddot{U}_{q,d}(\mathfrak{sl}_n)$, $U_q^{\mathrm{h}}(L\mathfrak{gl}_n), U_q^{\mathrm{v}}(\widehat{\mathfrak{gl}}_n) \subset \ddot{U}_{q,d}'(\mathfrak{sl}_n)$, and $U_q^{\mathrm{h}}(\widehat{\mathfrak{gl}}_n), U_q^{\mathrm{v}}(L\mathfrak{gl}_n) \subset {}'\ddot{U}_{q,d}(\mathfrak{sl}_n)$.

(b) The explicit formulas for the aforementioned images $\overline{\varpi}(h_{-r}^{\mathrm{v}})$, $\overline{\varpi}(h_{i,-r})$ expressing them via the original algebra generators are complicated and unknown. However, their shuffle realization is provided in the end of this chapter, see Theorem 3.19.

(c) We note that [20] uses an apparently different horizontal quantum affine $\mathfrak{gl}_n$ arising through the RTT realization [4] of $U_q^{\mathrm{aff}}(\mathfrak{gl}_n)$. However, it actually coincides with $\dot{U}_q^{\mathrm{h}}(\mathfrak{gl}_n)$ as follows from Section 3.2.5, see Corollary 3.1.

## 3.2 Shuffle algebra and its commutative subalgebras

In this section, we recall the shuffle algebra realization of $\ddot{U}_{q,d}^{>}$ established in [20]. Following [12], we also provide a shuffle realization of a family of its commutative subalgebras, reminiscent to the commutative subalgebra $\mathcal{A}$ of the quantum toroidal $\mathfrak{gl}_1$ from Section 2.1.6. These commutative subalgebras will be identified with the Bethe subalgebras in Section 3.5.4, while their particular limit recovers $\dot{U}_q^{\mathrm{h}}(\mathfrak{b}_n)$, the subalgebra generated by $\{\overline{\varpi}(h_{i,-r})\}_{i\in[n]}^{r>0}$ (that is, half of the horizontal Heisenberg).

### 3.2.1 Shuffle algebra realization

Consider an $\mathbb{N}^{[n]}$-graded $\mathbb{C}$-vector space

$$\mathbb{S} = \bigoplus_{\underline{k}=(k_0,\ldots,k_{n-1})\in\mathbb{N}^{[n]}} \mathbb{S}_{\underline{k}},$$

where $\mathbb{S}_{\underline{k}}$ consists of $\Sigma_{\underline{k}}$-symmetric rational functions in the variables $\{x_{i,r}\}_{i\in[n]}^{1\le r\le k_i}$. We fix a matrix of rational functions $(\zeta_{i,j}(z))_{i,j\in[n]} \in \mathrm{Mat}_{[n]\times[n]}(\mathbb{C}(z))$ via:

$$\zeta_{i,i}(z) = \frac{z - q^{-2}}{z-1}, \ \zeta_{i,i+1}(z) = \frac{d^{-1}z - q}{z-1}, \ \zeta_{i,i-1}(z) = \frac{z - qd^{-1}}{z-1}, \ \zeta_{i,j\notin\{i,i\pm1\}}(z) = 1 \, .$$

Similarly to (1.32), let us now introduce the bilinear shuffle product $\star$ on $\mathbb{S}$: given $F \in \mathbb{S}_{\underline{k}}$ and $G \in \mathbb{S}_{\underline{\ell}}$, define $F \star G \in \mathbb{S}_{\underline{k}+\underline{\ell}}$ via

$$(F \star G)(x_{0,1}, \ldots, x_{0,k_0+\ell_0}; \ldots; x_{n-1,1}, \ldots, x_{n-1,k_{n-1}+\ell_{n-1}}) := \frac{1}{\underline{k}! \cdot \underline{\ell}!} \times$$

$$\text{Sym}\left( F\left(\{x_{i,r}\}_{i\in[n]}^{1\leq r\leq k_i}\right) G\left(\{x_{i',r'}\}_{i'\in[n]}^{k_{i'}<r'\leq k_{i'}+\ell_{i'}}\right) \cdot \prod_{i\in[n]}^{i'\in[n]} \prod_{r\leq k_i}^{r'>k_{i'}} \zeta_{i,i'}\left(\frac{x_{i,r}}{x_{i',r'}}\right) \right). \qquad (3.26)$$

Here, for $\underline{m} = (m_0, m_1, \ldots, m_{n-1}) \in \mathbb{N}^{[n]}$, we set $\underline{m}! := \prod_{i\in[n]} m_i!$, while the *symmetrization* of $f \in \mathbb{C}(\{x_{i,1}, \ldots, x_{i,m_i}\}_{i\in[n]})$ is defined similarly to (1.33):

$$\text{Sym}(f)\left(\{x_{i,1}, \ldots, x_{i,m_i}\}_{i\in[n]}\right) := \sum_{(\sigma_0,\ldots,\sigma_{n-1})\in\Sigma_{\underline{m}}} f\left(\{x_{i,\sigma_i(1)}, \ldots, x_{i,\sigma_i(m_i)}\}_{i\in[n]}\right).$$

This endows $\mathbb{S}$ with a structure of an associative algebra with the unit $\mathbf{1} \in \mathbb{S}_{(0,\ldots,0)}$.

As before, we consider the subspace of $\mathbb{S}$ defined by the *pole* and *wheel* conditions:

- We say that $F \in \mathbb{S}_{\underline{k}}$ satisfies the *pole conditions* if

$$F = \frac{f(x_{0,1}, \ldots, x_{n-1,k_{n-1}})}{\prod_{i\in[n]} \prod_{r\leq k_{i+1}}^{r'\leq k_{i+1}} (x_{i,r} - x_{i+1,r'})} \ \text{ with } \ f \in \mathbb{C}\left[\{x_{i,r}^{\pm1}\}_{i\in[n]}^{1\leq r\leq k_i}\right]^{\Sigma_{\underline{k}}}. \qquad (3.27)$$

- We say that $F \in \mathbb{S}_{\underline{k}}$ satisfies the *wheel conditions* if

$$F\left(\{x_{i,r}\}\right) = 0 \ \text{ once } \ x_{i,r_1} = qd^\epsilon x_{i+\epsilon,s} = q^2 x_{i,r_2} \ \text{ for some } \ i, \epsilon, r_1, r_2, s, \qquad (3.28)$$

where $\epsilon \in \{\pm1\}, i \in [n], 1 \leq r_1 \neq r_2 \leq k_i, 1 \leq s \leq k_{\overline{i+\epsilon}}, \overline{i+\epsilon} := i+\epsilon \bmod n$.

Let $S_{\underline{k}} \subset \mathbb{S}_{\underline{k}}$ denote the subspace of all $F$ satisfying these two conditions and set

$$S := \bigoplus_{\underline{k}\in\mathbb{N}^{[n]}} S_{\underline{k}} \, .$$

Similarly to Lemma 1.2, it is straightforward to check that $S \subset \mathbb{S}$ is $\star$-closed:

**Lemma 3.1** *For any $F \in S_{\underline{k}}$ and $G \in S_{\underline{\ell}}$, we have $F \star G \in S_{\underline{k}+\underline{\ell}}$.*

The algebra $(S, \star)$ is called the $A_{n-1}^{(1)}$-type **shuffle algebra**. It is $\mathbb{N}^{[n]} \times \mathbb{Z}$-graded via:

$$S = \bigoplus_{(\underline{k},d)\in\mathbb{N}^{[n]}\times\mathbb{Z}} S_{\underline{k},d} \qquad \text{with} \qquad S_{\underline{k},d} := \left\{ F \in S_{\underline{k}} \, \middle| \, \text{tot.deg}(F) = d \right\}. \qquad (3.29)$$

Similarly to Proposition 1.3, the shuffle algebra $(S, \star)$ is linked to the "positive" subalgebra of the quantum toroidal $\mathfrak{sl}_n$, cf. (3.6), via:

**Proposition 3.3** *There exists an injective $\mathbb{C}$-algebra homomorphism*

$$\Psi \colon \ddot{U}^>_{q,d} \longrightarrow S \tag{3.30}$$

*such that $e_{i,r} \mapsto x^r_{i,1}$ for any $i \in [n], r \in \mathbb{Z}$.*

The following result was established in [20]:

**Theorem 3.6** *[20, Theorem 1.1] $\Psi \colon \ddot{U}^>_{q,d} \overset{\sim}{\longrightarrow} S$ of (3.30) is an isomorphism of $\mathbb{N}^{[n]} \times \mathbb{Z}$-graded $\mathbb{C}$-algebras.*

*Remark 3.6* We note that while the shuffle product of [20, §3.1] as well as the pole and wheel conditions of [20, §3.2] differ from (3.26) and (3.27, 3.28), the shuffle algebra $\mathcal{A}^+$ of *loc.cit.* is easily seen to be isomorphic to our $S$, see [20, Remark 3.4].

The proof of Theorem 3.6 is similar in spirit to its $\mathfrak{gl}_1$-counterpart (2.25). In analogy to Definition 2.2, the key tool of *loc.cit.* is the study of slope $\leq \mu$ subalgebras:

$$S^\mu = \bigoplus_{(\underline{k},d) \in \mathbb{N}^{[n]} \times \mathbb{Z}} S^\mu_{\underline{k},d} \quad \text{with} \quad S^\mu_{\underline{k},d} \quad \text{consisting of all} \quad F \in S_{\underline{k},d} \quad \text{for which}$$

$$\exists \lim_{\xi \to \infty} \frac{F(\dots, \xi \cdot x_{i,1}, \dots, \xi \cdot x_{i,\ell_i}, x_{i,\ell_i+1}, \dots, x_{i,k_i}, \dots)}{\xi^{\underline{\ell}|\mu}} \quad \forall \, 0 \leq \underline{\ell} \leq \underline{k} \,. \tag{3.31}$$

Similarly to Lemma 2.2, the subspace $S^\mu$ is easily seen to be a subalgebra of $S$. The key point is that each graded component $S^\mu_{\underline{k},d}$ is finite dimensional and we have:

**Proposition 3.4** *[20, Lemma 3.17] $\dim(S^\mu_{\underline{k},d}) \leq$ number of unordered collections*

$$L = \{(i_1, j_1, d_1), \dots, (i_t, j_t, d_t)\} \quad \text{such that}$$
$$[i_1, j_1 - 1] + \dots + [i_t, j_t - 1] = \underline{k}, \quad d_1 + \dots + d_t = d, \quad d_r \leq \mu(j_r - i_r), \tag{3.32}$$

*with $t \geq 1$, $i_r < j_r \in \mathbb{Z}$, $d_r \in \mathbb{Z}$, and the notation $[a,b] \in \mathbb{N}^{[n]}$ defined as in (3.35). Furthermore, we identify $(i_r, j_r, d_r) \sim (i_r - n, j_r - n, d_r)$ in the above counting.*

This result is proved similarly to Proposition 2.4 and uses the specialization maps $\phi_L$ of (3.38) utilized in our study of $\mathcal{A}(\underline{s})$ below (we disregard $d_r$'s in (3.38) as (3.32) forces all $d_r = 0$).

A more refined object from [20] is the following subalgebra of $S^\mu$:

$$B^\mu = \bigoplus_{\substack{\mu|\underline{k}| \in \mathbb{Z} \\ \underline{k} \in \mathbb{N}^{[n]}}} B^\mu_{\underline{k}} = \bigoplus_{\substack{\mu|\underline{k}| = d \in \mathbb{Z} \\ \underline{k} \in \mathbb{N}^{[n]}}} S^\mu_{\underline{k},d} \,. \tag{3.33}$$

Let $U^{\mathrm{DJ}>}_q(\widehat{\mathfrak{gl}}_n)$ be the "positive" subalgebra of the Drinfeld-Jimbo quantum affine $\mathfrak{gl}_n$, obtained from $U^{\mathrm{DJ}>}_q(\widehat{\mathfrak{sl}}_n)$ by adding central elements $\{P_k\}_{k \geq 1}$ (Heisenberg's half). For $\mu = \frac{b}{a} \in \mathbb{Q}$ with $\gcd(a,b) = 1$, the beautiful result [20, Lemma 3.15] identifies:

$$B^\mu \simeq U^{\mathrm{DJ}>}_q(\widehat{\mathfrak{gl}}_{\frac{n}{g}})^{\otimes g} \qquad \text{where} \qquad g = \gcd(n,a) \,. \tag{3.34}$$

### 3.2.2 Commutative subalgebras $\mathcal{A}(s_0, \ldots, s_{n-1})$

In this section, we introduce a natural counterpart of the commutative subalgebra $\mathcal{A}$ from Section 2.1.6, which now depends on $n-1$ additional parameters.

For any $0 \leq \underline{\ell} \leq \underline{k} \in \mathbb{N}^{[n]}$, $\xi \in \mathbb{C}^\times$, $F \in S_{\underline{k}}$, define $F_\xi^{\underline{\ell}} \in \mathbb{C}(x_{0,1}, \ldots, x_{n-1,k_{n-1}})$ via:

$$F_\xi^{\underline{\ell}} :=$$
$$F(\xi \cdot x_{0,1}, \ldots, \xi \cdot x_{0,\ell_0}, x_{0,\ell_0+1}, \ldots, x_{0,k_0}; \ldots; \xi \cdot x_{n-1,1}, \ldots, \xi \cdot x_{n-1,\ell_{n-1}}, \ldots, x_{n-1,k_{n-1}})$$

For any pair of integers $a \leq b$, we define the degree vector $\underline{\ell} := [a; b] \in \mathbb{N}^{[n]}$ via:

$$\underline{\ell} = (\ell_0, \ldots, \ell_{n-1}) \quad \text{with} \quad \ell_i = \#\left\{ c \in \mathbb{Z} \,\middle|\, a \leq c \leq b \text{ and } c \equiv i \pmod{n} \right\}. \quad (3.35)$$

For such a choice (3.35) of $\underline{\ell}$, we will denote $F_\xi^{\underline{\ell}}$ simply by $F_\xi^{(a,b)}$. We also define:

$$\partial^{(\infty;a,b)}F := \lim_{\xi \to \infty} F_\xi^{(a,b)} \quad \text{and} \quad \partial^{(0;a,b)}F := \lim_{\xi \to 0} F_\xi^{(a,b)} \quad \text{if they exist}. \quad (3.36)$$

**Definition 3.1** For any $\underline{s} = (s_0, \ldots, s_{n-1}) \in (\mathbb{C}^\times)^{[n]}$, consider an $\mathbb{N}^{[n]}$-graded subspace $\mathcal{A}(\underline{s}) \subset S$ whose degree $\underline{k} = (k_0, \ldots, k_{n-1})$ component is defined by:

$$\mathcal{A}(\underline{s})_{\underline{k}} := \left\{ F \in S_{\underline{k},0} \,\middle|\, \partial^{(\infty;a,b)}F = (s_a s_{a+1} \cdots s_b) \cdot \partial^{(0;a,b)}F \right.$$
$$\left. \text{for any } a \leq b \text{ such that } [a; b] \leq \underline{k} \right\},$$

where we again use the *cyclic* (that is, modulo $n$) notations $s_i := s_{i \bmod n}$.

A certain class of such elements is provided by the following result:

**Lemma 3.2** *For any $k \in \mathbb{N}$, $\mu \in \mathbb{C}$, and $\underline{s} \in (\mathbb{C}^\times)^{[n]}$, define $F_k^\mu(\underline{s}) \in S_{(k,\ldots,k)}$ by*

$$F_k^\mu(\underline{s}) :=$$
$$\frac{\prod_{i \in [n]} \prod_{1 \leq r \neq r' \leq k}(x_{i,r} - q^{-2}x_{i,r'}) \cdot \prod_{i \in [n]}(s_0 \cdots s_i \prod_{r=1}^k x_{i,r} - \mu \prod_{r=1}^k x_{i+1,r})}{\prod_{i \in [n]} \prod_{1 \leq r, r' \leq k}(x_{i,r} - x_{i+1,r'})},$$

*where we set $x_{n,r} := x_{0,r}$ as before. If $s_0 \cdots s_{n-1} = 1$, then $F_k^\mu(\underline{s}) \in \mathcal{A}(\underline{s})$.*

*Proof* It suffices to verify the claim for $\mu \neq 0, a = 0, b = nr + c$ with $0 \leq r < k$ and $0 \leq c \leq n-1$. Then $\ell_0 = \ldots = \ell_c = r+1, \ell_{c+1} = \ldots = \ell_{n-1} = r$. As $\xi \to \infty$, the function $F_k^\mu(\underline{s})_\xi^{(a,b)}$ grows at the speed $\xi^{\sum_{i \in [n]} \ell_i(\ell_{i+1} - \ell_i - 1) + \sum_{i \in [n]} \max\{\ell_i, \ell_{i+1}\}}$, and as $\xi \to 0$, the function $F_k^\mu(\underline{s})_\xi^{(a,b)}$ grows at the speed $\xi^{\sum_{i \in [n]} \ell_i(-\ell_{i+1} + \ell_i - 1) + \sum_{i \in [n]} \min\{\ell_i, \ell_{i+1}\}}$. For the above values of $\ell_i$, both powers of $\xi$ are zero and hence both limits $\partial^{(\infty;a,b)}F_k^\mu(\underline{s})$ and $\partial^{(0;a,b)}F_k^\mu(\underline{s})$ do exist. Moreover, for $\epsilon$ being 0 or $\infty$, we have $\partial^{(\epsilon;a,b)}F_k^\mu(\underline{s}) = (-1)^{\sum_{i \in [n]} \ell_i(\ell_i - \ell_{i-1})} q^{-2\sum_{i \in [n]} \ell_i(k - \ell_i)} \cdot G \cdot \prod_{i \in [n]} G_{\epsilon,i}$, where

$$G = \frac{\prod_{i \in [n]} \prod_{1 \le r \ne r' \le \ell_i} (x_{i,r} - q^{-2} x_{i,r'}) \cdot \prod_{i \in [n]} \prod_{\ell_i < r < r' \le k} (x_{i,r} - q^{-2} x_{i,r'})}{\prod_{i \in [n]} \prod_{\substack{1 \le r' \le \ell_{i+1} \\ 1 \le r \le \ell_i}} (x_{i,r} - x_{i+1,r'}) \cdot \prod_{i \in [n]} \prod_{\substack{\ell_{i+1} < r' \le k \\ \ell_i < r \le k}} (x_{i,r} - x_{i+1,r'})},$$

$$G_{\infty,i} = \prod_{r=1}^{\ell_i} x_{i,r}^{\ell_{i+1} + \ell_{i-1} - 2\ell_i} \cdot \begin{cases} s_0 \cdots s_i \prod_{r=1}^{k} x_{i,r} - \mu \prod_{r=1}^{k} x_{i+1,r} & \text{if } \ell_i = \ell_{i+1} \\ s_0 \cdots s_i \prod_{r=1}^{k} x_{i,r} & \text{if } \ell_i > \ell_{i+1}, \\ -\mu \prod_{r=1}^{k} x_{i+1,r} & \text{if } \ell_i < \ell_{i+1} \end{cases}$$

$$G_{0,i} = \prod_{r=\ell_i+1}^{k} x_{i,r}^{-\ell_{i+1} - \ell_{i-1} + 2\ell_i} \cdot \begin{cases} s_0 \cdots s_i \prod_{r=1}^{k} x_{i,r} - \mu \prod_{r=1}^{k} x_{i+1,r} & \text{if } \ell_i = \ell_{i+1} \\ -\mu \prod_{r=1}^{k} x_{i+1,r} & \text{if } \ell_i > \ell_{i+1}. \\ s_0 \cdots s_i \prod_{r=1}^{k} x_{i,r} & \text{if } \ell_i < \ell_{i+1} \end{cases}$$

The desired equality $\partial^{(\infty;a,b)} F_k^\mu(\underline{s}) = \prod_{i=0}^{c} s_i \cdot \partial^{(0;a,b)} F_k^\mu(\underline{s})$ follows immediately from these formulas, while the condition tot.deg$(F_k^\mu(\underline{s})) = 0$ is obvious.  □

We call $\underline{s} = (s_0, \dots, s_{n-1}) \in (\mathbb{C}^\times)^{[n]}$ satisfying $s_0 \cdots s_{n-1} = 1$ *generic* if

$$s_0^{a_0} \cdots s_{n-1}^{a_{n-1}} \in q^{\mathbb{Z}} \cdot d^{\mathbb{Z}} \ (a_0, \dots, a_{n-1} \in \mathbb{Z}) \iff a_0 = \dots = a_{n-1}. \tag{3.37}$$

For such $\underline{s} \in (\mathbb{C}^\times)^{[n]}$, the following description of $\mathcal{A}(\underline{s})$ was established in [12]:

**Theorem 3.7** *For a generic $\underline{s} = (s_0, \dots, s_{n-1}) \in (\mathbb{C}^\times)^{[n]}$ satisfying $s_0 \cdots s_{n-1} = 1$, the space $\mathcal{A}(\underline{s})$ is shuffle-generated by $\{F_k^\mu(\underline{s}) | k \ge 1, \mu \in \mathbb{C}\}$. Moreover, $\mathcal{A}(\underline{s})$ is a polynomial algebra in generators $\{F_k^{\mu_l}(\underline{s}) | k \ge 1, 1 \le l \le n\}$ for arbitrary pairwise distinct $\mu_1, \dots, \mu_n \in \mathbb{C}$. In particular, $\mathcal{A}(\underline{s})$ is a commutative subalgebra of $S$.*

The proof of this result is inspired by that of Theorem 2.4 but is more tedious and is presented in the next subsection.

### 3.2.3 Proof of Theorem 3.7

The proof of Theorem 3.7 will proceed in several steps. First, we will use an analogue of the *Gordon filtration* (2.29) (already used in the proof of Theorem 3.6 as presented in [20]) to obtain an upper bound on dimensions of $\mathcal{A}(\underline{s})_k$. Next, we will show that the subalgebra $\mathcal{A}'(\underline{s}) \subset S$, shuffle generated by all $F_k^\mu(\underline{s})$, is contained in $\mathcal{A}(\underline{s})$. We will use another filtration to argue that the dimension of $\mathcal{A}'(\underline{s})_k$ is at least as big as the upper bound for the dimension of $\mathcal{A}(\underline{s})_k$, implying the equality $\mathcal{A}'(\underline{s}) = \mathcal{A}(\underline{s})$. Finally, the commutativity of $F_k^\mu(\underline{s})$, and hence of $\mathcal{A}(\underline{s})$, will immediately follow from their Bethe algebra realization presented in Section 3.5.4, see Corollary 3.4.

**Lemma 3.3** *Consider a polynomial algebra $\mathcal{R} = \mathbb{C}[\{T_{i,m}\}_{1 \le i \le n}^{m \ge 1}]$, endowed with an $\mathbb{N}$-grading via $\deg(T_{i,m}) = m$. Then:*

*(a) For $\underline{k} = k\delta := (k, \dots, k) \in \mathbb{N}^{[n]}$, we have $\dim \mathcal{A}(\underline{s})_{\underline{k}} \le \dim \mathcal{R}_k$.*

*(b) For $\underline{k} \notin \{0, \delta, 2\delta, \dots\}$, we have $\mathcal{A}(\underline{s})_{\underline{k}} = 0$.*

*Proof* An unordered set $L$ of integer intervals $\{[a_1; b_1], \ldots, [a_r; b_r]\}$ is called a *partition* of $\underline{k} \in \mathbb{N}^{[n]}$, denoted by $L \vdash \underline{k}$, if $\underline{k} = [a_1; b_1] + \cdots + [a_r; b_r]$. We order the elements of $L$ so that $b_1 - a_1 \geq b_2 - a_2 \geq \cdots \geq b_r - a_r$. Two such sets $L$ and $L'$ are said to be equivalent if $|L| = |L'|$ and we can order their elements so that $b_i' - b_i = a_i' - a_i = nc_i$ for all $i$ and some $c_i \in \mathbb{Z}$. We note that the collection of $L \vdash \underline{k}$, up to this equivalence, is finite for any $\underline{k} \in \mathbb{N}^{[n]}$. Finally, we say $L' > L$ if there exists $s$, such that $b_s' - a_s' > b_s - a_s$ and $b_t' - a_t' = b_t - a_t$ for $1 \leq t < s$.

Similarly to (1.40), any partition $L \vdash \underline{k}$ gives rise to a **specialization map**

$$\phi_L \colon \mathcal{A}(\underline{s})_{\underline{k}} \longrightarrow \mathbb{C}[y_1^{\pm 1}, \ldots, y_r^{\pm 1}]. \tag{3.38}$$

To this end, we split the $x_{*,*}$-variables into $r$ groups, corresponding to the above intervals in $L$, and specialize those corresponding to the interval $[a_t; b_t]$ $(1 \leq t \leq r)$ in the natural order to:

$$(qd)^{-a_t} \cdot y_t, \ \ldots, \ (qd)^{-b_t} \cdot y_t$$

For $F = \frac{f(x_{0,1}, \ldots, x_{n-1,k_{n-1}})}{\prod_{i \in [n]} \prod_{1 \leq r \leq k_{i+1}}^{1 \leq r' \leq k_i} (x_{i,r} - x_{i+1,r'})} \in \mathcal{A}(\underline{s})_{\underline{k}}$, we define $\phi_L(F)$ as the corresponding specialization of $f$ (which is independent of our splitting of the $x_{*,*}$-variables since $f$ is symmetric). This gives rise to the filtration on $\mathcal{A}(\underline{s})_{\underline{k}}$ via:

$$\mathcal{A}(\underline{s})_{\underline{k}}^L := \bigcap_{L' > L} \mathrm{Ker}(\phi_{L'}). \tag{3.39}$$

We shall now consider the images $\phi_L(F)$ for $F \in \mathcal{A}(\underline{s})_{\underline{k}}^L$.

First, we note that the wheel conditions (3.28) for $F$ guarantee that $\phi_L(F)$, which is a Laurent polynomial in $\{y_1, \ldots, y_r\}$, vanishes under the following specializations:

(i) $(qd)^{-x'} y_{t'} = (q/d)(qd)^{-x} y_t$
   for any $1 \leq t < t' \leq r$, $a_t \leq x < b_t$, $a_{t'} \leq x' \leq b_{t'}$ with $x' \equiv x + 1 \bmod n$,

(ii) $(qd)^{-x'} y_{t'} = (d/q)(qd)^{-x} y_t$
   for any $1 \leq t < t' \leq r$, $a_t < x \leq b_t$, $a_{t'} \leq x' \leq b_{t'}$ with $x' \equiv x - 1 \bmod n$.

Second, the condition $\phi_{L'}(F) = 0$ for any $L' > L$ also implies that $\phi_L(F)$ vanishes under the following specializations:

(iii) $(qd)^{-x'} y_{t'} = (qd)^{-b_t-1} y_t$
   for any $1 \leq t < t' \leq r$, $a_{t'} \leq x' \leq b_{t'}$ with $x' \equiv b_t + 1 \bmod n$,

(iv) $(qd)^{-x'} y_{t'} = (qd)^{-a_t+1} \cdot y_t$
   for any $1 \leq t < t' \leq r$, $a_{t'} \leq x' \leq b_{t'}$ with $x' \equiv a_t - 1 \bmod n$.

Combining the vanishing conditions (i)–(iv) for $\phi_L(F)$ with $F \in \mathcal{A}(\underline{s})_{\underline{k}}^L$, we conclude that the specialization $\phi_L(F)$ is divisible by a polynomial $Q_L \in \mathbb{C}[y_1, \ldots, y_r]$, defined as a product of the linear terms in $y_t$'s arising from (i)–(iv), counted with the correct multiplicities.

For $1 \leq t \leq r$, we define $\underline{\ell}^t := [a_t; b_t] \in \mathbb{N}^{[n]}$. Then, the total degree of $Q_L$ equals

$$\text{tot.deg}(Q_L) = \sum_{1 \le t < t' \le r} \sum_{i \in [n]} (\ell_i^t \ell_{i+1}^{t'} + \ell_i^t \ell_{i-1}^{t'}) = \sum_{i \in [n]} k_i k_{i+1} - \sum_{1 \le t \le r} \sum_{i \in [n]} \ell_i^t \ell_{i+1}^t , \tag{3.40}$$

while the degree with respect to each variable $y_t$ $(1 \le t \le r)$ is given by:

$$\text{deg}_{y_t}(Q_L) = \sum_{i \in [n]} \left( \ell_i^t (k_{i-1} + k_{i+1}) - 2\ell_i^t \ell_{i+1}^t \right) . \tag{3.41}$$

On the other hand, for the images $\phi_L(F)$ with $F \in \mathcal{A}(\underline{s})_{\underline{k}}^L$, we have:

- Due to $F \in S_{\underline{k},0}$, the total degree of $\phi_L(F)$ is

$$\text{tot.deg}(\phi_L(F)) = \sum_{i \in [n]} k_i k_{i+1} . \tag{3.42}$$

- For each $1 \le t \le r$, the degree of $\phi_L(F)$ with respect to $y_t$ is bounded by

$$\text{deg}_{y_t}(\phi_L(F)) \le \sum_{i \in [n]} \left( \ell_i^t (k_{i-1} + k_{i+1}) - \ell_i^t \ell_{i+1}^t \right) \tag{3.43}$$

due to the existence of the limit $\partial^{(\infty; a_t, b_t)} F$ from (3.36), cf. (2.32).

Define $r_L := \phi_L(F)/Q_L \in \mathbb{C}[y_1^{\pm 1}, \ldots, y_r^{\pm 1}]$. Combining (3.40)–(3.43), we obtain:

$$\text{tot.deg}(r_L) = \sum_{1 \le t \le r} \sum_{i \in [n]} \ell_i^t \ell_{i+1}^t \quad \text{and} \quad \text{deg}_{y_t}(r_L) \le \sum_{i \in [n]} \ell_i^t \ell_{i+1}^t .$$

Thus, $r_L = \nu \cdot \prod_{t=1}^r y_t^{\sum_{i \in [n]} \ell_i^t \ell_{i+1}^t}$ with $\nu \in \mathbb{C}$, so that $\phi_L(F) = \nu \cdot \prod_{t=1}^r y_t^{\sum_{i \in [n]} \ell_i^t \ell_{i+1}^t} \cdot Q_L$.

On the other hand, applying the above specialization to the entire $F$ rather than its numerator $f$ of (3.27), we get $\phi_L(F)/Q$ with $Q \in \mathbb{C}[y_1, \ldots, y_r]$ given by

$$Q = \nu' \cdot \prod_{t=1}^r y_t^{\sum_{i \in [n]} \ell_i^t \ell_{i+1}^t} \cdot \prod_{1 \le t < t' \le r} \prod_{a_t \le x \le b_t} \prod_{a_{t'} \le x' \le b_{t'}}^{x' \equiv x \pm 1} \left( (qd)^{-x} y_t - (qd)^{-x'} y_{t'} \right)$$

for some $\nu' \in \mathbb{C}^\times$. The condition $F \in \mathcal{A}(\underline{s})$ implies:

$$\lim_{\xi \to \infty} \left( \frac{\phi_L(F)}{Q} \right)_{|y_t \mapsto \xi \cdot y_t} = s_{a_t} \cdots s_{b_t} \cdot \lim_{\xi \to 0} \left( \frac{\phi_L(F)}{Q} \right)_{|y_t \mapsto \xi \cdot y_t} \quad \forall \, 1 \le t \le r .$$

For $\nu \ne 0$, this forces $s_{a_t} \cdots s_{b_t} \in q^{\mathbb{Z}} \cdot d^{\mathbb{Z}}$. We thus get $b_t - a_t + 1 = nc_t$ for every $1 \le t \le r$ and some $c_t \in \mathbb{N}$, since $\underline{s}$ is generic (3.37). This implies part (b) of the lemma. Likewise, part (a) follows from the inequality $\dim \mathcal{A}(\underline{s})_{\underline{k}} \le \sum \dim \phi_L(\mathcal{A}(\underline{s})_{\underline{k}}^L)$ with the sum taken over all equivalence classes of $L \vdash \underline{k} = k\delta$. To this end, we note that each summand $\dim \phi_L(\mathcal{A}(\underline{s})_{\underline{k}}^L)$ equals $0$ or $1$, the latter being possible only if all $b_t - a_t + 1$ are divisible by $n$. Associating to such $L$ a monomial $\prod_{t=1}^r T_{i_t, m_t} \in \mathcal{R}$ with $i_t \equiv a_t$ and $m_t = \frac{b_t - a_t + 1}{n}$, we get the desired inequality $\dim \mathcal{A}(\underline{s})_{k\delta} \le \dim \mathcal{R}_k$. $\square$

**Lemma 3.4** *Let $\mathcal{A}'(\underline{s})$ denote the subalgebra of $S$ shuffle-generated by $\{F_k^\mu(\underline{s})\}_{k \geq 1}^{\mu \in \mathbb{C}}$. Then: $\mathcal{A}'(\underline{s}) \subseteq \mathcal{A}(\underline{s})$.*

*Proof* It suffices to show that $F_{k_1,\dots,k_N}^{\mu_1,\dots,\mu_N}(\underline{s}) := F_{k_1}^{\mu_1}(\underline{s}) \star \cdots \star F_{k_N}^{\mu_N}(\underline{s})$ are in $\mathcal{A}(\underline{s})$ for any $N$, $k_r \geq 1$, $\mu_r \in \mathbb{C}^\times$. The $N = 1$ case is due to Lemma 3.2, while the arguments for a general $N$ are similar. To this end, choose any $a \leq b$ such that $[a; b] \leq k\delta$, where $k := k_1 + \cdots + k_N$. We can further assume $a = 0$. Let us consider any summand from the symmetrization appearing in $F_{k_1,\dots,k_N}^{\mu_1,\dots,\mu_N}(\underline{s})$ with $\underline{\ell} := [a; b]$ variables being multiplied by $\xi$. We claim that as $\xi$ tends to $\infty$ or $0$ both limits (3.36) exist and differ by the factor $s_a \cdots s_b$.

To prove this, we consider the decomposition $\underline{\ell} = \underline{\ell}^1 + \cdots + \underline{\ell}^N$ with $\underline{\ell}^t$ denoting the number of $x_{*,*}$-variables that are multiplied by $\xi$ and get placed into $F_{k_t}^{\mu_t}(\underline{s})$ in the given summand. As noticed in the proof of Lemma 3.2, the function $F_{k_t}^{\mu_t}(\underline{s})\frac{\underline{\ell}^t}{\xi}$ grows at the speed $\xi^{\sum_{i\in[n]}\ell_i^t(\ell_{i+1}^t-\ell_i^t-1)+\sum_{i\in[n]}\max\{\ell_i^t,\ell_{i+1}^t\}}$ as $\xi \to \infty$ and at the speed $\xi^{\sum_{i\in[n]}\ell_i^t(-\ell_{i+1}^t+\ell_i^t-1)+\sum_{i\in[n]}\min\{\ell_i^t,\ell_{i+1}^t\}}$ as $\xi \to 0$. To bound these powers, we note that $(a-b)(a-b-1) \geq 0$ for any $a, b \in \mathbb{Z}$, thus implying

$$\min\{a,b\} + \frac{a^2 + b^2 - a - b}{2} \geq ab \qquad (3.44)$$

with the equality taking place iff $a - b \in \{-1, 0, 1\}$. Thus, for $1 \leq t \leq N$ we have:

$$\sum_{i\in[n]} \ell_i^t(\ell_{i+1}^t - \ell_i^t - 1) + \sum_{i\in[n]} \max\left\{\ell_i^t, \ell_{i+1}^t\right\} \leq 0,$$

$$\sum_{i\in[n]} \ell_i^t(-\ell_{i+1}^t + \ell_i^t - 1) + \sum_{i\in[n]} \min\left\{\ell_i^t, \ell_{i+1}^t\right\} \geq 0,$$

with the equalities only when $\ell_i^t - \ell_{i+1}^t \in \{-1, 0, 1\}$ for any $i \in [n]$. As the limits of

$$\zeta_{i,j}\left(\frac{\xi \cdot x}{y}\right), \quad \zeta_{i,j}\left(\frac{x}{\xi \cdot y}\right), \quad \zeta_{i,j}\left(\frac{\xi \cdot x}{\xi \cdot y}\right) \quad \text{for} \quad \xi \to 0, \infty$$

exist for any $i, j \in [n]$, the $\xi \to 0, \infty$ limits of the corresponding summand in the symmetrization also exist and vanish if $|\ell_i^t - \ell_{i+1}^t| > 1$ for some $1 \leq t \leq N$, $i \in [n]$.

Finally, if $|\ell_i^t - \ell_{i+1}^t| \leq 1$ for any $i \in [n]$ and $1 \leq t \leq N$, the formulas from the proof of Lemma 3.2 imply that the ratio of the limits as $\xi$ tends to $\infty$ and $0$ equals

$$\prod_{1\leq t\leq N} \prod_{i\in[n]} \left(\frac{s_0 \cdots s_i}{-\mu_t}\right)^{\ell_i^t - \ell_{i+1}^t} = \prod_{i\in[n]} (s_0 \cdots s_i)^{\ell_i - \ell_{i+1}} = s_a \cdots s_b \,.$$

This completes the proof of $F_{k_1,\dots,k_N}^{\mu_1,\dots,\mu_N}(\underline{s}) \in \mathcal{A}(\underline{s})$. $\qquad\square$

**Lemma 3.5** $\dim \mathcal{A}'(\underline{s})_{k\delta} \geq \dim \mathcal{R}_k$ *for any $k \geq 1$.*

*Proof* Let us choose any pairwise distinct $\mu_1, \dots, \mu_n \in \mathbb{C}$ and consider a subspace $\mathcal{A}''(\underline{s})$ of $\mathcal{A}'(\underline{s})$ spanned by $F_{\underline{\pi}}^{\underline{\mu}}(\underline{s}) = F_{k_1,\dots,k_N}^{\mu_{i_1},\dots,\mu_{i_N}}(\underline{s})$ with $N \geq 0$, $1 \leq i_1, \dots, i_N \leq n$, and $\pi = \{k_1 \geq \cdots \geq k_N > 0\}$. It suffices to show that

$$\dim \mathcal{A}''(\underline{s})_{k\delta} \geq \dim \mathcal{R}_k . \tag{3.45}$$

To prove (3.45), for any Young diagram $\lambda = (\lambda_1, \ldots, \lambda_l)$, we introduce the *vertical specialization map*

$$\varphi_\lambda \colon S_{|\lambda| \cdot \delta} \longrightarrow \mathbb{C}\left(\{y_{i,t}\}_{i \in [n]}^{1 \leq t \leq l}\right) \tag{3.46}$$

by specializing the $x_{*,*}$-variables as follows:

$$x_{i, \lambda_1 + \cdots + \lambda_{t-1} + a} \mapsto q^{2a} y_{i,t} \quad \text{for any} \quad 1 \leq t \leq l,\ 1 \leq a \leq \lambda_t,\ i \in [n] .$$

It is clear that for any $\pi = \{k_1 \geq \cdots \geq k_N > 0\}$ with $|\pi| = k_1 + \cdots + k_N = k = |\lambda|$ and $\lambda > \pi'$ (where $>$ denotes the dominance order on size $k$ Young diagrams and $\pi'$ denotes the transposed Young diagram, as before), we have (cf. Lemmas 1.5, 2.3):

$$\varphi_\lambda(F_\pi^\mu(\underline{s})) = 0 \quad \text{for any} \quad \underline{\mu} \in \mathbb{C}^N . \tag{3.47}$$

Therefore, it suffices to prove (cf. Lemma 1.6):

$$\sum_{\pi \vdash k} \dim \varphi_{\pi'}\left(\mathrm{span}\left\{F_\pi^{\mu_{i_1}, \ldots, \mu_{i_N}}(\underline{s}) \,\Big|\, 1 \leq i_1, \ldots, i_N \leq n\right\}\right) \geq \dim \mathcal{R}_k . \tag{3.48}$$

To this end, we first consider the case $k_1 = \ldots = k_N \Rightarrow k = N k_1 = |\pi|$. Then

$$\varphi_{\pi'}\left(F_\pi^\mu(\underline{s})\right) = Z \cdot \prod_{r=1}^N \prod_{i \in [n]} \left(s_0 \cdots s_i \prod_{t=1}^{k_1} y_{i,t} - \mu_{i_r} \prod_{t=1}^{k_1} y_{i+1,t}\right) \tag{3.49}$$

for a certain nonzero common factor $Z$. Define $Y_i := y_{i,1} \cdots y_{i,k_1}$ for $i \in [n]$. Since

$$f_r(Y_1, \ldots, Y_n) := \prod_{i \in [n]} (s_0 \cdots s_i Y_i - \mu_r Y_{i+1}), \quad 1 \leq r \leq n ,$$

are algebraically independent, we see that $\dim \mathrm{Span}\{\prod_{r=1}^N f_{i_r}(Y_1, \ldots, Y_n)\}$ is $\geq$ the number of degree $N$ monomials in $\{T_{i,k_1}\}_{i=1}^n$. Generalizing this argument to a general $\pi$ implies the stated inequality (3.48). □

By Lemmas 3.3–3.5, the subspace $\mathcal{A}(\underline{s})$ is shuffle-generated by $\{F_k^\mu(\underline{s})\}_{k \geq 1}^{\mu \in \mathbb{C}}$ and thus it is a subalgebra of $S$. Moreover, it has the prescribed dimensions of each graded component, with the equality taking place in Lemma 3.3(a). The following result will be derived in Section 3.5.4 from the transfer matrix realization of $F_k^\mu(\underline{s})$:

**Lemma 3.6** *The elements* $\{F_k^\mu(\underline{s})\}_{k \geq 1}^{\mu \in \mathbb{C}}$ *pairwise commute.*

The above Lemmas 3.3–3.6 complete our proof of Theorem 3.7.

*Remark 3.7* We note that the proof of Lemma 3.3 implies that $\mathcal{A}(\underline{s}) = \mathbb{C}$ for any $\underline{s} = (s_0, \ldots, s_{n-1}) \in (\mathbb{C}^\times)^{[n]}$ such that $\prod_{i \in [n]} s_i^{a_i} \notin q^{\mathbb{Z}} \cdot d^{\mathbb{Z}}$ unless all integers $a_i = 0$.

### 3.2.4  Special limit $\mathcal{A}^h$

Let us consider a special limit $\mathcal{A}^h$ of $\mathcal{A}(\underline{s})$, corresponding to the case $s_1, \ldots, s_{n-1} \to 0$ and $s_0 = 1/(s_1 \cdots s_{n-1})$. Explicitly, we define $\mathcal{A}^h$ as follows:

**Definition 3.2** For $\underline{k} \in \mathbb{N}^{[n]}$, let $\mathcal{A}^h_{\underline{k}} \subset S_{\underline{k}}$ be the subspace of all $F \in S_{\underline{k},0}$ such that:

- $\partial^{(0;a,b)} F$ exists for any $a \leq b$ satisfying $\underline{\ell} := [a;b] \leq \underline{k}$, and furthermore vanishes whenever $\ell_0 > \min\{\ell_i\}_{i \in [n]}$.
- $\partial^{(\infty;a,b)} F$ exists for any $a \leq b$ satisfying $\underline{\ell} := [a;b] \leq \underline{k}$, and furthermore vanishes whenever $\ell_0 < \max\{\ell_i\}_{i \in [n]}$.
- $\partial^{(0;a,b)} F = \partial^{(\infty;a,b)} F$ for any $a \leq b$ satisfying $\ell\delta = [a;b] \leq \underline{k}$.

We set $\mathcal{A}^h := \bigoplus_{\underline{k} \in \mathbb{N}^{[n]}} \mathcal{A}^h_{\underline{k}}$.

For any $p \in [n]$ and $k \in \mathbb{N}$, consider the elements $\Gamma^0_{p;k} \in S_{k\delta,0}$ defined via:

$$\Gamma^0_{p;k} := \frac{\prod_{i \in [n]} \prod_{1 \leq r \neq r' \leq k} (x_{i,r} - q^{-2} x_{i,r'}) \cdot \prod_{i \in [n]} \prod_{r=1}^{k} x_{i,r}}{\prod_{i \in [n]} \prod_{1 \leq r, r' \leq k} (x_{i,r} - x_{i+1,r'})} \cdot \prod_{r=1}^{k} \frac{x_{0,r}}{x_{p,r}} \, . \quad (3.50)$$

*Remark 3.8* The elements (3.50) arise as the leading terms of the $\mu^{n-p}$-coefficients in $F^{\mu}_k(\underline{s})$ from Lemma 3.2 in the aforementioned limit $s_1, s_2, \ldots, s_{n-1} \to 0$.

It is straightforward to see that $\Gamma^0_{p;k} \in \mathcal{A}^h$, cf. Lemma 3.2. In fact, we have the following description of $\mathcal{A}^h$ (proved completely analogously to Theorem 3.7):

**Theorem 3.8** *The subspace $\mathcal{A}^h$ is a subalgebra of $S$. Moreover, $\mathcal{A}^h$ is a polynomial algebra in the generators $\{\Gamma^0_{p;k}\}^{k \geq 1}_{p \in [n]}$ of (3.50). In particular, $\mathcal{A}^h$ is a commutative subalgebra of $S$.*

Combining this result with Theorem 3.19 from the end of this chapter, we obtain:

**Theorem 3.9** *The restriction of the algebra isomorphism $\Psi \colon \ddot{U}^{>}_{q,d} \xrightarrow{\sim} S$ from Theorem 3.6 gives rise to an algebra isomorphism $\Psi \colon \dot{U}^{h>}_q(\mathfrak{h}_n) \xrightarrow{\sim} \mathcal{A}^h$.*

*Remark 3.9* We note that $\mathcal{A}^h = \bigoplus_{k \in \mathbb{N}} \mathcal{A}^h_{k\delta}$ and $\mathcal{A}^h_{k\delta} \subset S_{k\delta,0}$ can be described only by the first condition of Definition 3.2 (the other two being implied).

### 3.2.5  Identification of two extra horizontal Heisenbergs

Recall the *slope filtration* $S = \bigcup_{\mu \in \mathbb{R}} S^{\mu}$ of (3.31). For a slope $\mu = 0$, consider the subalgebra $B^0 = \bigoplus_{\underline{k} \in \mathbb{N}^{[n]}} B^0_{\underline{k}}$ with $B^0_{\underline{k}} := S^0 \cap S_{\underline{k},0}$, see (3.33). Explicitly, we have:

$$F \in B^0_{\underline{k}} \iff F \in S_{\underline{k},0} \quad \text{and} \quad \exists \lim_{\underline{\xi} \to \infty} F_{\underline{\xi}}^{\underline{\ell}} \quad \forall \, 0 \leq \underline{\ell} \leq \underline{k} \, . \quad (3.51)$$

Let $\dot{U}_q^{\mathrm{h}>}(\mathfrak{sl}_n)$ be the subalgebra of $\dot{U}_q^{\mathrm{h}}(\mathfrak{sl}_n)$ generated by $\{e_{i,0}\}_{i\in[n]}$. In [20] it was enlarged to $\dot{U}_q^{\mathrm{h}>}(\mathfrak{gl}_n)$, isomorphic to the "positive" subalgebra $U_q^{\mathrm{DJ}>}(\widehat{\mathfrak{gl}}_n)$ of the quantum group $U_q^{\mathrm{DJ}}(\widehat{\mathfrak{gl}}_n)$, by adjoining a family of commuting elements $\{X_k\}_{k\geq 1}$ (the subalgebra generated by $X_k$ coincides with the subalgebra generated by $h_k^{\mathrm{h}} = \overline{\varpi}(h_{-k}^{\vee})$ from Section 3.1.5, due to Corollary 3.1). As the key ingredient in the proof of Theorem 3.6, the following result was established in [20], cf. (3.34):

**Proposition 3.5** [20] (a) *The restriction of the algebra embedding* $\Psi \colon \ddot{U}_{q,d}^{>} \hookrightarrow S$ *of Proposition 3.3 gives rise to the algebra isomorphism* $\Psi \colon \dot{U}_q^{\mathrm{h}>}(\mathfrak{gl}_n) \xrightarrow{\sim} B^0$.

(b) *The images* $P_k := \Psi(X_k)$ *are uniquely (up to a constant) characterized by:*

$$\lim_{\xi\to\infty}(P_k)_{\overline{\xi}}^{\underline{\ell}} = 0 \quad \forall\, 0 < \underline{\ell} < k\delta \qquad \text{and} \qquad P_k \text{ is } \mathbb{Z}/n\mathbb{Z} - \text{invariant}, \qquad (3.52)$$

*where the cyclic group* $\mathbb{Z}/n\mathbb{Z}$ *acts on* $S_{k\delta,0}$ *by permuting* $x_{i,r} \mapsto x_{i+1,r}$ *for all* $i, r$.

As an immediate corollary, we get the following result:

**Theorem 3.10** $\Psi^{-1}(\mathcal{A}(\underline{s})) \subset \dot{U}_q^{\mathrm{h}>}(\mathfrak{gl}_n)$ *for generic* $\{s_i\}_{i\in[n]}$ *such that* $s_0\cdots s_{n-1} = 1$.

*Proof* By Theorem 3.7 and Proposition 3.5(a), it suffices to show that $F_k^{\mu}(\underline{s}) \in B^0$. The latter is equivalent to the existence of limits $\lim\limits_{\xi\to\infty}(F_k^{\mu}(\underline{s}))_{\overline{\xi}}^{\underline{\ell}}$ for all $0 \leq \underline{\ell} \leq k\delta$. As $\xi \to \infty$, the function $(F_k^{\mu}(\underline{s}))_{\overline{\xi}}^{\underline{\ell}}$ grows at the speed $\xi^{\sum_{i\in[n]}\ell_i(\ell_{i+1}-\ell_i+1)-\sum_{i\in[n]}\min\{\ell_i,\ell_{i+1}\}}$ (see the proof of Lemma 3.2). Since $\sum_{i\in[n]}\ell_i(\ell_{i+1}-\ell_i+1) - \sum_{i\in[n]}\min\{\ell_i,\ell_{i+1}\} = \sum_{i\in[n]}\left(\ell_i\ell_{i+1} - \frac{\ell_i^2+\ell_{i+1}^2-\ell_i-\ell_{i+1}}{2} - \min\{\ell_i,\ell_{i+1}\}\right)$ and each summand is nonpositive by (3.44), so is the aforementioned power of $\xi$. Hence, the limits $\lim\limits_{\xi\to\infty}(F_k^{\mu}(\underline{s}))_{\overline{\xi}}^{\underline{\ell}}$ exist and are finite for any $0 \leq \underline{\ell} \leq k\delta$. This completes the proof. $\qquad\square$

We conclude this section by providing explicit formulas for $P_k$. To this end, define:

$$F_k := \frac{\prod_{i\in[n]}\prod_{1\leq r\neq r'\leq k}(q^{-1}x_{i,r}-qx_{i,r'})\cdot\prod_{i\in[n]}\prod_{1\leq r\leq k}x_{i,r}}{\prod_{i\in[n]}\prod_{1\leq r,r'\leq k}(x_{i+1,r'}-x_{i,r})} \in S_{k\delta,0} \quad (3.53)$$

for $k > 0$, and $F_0 := \mathbf{1} \in S_{(0,\dots,0),0}$. We note that $F_k$ is a multiple of $F_k^0(\underline{s}) \in \mathcal{A}(\underline{s})$ whenever $\prod_{i\in[n]}s_i = 1$. Furthermore, $F_k$ is clearly $\mathbb{Z}/n\mathbb{Z}$-invariant. We also define

$$L_k \in S_{k\delta,0} \qquad \text{via} \qquad \exp\left(\sum_{k\geq 1}L_k z^{-k}\right) = \sum_{k\geq 0}F_k z^{-k}. \qquad (3.54)$$

The relevant properties of these elements are stated in the next theorem:

**Theorem 3.11** (a) *For* $\underline{\ell} \notin \{0, \delta, 2\delta, \dots, k\delta\}$, *we have* $\lim\limits_{\xi\to\infty}(F_k)_{\overline{\xi}}^{\underline{\ell}} = 0$.

(b) *For any* $0 \leq \ell \leq k$, *we have* $\lim\limits_{\xi\to\infty}(F_k)_{\overline{\xi}}^{\ell\delta} = F_\ell(\{x_{i,r}\}_{i\in[n]}^{r\leq\ell})\cdot F_{k-\ell}(\{x_{i,r}\}_{i\in[n]}^{r>\ell})$.

(c) *For any* $0 < \underline{\ell} < k\delta$, *we have* $\lim\limits_{\xi\to\infty}(L_k)_{\overline{\xi}}^{\underline{\ell}} = 0$.

*Proof* For any $0 \leq \underline{\ell} \leq k\delta$, the function $(F_k)\frac{\ell}{\xi}$ grows at the speed $\xi^{\sum_{i\in[n]}\ell_i(\ell_{i+1}-\ell_i)}$ as $\xi \to \infty$. Note that $\sum_{i\in[n]}\ell_i(\ell_{i+1}-\ell_i) = -\frac{1}{2}\sum_{i\in[n]}(\ell_i-\ell_{i+1})^2 \leq 0$ and the equality takes place iff $\ell_0 = \ldots = \ell_{n-1} \Leftrightarrow \underline{\ell} \in \{0, \delta, 2\delta, \ldots, k\delta\}$. This implies part (a). Part (b) is straightforward. Finally, part (c) follows directly from parts (a, b), due to the *logarithmical* definition (3.54) of $L_k$. □

The shuffle elements $L_k$'s are clearly $\mathbb{Z}/n\mathbb{Z}$-invariant as so are $F_k$'s. Combining this with Theorem 3.11(c), the characterization from Proposition 3.5(b), and Theorem 3.19 from the end of this chapter, we obtain:

**Corollary 3.1** *$P_k$ is a nonzero multiple of $L_k$ from (3.54). The subalgebra generated by $\{X_k\}$ coincides with the subalgebra generated by $h_k^{\mathrm{h}} = \overline{\varpi}(h_{-k}^{\vee})$ from Section 3.1.5. In particular, the horizontal quantum affine $\mathfrak{gl}_n$ of [9] and [20] actually coincide.*

## 3.3 Representations of quantum toroidal $\mathfrak{sl}_n$

In this section, following [9, 22, 23], we provide three different constructions of the modules over quantum toroidal algebras of $\mathfrak{sl}_n$ that are reminiscent of the corresponding constructions for the quantum toroidal algebras of $\mathfrak{gl}_1$ from Chapter 2. In fact, these will be actually modules of either $'\ddot{U}_{q,d}(\mathfrak{sl}_n)$ or $'\ddot{U}'_{q,d}(\mathfrak{sl}_n)$, as no $\ddot{U}_{q,d}(\mathfrak{sl}_n)$-modules are presently known (in the presence of the degree operators $q^{d_1}, q^{d_2}$) and in all known $U_{q,d}^{\mathrm{tor}}(\mathfrak{sl}_n)$-modules one of the two central elements $c, c'$ acts trivially.

### 3.3.1 Vector, Fock, and Macmahon modules

The general theory of the category $O^-$ and the lowest weight modules for $'\ddot{U}_{q,d}(\mathfrak{sl}_n)$ is similar to that for $'\ddot{U}_{q_1,q_2,q_3}(\mathfrak{gl}_1)$ from Section 2.2.1. The category $O^-$ is defined as in Definition 2.3 with the algebra $'\ddot{U}_{q,d}(\mathfrak{sl}_n)$ being $\mathbb{Z}$-graded via the adjoint action $x \mapsto q^{d_2}xq^{-d_2}$, cf. (T0.4). Likewise, the lowest weight modules are defined as in Definition 2.4, but with the lowest weight being now encoded by a nonzero constant $\nu \in \mathbb{C}^\times$ and a tuple of series $\boldsymbol{\phi} = \{\phi_i^{\pm}(z)\}_{i\in[n]}$, $\phi_i^{\pm}(z) \in \mathbb{C}[[z^{\mp 1}]]$, so that on the lowest weight vector $\mathrm{v}_0$ we have:

$$q^{d_2}\mathrm{v}_0 = \nu \cdot \mathrm{v}_0, \quad \psi_i^{\pm}(z)\mathrm{v}_0 = \phi_i^{\pm}(z) \cdot \mathrm{v}_0, \quad f_i(z)\mathrm{v}_0 = 0 \quad \text{for} \quad i \in [n].$$

Let $V_{\boldsymbol{\phi},\nu}$ be the corresponding irreducible lowest weight $'\ddot{U}_{q,d}(\mathfrak{sl}_n)$-module. The following result is completely analogous to Proposition 2.7:

**Proposition 3.6** *$V_{\boldsymbol{\phi},\nu}$ is in $O^-$ if and only if there are rational functions $\phi_i(z) \in \mathbb{C}(z)$ such that $\phi_i(0)\phi_i(\infty) = 1$ and $\phi_i^{\pm}(z) = \phi_i(z)^{\pm}$–the expansions of $\phi_i(z)$ in $z^{\mp 1}$.*

Similarly to the observation from the end of Section 2.2.1, all the above admits natural opposite counterparts: the category $O^+$ and the notion of the highest weight modules. Then, the analogue of Proposition 3.6 still holds. Moreover, we have natural equivalences $O^{\pm} \to O^{\mp}$ given by twists with the automorphisms $\vartheta_k$ of $\,'\ddot{U}_{q,d}(\mathfrak{sl}_n)$:

$$\vartheta_k: e_i(z) \mapsto f_{i'}(1/z),\ f_i(z) \mapsto e_{i'}(1/z),\ \psi_i^{\pm}(z) \mapsto \psi_{i'}^{\mp}(1/z),\ q^{d_2} \mapsto q^{-d_2} \quad (3.55)$$

for $k \in [n]$, where $i' := 2k - i \bmod n \in [n]$. We note that $\vartheta_k$ are just *cyclic rotations* of $\vartheta_0$, that is, $\vartheta_k = \varrho^k \circ \vartheta_0 \circ \varrho^{-k}$, where $\varrho$ is the following automorphism of $\,'\ddot{U}_{q,d}(\mathfrak{sl}_n)$:

$$\varrho: e_i(z) \mapsto e_{i+1}(z),\ f_i(z) \mapsto f_{i+1}(z),\ \psi_i^{\pm}(z) \mapsto \psi_{i+1}^{\pm}(z),\ q^{d_2} \mapsto q^{d_2}. \quad (3.56)$$

We shall now describe some of such modules combinatorially, following [9]. All of them depend on two parameters: $p \in [n]$, $u \in \mathbb{C}^{\times}$. Following (2.59, 2.79), define:

$$\phi(t) := \frac{q^{-1}t - q}{t - 1} \qquad \text{and} \qquad \phi^K(t) := \frac{K^{-1}t - K}{t - 1}. \quad (3.57)$$

We shall assume that $q$ and $d$ are *generic* (equivalently, $q_1$ and $q_3$ are generic (2.39)):

$$q^a d^b = 1\ (a, b \in \mathbb{Z}) \iff a = b = 0. \quad (\text{G})$$

For integers $a$ and $b$, we shall write $a \equiv b$ if $a - b$ is divisible by $n$, use $\overline{a} \in [n]$ to denote the residue $a \bmod n$, and finally set $\bar{\delta}_{ab} := \delta_{\overline{ab}}$ (vanishing unless $a \equiv b$).

- *Vector representations and their tensor products*

The "building block" of all the constructions below is the family of *vector representations* $\{V^{(p)}(u)\}_{u \in \mathbb{C}^{\times}}^{p \in [n]}$ with a basis parametrized by $\mathbb{Z}$, see [9, Lemma 3.1]:

**Proposition 3.7 (Vector representations)** *Let* $V^{(p)}(u)$ *be a $\mathbb{C}$-vector space with a basis* $\{[u]_r\}_{r \in \mathbb{Z}}$. *The following formulas define a* $\,'\ddot{U}_{q,d}(\mathfrak{sl}_n)$-*action on it:*

$$\begin{aligned}
e_i(z)[u]_r &= \begin{cases} \delta(q_1^r u/z) \cdot [u]_{r+1} & \text{if } i + r + 1 \equiv p \\ 0 & \text{otherwise} \end{cases}, \\[4pt]
f_i(z)[u]_r &= \begin{cases} \delta(q_1^{r-1} u/z) \cdot [u]_{r-1} & \text{if } i + r \equiv p \\ 0 & \text{otherwise} \end{cases}, \\[4pt]
q^{\pm d_2}[u]_r &= q^{\pm r} \cdot [u]_r, \\[4pt]
\psi_i^{\pm}(z)[u]_r &= \begin{cases} \phi(q_1^{r-1}u/z)^{\pm} \cdot [u]_r & \text{if } i + r \equiv p \\ \phi(z/q_1^r u)^{\pm} \cdot [u]_r & \text{if } i + r + 1 \equiv p \\ [u]_r & \text{otherwise} \end{cases}.
\end{aligned} \quad (3.58)$$

*Proof* The proof is completely analogous to its $\mathfrak{gl}_1$-counterpart from Section 2.2.2. $\square$

Furthermore, for any $p_1, \ldots, p_N \in [n]$ and $u_1, \ldots, u_N \in \mathbb{C}^{\times}$ satisfying $\frac{u_i}{u_j} \notin q_1^{\mathbb{Z}}$ for $i \neq j$, the analogue of (2.45)–(2.48) gives rise to a natural $\,'\ddot{U}_{q,d}(\mathfrak{sl}_n)$-action on $V^{(p_1)}(u_1) \otimes \cdots \otimes V^{(p_N)}(u_N)$. In particular, if $q$ and $d$ are generic (G), then we get a $\,'\ddot{U}_{q,d}(\mathfrak{sl}_n)$-action on the tensor products (cf. (2.52))

$$V^{(p)N}(u) := V^{(p)}(u) \otimes V^{(p)}(q_2^{-1}u) \otimes V^{(p)}(q_2^{-2}u) \otimes \cdots \otimes V^{(p)}(q_2^{-N+1}u),$$

and $W^{(p)N}(u) \subset V^{(p)N}(u)$, spanned by $\{|\lambda\rangle_u\}_{\lambda \in \mathcal{P}^N}$, is a $'\ddot{U}_{q,d}(\mathfrak{sl}_n)$-submodule.

- *Fock representations and their tensor products*

Similarly to the $\mathfrak{gl}_1$-counterpart, the above action of $'\ddot{U}_{q,d}(\mathfrak{sl}_n)$ on $W^{(p)N}(u)$ can be used to construct an action on the inductive limit of the subspaces $W^{(p)N,+}(u)$ spanned by $\{|\lambda\rangle_u\}_{\lambda \in \mathcal{P}^{N,+}}$, see [9, Proposition 3.3]. The resulting Fock modules $\{\mathcal{F}^{(p)}(u)\}_{u \in \mathbb{C}^\times}^{p \in [n]}$ have bases $\{|\lambda\rangle\}_{\lambda \in \mathcal{P}^+}$, see (2.58), and are explicitly described via:

**Proposition 3.8** *Assume that $q$ and $d$ are generic* (G). *Then, the $'\ddot{U}_{q,d}(\mathfrak{sl}_n)$-action on $\mathcal{F}^{(p)}(u)$ is given by the following explicit matrix coefficients:*

$\langle \lambda + 1_l | e_i(z) | \lambda \rangle =$

$$\bar{\delta}_{p+l-\lambda_l, i+1} \prod_{1 \leq k < l}^{p+k-\lambda_k \equiv i} \phi\left(q_1^{\lambda_k - \lambda_l - 1} q_3^{k-l}\right) \prod_{1 \leq k < l}^{p+k-\lambda_k \equiv i+1} \phi\left(q_1^{\lambda_l - \lambda_k} q_3^{l-k}\right) \delta\left(q_1^{\lambda_l} q_3^{l-1} u/z\right),$$

$\langle \lambda | f_i(z) | \lambda + 1_l \rangle =$

$$\bar{\delta}_{p+l-\lambda_l, i+1} \prod_{k > l}^{p+k-\lambda_k \equiv i} \phi\left(q_1^{\lambda_k - \lambda_l - 1} q_3^{k-l}\right) \prod_{k > l}^{p+k-\lambda_k \equiv i+1} \phi\left(q_1^{\lambda_l - \lambda_k} q_3^{l-k}\right) \delta\left(q_1^{\lambda_l} q_3^{l-1} u/z\right),$$

$$\langle \lambda | \psi_i^\pm(z) | \lambda \rangle = \left( \prod_{k \geq 1}^{p+k-\lambda_k \equiv i} \phi\left(q_1^{\lambda_k - 1} q_3^{k-1} u/z\right) \prod_{k \geq 1}^{p+k-\lambda_k \equiv i+1} \phi\left(q_1^{\lambda_k - 1} q_3^{k-2} u/z\right)^{-1} \right)^\pm,$$

$$\langle \lambda | q^{\pm d_2} | \lambda \rangle = q^{\pm |\lambda|},$$

*and all other matrix coefficients vanishing. It is an irreducible lowest weight module in the category $\mathcal{O}^-$ with the lowest weight $(1; \phi(z/u)^{\delta_{ip}})_{i \in [n]}$, and $|\emptyset\rangle$ is its lowest weight vector.*

*Proof* The proof is completely analogous to its $\mathfrak{gl}_1$-counterpart from Section 2.2.2. □

**Definition 3.3** For $\underline{c} \in (\mathbb{C}^\times)^{[n]}$, let $\tau_{u,\underline{c}}^p$ denote the twist of this representation by the algebra automorphism $\chi_{p,\underline{c}}$ of $'\ddot{U}_{q,d}(\mathfrak{sl}_n)$ defined via (cf. Remark 3.13):

$$e_{i,r} \mapsto c_i e_{i,r}, \quad f_{i,r} \mapsto c_i^{-1} f_{i,r}, \quad \psi_{i,r} \mapsto \psi_{i,r}, \quad q^{d_2} \mapsto q^{-\frac{p(n-p)}{2}} \cdot q^{d_2} \ \forall i \in [n], r \in \mathbb{Z}.$$

We call a collection $\{(p_k, u_k, \underline{c}_k)\}_{k=1}^N \subset [n] \times \mathbb{C}^\times \times (\mathbb{C}^\times)^{[n]}$ *generic* if

$$\forall \ 1 \leq s' < s \leq N \ \ \nexists a, b \in \mathbb{Z} \ \text{ s.t. } \ b - a \equiv p_{s'} - p_s \ \text{ and } \ u_s = u_{s'} q_1^{-a} q_3^{-b}. \quad (3.59)$$

We have the following simple result (see [9, Lemma 4.1]):

**Lemma 3.7** *For a generic collection $\{(p_k, u_k, \underline{c}_k)\}_{k=1}^N$, an analogue of (2.45)–(2.48) defines a $'\ddot{U}_{q,d}(\mathfrak{sl}_n)$-action on $\tau_{u_1,\underline{c}_1}^{p_1} \otimes \cdots \otimes \tau_{u_N,\underline{c}_N}^{p_N}$. Moreover, it is an irreducible lowest weight $'\ddot{U}_{q,d}(\mathfrak{sl}_n)$-module in $\mathcal{O}^-$, and $|\emptyset\rangle \otimes \cdots \otimes |\emptyset\rangle$ is its lowest weight vector.*

*Proof* The proof is completely analogous to its $\mathfrak{gl}_1$-counterpart from Section 2.2.2. In particular, the irreducibility of this module follows from the fact that $\{\psi_{i,r}\}_{i\in[n]}^{r\in\mathbb{Z}}$ act diagonally in the basis $|\lambda^1\rangle\otimes\cdots\otimes|\lambda^N\rangle$ with a simple joint spectrum. $\square$

● *Macmahon representations*

Finally, as in the $\mathfrak{gl}_1$-case, one can further construct a $'\ddot{U}_{q,d}(\mathfrak{sl}_n)$-action on the subspace of $\mathcal{F}^{(p)}(u)\otimes\mathcal{F}^{(p)}(q_2u)\otimes\cdots\otimes\mathcal{F}^{(p)}(q_2^{N-1}u)$ spanned by $|\lambda^1\rangle\otimes\cdots\otimes|\lambda^N\rangle$ satisfying $\lambda_l^k\geq\lambda_l^{k+1}$ for $1\leq k<N,l\geq 1$. Considering the $N\to\infty$ limit again, one further obtains a one-parametric family of $'\ddot{U}_{q,d}(\mathfrak{sl}_n)$-actions on the vector space with a basis labeled by plane partitions $\bar\lambda$ of (2.66), depending on a parameter $K\in\mathbb{C}^\times$. The resulting modules $\mathcal{M}^{(p)}(u,K)$ are called *(vacuum) Macmahon modules*. As we will not presently need explicit formulas for the action, we refer the interested reader to [9, §4.2] for more details. We shall only need the following result, cf. (3.57):

**Proposition 3.9** *[9, Theorem 4.3] Let $p\in[n]$, $u\in\mathbb{C}^\times$, and assume that $q$ and $d$ are generic (G) and $K^2\notin q^\mathbb{Z}d^\mathbb{Z}$. Then, the Macmahon module $\mathcal{M}^{(p)}(u,K)$ can be invariantly described as the lowest weight $'\ddot{U}_{q,d}(\mathfrak{sl}_n)$-module in the category $O^-$ with the lowest weight $(1;\phi^K(z/u)^{\delta_{ip}})_{i\in[n]}$. In the above realization, the empty plane partition $\bar\emptyset=(\emptyset,\emptyset,\ldots)$ is its lowest weight vector.*

## 3.3.2 Vertex representations

In this section, we recall an important family of vertex-type $\ddot{U}'_{q,d}(\mathfrak{sl}_n)$-modules from [22], generalizing the construction of Frenkel-Jing [13] for quantum affine algebras.

To this end, we consider a Heisenberg algebra $\mathfrak{A}_n$ generated by a central element $H_0$ and $\{H_{i,k}\}_{i\in[n]}^{k\neq 0}$ with the following defining relations, cf. (T1′):

$$[H_{i,k},H_{j,l}]=d^{-km_{ij}}\frac{[k]_q\cdot[ka_{ij}]_q}{k}\delta_{k,-l}\cdot H_0\,. \tag{3.60}$$

Let $\mathfrak{A}_n^{\geq}$ be the subalgebra generated by $\{H_{i,k}\}_{i\in[n]}^{k>0}\cup\{H_0\}$, and $\mathbb{C}\mathrm{v}_0$ be the $\mathfrak{A}_n^{\geq}$-module with $H_{i,k}\mathrm{v}_0=0$ and $H_0\mathrm{v}_0=\mathrm{v}_0$. The induced representation $F_n:=\mathrm{Ind}_{\mathfrak{A}_n^{\geq}}^{\mathfrak{A}_n}(\mathbb{C}\mathrm{v}_0)$ is called the *Fock representation* of $\mathfrak{A}_n$.

We denote the simple roots of $\mathfrak{sl}_n$ by $\{\bar\alpha_i\}_{i=1}^{n-1}$, the fundamental weights of $\mathfrak{sl}_n$ by $\{\bar\Lambda_i\}_{i=1}^{n-1}$, the simple coroots of $\mathfrak{sl}_n$ by $\{\bar h_i\}_{i=1}^{n-1}$. Let $\bar Q:=\bigoplus_{i=1}^{n-1}\mathbb{Z}\bar\alpha_i$ be the root lattice of $\mathfrak{sl}_n$, $\bar P:=\bigoplus_{i=1}^{n-1}\mathbb{Z}\bar\Lambda_i=\bigoplus_{i=2}^{n-1}\mathbb{Z}\bar\alpha_i\oplus\mathbb{Z}\bar\Lambda_{n-1}$ be the weight lattice of $\mathfrak{sl}_n$. We set

$$\bar\alpha_0:=-\sum_{1\leq i\leq n-1}\bar\alpha_i\in\bar Q,\qquad \bar\Lambda_0:=0\in\bar P,\qquad \bar h_0:=-\sum_{1\leq i\leq n-1}\bar h_i\,.$$

Let $\mathbb{C}\{\bar P\}$ be the $\mathbb{C}$-algebra generated by $e^{\bar\alpha_2},\ldots,e^{\bar\alpha_{n-1}},e^{\bar\Lambda_{n-1}}$ subject to the relations:

$$e^{\bar\alpha_i}\cdot e^{\bar\alpha_j}=(-1)^{\langle\bar h_i,\bar\alpha_j\rangle}e^{\bar\alpha_j}\cdot e^{\bar\alpha_i}\,,\qquad e^{\bar\alpha_i}\cdot e^{\bar\Lambda_{n-1}}=(-1)^{\delta_{i,n-1}}e^{\bar\Lambda_{n-1}}\cdot e^{\bar\alpha_i}\,.$$

For $\alpha = \sum_{i=2}^{n-1} m_i \bar\alpha_i + m_n \bar\Lambda_{n-1}$, we define $e^{\bar\alpha} \in \mathbb{C}\{\bar P\}$ via:

$$e^{\bar\alpha} := (e^{\bar\alpha_2})^{m_2} \cdots (e^{\bar\alpha_{n-1}})^{m_{n-1}} (e^{\bar\Lambda_{n-1}})^{m_n}.$$

Let $\mathbb{C}\{\bar Q\}$ be the subalgebra of $\mathbb{C}\{\bar P\}$ generated by $\{e^{\bar\alpha_i}\}_{i=1}^{n-1}$.

For every $0 \le p \le n-1$, consider the subspace

$$W(p)_n := F_n \otimes \mathbb{C}\{\bar Q\} e^{\bar\Lambda_p} \subset F_n \otimes \mathbb{C}\{\bar P\}.$$

We define the operators $e^{\bar\alpha}, \partial_{\bar\alpha_i}, H_{i,l}, z^{H_{i,0}}, \mathrm{d}$ acting on $W(p)_n$, which assign to every

$$\mathrm{v} \otimes e^{\bar\beta} = \left(H_{i_1,-k_1} \cdots H_{i_N,-k_N} \mathrm{v}_0\right) \otimes e^{\sum_{j=1}^{n-1} m_j \bar\alpha_j + \bar\Lambda_p} \in W(p)_n \qquad (3.61)$$

the following values:

$$
\begin{aligned}
e^{\bar\alpha}(\mathrm{v} \otimes e^{\bar\beta}) &:= \mathrm{v} \otimes e^{\bar\alpha} e^{\bar\beta}, \\
\partial_{\bar\alpha_i}(\mathrm{v} \otimes e^{\bar\beta}) &:= \langle \bar h_i, \bar\beta \rangle \cdot \mathrm{v} \otimes e^{\bar\beta}, \\
H_{i,l}(\mathrm{v} \otimes e^{\bar\beta}) &:= (H_{i,l}\mathrm{v}) \otimes e^{\bar\beta}, \\
z^{H_{i,0}}(\mathrm{v} \otimes e^{\bar\beta}) &:= z^{\langle \bar h_i, \bar\beta\rangle} d^{\frac{1}{2} \sum_{j=1}^{n-1} \langle \bar h_i, m_j \bar\alpha_j\rangle m_{ij}} \cdot \mathrm{v} \otimes e^{\bar\beta}, \\
\mathrm{d}(\mathrm{v} \otimes e^{\bar\beta}) &:= \left(\sum k_i + \frac{(\bar\beta,\bar\beta) - (\bar\Lambda_p,\bar\Lambda_p)}{2}\right) \cdot \mathrm{v} \otimes e^{\bar\beta}.
\end{aligned}
\qquad (3.62)
$$

The following result provides a natural structure of a $\ddot U'_{q,d}(\mathfrak{sl}_n)$-module on $W(p)_n$:

**Proposition 3.10** *[22, Proposition 3.2.2] For $\underline c = (c_0, \ldots, c_{n-1}) \in (\mathbb{C}^\times)^{[n]}, u \in \mathbb{C}^\times,$ $0 \le p \le n-1$, the following formulas define an action of $\ddot U'_{q,d}(\mathfrak{sl}_n)$ on $W(p)_n$:*

$$\rho^p_{u,\underline c}(e_i(z)) =$$
$$c_i \cdot \exp\left(\sum_{k>0} \frac{q^{-k/2}}{[k]_q} H_{i,-k} \left(\frac{z}{u}\right)^k\right) \exp\left(-\sum_{k>0} \frac{q^{-k/2}}{[k]_q} H_{i,k} \left(\frac{z}{u}\right)^{-k}\right) \cdot e^{\bar\alpha_i} \left(\frac{z}{u}\right)^{1+H_{i,0}},$$

$$\rho^p_{u,\underline c}(f_i(z)) =$$
$$\frac{(-1)^{n\delta_{i0}}}{c_i} \cdot \exp\left(-\sum_{k>0} \frac{q^{k/2}}{[k]_q} H_{i,-k} \left(\frac{z}{u}\right)^k\right) \exp\left(\sum_{k>0} \frac{q^{k/2}}{[k]_q} H_{i,k} \left(\frac{z}{u}\right)^{-k}\right) \cdot e^{-\bar\alpha_i} \left(\frac{z}{u}\right)^{1-H_{i,0}},$$

$$\rho^p_{u,\underline c}(\psi_i^\pm(z)) = \exp\left(\pm(q - q^{-1}) \sum_{k>0} H_{i,\pm k} \left(\frac{z}{u}\right)^{\mp k}\right) \cdot q^{\pm\partial_{\bar\alpha_i}},$$

$$\rho^p_{u,\underline c}(\gamma^{\pm 1/2}) = q^{\pm 1/2}, \qquad \rho^p_{u,\underline c}(q^{\pm d_1}) = q^{\pm \mathrm{d}}.$$

*Remark 3.10* For $q$ and $d$ generic (G), the $\ddot U'_{q,d}(\mathfrak{sl}_n)$-module $\rho^p_{u,\underline c}$ is irreducible since so its restriction to $\dot U^v_q(\mathfrak{sl}_n)$, recovering the representation of [13], as proved in [2].

### 3.3.3 Shuffle bimodules

Following [23] (and generalizing the $\mathfrak{gl}_1$-counterparts from Section 2.2.4), let us introduce three families of $S$-bimodules which also become $'\ddot{U}_{q,d}(\mathfrak{sl}_n)$-modules.

- *Shuffle bimodules $S_{1,p}(u)$*

  For $u \in \mathbb{C}^\times$ and $0 \leq p \leq n-1$, consider an $\mathbb{N}^{[n]}$-graded $\mathbb{C}$-vector space

  $$S_{1,p}(u) = \bigoplus_{\underline{k}=(k_0,\dots,k_{n-1}) \in \mathbb{N}^{[n]}} S_{1,p}(u)_{\underline{k}},$$

  where the degree $\underline{k}$ component $S_{1,p}(u)_{\underline{k}}$ consists of $\Sigma_{\underline{k}}$-symmetric rational functions $F$ in the variables $\{x_{i,r}\}_{i \in [n]}^{1 \leq r \leq k_i}$ satisfying the following three conditions:

(i) *Pole conditions* on $F$:

$$F = \frac{f(x_{0,1},\dots,x_{n-1,k_{n-1}})}{\prod_{i \in [n]} \prod_{r \leq k_i}^{r' \leq k_{i+1}} (x_{i,r} - x_{i+1,r'}) \cdot \prod_{r=1}^{k_p} (x_{p,r} - u)}, \quad f \in \mathbb{C}\left[\{x_{i,r}^{\pm 1}\}_{i \in [n]}^{1 \leq r \leq k_i}\right]^{\Sigma_{\underline{k}}}.$$

(ii) *First kind wheel conditions* on $F$:

$$F\left(\{x_{i,r}\}\right) = 0 \quad \text{once} \quad x_{i,r_1} = qd^\epsilon x_{i+\epsilon,s} = q^2 x_{i,r_2} \quad \text{for some} \quad \epsilon, i, r_1, r_2, s.$$

where $\epsilon \in \{\pm 1\}$, $i \in [n]$, $1 \leq r_1 \neq r_2 \leq k_i$, and $1 \leq s \leq k_{\overline{i+\epsilon}}$.

(iii) *Second kind wheel conditions* on $F$ (with $f$ from (i)):

$$f\left(\{x_{i,r}\}\right) = 0 \quad \text{once} \quad x_{p,r_1} = u \quad \text{and} \quad x_{p,r_2} = q^2 u \quad \text{for some} \quad 1 \leq r_1 \neq r_2 \leq k_p.$$

Fix $\underline{c} = (c_0, c_1, \dots, c_{n-1}) \in (\mathbb{C}^\times)^{[n]}$. Given $F \in S_{\underline{k}}$ and $G \in S_{1,p}(u)_{\underline{\ell}}$, we define $F \star G, G \star F \in S_{1,p}(u)_{\underline{k}+\underline{\ell}}$ via (cf. (2.77, 2.78) and (3.26)):

$$(F \star G)(x_{0,1}, \dots, x_{0,k_0+\ell_0}; \dots; x_{n-1,1}, \dots, x_{n-1,k_{n-1}+\ell_{n-1}}) := \frac{1}{\underline{k}! \cdot \underline{\ell}!} \cdot \prod_{i \in [n]} c_i^{k_i} \times$$

$$\text{Sym}\left(F\left(\{x_{i,r}\}_{i \in [n]}^{r \leq k_i}\right) G\left(\{x_{i',r'}\}_{i' \in [n]}^{r' > k_{i'}}\right) \prod_{i \in [n]}^{i' \in [n]} \prod_{r \leq k_i}^{r' > k_{i'}} \zeta_{i,i'}\left(\frac{x_{i,r}}{x_{i',r'}}\right) \prod_{r=1}^{k_p} \phi\left(\frac{x_{p,r}}{u}\right)\right) \tag{3.63}$$

and

$$(G \star F)(x_{0,1}, \dots, x_{0,k_0+\ell_0}; \dots; x_{n-1,1}, \dots, x_{n-1,k_{n-1}+\ell_{n-1}}) := \frac{1}{\underline{k}! \cdot \underline{\ell}!} \times$$

$$\text{Sym}\left(G\left(\{x_{i,r}\}_{i \in [n]}^{r \leq \ell_i}\right) F\left(\{x_{i',r'}\}_{i' \in [n]}^{r' > \ell_{i'}}\right) \prod_{i \in [n]}^{i' \in [n]} \prod_{r \leq \ell_i}^{r' > \ell_{i'}} \zeta_{i,i'}\left(\frac{x_{i,r}}{x_{i',r'}}\right)\right). \tag{3.64}$$

These formulas clearly endow $S_{1,p}(u)$ with an $S$-bimodule structure. Identifying $S$ with $\,^{\prime}\ddot{U}^{>}_{q,d}$ via Theorem 3.6, we thus get two commuting $\,^{\prime}\ddot{U}^{>}_{q,d}$-actions on $S_{1,p}(u)$.

In fact, similarly to Proposition 2.12, the left action can be extended to an action of the entire quantum toroidal algebra $\,^{\prime}\ddot{U}_{q,d}(\mathfrak{sl}_n)$:

**Proposition 3.11** *The following formulas define a* $\,^{\prime}\ddot{U}_{q,d}(\mathfrak{sl}_n)$-*action on* $S_{1,p}(u)$:

$$\pi^p_{u,\underline{c}}(q^{d_2})G = q^{-\frac{p(n-p)}{2}+|\underline{k}|}\cdot G\,,$$

$$\pi^p_{u,\underline{c}}(e_{i,k})G = x^k_{i,1}\star G\,,$$

$$\pi^p_{u,\underline{c}}(h_{i,0})G = (2k_i - k_{i-1} - k_{i+1} - \delta_{ip})\cdot G\,,$$

$$\pi^p_{u,\underline{c}}(h_{i,l})G = \left(\frac{1}{l}\sum_{i'\in[n]}\sum_{r'=1}^{k_{i'}}[la_{ii'}]_q\, d^{-lm_{ii'}}x^l_{i',r'} - \delta_{ip}\frac{[l]_q}{l}\,q^l u^l\right)\cdot G \ \text{ for } \ l\neq 0\,,$$

$$\pi^p_{u,\underline{c}}(f_{i,k})G = \frac{k_i c_i^{-1}}{q^{-1}-q}\left(\underset{z=0}{\mathrm{Res}}+\underset{z=\infty}{\mathrm{Res}}\right)\left\{\frac{z^k G(\{x_{i',r'}\}_{|x_{i,k_i}\mapsto z})}{\prod_{i'}\prod_{r'=1}^{k_{i'}-\delta_{ii'}}\zeta_{i',i}(\frac{x_{i',r'}}{z})}\,\frac{dz}{z}\right\}\,.$$

*Here,* $k\in\mathbb{Z}$, $\underline{c} = (c_0,\dots,c_{n-1})\in(\mathbb{C}^{\times})^{[n]}$, $G\in S_{1,p}(u)_{\underline{k}}$, *and* $|\underline{k}| = \sum_{i\in[n]}k_i$.

The factor $q^{-\frac{p(n-p)}{2}}$ in $\pi^p_{u,\underline{c}}(q^{d_2})$ is solely needed for Section 3.4, see Remark 3.13.

*Remark 3.11* The formulas of Proposition 3.11 can be equivalently written as follows:

$$\begin{aligned}\pi^p_{u,\underline{c}}(q^{d_2})G &= q^{-\frac{p(n-p)}{2}+|\underline{k}|}\cdot G\,,\\ \pi^p_{u,\underline{c}}(e_i(z))G &= \delta(x_{i,1}/z)\star G\,,\end{aligned}\tag{3.65}$$

$$\pi^p_{u,\underline{c}}(\psi^{\pm}_i(z))G =$$
$$\left(\prod_{r=1}^{k_i}\frac{q^2 z - x_{i,r}}{z - q^2 x_{i,r}}\cdot\prod_{r=1}^{k_{i+1}}\frac{z - qdx_{i+1,r}}{qz - dx_{i+1,r}}\cdot\prod_{r=1}^{k_{i-1}}\frac{dz - qx_{i-1,r}}{qdz - x_{i-1,r}}\cdot\phi\left(\frac{z}{u}\right)^{\delta_{ip}}\right)^{\pm}\cdot G\,,\tag{3.66}$$

$$\pi^p_{u,\underline{c}}(f_i(z))G = \frac{k_i c_i^{-1}}{q^{-1}-q}\times$$
$$\left\{\left(\frac{G(\{x_{i',r'}\}_{|x_{i,k_i}\mapsto z})}{\prod_{i'}\prod_{r'=1}^{k_{i'}-\delta_{ii'}}\zeta_{i',i}(\frac{x_{i',r'}}{z})}\right)^{+} - \left(\frac{G(\{x_{i',r'}\}_{|x_{i,k_i}\mapsto z})}{\prod_{i'}\prod_{r'=1}^{k_{i'}-\delta_{ii'}}\zeta_{i',i}(\frac{x_{i',r'}}{z})}\right)^{-}\right\}\,,\tag{3.67}$$

where $\gamma(z)^{\pm}$ denotes the expansion of a rational function $\gamma(z)$ in $z^{\mp 1}$, respectively.

*Proof* We need to verify the compatibility of the above assignment $\pi^p_{u,\underline{c}}$ with the defining relations (T0.1)–(T8). The only nontrivial of those are (T3, T4, T6, T8). To check (T3, T6), we use formulas (3.66, 3.67) together with an obvious identity $\frac{\zeta_{i,j}(z/w)}{\zeta_{j,i}(w/z)} = g_{a_{ij}}(d^{m_{ij}}z/w)$ for any $i,j\in[n]$. The verification of (T8) boils down to:

$$\mathrm{Sym}_{z_1,z_2}\left(\frac{\zeta_{i,i}(z_1/z_2)^{-1}}{\zeta_{i,i\pm1}(z_1/w)\zeta_{i,i\pm1}(z_2/w)} - \frac{(q+q^{-1})\zeta_{i,i}(z_1/z_2)^{-1}}{\zeta_{i,i\pm1}(z_1/w)\zeta_{i\pm1,i}(w/z_2)} + \right.$$
$$\left. \frac{\zeta_{i,i}(z_1/z_2)^{-1}}{\zeta_{i\pm1,i}(w/z_1)\zeta_{i\pm1,i}(w/z_2)}\right) = 0\,.$$

Finally, to verify (T4) we note that the $k_i + 1 - \delta_{ij}$ different summands from the symmetrization appearing in $\pi_{u,\underline{c}}^{p}(e_i(z))\pi_{u,\underline{c}}^{p}(f_j(w))G$ cancel the $k_i + 1 - \delta_{ij}$ terms (out of $k_i + 1$) from the symmetrization appearing in $\pi_{u,\underline{c}}^{p}(f_j(w))\pi_{u,\underline{c}}^{p}(e_i(z))G$.   □

● *Shuffle bimodules $S(\underline{u})$*

The above construction admits a *higher rank* counterpart. For $\underline{m} \in \mathbb{N}^{[n]}$, consider

$$\underline{u} = (u_{0,1},\ldots,u_{0,m_0};\ldots;u_{n-1,1},\ldots,u_{n-1,m_{n-1}}) \quad \text{with} \quad u_{i,s} \in \mathbb{C}^{\times}.$$

Define $S(\underline{u}) = \bigoplus_{\underline{k}\in\mathbb{N}^{[n]}} S(\underline{u})_{\underline{k}}$ similarly to $S_{1,p}(u)$ with the following modifications:

(i′) *Pole conditions* for $F$ from a degree $\underline{k}$ component should read as follows:

$$F = \frac{f(x_{0,1},\ldots,x_{n-1,k_{n-1}})}{\prod_{i\in[n]}\prod_{r'\le k_{i+1}}^{r'\le k_{i+1}}(x_{i,r}-x_{i+1,r'})\cdot\prod_{i\in[n]}\prod_{s=1}^{m_i}\prod_{r=1}^{k_i}(x_{i,r}-u_{i,s})}$$

with $f \in \mathbb{C}\left[\{x_{i,r}^{\pm1}\}_{i\in[n]}^{1\le r\le k_i}\right]^{\Sigma_{\underline{k}}}$.

(iii′) *Second kind wheel conditions* for such $F$ should read as follows:

$$f\left(\{x_{i,r}\}\right) = 0 \quad \text{once} \quad x_{i,r_1} = u_{i,s} \quad \text{and} \quad x_{i,r_2} = q^2 u_{i,s} \quad \text{for some} \quad i\,,s\,,r_1 \ne r_2\,.$$

We endow $S(\underline{u})$ with an $S$-bimodule structure by the formulas (3.63, 3.64) with

$$\prod_{1\le r\le k_p}\phi(x_{p,r}/u) \rightsquigarrow \prod_{i\in[n]}\prod_{1\le s\le m_i}\prod_{1\le r\le k_i}\phi(x_{i,r}/u_{i,s})\,.$$

The resulting left action of ${}'\ddot{U}_{q,d}^{>} \simeq S$ on $S(\underline{u})$ can be extended to the action of the entire quantum toroidal ${}'\ddot{U}_{q,d}(\mathfrak{sl}_n)$, denoted by $\pi_{\underline{u},\underline{c}}$. Explicitly, the latter is defined by the formulas (3.65)–(3.67) with the following two modifications:

$$\phi(z/u)^{\delta_{ip}} \rightsquigarrow \prod_{1\le s\le m_i}\phi(z/u_{i,s})\,, \qquad q^{-\frac{p(n-p)}{2}} \rightsquigarrow q^{-\sum_{p=0}^{n-1}m_p\cdot\frac{p(n-p)}{2}}\,.$$

● *Shuffle bimodules $S_{1,p}^{K}(u)$ and $S^{\underline{K}}(\underline{u})$*

Another natural generalization of the above $S_{1,p}(u)$ is provided by the $S$-bimodules $S_{1,p}^{K}(u)$. As a vector space, $S_{1,p}^{K}(u)$ is defined similarly to $S_{1,p}(u)$ but without imposing the second kind wheel conditions. The $S$-bimodule structure on $S_{1,p}^{K}(u)$ is defined by the formulas (3.63, 3.64) with the following modification in (3.63):

$$\phi(t) \rightsquigarrow \phi^{K}(t) = (K^{-1}\cdot t - K)/(t-1)\,.$$

The resulting left action of ${}'\ddot{U}^{>}_{q,d} \simeq S$ on $S^K_{1,p}(u)$ can be extended to the ${}'\ddot{U}_{q,d}(\mathfrak{sl}_n)$-action, denoted by $\pi^{p,K}_{u,\underline{c}}$, defined by (3.65)–(3.67) but with $\phi$ being replaced by $\phi^K$.

Inspired by the above construction of $S(\underline{u})$, it is clear how to define the *higher rank* counterparts $S^{\underline{K}}(\underline{u})$ of $S^K_{1,p}(u)$, equip them with the $S$-bimodule structure, and extend the resulting left action of ${}'\ddot{U}^{>}_{q,d} \simeq S$ to the action $\pi^K_{\underline{u},\underline{c}}$ of ${}'\ddot{U}_{q,d}(\mathfrak{sl}_n)$ on $S^{\underline{K}}(\underline{u})$.

## 3.4 Identification of representations

In this section, we present explicit relation between three different families of modules from Section 3.3 (assuming that $q$ and $d$ are generic in the sense of (G)).

### 3.4.1 Shuffle realization of Fock modules

For $0 \leq p \leq n-1$, $u \in \mathbb{C}^{\times}$, $\underline{c} \in (\mathbb{C}^{\times})^{[n]}$, recall the action $\pi^p_{u,\underline{c}}$ of ${}'\ddot{U}_{q,d}(\mathfrak{sl}_n)$ on $S_{1,p}(u)$ from Proposition 3.11. We define

$$S' := \bigoplus_{\underline{k} \neq (0,\ldots,0)} S_{\underline{k}} \subset S$$

and consider the following subspace

$$J_0 := S_{1,p}(u) \star S' = \mathrm{span}_{\mathbb{C}}\left\{ G \star F \mid G \in S_{1,p}(u), F \in S' \right\} \subset S_{1,p}(u). \qquad (3.68)$$

The following result is straightforward, cf. Proposition 2.13(a):

**Lemma 3.8** $J_0$ *of* (3.68) *is a submodule of the* ${}'\ddot{U}_{q,d}(\mathfrak{sl}_n)$-action $\pi^p_{u,\underline{c}}$ *on* $S_{1,p}(u)$.

*Proof* The invariance of $J_0$ with respect to the left action of ${}'\ddot{U}^{>}_{q,d} \simeq S$ is obvious as it coincides with (3.63) and hence commutes with the right action of $S$ from (3.64). The invariance of $J_0$ with respect to the left ${}'\ddot{U}^0_{q,d}$-action follows from the equality

$$\pi^p_{u,\underline{c}}(h_{i,l})(G \star F) = \left(\pi^p_{u,\underline{c}}(h_{i,l})G\right) \star F + G \star \left(\sum_{j \in [n]} \sum_{s=1}^{k_j} \frac{[la_{ij}]_q d^{-lm_{ij}}}{l} x^l_{j,s} \cdot F\right)$$

for any $i \in [n]$, $l \neq 0$, and $G \in S_{1,p}(u)_{\underline{\ell}}$, $F \in S'_{\underline{k}}$ with $\underline{k}, \underline{\ell} \in \mathbb{N}^{[n]}$. Finally, for $i \in [n]$, $l \in \mathbb{Z}$, and $F, G$ as above, using the formula (3.67) one can easily see that

$$\pi^p_{u,\underline{c}}(f_{i,l})(G \star F) = G' \star F + G \star F'$$

for some $G' \in S_{1,p}(u)_{\underline{\ell}-1_i}$ and $F' \in S_{\underline{k}-1_i} \cap S'$, with the two terms corresponding to cases when the $x_{*,*}$-variable specialized to $z$ in (3.67) is placed either in $G$ or $F$. $\qquad \square$

Let $\bar{\pi}_{u,\underline{c}}^{p}$ denote the corresponding quotient representation of $'\ddot{U}_{q,d}(\mathfrak{sl}_n)$ on

$$\bar{S}_{1,p}(u) := S_{1,p}(u)/J_0 \,. \tag{3.69}$$

The following result is analogous to Proposition 2.13(b):

**Theorem 3.12** *For any $p, u, \underline{c}$ as above, we have a $'\ddot{U}_{q,d}(\mathfrak{sl}_n)$-module isomorphism*

$$\bar{\pi}_{u,\underline{c}}^{p} \simeq \tau_{u,\underline{c}}^{p} \,. \tag{3.70}$$

*Proof* By Proposition 3.8, $\tau_{u,\underline{c}}^{p}$ is an irreducible $'\ddot{U}_{q,d}(\mathfrak{sl}_n)$-module. Moreover, both $\bar{\mathbf{1}}_u \in \bar{S}_{1,p}(u)$ (the image of $1 \in S_{1,p}(u)_{(0,\ldots,0)}$) and $|\emptyset\rangle \in \mathcal{F}^{(p)}(u)$ are the lowest weight vectors of the same lowest weight. Hence, it suffices to compare dimensions of the $\mathbb{N}^{[n]}$-graded components of the quotient space $\bar{S}_{1,p}(u)$ from (3.69):

$$\sum_{\underline{k}\in\mathbb{N}^{[n]}\,:\,|\underline{k}|=m} \dim \bar{S}_{1,p}(u)_{\underline{k}} = p(m) := \#\Big\{\text{size } m \text{ Young diagrams}\Big\} \quad \forall\, m \in \mathbb{N}. \quad (\dagger)$$

*Descending filtration*

To prove $(\dagger)$ we equip $S_{1,p}^{m}(u) := \bigoplus_{|\underline{k}|=m} S_{1,p}(u)_{\underline{k}}$ with a filtration $\{S_{1,p}^{m,\lambda}(u)\}_\lambda$ labeled by Young diagrams $\lambda$ of size $\leq m$, with specialization maps $\rho_\lambda$ defined below:

$$S_{1,p}^{m,\lambda}(u) := \bigcap_{\mu > \lambda} \mathrm{Ker}(\rho_\mu) \subset S_{1,p}^{m}(u), \tag{3.71}$$

cf. (2.29, 3.39), where $>$ is the lexicographical order on size $\leq m$ Young diagrams.

Consider an $[n]$-coloring of the boxes of a Young diagram $\lambda$ by assigning color $c(\square) := \overline{p - a + b} \in [n]$ to a box $\square = (a, b) \in \lambda$ located at the $b$-th row and $a$-th column (with $1 \leq b \leq \lambda_1'$ and $1 \leq a \leq \lambda_b$). We define the $\mathbb{N}^{[n]}$-degree of $\lambda$ via:

$$\underline{k}^{\lambda} := (k_0^\lambda, \ldots, k_{n-1}^\lambda) \in \mathbb{N}^{[n]} \qquad \text{with} \qquad k_i^\lambda = \#\{\square \in \lambda \,|\, c(\square) = i\}\,.$$

*Remark 3.12* Let $\tau_u^p$ denote $\tau_{u,(1,\ldots,1)}^p$. Then, the assignment $|\lambda\rangle \mapsto \prod_{\square\in\lambda} c_{c(\square)} \cdot |\lambda\rangle$ gives rise to an isomorphism of $'\ddot{U}_{q,d}(\mathfrak{sl}_n)$-modules $\tau_u^p \xrightarrow{\sim} \tau_{u,\underline{c}}^p$ for any $\underline{c} \in (\mathbb{C}^\times)^{[n]}$. $\square$

Let us fill each box $\square = (a, b) \in \lambda$ with $q_1^{a-1} q_3^{b-1} u$. For $F \in S_{1,p}(u)_{\underline{k}}$, the **specialization** $\rho_\lambda(F)$ shall be obtained by specializing $\underline{k}^\lambda$ variables of $F$ to the corresponding entries of $\lambda$, while taking care of the zeroes in the numerator and denominator.

*Specialization maps $\rho_\lambda$*

For $\underline{k} \in \mathbb{N}^{[n]}$ and a Young diagram $\lambda$, define $\underline{\ell} := \underline{k} - \underline{k}^\lambda \in \mathbb{Z}^{[n]}$. If $\underline{\ell} \notin \mathbb{N}^{[n]}$, we set $\rho_\lambda(F) = 0$ for any $F \in S_{1,p}(u)_{\underline{k}}$. Otherwise, we define $\rho_\lambda(F)$ via (3.72) as follows:

• First, we pick the corner box $\square = (1, 1) \in \lambda$ of color $p$ and specialize $x_{p,k_p} \mapsto u$. Since $F$ has the first order pole at $x_{p,k_p} = u$, the following is well-defined:

$$\rho_\lambda^{(1)}(F) := \big[(x_{p,k_p} - u) \cdot F\big]_{|x_{p,k_p}\mapsto u} \,.$$

• Next, we specialize more variables to the entries of the remaining boxes from the first row and the first column. For every box $(a+1, 1) \in \lambda$ $(0 < a < \lambda_1)$ of color $\overline{p-a}$, we choose an unspecified yet $x_{\overline{p-a},*}$-variable and set it to $q_1^a u$. Likewise, for every box $(1, b+1) \in \lambda$ $(0 < b < \lambda_1')$ of color $\overline{p+b}$, we choose an unspecified yet $x_{\overline{p+b},*}$-variable and set it to $q_3^b u$. We perform this procedure step-by-step moving from $(1,1)$ to the right and then from $(1,1)$ upwards. Let $\rho_\lambda^{(\lambda_1 + \lambda_1' - 1)}(F)$ denote the resulting overall specialization of $F$.

• If $(2,2) \notin \lambda$, then we set $\rho_\lambda(F) := \rho_\lambda^{(\lambda_1 + \lambda_1' - 1)}(F)$. Otherwise, we pick another variable of the $p$-th family, say $x_{p,k_p-1}$, and specialize it to $q_1 q_3 u$, but first cancel the zero of $\rho_\lambda^{(\lambda_1 + \lambda_1' - 1)}(F)$ due to the first kind wheel conditions:

$$\rho_\lambda^{(\lambda_1 + \lambda_1')}(F) := \left[ \frac{1}{x_{p,k_p-1} - q_1 q_3 u} \cdot \rho_\lambda^{(\lambda_1 + \lambda_1' - 1)}(F) \right]_{|x_{p,k_p-1} \mapsto q_1 q_3 u} .$$

• Next, we move step-by-step from $(2,2)$ to the right and then from $(2,2)$ upwards, on each step specializing the corresponding $x_{*,*}$-variable to the prescribed entry of $\square \in \lambda$ located in the second row or column, but first eliminating order 1 zeroes as above (due to the first kind wheel conditions).

• Having performed this procedure $|\lambda|$ times, we get $\rho_\lambda^{(|\lambda|)}(F) \in \mathbb{C}(\{x_{i,r}\}_{i \in [n]}^{1 \le r \le \ell_i})$. Then, we finally set

$$\rho_\lambda(F) := \rho_\lambda^{(|\lambda|)}(F) . \tag{3.72}$$

*Key properties of $\rho_\lambda$*

Let us now identify the key properties of the specializations maps $\rho_\lambda$ similar to those of $\phi_{\underline{d}}$ from Lemmas 1.5–1.6 and of $\phi_L$ from Lemma 3.3. To do so, we define

$$Q_{\underline{\ell},\lambda} := \prod_{r=1}^{\ell_{p-\lambda_1}} (x_{\overline{p-\lambda_1},r} - q_1^{\lambda_1} u) \cdot \prod_{\substack{\lambda_{b+1} < \lambda_b \\ b \ge 1}}^{} \prod_{r=1}^{\ell_{p-\lambda_{b+1}+b}} (x_{\overline{p-\lambda_{b+1}+b},r} - q_1^{\lambda_{b+1}} q_3^b u)$$

and

$$G_{\underline{\ell},\lambda} := \prod_{r=1}^{\ell_p} \frac{x_{p,r} - q^2 u}{x_{p,r} - u} \times$$

$$\frac{\prod_{\square=(a,b)\in X_\lambda^+} \prod_{r=1}^{\ell_{c(\square)}} (x_{c(\square),r} - q_1^{a-1} q_3^{b-1} u) \prod_{\square=(a,b)\in X_\lambda^-} \prod_{r=1}^{\ell_{c(\square)}} (x_{c(\square),r} - q_1^{a-1} q_3^{b-1} u)}{\prod_{\square=(a,b)\in\lambda} \left\{ \prod_{r=1}^{\ell_{c(\square)-1}} (x_{c(\square)-1,r} - q_1^{a-1} q_3^{b-1} u) \prod_{r=1}^{\ell_{c(\square)+1}} (x_{c(\square)+1,r} - q_1^{a-1} q_3^{b-1} u) \right\}}$$

Here, the set $X_\lambda^+ \subset \mathbb{Z}^2$ consists of those $(a,b) \in \mathbb{Z}^2$ such that $(a+1,b),(a+1,b+1) \in \lambda$ or $(a,b+1),(a+1,b+1) \in \lambda$, while $X_\lambda^- \subset \mathbb{Z}^2$ consists of those $(a,b) \in \mathbb{Z}^2$ such that $(a-1,b),(a-1,b-1) \in \lambda$ or $(a,b-1),(a-1,b-1) \in \lambda$. Then, we have:

**Lemma 3.9** (a) $\rho_\lambda(S_{1,p}(u)_{\underline{k}^\lambda + \underline{\ell}}) \subset S_{\underline{\ell}} \cdot G_{\underline{\ell},\lambda} := \{F' \cdot G_{\underline{\ell},\lambda} \mid F' \in S_{\underline{\ell}}\}$.

(b) $\rho_\lambda(S_{1,p}^{|\lambda|+|\underline{\ell}|,\lambda}(u)_{\underline{k}^\lambda + \underline{\ell}}) \subset S_{\underline{\ell}} \cdot G_{\underline{\ell},\lambda} Q_{\underline{\ell},\lambda} := \{F' \cdot G_{\underline{\ell},\lambda} Q_{\underline{\ell},\lambda} \mid F' \in S_{\underline{\ell}}\}$.    □

*Proof* Part (a) is a consequence of the first and second kind wheel conditions, while the extra divisibility by $Q_{\underline{\ell},\lambda}$ in part (b) is due to the definition of the filtration (3.71). □

**Lemma 3.10** *(a) If $\underline{k} - \underline{k}^\lambda \notin \mathbb{N}^{[n]}$, then $\rho_\lambda(S_{1,p}(u)_{\underline{k}} \star S_{\underline{\ell}}) = 0$ for any $\underline{\ell} \in \mathbb{N}^{[n]}$.*

*(b) $\rho_\lambda(S_{1,p}^{|\lambda|+|\underline{\ell}|,\lambda}(u)_{\underline{k}^\lambda+\underline{\ell}}) = \rho_\lambda((S_{1,p}(u)_{\underline{k}^\lambda} \star S_{\underline{\ell}}) \cap S_{1,p}^{|\lambda|+|\underline{\ell}|,\lambda}(u)_{\underline{k}^\lambda+\underline{\ell}})$ for any $\underline{\ell} \in \mathbb{N}^{[n]}$.* □

*Proof* (a) For $F_1 \in S_{1,p}(u)_{\underline{k}}$ and $F_2 \in S_{\underline{\ell}}$, let us evaluate the $\rho_\lambda$-specialization of any summand in $F_1 \star F_2$ from (3.64). In what follows, we say "$q_1^a q_3^b u$ gets into $F_2$" in the chosen summand if the $x_{*,*}$-variable which is specialized to $q_1^a q_3^b u$ is plugged to $F_2$ rather than $F_1$. If $u$ gets into $F_2$, we obtain zero once we apply $\rho_\lambda^{(1)}$. A simple inductive argument shows that if at least one of the variables $\{q_1^a u\}_{a=1}^{\lambda_1-1} \cup \{q_3^b u\}_{b=1}^{\lambda'_1-1}$ gets into $F_2$, we also obtain zero after applying $\rho_\lambda^{(a+1)}$ or $\rho_\lambda^{(\lambda_1+b)}$ since we get a vanishing $\zeta$-factor. If $q_1 q_3 u$ gets into $F_2$, but all the entries from the first hook of $\lambda$ get into $F_1$, then there are two vanishing $\zeta$-factors, and so we get zero after applying $\rho_\lambda^{(\lambda_1+\lambda'_1)}$, and so on. However, all the specialized variables cannot get into $F_1$ as $\underline{k} - \underline{k}^\lambda \notin \mathbb{N}^{[n]}$. Thus, the $\rho_\lambda$-specialization of this summand is zero, and so $\rho_\lambda(F_1 \star F_2) = 0$.

(b) For $F_1 \in S_{1,p}(u)_{\underline{k}^\lambda}$ and $F_2 \in S_{\underline{\ell}}$, the specialization $\rho_\lambda(F_1 \star F_2)$ is a sum of $\rho_\lambda$-specializations of all summands from the symmetrization appearing in $F_1 \star F_2$. Following the above argument used to prove part (a), only one such specialization is nonzero and we have:

$$\rho_\lambda(F_1 \star F_2) = \rho_\lambda(F_1) \cdot F_2(\{x_{i,r}\}_{i \in [n]}^{1 \le r \le \ell_i}) \cdot P_{\underline{\ell},\lambda} \qquad (3.73)$$

with the factor $P_{\underline{\ell},\lambda}$ denoting the product of the corresponding $\zeta$-factors:

$$P_{\underline{\ell},\lambda} = \prod_{\square=(a,b) \in \lambda} \prod_{i \in [n]} \prod_{1 \le r \le \ell_i} \zeta_{c(\square),i}(q_1^{a-1} q_3^{b-1} u / x_{i,r}).$$

Comparing (3.73) to Lemma 3.9(b) and noting that $P_{\underline{\ell},\lambda}$ is a nonzero multiple of $G_{\underline{\ell},\lambda} Q_{\underline{\ell},\lambda}$, it remains to find $F_1 \in S_{1,p}(u)_{\underline{k}^\lambda}$ such that $\rho_\lambda(F_1) \ne 0$ and $\rho_\mu(F_1 \star F_2) = 0$ for any $F_2 \in S_{\underline{\ell}}$, $\mu > \lambda$. Such an example is $F_1 = K_{\underline{\lambda}'_1} \star \cdots \star K_{\underline{\lambda}'_t} \cdot \prod_{r=1}^{k_p^\lambda}(x_{p,r} - u)^{-1}$, where $\underline{\lambda}'_r \in \mathbb{N}^{[n]}$ encodes the number of color $j \in [n]$ boxes in the $r$-th column of $\lambda$ and $K_{\underline{m}} = \prod_{i \in [n]} \prod_{1 \le r \ne r' \le m_i}(x_{i,r} - q^{-2}x_{i,r'}) \cdot \prod_{i \in [n]} \prod_{1 \le r \le m_i}^{1 \le r' \le m_{i+1}} \frac{x_{i,r} - q_1 x_{i+1,r'}}{x_{i,r} - x_{i+1,r'}}$. □

### *Proof of* (†)

Let us now establish (†), thus completing the proof of Theorem 3.12. Note that:

$$\dim S_{1,p}^m(u) = \sum_{\lambda: |\lambda| \le m} \dim \mathrm{gr}_\lambda(S_{1,p}^m(u)), \qquad \dim \bar{S}_{1,p}^m(u) = \sum_{\lambda: |\lambda| \le m} \dim \mathrm{gr}_\lambda(\bar{S}_{1,p}^m(u)),$$

with the filtration $\{\bar{S}_{1,p}^{m,\lambda}(u)\}_\lambda$ on $\bar{S}_{1,p}^m(u) = \bigoplus_{|\underline{k}|=m} \bar{S}_{1,p}(u)_{\underline{k}}$ induced from $\{S_{1,p}^{m,\lambda}(u)\}_\lambda$ of (3.71). As $\rho_\lambda$-specialization identifies $\mathrm{gr}_\lambda(S_{1,p}^m(u))$ with $\rho_\lambda(S_{1,p}^{m,\lambda}(u))$, Lemma 3.10 and Lemma 3.9(b) imply that $\mathrm{gr}_\lambda(\bar{S}_{1,p}^m(u))$ is zero if $|\lambda| < m$ and is one-dimensional if $|\lambda| = m$ (note that it is nonzero due to the above example of $F_1$). This proves (†). □

### 3.4.2 Generalization for tensor products of Fock and Macmahon

The result of Theorem 3.12 can be generalized in both directions mentioned in the end of Section 3.3.3. Recall the $'\check{U}_{q,d}(\mathfrak{sl}_n)$-action $\pi_{\underline{u},\underline{c}}$ on the space $S(\underline{u})$, which preserves the subspace $S(\underline{u}) \star S'$, cf. Lemma 3.8. Let $\bar{\pi}_{\underline{u},\underline{c}}$ denote the corresponding quotient representation of $'\check{U}_{q,d}(\mathfrak{sl}_n)$ on $\bar{S}(\underline{u}) := S(\underline{u})/(S(\underline{u}) \star S')$. We call a tuple $\underline{u} = \{u_{i,s}\}_{i \in [n]}^{1 \le s \le m_i}$ *generic* if $\{(i, u_{i,s}, (1, \dots, 1))\}$ is generic in the sense of (3.59).

**Theorem 3.13** *We have a $'\check{U}_{q,d}(\mathfrak{sl}_n)$-module isomorphism*

$$\bar{\pi}_{\underline{u},\underline{c}} \simeq \bigotimes_{i=0}^{n-1} \bigotimes_{s=1}^{m_i} \tau_{u_{i,s}}^i$$

*for any generic tuple $\underline{u} = \{u_{i,s}\}_{i \in [n]}^{1 \le s \le m_i}$ (using notation $\tau_u^p$ from Remark 3.12).*

*Proof* The proof proceeds along the same lines as in the above case $\sum m_i = 1$. According to Lemma 3.7, $\bigotimes_{i=0}^{n-1} \bigotimes_{s=1}^{m_i} \tau_{u_{i,s}}^i$ is a well-defined, irreducible, lowest weight $'\check{U}_{q,d}(\mathfrak{sl}_n)$-module with the lowest weight vector $|\emptyset\rangle_{\underline{u}} := \otimes_{i=0}^{n-1} \otimes_{s=1}^{m_i} |\emptyset\rangle$. On the other hand, the vector $\bar{\mathbf{1}}_{\underline{u}} \in \bar{S}(\underline{u})$ is the lowest weight vector of the same lowest weight as $|\emptyset\rangle_{\underline{u}}$. Hence, it suffices to compare the dimensions of their graded components. This can be accomplished exactly as in the proof of Theorem 3.12 by using the specialization maps $\rho_{\underline{\lambda}}$ labeled by tuples $\underline{\lambda} = \{\lambda^{(0,1)}, \dots, \lambda^{(n-1,m_{n-1})}\}$ of Young diagrams. The maps $\rho_{\underline{\lambda}}$ are defined alike $\rho_\lambda$, but now for a box $\square = (a,b) \in \lambda^{(i,s)}$ we fill it with $q_1^{a-1} q_3^{b-1} u_{i,s}$ and color it into $c(\square) := \overline{i-a+b} \in [n]$. ∎

**Corollary 3.2** *$\bar{\pi}_{\underline{u},\underline{c}}$ is an irreducible $'\check{U}_{q,d}(\mathfrak{sl}_n)$-module for generic $\underline{u}$.*

Another generalization of Theorem 3.12 establishes isomorphisms between $\bar{\pi}_{\underline{u},\underline{c}}^K$ and tensor products of vacuum Macmahon modules for generic parameters. To simplify our exposition, we illustrate this in the simplest case of $\pi_{u,\underline{c}}^{p,K}$ for generic $K$. Let $\bar{\pi}_{u,\underline{c}}^{p,K}$ denote the corresponding quotient representation of $'\check{U}_{q,d}(\mathfrak{sl}_n)$ on

$$\bar{S}_{1,p}^K(u) := S_{1,p}^K(u)/(S_{1,p}^K(u) \star S').$$

**Theorem 3.14** *For any $p, \underline{c}$, generic $K$, we have a $'\check{U}_{q,d}(\mathfrak{sl}_n)$-module isomorphism*

$$\bar{S}_{1,p}^K(u) \simeq \mathcal{M}^{(p)}(u,K). \tag{3.74}$$

*Proof* The proof is similar to that of Theorem 3.12, but with the specialization maps and hence filtration pieces parametrized by plane partitions $\bar{\lambda} = (\lambda^1, \lambda^2, \dots)$ of (2.66). For such $\bar{\lambda}$, the specialization map $\rho_{\bar{\lambda}}$ is defined alike $\rho_\lambda$ but with a box $\square = (a,b) \in \lambda^c$ filled with $q_1^{a-1} q_2^{c-1} q_3^{b-1} u$. Let us point out that while the only place where the second kind wheel conditions were used in the proof of Theorem 3.12 was the presence of the factor $\prod_{r=1}^{\ell_p} (x_{p,r} - q^2 u)$ in $G_{\underline{\ell},\lambda}$, its absence now is compensated by a change of $Q_{\underline{\ell},\lambda}$–the factor keeping track of the filtration depth. ∎

### 3.4.3 Twisting vertex representations by Miki's isomorphism

For an action $\rho$ of an algebra $B$ on a vector space $V$ and an algebra homomorphism $\sigma \colon A \to B$, let $\rho^\sigma$ denote as before the corresponding "twisted" action of $A$ on $V$:

$$\rho^\sigma(x) = \rho(\sigma(x)) \quad \text{for any} \quad x \in A \,. \tag{3.75}$$

Recall that Fock ${}'\ddot{U}_{q,d}(\mathfrak{sl}_n)$-module $(\mathcal{F}^{(p)}(u), \tau^p_{1/q_2 u})$ is the lowest weight module in the category $O^-$. As pointed out after Proposition 3.6, the twist $(\mathcal{F}^{(p)}(u), \tau^{p,\vartheta_p}_{1/q_2 u})$ with $\vartheta_p \colon {}'\ddot{U}_{q,d}(\mathfrak{sl}_n) \to {}'\ddot{U}_{q,d}(\mathfrak{sl}_n)$ defined in (3.55) is then the highest weight ${}'\ddot{U}_{q,d}(\mathfrak{sl}_n)$-module in the category $O^+$ with the highest weight $(1; (\phi(z/u)^{-1})^{\delta_{ip}})_{i \in [n]}$. We also consider the ${}'\ddot{U}_{q,d}(\mathfrak{sl}_n)$-module $\rho^{p,\varpi}_{v,\underline{c}}$ obtained by twisting the $\rho^p_{v,\underline{c}}$-action of $\ddot{U}'_{q,d}(\mathfrak{sl}_n)$ on $W(p)_n$ by the isomorphism $\varpi \colon {}'\ddot{U}_{q,d}(\mathfrak{sl}_n) \to \ddot{U}'_{q,d}(\mathfrak{sl}_n)$ of Theorem 3.5.

Then, we have the following result:

**Theorem 3.15** *For any* $0 \le p \le n-1$, $v \in \mathbb{C}^\times$, $\underline{c} = (c_0, c_1, \dots, c_{n-1}) \in (\mathbb{C}^\times)^{[n]}$, *set* $u := (-1)^{\frac{(n-2)(n-3)}{2}} \frac{q^{-1} d^{-n/2}}{c_0 \cdots c_{n-1}}$. *Then, we have a* ${}'\ddot{U}_{q,d}(\mathfrak{sl}_n)$-*module isomorphism:*

$$\rho^{p,\varpi}_{v,\underline{c}} \simeq \tau^{p,\vartheta_p}_{1/q_2 u} \,. \tag{3.76}$$

*Proof* The proof is inspired by that of Theorem 2.5 and will proceed in three steps. First, we verify that $v_0 \otimes e^{\tilde{\Lambda}_p} \in \rho^{p,\varpi}_{v,\underline{c}}$ and $|0\rangle \in \tau^{p,\vartheta_p}_{1/q_2 u}$ have the same eigenvalues with respect to the action of the *finite Cartan subalgebra* $\mathbb{C}[\psi^{\pm 1}_{0,0}, \dots, \psi^{\pm 1}_{n-1,0}, q^{\pm d_2}]$. Second, we show that both vectors are annihilated by all $\{e_{i,r}\}^{r \in \mathbb{Z}}_{i \in [n]}$. Finally, we prove that they have the same eigenvalues with respect to all $\{\psi_{i,l}\}^{l \in \mathbb{Z}}_{i \in [n]}$.

*Step 1: Comparing* $\mathbb{C}[\psi^{\pm 1}_{0,0}, \dots, \psi^{\pm 1}_{n-1,0}, q^{\pm d_2}]$-*weights.*

By Proposition 3.10, the generators $\psi_{i,0}, \gamma, q^{d_1}$ act on $v_0 \otimes e^{\tilde{\Lambda}_p}$ via the multiplication by $q^{\langle \tilde{h}_i, \tilde{\Lambda}_p \rangle}, q, 1$, respectively. Evoking Proposition 3.2(a), we thus obtain:

$$\rho^{p,\varpi}_{v,\underline{c}}(q^{d_2}) \, v_0 \otimes e^{\tilde{\Lambda}_p} = q^{\frac{p(n-p)}{2}} \cdot v_0 \otimes e^{\tilde{\Lambda}_p} \,,$$
$$\rho^{p,\varpi}_{v,\underline{c}}(\psi_{i,0}) \, v_0 \otimes e^{\tilde{\Lambda}_p} = q^{\delta_{ip}} \cdot v_0 \otimes e^{\tilde{\Lambda}_p} \quad \forall \, i \in [n] \,.$$

We also have: $\tau^p_u(q^{d_2})|0\rangle = q^{-\frac{p(n-p)}{2}} \cdot |0\rangle$, $\tau^p_u(\psi_{i,0})|0\rangle = q^{-\delta_{ip}} \cdot |0\rangle$. Thus $|0\rangle \in \tau^{p,\vartheta_p}_{1/q_2 u}$ and $v_0 \otimes e^{\tilde{\Lambda}_p} \in \rho^{p,\varpi}_{v,\underline{c}}$ have the same weights with respect to $\mathbb{C}[\psi^{\pm 1}_{0,0}, \dots, \psi^{\pm 1}_{n-1,0}, q^{\pm d_2}]$.

*Remark 3.13* This explains the power of $q$ in the formulas for $\tau^p_{u,\underline{c}}(q^{d_2})$, $\pi^p_{u,\underline{c}}(q^{d_2})$. □

*Step 2: Verifying* ${}'\ddot{U}^>_{q,d}$-*annihilation.*

First, let us prove that $\rho^{p,\varpi}_{v,\underline{c}}(e_{i,0}) v_0 \otimes e^{\tilde{\Lambda}_p} = 0$ for any $i \in [n]$. For $i \ne 0$, this is clear as $\langle \tilde{h}_i, \tilde{\Lambda}_p \rangle + 1 > 0$, while $H_{j,k} v_0 = 0$ for all $j \in [n], k > 0$. For $i = 0$, we have $\varpi(e_{0,0}) = d \cdot \gamma \psi_{0,0} \cdot [\cdots [f_{1,1}, f_{2,0}]_q, \cdots, f_{n-1,0}]_q$ by Proposition 3.2(a), and so the equality $\rho^{p,\varpi}_{v,\underline{c}}(e_{0,0}) v_0 \otimes e^{\tilde{\Lambda}_p} = 0$ follows from the next result:

**Lemma 3.11** $\rho^p_{v,\underline{c}}([\cdots[f_{1,1},f_{2,0}]_q,\cdots,f_{n-1,0}]_q)\,v_0\otimes e^{\tilde{\Lambda}_p}=0.$ □

*Proof* Let us show that any summand $f_{i_{n-1},r_{n-1}}\cdots f_{i_1,r_1}$ from the multicommutator $[\cdots[f_{1,1},f_{2,0}]_q,\cdots,f_{n-1,0}]_q$ (thus, $\{i_1,\ldots,i_{n-1}\}=[n]^\times$ and $r_k=\delta_{i_k 1}$) acts trivially on $v_0\otimes e^{\tilde{\Lambda}_p}$. Since $-\langle\bar{h}_i,\bar{\Lambda}_p\rangle+1>0$ for $i\neq p$, we see that $\rho^p_{v,\underline{c}}(f_{i_1,r_1})v_0\otimes e^{\tilde{\Lambda}_p}=0$ unless $i_1=p\neq 1$. For $i_1=p\neq 1$, we get $\rho^p_{v,\underline{c}}(f_{i_1,r_1})v_0\otimes e^{\tilde{\Lambda}_p}=\pm c_{i_1}^{-1}v_0\otimes e^{\tilde{\Lambda}_p^{(1)}}$ with $\bar{\Lambda}_p^{(1)}:=\bar{\Lambda}_p-\bar{\alpha}_{i_1}$. The key property of this weight is $-\langle\bar{h}_i,\bar{\Lambda}_p^{(1)}\rangle+1\geq 0$ with an equality only for $i=p\pm 1$. In particular, $\rho^p_{v,\underline{c}}(f_{i_2,r_2})v_0\otimes e^{\tilde{\Lambda}_p^{(1)}}=0$ unless $i_2=p\pm 1\neq 1$. In the latter cases, the result is a multiple of $v_0\otimes e^{\tilde{\Lambda}_p^{(2)}}$ with $\bar{\Lambda}_p^{(2)}:=\bar{\Lambda}_p^{(1)}-\bar{\alpha}_{i_2}$ satisfying a similar property. Iterating the argument, we finally get to the $k$-th place with $i_k=1$ and $r_k=1$. As $-\langle\bar{h}_1,\bar{\Lambda}_p^{(k-1)}\rangle+1\geq 0$, we have $\rho^p_{v,\underline{c}}(f_{i_k,r_k}\cdots f_{i_1,r_1})v_0\otimes e^{\tilde{\Lambda}_p}=0.$ □

This completes the proof of

$$\rho^{p,\varpi}_{v,\underline{c}}(e_{i,0})\,v_0\otimes e^{\tilde{\Lambda}_p}=0\qquad\forall\,i\in[n]\,.$$

According to the formulas (3.78, 3.79) below, we have $\rho^{p,\varpi}_{v,\underline{c}}(h_{j,\pm 1})v_0\otimes e^{\tilde{\Lambda}_p}=0$ for $j\neq p$. Combining this with $[h_{j,\pm 1},e_{i,k}]=d^{\mp m_{ji}}\gamma^{-1/2}[a_{ji}]_q\cdot e_{i,k\pm 1}$ of (T5′), we get:

$$\rho^{p,\varpi}_{v,\underline{c}}(e_{i,k})\,v_0\otimes e^{\tilde{\Lambda}_p}=0\qquad\forall\,i\in[n],k\in\mathbb{Z}\,.$$

On the other hand, $|\emptyset\rangle\in\tau^{p,\vartheta_p}_{1/q_2u}$ is the highest weight vector and so

$$\tau^{p,\vartheta_p}_{1/q_2u}(e_{i,k})|\emptyset\rangle=0\qquad\forall\,i\in[n],k\in\mathbb{Z}\,.$$

<u>*Step 3: Comparing* $'\ddot{U}^0_{q,d}$-*weights.*</u>

Let us now prove that both $v_0\otimes e^{\tilde{\Lambda}_p}\in\rho^{p,\varpi}_{v,\underline{c}}$ and $|\emptyset\rangle\in\tau^{p,\vartheta_p}_{1/q_2u}$ are eigenvectors of all Cartan generators $\{\psi_{i,l}\}^{l\neq 0}_{i\in[n]}$ with the same eigenvalues. By the construction:

$$\tau^{p,\vartheta_p}_{1/q_2u}(\psi_{i,\pm r})|\emptyset\rangle=\pm\delta_{ip}(q-q^{-1})(q^2u)^{\pm r}|\emptyset\rangle\qquad\forall\,i\in[n],r>0\,.$$

Therefore, it remains to prove the following formula, cf. (2.69):

$$\rho^{p,\varpi}_{v,\underline{c}}(\psi_{i,\pm r})\,v_0\otimes e^{\tilde{\Lambda}_p}=\pm\delta_{ip}(q-q^{-1})(q^2u)^{\pm r}v_0\otimes e^{\tilde{\Lambda}_p}\quad\forall\,i\in[n],\,r>0\,.\quad(\ddagger)$$

The proof of $(\ddagger)$ is based on the following technical result:

**Lemma 3.12** *We have the following equalities:*

$$\rho^{p,\varpi}_{v,\underline{c}}(f_{i,0})\,v_0\otimes e^{\tilde{\Lambda}_p}=\delta_{ip}c_p^{-1}\varepsilon^{\delta_{p0}}\cdot v_0\otimes e^{-\bar{\alpha}_p}e^{\tilde{\Lambda}_p}\,,\tag{3.77}$$

$$\rho^{p,\varpi}_{v,\underline{c}}(h_{i,-1})\,v_0\otimes e^{\tilde{\Lambda}_p}=\delta_{ip}q^{-1}u^{-1}\cdot v_0\otimes e^{\tilde{\Lambda}_p}\,,\tag{3.78}$$

$$\rho^{p,\varpi}_{v,\underline{c}}(h_{i,1})\,v_0\otimes e^{\tilde{\Lambda}_p}=\delta_{ip}qu\cdot v_0\otimes e^{\tilde{\Lambda}_p}\,,\tag{3.79}$$

$$\rho_{v,\underline{c}}^{p,\varpi}(h_{p,-1})\, v_0 \otimes e^{-\bar{\alpha}_p}\, e^{\bar{\Lambda}_p} = -q^{-3} u^{-1} \cdot v_0 \otimes e^{-\bar{\alpha}_p}\, e^{\bar{\Lambda}_p}\,, \qquad (3.80)$$

$$\rho_{v,\underline{c}}^{p,\varpi}(h_{p,1})\, v_0 \otimes e^{-\bar{\alpha}_p}\, e^{\bar{\Lambda}_p} = -q^3 u \cdot v_0 \otimes e^{-\bar{\alpha}_p}\, e^{\bar{\Lambda}_p}\,, \qquad (3.81)$$

*where the constants $u$ and $\mathfrak{e}$ are defined through $\mathfrak{c} = \prod_{j \in [n]} c_j$ via:*

$$u := (-1)^{(n-2)(n-3)/2} \mathfrak{c}^{-1} q^{-1} d^{-n/2} \quad \text{and} \quad \mathfrak{e} := (-1)^{(n-2)(n-3)/2} \mathfrak{c} v^{-1} q^{-1} d^{-n/2}\,.$$

*Proof* This lemma is proved by rather tedious computations involving the formulas from Propositions 3.2 and 3.10 (thus, it can be skipped at the first reading).

• For $i \neq 0$, we have $\varpi(f_{i,0}) = f_{i,0}$ and $-\langle \bar{h}_i, \bar{\Lambda}_p \rangle + 1 = 1 - \delta_{ip} \geq 0$, so that

$$\rho_{v,\underline{c}}^{p,\varpi}(f_{i,0})\, v_0 \otimes e^{\bar{\Lambda}_p} = \delta_{ip} c_p^{-1} \cdot v_0 \otimes e^{-\bar{\alpha}_p}\, e^{\bar{\Lambda}_p}$$

thus establishing (3.77) for $i \neq 0$. For $i = 0$, we apply Proposition 3.2(a) to get:

$$\rho_{v,\underline{c}}^{p,\varpi}(f_{0,0})\, v_0 \otimes e^{\bar{\Lambda}_p} = q^{-\delta_{p0}} d^{-1} \rho_{v,\underline{c}}^p([e_{n-1,0}, \cdots, [e_{2,0}, e_{1,-1}]_{q^{-1}} \cdots ]_{q^{-1}})\, v_0 \otimes e^{\bar{\Lambda}_p}\,.$$

Expanding this multicommutator and using $\rho_{v,\underline{c}}^p(e_{j,0}) v_0 \otimes e^{\bar{\Lambda}_p} = 0$ for $j \neq 0$, we get:

$$\rho_{v,\underline{c}}^p([e_{n-1,0}, \cdots, [e_{2,0}, e_{1,-1}]_{q^{-1}} \cdots ]_{q^{-1}})\, v_0 \otimes e^{\bar{\Lambda}_p} =$$
$$\rho_{v,\underline{c}}^p(e_{n-1,0}) \cdots \rho_{v,\underline{c}}^p(e_{2,0}) \rho_{v,\underline{c}}^p(e_{1,-1})\, v_0 \otimes e^{\bar{\Lambda}_p}\,.$$

For $p \neq 0$, one shows that $\rho_{v,\underline{c}}^p(e_{p,0}) \cdots \rho_{v,\underline{c}}^p(e_{2,0}) \rho_{v,\underline{c}}^p(e_{1,-1}) v_0 \otimes e^{\bar{\Lambda}_p} = 0$ by applying the same argument as above, while for $p = 0$ we obtain:

$$\rho_{v,\underline{c}}^p(e_{n-1,0}) \cdots \rho_{v,\underline{c}}^p(e_{2,0}) \rho_{v,\underline{c}}^p(e_{1,-1})\, v_0 \otimes e^0 =$$
$$v^{-1} d^{(2-n)/2} (c_1 \cdots c_{n-1}) \cdot v_0 \otimes (e^{\bar{\alpha}_{n-1}} \cdots e^{\bar{\alpha}_1})\,,$$

with the power of $d$ arising from $(z/v)^{H_{*,0}}$. Since $e^{\bar{\alpha}_{n-1}} \cdots e^{\bar{\alpha}_1} = (-1)^{(n-2)(n-3)/2} e^{-\bar{\alpha}_0}$, we finally get:

$$\rho_{v,\underline{c}}^{0,\varpi}(f_{0,0})\, v_0 \otimes e^0 = (-1)^{(n-2)(n-3)/2} v^{-1} q^{-1} d^{-n/2} \mathfrak{c} c_0^{-1} \cdot v_0 \otimes e^{-\bar{\alpha}_0}\,.$$

Let us now establish (3.78)–(3.81) for $p \neq 0$ (leaving $p = 0$ case to the reader).

• Combining the formula for $\varpi(h_{0,-1})$ from Proposition 3.2(c) with

$$\rho_{v,\underline{c}}^p(e_{0,1})\, v_0 \otimes e^{\bar{\Lambda}_p} = \rho_{v,\underline{c}}^p(e_{2,0})\, v_0 \otimes e^{\bar{\Lambda}_p} = \cdots = \rho_{v,\underline{c}}^p(e_{n-1,0})\, v_0 \otimes e^{\bar{\Lambda}_p} = 0\,,$$

we obtain:

$$\rho_{v,\underline{c}}^{p,\varpi}(h_{0,-1})\, v_0 \otimes e^{\bar{\Lambda}_p} =$$
$$(-1)^n d^{n-1} \cdot \rho_{v,\underline{c}}^p(e_{0,1}) \rho_{v,\underline{c}}^p(e_{n-1,0}) \cdots \rho_{v,\underline{c}}^p(e_{2,0}) \rho_{v,\underline{c}}^p(e_{1,-1})\, v_0 \otimes e^{\bar{\Lambda}_p}\,. \qquad (3.82)$$

Since $\rho^p_{v,\underline{c}}(e_{p-1,0}) \cdots \rho^p_{v,\underline{c}}(e_{1,-1}) v_0 \otimes e^{\bar{\Lambda}_p}$ is a multiple of $v_0 \otimes e^{\bar{\Lambda}_p + \bar{\alpha}_1 + \cdots + \bar{\alpha}_{p-1}}$ and $\langle \bar{h}_p, \bar{\Lambda}_p + \bar{\alpha}_1 + \cdots + \bar{\alpha}_{p-1} \rangle + 1 > 0$, the right-hand side of (3.82) vanishes. Hence, so does the left-hand side, that is, $\rho^{p,\varpi}_{v,\underline{c}}(h_{0,-1}) v_0 \otimes e^{\bar{\Lambda}_p} = 0$ for $p \neq 0$.

For $i \neq 0$, the formula for $\varpi(h_{i,-1})$ and $\rho^p_{v,\underline{c}}(e_{j,0}) v_0 \otimes e^{\bar{\Lambda}_p} = 0$ ($j \neq 0$) imply:

$$\rho^{p,\varpi}_{v,\underline{c}}(h_{i,-1}) v_0 \otimes e^{\bar{\Lambda}_p} =$$
$$(-1)^{i+1} d^i \rho^p_{v,\underline{c}}(e_{i,0}) \cdots \rho^p_{v,\underline{c}}(e_{1,0}) \rho^p_{v,\underline{c}}(e_{i+1,0}) \cdots \rho^p_{v,\underline{c}}(e_{n-1,0}) \rho^p_{v,\underline{c}}(e_{0,0}) v_0 \otimes e^{\bar{\Lambda}_p} .$$

If $0 < i < p$, then $\langle \bar{h}_p, \bar{\Lambda}_p + \bar{\alpha}_0 + \sum_{j=p+1}^{n-1} \bar{\alpha}_j \rangle + 1 > 0$, and thus we get:

$$\rho^p_{v,\underline{c}}(e_{p,0}) \cdots \rho^p_{v,\underline{c}}(e_{n-1,0}) \rho^p_{v,\underline{c}}(e_{0,0}) v_0 \otimes e^{\bar{\Lambda}_p} = 0 \implies \rho^{p,\varpi}_{v,\underline{c}}(h_{i,-1}) v_0 \otimes e^{\bar{\Lambda}_p} = 0 .$$

Likewise, if $p < i < n$, then $\langle \bar{h}_p, \bar{\Lambda}_p + \bar{\alpha}_0 + \sum_{j=i+1}^{n-1} \bar{\alpha}_j + \sum_{j=1}^{p-1} \bar{\alpha}_j \rangle + 1 > 0$ and hence $\rho^{p,\varpi}_{v,\underline{c}}(h_{i,-1}) v_0 \otimes e^{\bar{\Lambda}_p} = 0$. Finally, for $i = p$:

$$\rho^{p,\varpi}_{v,\underline{c}}(h_{p,-1}) v_0 \otimes e^{\bar{\Lambda}_p} = (-1)^{p+1} d^{n/2} \mathfrak{c} \cdot v_0 \otimes (e^{\bar{\alpha}_p} \cdots e^{\bar{\alpha}_1} e^{\bar{\alpha}_{p+1}} \cdots e^{\bar{\alpha}_{n-1}} e^{\bar{\alpha}_0} e^{\bar{\Lambda}_p}) =$$
$$(-1)^{(n-2)(n-3)/2} d^{n/2} \mathfrak{c} \cdot v_0 \otimes e^{\bar{\Lambda}_p} ,$$

due to $e^{\bar{\alpha}_p} \cdots e^{\bar{\alpha}_1} e^{\bar{\alpha}_{p+1}} \cdots e^{\bar{\alpha}_{n-1}} e^{\bar{\alpha}_0} = (-1)^{n(n-1)/2+p} e^0$. This establishes (3.78).

- By similar arguments we have $\rho^{p,\varpi}_{v,\underline{c}}(h_{i,1}) v_0 \otimes e^{\bar{\Lambda}_p} = 0$ if $i \neq p$. For $i = p$, we get:

$$\rho^{p,\varpi}_{v,\underline{c}}(h_{p,1}) v_0 \otimes e^{\bar{\Lambda}_p} = (-1)^{n+p+1} d^{-p} \mathfrak{c}^{-1} d^{p-n/2} \times$$
$$v_0 \otimes (e^{-\bar{\alpha}_0} e^{-\bar{\alpha}_{n-1}} \cdots e^{-\bar{\alpha}_{p+1}} e^{-\bar{\alpha}_1} \cdots e^{-\bar{\alpha}_p} e^{\bar{\Lambda}_p}) = (-1)^{(n-2)(n-3)/2} d^{-n/2} \mathfrak{c}^{-1} \cdot v_0 \otimes e^{\bar{\Lambda}_p} ,$$

due to $e^{-\bar{\alpha}_0} e^{-\bar{\alpha}_{n-1}} \cdots e^{-\bar{\alpha}_{p+1}} e^{-\bar{\alpha}_1} \cdots e^{-\bar{\alpha}_p} = (-1)^{n(n+1)/2+p} e^0$. This establishes (3.79).

- Again, only one summand of the corresponding multicommutator acts nontrivially:

$$\rho^{p,\varpi}_{v,\underline{c}}(h_{p,-1}) v_0 \otimes e^{-\bar{\alpha}_p} e^{\bar{\Lambda}_p} = (-1)^{p+1} d^p (-q^{-2}) \times$$
$$\rho^p_{v,\underline{c}}(e_{p-1,0}) \cdots \rho^p_{v,\underline{c}}(e_{1,0}) \rho^p_{v,\underline{c}}(e_{p+1,0}) \cdots \rho^p_{v,\underline{c}}(e_{0,0}) \rho^p_{v,\underline{c}}(e_{p,0}) v_0 \otimes e^{-\bar{\alpha}_p} e^{\bar{\Lambda}_p} =$$
$$(-1)^{1+(n-2)(n-3)/2} d^{n/2} q^{-2} \mathfrak{c} \cdot v_0 \otimes e^{-\bar{\alpha}_p} e^{\bar{\Lambda}_p} .$$

This establishes (3.80) for $p \neq 0$.

- Again, only one summand of the corresponding multicommutator acts nontrivially:

$$\rho^{p,\varpi}_{v,\underline{c}}(h_{p,1}) v_0 \otimes e^{-\bar{\alpha}_p} e^{\bar{\Lambda}_p} = (-1)^{p+1} d^{-p} (-q^2) \times$$
$$\rho^p_{v,\underline{c}}(f_{p,0}) \rho^p_{v,\underline{c}}(f_{0,0}) \cdots \rho^p_{v,\underline{c}}(f_{p+1,0}) \rho^p_{v,\underline{c}}(f_{1,0}) \cdots \rho^p_{v,\underline{c}}(f_{p-1,0}) v_0 \otimes e^{-\bar{\alpha}_p} e^{\bar{\Lambda}_p} =$$
$$(-1)^{1+(n-2)(n-3)/2} d^{-n/2} q^2 \mathfrak{c}^{-1} \cdot v_0 \otimes e^{-\bar{\alpha}_p} e^{\bar{\Lambda}_p} .$$

This completes our proof of (3.81) for $p \neq 0$. □

Let us finally prove (‡). Note that $\rho_{v,\underline{c}}^{p,\varpi}(e_{p,0})v_0 \otimes e^{-\bar{\alpha}_p}e^{\bar{\Lambda}_p} = c_p e^{-\delta_{p0}} \cdot v_0 \otimes e^{\bar{\Lambda}_p}$. Combining this with the identity $[h_{p,\pm1}, e_{p,\pm r}] = \gamma^{-1/2}(q + q^{-1})e_{p,\pm(r+1)}$ of (T5′) and the equalities (3.78)–(3.81) of Lemma 3.12, we obtain by induction on $r$:

$$\rho_{v,\underline{c}}^{p,\varpi}(e_{p,\pm r})\, v_0 \otimes e^{-\bar{\alpha}_p}e^{\bar{\Lambda}_p} = c_p e^{-\delta_{p0}}(q^2 u)^{\pm r} \cdot v_0 \otimes e^{\bar{\Lambda}_p} \qquad \forall\, r > 0\,.$$

On the other hand, according to (T4), we also have:

$$\rho_{v,\underline{c}}^{p,\varpi}(\psi_{i,\pm r}) = \pm(q - q^{-1})[\rho_{v,\underline{c}}^{p,\varpi}(e_{i,\pm r}), \rho_{v,\underline{c}}^{p,\varpi}(f_{i,0})] \qquad \forall\, r > 0\,.$$

Since $\rho_{v,\underline{c}}^{p,\varpi}(e_{i,\pm r})v_0 \otimes e^{\bar{\Lambda}_p} = 0 = \rho_{v,\underline{c}}^{p,\varpi}(f_{i,0})v_0 \otimes e^{\bar{\Lambda}_p}$ for $i \neq p$, see (3.77), we get $\rho_{v,\underline{c}}^{p,\varpi}(\psi_{i,\pm r})v_0 \otimes e^{\bar{\Lambda}_p} = 0$ if $i \neq p$. For $i = p$, the equality (‡) follows now from

$$\rho_{v,\underline{c}}^{p,\varpi}(\psi_{p,\pm r})\, v_0 \otimes e^{\bar{\Lambda}_p} = \pm(q - q^{-1})\rho_{v,\underline{c}}^{p,\varpi}(e_{p,\pm r})\rho_{v,\underline{c}}^{p,\varpi}(f_{p,0})\, v_0 \otimes e^{\bar{\Lambda}_p} =$$
$$\pm (q - q^{-1})(q^2 u)^{\pm r} \cdot v_0 \otimes e^{\bar{\Lambda}_p}\,.$$

### Step 4: Completion of the proof

According to Proposition 3.8 and Remark 3.10 both $\rho_{v,\underline{c}}^p$ and $\tau_{1/q_2 u}^p$ are irreducible (assuming $q, d$ are generic (G)), hence, so are $'\ddot{U}_{q,d}(\mathfrak{sl}_n)$-modules $\rho_{v,\underline{c}}^{p,\varpi}$ and $\tau_{1/q_2 u}^{p,\vartheta_p}$. Moreover, $v_0 \otimes e^{\bar{\Lambda}_p} \in \rho_{v,\underline{c}}^{p,\varpi}$ and $|\emptyset\rangle \in \tau_{1/q_2 u}^{p,\vartheta_p}$ are the highest weight vectors of the same highest weight, due to Steps 1–3 above. This completes the proof of Theorem 3.15.□

**Corollary 3.3** *For any* $0 \leq p \leq n - 1$, $v, v' \in \mathbb{C}^\times$, *and* $\underline{c}, \underline{c}' \in (\mathbb{C}^\times)^{[n]}$ *satisfying* $\prod_{i \in [n]} c_i = \prod_{i \in [n]} c'_i$, *we have a* $\ddot{U}'_{q,d}(\mathfrak{sl}_n)$-*module isomorphism:* $\rho_{v,\underline{c}}^p \simeq \rho_{v',\underline{c}'}^p$.

*Remark 3.14* Analogously to Theorem 3.15, one can show that $\tau_u^p \simeq (\rho_{v,\underline{c}}^{p,*})^\varpi$, similarly to Remark 2.13. Here, $(\rho_{v,\underline{c}}^{p,*})^\varpi$ is the $\varpi$-twist of the $\ddot{U}'_{q,d}(\mathfrak{sl}_n)$-module $\rho_{v,\underline{c}}^{p,*}$, dual to $\rho_{v,\underline{c}}^p$ with respect to the antipode (H3). Consider the algebra anti-automorphism $\flat\colon \mathfrak{A}_n \to \mathfrak{A}_n$ defined by $H_0^\flat = H_0$, $H_{i,-k}^\flat = \sum_{j \in [n]} d^{km_{ij}}\frac{[k]_q[ka_{ij}]_q}{k} \cdot H_{j,k}$, $H_{i,k}^\flat = \widetilde{H}_{i,-k}$ for $i \in [n], k > 0$ with $\widetilde{H}_{i,-k} \in \text{span}_{\mathbb{C}}\{H_{j,-k}\}_{j \in [n]}$ such that $[H_{j,l}, \widetilde{H}_{i,-k}] = \delta_{ij}\delta_{kl}H_0$ for any $i, j \in [n]$ and $k, l > 0$. We define the bilinear form $(\cdot, \cdot)$ on $F_n \otimes \mathbb{C}\{\tilde{P}\}$ via

$$(P \otimes e^{\bar{\beta}}, Q \otimes e^{\bar{\gamma}}) = \delta_{\bar{\beta},\bar{\gamma}} \cdot \langle v_0 \mid P^\flat Q \mid v_0 \rangle\,.$$

Then, the dual vector space $W(p)_n^*$ can be naturally identified with $W(p)_n$, and thus the dual representation $\rho_{v,\underline{c}}^{p,*}$ can be naturally realized on $W(p)_n$ by vertex-type formulas: $\rho_{v,\underline{c}}^{p,*}(e_i(z)) = -(-1)^{n\delta_{i0}}c_i \cdot \exp(-\sum_{k>0} \frac{q^{k/2}}{[k]_q}H_{i,k}^\flat (\frac{z}{v})^{-k}) \cdot \exp(\sum_{k>0} \frac{q^{k/2}}{[k]_q}H_{i,-k}^\flat (\frac{z}{v})^k) \cdot e^{-\bar{\alpha}_i}(\frac{z}{qv})^{H_{i,0}-1}q^{\partial_{\bar{\alpha}_i}}$, $\rho_{v,\underline{c}}^{p,*}(\psi_i^\pm(z)) = \exp(\mp(q - q^{-1})\sum_{k>0} H_{i,\pm k}^\flat (\frac{z}{v})^{\mp k}) \cdot q^{\mp\partial_{\bar{\alpha}_i}}$, $\rho_{v,\underline{c}}^{p,*}(f_i(z)) = -\frac{1}{c_i q^2} \cdot \exp(\sum_{k>0} \frac{q^{-k/2}}{[k]_q}H_{i,k}^\flat (\frac{z}{v})^{-k}) \cdot \exp(-\sum_{k>0} \frac{q^{-k/2}}{[k]_q}H_{i,-k}^\flat (\frac{z}{v})^k) \cdot e^{\bar{\alpha}_i}(\frac{z}{qv})^{-H_{i,0}-1}q^{-\partial_{\bar{\alpha}_i}}$, $\rho_{v,\underline{c}}^{p,*}(\gamma^{\pm1/2}) = q^{\pm1/2}$, and $\rho_{v,\underline{c}}^{p,*}(q^{\pm d_1}) = q^{\mp d}$.

## 3.5 Bethe algebra realization of $\mathcal{A}(s_0,\ldots,s_{n-1})$

In this section, we finally prove the commutativity of the subalgebras $\mathcal{A}(\underline{s})$ from Section 3.2.2 (and their one-parameter deformations) for generic $\underline{s} = (s_0,\ldots,s_{n-1})$ satisfying $\prod_{i\in[n]} s_i = 1$ by providing a transfer matrix realization of their generators $F_k^\mu(\underline{s})$. To this end, we study functionals on $\ddot{U}'^{\leq}_{q,d}$ (arising as linear combinations of matrix coefficients), and realize them as pairings with explicit elements of $\ddot{U}'^{\geq}_{q,d}$ via $\varphi'$ from Theorem 3.4(c). This can be interpreted via the universal $R$-matrix of Theorem 3.2, thus relating $\mathcal{A}(\underline{s})$ to the standard Bethe subalgebras of $U_q(L\mathfrak{gl}_n)$.

### 3.5.1 Trace functionals

In this section, we introduce and evaluate three functionals on $\ddot{U}'^{\leq}_{q,d}$. The corresponding computations use the standard technique of the *normal ordering* (which is also how Proposition 3.10 is proved, see [22]). Explicitly, it is a procedure that for any given unordered product $P$ of terms $\exp(yH_{i,-k})$, $\exp(xH_{i,k})$, $e^{\tilde{\alpha}_i}$, $z^{H_{i,0}}$, $q^{\partial_{\tilde{\alpha}_i}}$ (with $i\in[n], k>0$) reorders them to $:P:$ by placing specifically in the order indicated above (the relative order of $\{\exp(yH_{i,-k})\}$ as well as $\{\exp(xH_{i,k})\}$ is not important). When doing so, we note that most of the terms pairwise commute except for:

$$\exp(xH_{i,k})\exp(yH_{j,-l}) =$$
$$\exp\left(\delta_{kl}\frac{xyd^{-km_{ij}}[k]_q[ka_{ij}]_q}{k}\right)\cdot\exp(yH_{j,-l})\exp(xH_{i,k}) \quad (3.83)$$

(due to the defining relations (3.60)) as well as

$$\begin{aligned}
e^{\tilde{\alpha}_i}e^{\tilde{\alpha}_j} &= (-1)^{\langle\bar{h}_i,\bar{\alpha}_j\rangle}\cdot e^{\tilde{\alpha}_j}e^{\tilde{\alpha}_i}, \\
q^{\partial_{\tilde{\alpha}_i}}e^{\tilde{\alpha}_j} &= q^{\langle\bar{h}_i,\bar{\alpha}_j\rangle}\cdot e^{\tilde{\alpha}_j}q^{\partial_{\tilde{\alpha}_i}}, \\
z^{H_{i,0}}e^{\tilde{\alpha}_j} &= z^{\langle\bar{h}_i,\bar{\alpha}_j\rangle}d^{m_{ij}\langle\bar{h}_i,\bar{\alpha}_j\rangle/2}\cdot e^{\tilde{\alpha}_j}z^{H_{i,0}}.
\end{aligned} \quad (3.84)$$

- *Top matrix coefficient*

Consider the functional

$$\phi^0_{p,\underline{c}}:\ddot{U}'^{\leq}_{q,d}\longrightarrow\mathbb{C} \quad\text{defined by}\quad \phi^0_{p,\underline{c}}(y):=\left\langle \mathrm{v}_0\otimes e^{\bar{\Lambda}_p}\middle|\rho^p_{1,\underline{c}}(y)\middle|\mathrm{v}_0\otimes e^{\bar{\Lambda}_p}\right\rangle \quad (3.85)$$

evaluating the *top-to-top* matrix coefficient in the vertex $\ddot{U}'^{\leq}_{q,d}$-representation $\rho^p_{1,\underline{c}}$, cf. Corollary 3.3. Since $\rho^p_{1,\underline{c}}(h_{i,k})\mathrm{v}_0\otimes e^{\bar{\Lambda}_p} = 0$ for any $i\in[n], k>0$, the functional (3.85) is determined by its values at

$$y = f_{i_1,k_1}f_{i_2,k_2}\cdots f_{i_m,k_m}\cdot\psi^{r_0}_{0,0}\cdots\psi^{r_{n-1}}_{n-1,0}\cdot(\gamma^{1/2})^a(q^{d_1})^b$$

with $a, b \in \mathbb{Z}$, $\underline{r} := (r_0, \ldots, r_{n-1}) \in \mathbb{Z}^{[n]}$, $m \in \mathbb{N}$, $k_1, \ldots, k_m \in \mathbb{Z}$, and $i_1, \ldots, i_m \in [n]$ satisfying $\sum_{s=1}^m \bar{\alpha}_{i_s} = 0 \in \bar{Q}$. The latter is equivalent to $\{i_1, \ldots, i_m\}$ containing an equal number of each of the indices $\{0, \ldots, n-1\}$. Evoking the quadratic relation (T3), we thus conclude that $\phi^0_{p,\underline{c}}$ is determined by the following generating series:

$$\phi^0_{p,\underline{c};N,\underline{r},a,b}(z_{0,1}, \ldots, z_{n-1,N}) :=$$
$$\phi^0_{p,\underline{c}}\left(\prod_{k=1}^N \big(f_0(z_{0,k}) \cdots f_{n-1}(z_{n-1,k})\big) \cdot \prod_{i \in [n]} \psi^{r_i}_{i,0} \cdot \gamma^{a/2} q^{bd_1}\right) \quad (3.86)$$

for any $N \in \mathbb{N}$, $a, b \in \mathbb{Z}$, $\underline{r} \in \mathbb{Z}^{[n]}$, with the above $z_{*,*}$-variables ordered as follows:

$$z_{0,1}, \ldots, z_{n-1,1}, z_{0,2}, \ldots, z_{n-1,2}, \ldots, z_{0,N}, \ldots, z_{n-1,N}. \quad (3.87)$$

**Proposition 3.12** *Set $\mathfrak{c} = c_0 \cdots c_{n-1}$ and $\varepsilon = (-1)^{(n-1)(n+2)/2}$. We have:*

$$\phi^0_{p,\underline{c};N,\underline{r},a,b}(z_{0,1}, \ldots, z_{n-1,N}) = (\varepsilon \mathfrak{c})^{-N} q^{a/2+r_p-r_0} d^{\frac{N(n-2)}{2}} \cdot \prod_{k=1}^N \frac{z_{0,k}}{z_{p,k}} \times$$

$$\frac{\prod_{i \in [n]} \prod_{1 \le k < k' \le N}(z_{i,k} - z_{i,k'})(z_{i,k} - q^2 z_{i,k'}) \cdot \prod_{i \in [n]} \prod_{k=1}^N z_{i,k}}{\prod_{i \in [n]} \prod_{1 \le k \le k' \le N}(z_{i,k} - q d z_{i+1,k'}) \cdot \prod_{i \in [n]} \prod_{1 \le k < k' \le N}(z_{i,k} - q d^{-1} z_{i-1,k'})},$$

*where actually $k < k'$ for $i = n-1$ in the first product and $k \le k'$ for $i = 0$ in the second product from the denominator (that is, $z_{i,k}$ is to the left of $z_{i\pm 1,k'}$ in (3.87)).*

*Proof* This straightforward computation follows by normally ordering the product $\prod_{k=1}^N \big(f_0(z_{0,k}) \cdots f_{n-1}(z_{n-1,k})\big)$ using the equalities

$$\exp\left(\sum_{k>0} \frac{q^k d^{\pm k}}{k}\left(\frac{y}{x}\right)^k\right) = \left(1 - \frac{q d^{\pm 1} y}{x}\right)^{-1}, \quad (3.88)$$

$$\exp\left(-\sum_{k>0} \frac{q^k [2k]_q}{k [k]_q}\left(\frac{y}{x}\right)^k\right) = \left(1 - \frac{y}{x}\right)\left(1 - \frac{q^2 y}{x}\right), \quad (3.89)$$

and the commutation rules (3.83, 3.84). We note that $\varepsilon$ appears due to the equality

$$(-1)^n e^{-\bar{\alpha}_0} e^{-\bar{\alpha}_1} \cdots e^{-\bar{\alpha}_{n-1}} = \varepsilon e^0. \quad (3.90)$$

We leave the remaining details to the interested reader. $\qquad\square$

● *Top level graded trace*

Recall the degree operator d of (3.62) acting diagonally in the natural basis (3.61) of $W(p)_n$ with all eigenvalues in $\mathbb{N}$. Let $M(p)_n := \mathrm{Ker}(\mathrm{d}) \subset W(p)_n$ denote its kernel. Finally, let $U_q(\mathfrak{sl}_n)$ denote the subalgebra of $\ddot{U}'_{q,d}(\mathfrak{sl}_n)$ generated by $\{e_{i,0}, f_{i,0}, \psi^{\pm 1}_{i,0}\}_{i \in [n]^\times}$, which is isomorphic to $U^{\mathrm{DJ}}_q(\mathfrak{sl}_n)$–the Drinfeld-Jimbo quantum group of $\mathfrak{sl}_n$. We leave the following result as a simple exercise to the reader:

**Lemma 3.13** *(a) The subspace $M(p)_n$ is $U_q(\mathfrak{sl}_n)$-invariant and is isomorphic to $L_q(\bar{\Lambda}_p)$, the irreducible highest weight $U_q^{DJ}(\mathfrak{sl}_n)$-module of the highest weight $\bar{\Lambda}_p$.*

*(b) For any $\bar{\sigma} = \{1 \le \sigma_1 < \cdots < \sigma_p \le n\}$, let $\bar{\Lambda}_p^{\bar{\sigma}}$ be the $\mathfrak{sl}_n$-weight having entries $1 - \frac{p}{n}$ at the places $\{\sigma_i\}_{i=1}^p$ and $-\frac{p}{n}$ elsewhere. Then $\{v_0 \otimes e^{\bar{\Lambda}_p^{\bar{\sigma}}}\}_{\bar{\sigma}}$ is a basis of $M(p)_n$.*

Furthermore, we also define the degree operators $d_1, \ldots, d_{n-1}$ acting on $W(p)_n$ via:

$$d_r \left( v \otimes e^{\sum_{j=1}^{n-1} m_j \bar{\alpha}_j + \bar{\Lambda}_p} \right) = -m_r \cdot v \otimes e^{\sum_{j=1}^{n-1} m_j \bar{\alpha}_j + \bar{\Lambda}_p} \qquad \forall\, v \in F_n \,.$$

For $\underline{u} = (u_1, \ldots, u_{n-1}) \in (\mathbb{C}^\times)^{n-1}$, consider the functional

$$\phi_{\bar{p},\underline{c}}^{\underline{u}} : \ddot{U}_{q,d}^{\prime\le} \longrightarrow \mathbb{C}, \quad y \mapsto \sum_{\bar{\sigma}} \left\langle v_0 \otimes e^{\bar{\Lambda}_p^{\bar{\sigma}}} \middle| \rho_{1,\underline{c}}^p(y) u_1^{d_1} \cdots u_{n-1}^{d_{n-1}} \middle| v_0 \otimes e^{\bar{\Lambda}_p^{\bar{\sigma}}} \right\rangle$$

evaluating the $\bar{Q}$-graded trace of the $\ddot{U}_{q,d}^{\prime\le}$-action on the subspace (not a submodule!) $M(p)_n$. We note that $\rho_{1,\underline{c}}^p(h_{i,k}) v_0 \otimes e^{\bar{\Lambda}_p^{\bar{\sigma}}} = 0$ for $i \in [n], k > 0$. Therefore, similarly to the previous computation, the functional $\phi_{\bar{p},\underline{c}}^{\underline{u}}$ is determined by the generating series analogous to (3.86):

$$\phi_{p,\underline{c};N,\underline{r},a,b}^{\underline{u}}(z_{0,1}, \ldots, z_{n-1,N}) :=$$

$$\phi_{\bar{p},\underline{c}}^{\underline{u}} \left( \prod_{k=1}^N \big(f_0(z_{0,k}) \cdots f_{n-1}(z_{n-1,k})\big) \cdot \prod_{i \in [n]} \psi_{i,0}^{r_i} \cdot \gamma^{a/2} q^{b d_1} \right). \quad (3.91)$$

This series can be evaluated similarly to Proposition 3.12 by normally ordering the product $\prod_{k=1}^N (f_0(z_{0,k}) \cdots f_{n-1}(z_{n-1,k}))$ and using (3.83, 3.84) as well as (3.88, 3.89):

**Proposition 3.13** *Set $\mathfrak{c} = c_0 \cdots c_{n-1}$ and $\varepsilon = (-1)^{(n-1)(n+2)/2}$. We have:*

$$\phi_{\bar{p},\underline{c};N,\underline{r},a,b}^{\underline{u}}(z_{0,1}, \ldots, z_{n-1,N}) = (\varepsilon \mathfrak{c})^{-N} q^{a/2} d^{\frac{N(n-2)}{2}} \times$$

$$\frac{\prod_{i \in [n]} \prod_{1 \le k < k' \le N} (z_{i,k} - z_{i,k'})(z_{i,k} - q^2 z_{i,k'})}{\prod_{i \in [n]} \prod_{1 \le k \le k' \le N} (z_{i,k} - q d z_{i+1,k'}) \cdot \prod_{i \in [n]} \prod_{1 \le k < k' \le N} (z_{i,k} - q d^{-1} z_{i-1,k'})} \times$$

$$(-1)^p \prod_{j=1}^p \frac{1}{u_1 \cdots u_{j-1}} \cdot [\mu^p] \left\{ \prod_{i \in [n]} \left( \prod_{k=1}^N z_{i+1,k} - \mu \cdot u_1 \cdots u_i q^{r_{i+1} - r_i} \prod_{k=1}^N z_{i,k} \right) \right\},$$

*where the notation $[\mu^p]\{\cdots\}$ denotes the coefficient of $\mu^p$ in $\{\cdots\}$ and we use the same conventions on the products in the denominator as in Proposition 3.12.*

- *Full graded trace*

Finally, let us consider the functional

$$\phi_{\bar{p},\underline{c}}^{\underline{u},t} : \ddot{U}_{q,d}^{\prime\le} \longrightarrow \mathbb{C}[[t]] \quad \text{defined by} \quad \phi_{\bar{p},\underline{c}}^{\underline{u},t}(y) := \mathrm{Tr}_{W(p)_n} \left( \rho_{1,\underline{c}}^p(y) u_1^{d_1} \cdots u_{n-1}^{d_{n-1}} t^d \right)$$

evaluating the $\bar{Q} \times \mathbb{N}$-*graded trace* of the $\acute{U}'^{\leq}_{q,d}$-action $\rho^p_{1,\underline{\mathbf{c}}}$ on the space $W(p)_n$. Evoking the defining relations (T3, T6) and the $\bar{Q}$-grading of both $\acute{U}'^{\leq}_{q,d}$ and $W(p)_n$, we note that $\phi^{u,t}_{p,\underline{\mathbf{c}}}$ is determined by the following generating series:

$$\phi^{u,t}_{p,\underline{\mathbf{c}};N,\underline{\ell},\underline{r},a,b}(z_{0,1},\ldots,z_{n-1,N}; w_{0,1},\ldots,w_{0,\ell_0},\ldots,w_{n-1,1},\ldots,w_{n-1,\ell_{n-1}}) :=$$

$$\phi^{u,t}_{p,\underline{\mathbf{c}}}\left(\prod_{k=1}^{N}\left(f_0(z_{0,k})\cdots f_{n-1}(z_{n-1,k})\right)\cdot \prod_{i\in[n]}\prod_{k=1}^{\ell_i}\bar{\psi}^+_i(w_{i,k})\cdot \prod_{i\in[n]}\psi^{r_i}_{i,0}\cdot \gamma^{a/2}q^{bd_1}\right)$$

(3.92)

for any $N\in\mathbb{N}, a,b\in\mathbb{Z}, \underline{r}\in\mathbb{Z}^{[n]}, \underline{\ell}\in\mathbb{N}^{[n]}$ and $\bar{\psi}^+_i(z)=\psi^{-1}_{i,0}\psi^+_i(z)$ as in (3.4). The explicit formula for (3.92) involves the $t$-Pochhammer symbol $(z;t)_\infty$ defined via:

$$(z;t)_\infty := \prod_{k\geq 0}(1-zt^k).$$

**Theorem 3.16** *Set* $\mathfrak{c}=c_0\cdots c_{n-1}$ *and* $\varepsilon=(-1)^{(n-1)(n+2)/2}$. *We have:*

$$\phi^{u,t}_{p,\underline{\mathbf{c}};N,\underline{\ell},\underline{r},a,b}(z_{0,1},\ldots,z_{n-1,N}; w_{0,1},\ldots,w_{n-1,\ell_{n-1}}) = (\varepsilon\mathfrak{c})^{-N}q^{\frac{a}{2}}d^{\frac{N(n-2)}{2}}\times$$

$$\frac{\prod_{i\in[n]}\prod_{1\leq k<k'\leq N}(z_{i,k}-z_{i,k'})(z_{i,k}-q^2 z_{i,k'})\cdot \prod_{i\in[n]}\prod_{k=1}^{N}z_{i,k}}{\prod_{i\in[n]}\prod_{1\leq k\leq k'\leq N}(z_{i,k}-qdz_{i+1,k'})\cdot \prod_{i\in[n]}\prod_{1\leq k<k'\leq N}(z_{i,k}-qd^{-1}z_{i-1,k'})}\times$$

$$q^{r_p-r_0}\cdot\prod_{k=1}^{N}\frac{z_{0,k}}{z_{p,k}}\cdot\theta(\mathbf{y};\bar{\Omega})\times$$

$$\prod_{i\in[n]}\prod_{A=1}^{N}\prod_{B=1}^{\ell_i}\frac{(Tq^2\frac{q^{1/2}z_{i,A}}{w_{i,B}};T)_\infty\cdot(Tq^{-1}d\frac{q^{1/2}z_{i+1,A}}{w_{i,B}};T)_\infty\cdot(Tq^{-1}d^{-1}\frac{q^{1/2}z_{i-1,A}}{w_{i,B}};T)_\infty}{(Tq^{-2}\frac{q^{1/2}z_{i,A}}{w_{i,B}};T)_\infty\cdot(Tqd\frac{q^{1/2}z_{i,A}}{w_{i,B}};T)_\infty\cdot(Tqd^{-1}\frac{q^{1/2}z_{i-1,A}}{w_{i,B}};T)_\infty}\times$$

$$\frac{1}{(T;T)^n_\infty}\cdot\prod_{i\in[n]}\prod_{A,B=1}^{N}\frac{(T\cdot\frac{z_{i,A}}{z_{i,B}};T)_\infty\cdot(Tq^2\frac{z_{i,A}}{z_{i,B}};T)_\infty}{(Tqd\frac{z_{i+1,A}}{z_{i,B}};T)_\infty\cdot(Tqd^{-1}\frac{z_{i-1,A}}{z_{i,B}};T)_\infty},$$

($\diamond$)

*where* $T=tq^b$ *and* $\theta(\mathbf{y};\bar{\Omega})=\sum_{\mathbf{n}\in\mathbb{Z}^{n-1}}\exp(2\pi\sqrt{-1}(\frac{1}{2}\mathbf{n}\bar{\Omega}\mathbf{n}'+\mathbf{n}\mathbf{y}'))$ *is the Riemann theta function with* $'$ *denoting the transposition,* $\bar{\Omega}=\frac{1}{2\pi\sqrt{-1}}\cdot(a_{ij}\ln(T))^{n-1}_{i,j=1}$, *and*

$$\mathbf{y}=(y_1,\ldots,y_{n-1}) \text{ with } y_i=\frac{1}{2\pi\sqrt{-1}}\ln\left(u^{-1}_i T^{\delta_{pi}}q^{2r_i-r_{i-1}-r_{i+1}}\prod_{k=1}^{N}\frac{z_{i-1,k}z_{i+1,k}}{z^2_{i,k}}\right).$$

This theorem is proved similarly to Propositions 3.12 and 3.13 above, but it does require an additional step: a calculation of the trace of the normally ordered product, which gains nontrivial contributions from both factors $F_n$ and $\mathbb{C}\{\bar{Q}\}e^{\Lambda_p}$ in $W(p)_n$. To carry out this computation, let us record the following simple result:

**Lemma 3.14** *Let $\mathfrak{a}$ be a Heisenberg algebra generated by a central element $a_0$ and $\{a_k\}_{k \neq 0}$ subject to $[a_k, a_l] = \delta_{k,-l} \kappa_k a_0$ for some constants $\kappa_k$, and let $\mathfrak{a}^{\geq}$ be the subalgebra generated by $\{a_k\}_{k \geq 0}$. Consider the Fock $\mathfrak{a}$-module $F := \mathrm{Ind}_{\mathfrak{a}^{\geq}}^{\mathfrak{a}}(\mathbb{C} v_0)$ with $a_0(v_0) = v_0$, and define $\mathrm{d} \in \mathrm{End}(F)$ via $[\mathrm{d}, a_k] = -k a_k$, $\mathrm{d}(v_0) = 0$. Then:*

$$\mathrm{Tr}_F \left\{ \exp \left( \sum_{r \geq 1} x_r a_{-r} \right) \exp \left( \sum_{r \geq 1} y_r a_r \right) \cdot t^{\mathrm{d}} \right\} = \frac{1}{(t; t)_\infty} \cdot \exp \left( \sum_{r \geq 1} \frac{x_r y_r \kappa_r t^r}{1 - t^r} \right) \quad (3.93)$$

*for any constants $\{x_r, y_r\}_{r \geq 1}$.*

*Proof* The proof is based on the following simple formula for the matrix coefficients:

$$\left\langle a_{-r}^l v_0 \middle| a_{-r}^k a_r^{k'} \middle| a_{-r}^l v_0 \right\rangle = \delta_{kk'} l(l-1) \cdots (l-k+1) \kappa_r^k \ . \quad (3.94)$$

Therefore, we get:

$$\mathrm{Tr}_F \left( \exp \left( \sum_{r \geq 1} x_r a_{-r} \right) \exp \left( \sum_{r \geq 1} y_r a_r \right) t^{\mathrm{d}} \right) = \sum_{k_1, k_2, \ldots \geq 0} \mathrm{Tr}_F \left( \prod_{r \geq 1} \frac{(x_r y_r)^{k_r}}{(k_r!)^2} a_{-r}^{k_r} a_r^{k_r} t^{\mathrm{d}} \right)$$

$$= \prod_{r \geq 1} \left\{ \sum_{k_r \geq 0} \sum_{l_r \geq k_r} \frac{(x_r y_r)^{k_r}}{(k_r!)^2} \cdot \frac{l_r! \cdot \kappa_r^{k_r}}{(l_r - k_r)!} \cdot t^{r l_r} \right\} = \prod_{r \geq 1} \left\{ \sum_{k_r \geq 0} \frac{(x_r y_r \kappa_r t^r)^{k_r}}{k_r! (1 - t^r)^{k_r}} \cdot \frac{1}{1 - t^r} \right\}$$

with the formula (3.94) used in the second equality. This implies (3.93). $\qquad \square$

*Proof of Theorem 3.16* First, we use the commutation rules (3.83, 3.84) to normally order the factors in $\prod_{k=1}^N \left( f_0(z_{0,k}) \cdots f_{n-1}(z_{n-1,k}) \right) \cdot \prod_{i \in [n]} \prod_{k=1}^{\ell_i} \tilde{\psi}_i^+(w_{i,k}) \cdot \prod_{i \in [n]} \psi_{i,0}^{r_i}$. When doing so, the overall factors we gain recover exactly the product of factors from the first two lines in the right-hand side of $(\diamond)$. On the other hand, the $\bar{Q} \times \mathbb{N}$-graded trace of the resulting normally ordered product splits as $\mathrm{Tr}_1 \cdot \mathrm{Tr}_2$ with

$$\mathrm{Tr}_1 = \mathrm{Tr}_{\mathbb{C}\{\bar{Q}\} e^{\bar{\Lambda}_p}} \left( q^{\sum_{i \in [n]} r_i \partial_{\bar{\alpha}_i}} \cdot \prod_{i \in [n]} \prod_{k=1}^N z_{i,k}^{-H_{i,0}} \cdot \prod_{i=1}^{n-1} u_i^{\mathrm{d}_i} \cdot (t q^b)^{\mathrm{d}^{(2)}} \right),$$

$$\mathrm{Tr}_2 = \mathrm{Tr}_{F_n} \left( \exp \left( \sum_{i \in [n]} \sum_{r > 0} u_{i,r} H_{i,-r} \right) \cdot \exp \left( \sum_{i \in [n]} \sum_{r > 0} (v_{i,r}^{(1)} + v_{i,r}^{(2)}) H_{i,r} \right) \cdot (t q^b)^{\mathrm{d}^{(1)}} \right),$$

where

$$u_{i,r} = \frac{-q^{r/2}}{[r]_q} \sum_{1 \leq k \leq N} z_{i,k}^r, \quad v_{i,r}^{(1)} = \frac{q^{r/2}}{[r]_q} \sum_{1 \leq k \leq N} z_{i,k}^{-r}, \quad v_{i,r}^{(2)} = (q - q^{-1}) \sum_{1 \leq k \leq \ell_i} w_{i,k}^{-r},$$

and the operators $\mathrm{d}^{(1)} \in \mathrm{End}(F_n)$, $\mathrm{d}^{(2)} \in \mathrm{End}(\mathbb{C}\{\bar{Q}\} e^{\bar{\Lambda}_p})$ are defined via:

$$\mathrm{d}^{(1)}(H_{i_1, -k_1} \cdots H_{i_l, -k_l} v_0) = (k_1 + \cdots + k_l) \cdot H_{i_1, -k_1} \cdots H_{i_l, -k_l} v_0,$$

$$\mathrm{d}^{(2)}(e^{\bar{\beta}}) = \tfrac{1}{2}((\bar{\beta}, \bar{\beta}) - (\bar{\Lambda}_p, \bar{\Lambda}_p)) \cdot e^{\bar{\beta}}.$$

The computation of $\mathrm{Tr}_1$ is straightforward and recovers precisely the third line in the right-hand side of ($\diamond$). To evaluate $\mathrm{Tr}_2$, we pick a basis $\{\widetilde{H}_{i,-k}\}_{i \in [n]}$ of $\mathrm{span}_{\mathbb{C}}\{H_{j,-k}\}_{j \in [n]}$ for any $k > 0$ such that

$$[H_{j,l}, \widetilde{H}_{i,-k}] = \delta_{ij}\delta_{kl}H_0 \qquad \forall\, i, j \in [n]\,,\ k, l > 0\,, \tag{3.95}$$

as in Remark 3.14. In particular, the elements $\{H_{i,k}, \widetilde{H}_{i,-k}, H_0\}_{k>0}$ generate a Heisenberg algebra $\mathfrak{a}_i$ for any $i \in [n]$, and moreover $\mathfrak{a}_i$ commutes with $\mathfrak{a}_j$ for any $i \neq j \in [n]$.

According to (3.60), we have

$$\sum_{i \in [n]}\sum_{r>0} u_{i,r}H_{i,-r} = \sum_{i \in [n]}\sum_{r>0} \widetilde{u}_{i,r}\widetilde{H}_{i,-r} \quad \text{with} \quad \widetilde{u}_{i,r} = \sum_{j \in [n]} d^{-rm_{ij}}\frac{[r]_q[ra_{ij}]_q}{r} \cdot u_{j,r}\,.$$

Using the pairwise commutativity of $\mathfrak{a}_i$ and $\mathfrak{a}_j$ for $i \neq j$ one can rewrite $\mathrm{Tr}_2$ as a product of the corresponding traces over the Fock modules of each $\mathfrak{a}_i$. Applying Lemma 3.14 to each of these, it is now straightforward to see that $\mathrm{Tr}_2$ recovers exactly the product of factors from the last two lines in the right-hand side of ($\diamond$).

This completes the proof of Theorem 3.16.                                        $\square$

### 3.5.2 Functionals via pairing

Recall the Hopf pairing $\varphi' \colon \ddot{U}_{q,d}^{'\geq} \times \ddot{U}_{q,d}^{'\leq} \to \mathbb{C}$ from Theorem 3.4. As $\varphi'$ is non-degenerate, there exist unique elements $X_{p,\underline{c}}^0, X_{p,\underline{c}}^{\underline{u}} \in \ddot{U}_{q,d}^{'\geq,\wedge}$ and $X_{p,\underline{c}}^{\underline{u},t} \in \ddot{U}_{q,d}^{'\geq,\wedge}[[t]]$ (with $\wedge$ denoting an appropriate completion) such that

$$\phi_{p,\underline{c}}^0(y) = \varphi'(X_{p,\underline{c}}^0, y), \quad \phi_{p,\underline{c}}^{\underline{u}}(y) = \varphi'(X_{p,\underline{c}}^{\underline{u}}, y), \quad \phi_{p,\underline{c}}^{\underline{u},t}(y) = \varphi'(X_{p,\underline{c}}^{\underline{u},t}, y) \tag{3.96}$$

for any $y \in \ddot{U}_{q,d}^{'\leq}$. We shall now obtain explicit shuffle formulas for these elements.

To this end, we extend the isomorphism $\Psi$ of Theorem 3.6 to the isomorphism

$$\Psi^{\geq} \colon \ddot{U}_{q,d}^{'\geq} \xrightarrow{\sim} S^{\geq}\,. \tag{3.97}$$

Here, the *extended shuffle algebra* $S^{\geq}$ is obtained from $S$ by adjoining generators $\{\psi_{i,-s}, \psi_{i,0}^{-1}, \gamma^{\pm 1/2}, q^{\pm d_1}\}_{i \in [n]}^{s \in \mathbb{N}}$ with defining relations compatible with those for $\ddot{U}_{q,d}^{'\geq}$. In particular, we have:

$$q^{d_1}F q^{-d_1} = q^{-d} \cdot F\,, \tag{3.98}$$

$$\psi_i^-(z) \star F = \left[ F\left(\{x_{j,r}\}_{j \in [n]}^{1 \leq r \leq k_j}\right) \cdot \prod_{j \in [n]}^{r \leq k_j} \frac{\zeta_{i,j}(z/x_{j,r})}{\zeta_{j,i}(x_{j,r}/z)} \right] \star \psi_i^-(z) \tag{3.99}$$

for any $F \in S_{\underline{k},d}$, where we set $\psi_i^-(z) := \sum_{s \geq 0} \psi_{i,-s}^- z^s$, $\star$ denotes the multiplication in $S^{\geq}$, and the $\zeta$-factors in (3.99) are all expanded in the nonnegative powers of $z$.

Let $\Gamma_{p,\underline{c}}^0, \Gamma_{p,\underline{c}}^{\underline{u}}, \Gamma_{p,\underline{c}}^{\underline{u},t}$ denote the images of $X_{p,\underline{c}}^0, X_{p,\underline{c}}^{\underline{u}}, X_{p,\underline{c}}^{\underline{u},t}$ under $\Psi^{\geq}$ of (3.97):

$$\Gamma^0_{p,\underline{c}} := \Psi^\geq(X^0_{p,\underline{c}}), \quad \Gamma^{\underline{u}}_{p,\underline{c}} := \Psi^\geq(X^{\underline{u}}_{p,\underline{c}}), \quad \Gamma^{\underline{u},t}_{p,\underline{c}} := \Psi^\geq(X^{\underline{u},t}_{p,\underline{c}}). \tag{3.100}$$

The following is the main result of this section ($\delta = (1, 1, \ldots, 1) \in \mathbb{N}^{[n]}$ as before):

**Theorem 3.17** *(a) We have*

$$\Gamma^0_{p,\underline{c}} = \sum_{N \geq 0} (\varepsilon c)^{-N} (1 - q^{-2})^{nN} (-q^n d^{-n/2})^{N^2} \cdot \Gamma^0_{p;N} \cdot q^{\tilde{\Lambda}_p} q^{-d_1}$$

*with $\Gamma^0_{p;N} \in S_{N\delta}$ defined in (3.50), and $c = c_0 \cdots c_{n-1}$, $\varepsilon = (-1)^{(n-1)(n+2)/2}$ as before.*

*(b) We have*

$$\Gamma^{\underline{u}}_{p,\underline{c}} = \sum_{N \geq 0} (\varepsilon c)^{-N} (1 - q^{-2})^{nN} (-q^n d^{-n/2})^{N^2} \cdot \Gamma^{\underline{u}}_{p;N} \cdot q^{-d_1}$$

*with $\Gamma^{\underline{u}}_{p;N} \in S^\geq_{N\delta}$ given by:*

$$\Gamma^{\bar{u}}_{p;N} = \prod_{j=1}^p \frac{1}{u_1 \cdots u_{j-1}} \cdot \frac{\prod_{i \in [n]} \prod_{1 \leq r \neq r' \leq N} (x_{i,r} - q^{-2} x_{i,r'})}{\prod_{i \in [n]} \prod_{1 \leq r, r' \leq N} (x_{i,r} - x_{i+1,r'})} \times$$

$$(-1)^p [\mu^p] \left\{ \prod_{i \in [n]} \left( \prod_{r=1}^N x_{i+1,r} - \mu \cdot u_1 \cdots u_i \prod_{r=1}^N x_{i,r} \cdot q^{\tilde{\Lambda}_{i+1} - \tilde{\Lambda}_i} \right) \right\},$$

*where in the last product we take all $x_{i,r}$'s to the left and all $q^{\tilde{\Lambda}_{i+1} - \tilde{\Lambda}_i}$ to the right.*

*(c) We have*

$$\Gamma^{\underline{u},t}_{p,\underline{c}} = \sum_{N \geq 0} (\varepsilon c)^{-N} (1 - q^{-2})^{nN} (-q^n d^{-n/2})^{N^2} \cdot \Gamma^{\underline{u},t}_{p;N} \cdot q^{\tilde{\Lambda}_p} q^{-d_1}$$

*with $\Gamma^{\underline{u},t}_{p;N} \in (S^\geq_{N\delta})^\wedge[[t]]$ given by:*

$$\Gamma^{\underline{u},t}_{p;N} = \frac{\prod_{i \in [n]} \prod_{1 \leq r \neq r' \leq N} (x_{i,r} - q^{-2} x_{i,r'}) \cdot \prod_{i \in [n]} \prod_{r=1}^N x_{i,r}}{\prod_{i \in [n]} \prod_{1 \leq r, r' \leq N} (x_{i,r} - x_{i+1,r'})} \cdot \prod_{r=1}^N \frac{x_{0,r}}{x_{p,r}} \cdot \theta(\mathbf{x}; \tilde{\Omega}) \times$$

$$\frac{1}{(\bar{t}; \bar{t})^n_\infty} \cdot \prod_{i \in [n]} \prod_{\mathsf{A},\mathsf{B}=1}^N \frac{(\bar{t} \frac{x_{i,\mathsf{A}}}{x_{i,\mathsf{B}}}; \bar{t})_\infty \cdot (\bar{t} q^2 \frac{x_{i,\mathsf{A}}}{x_{i,\mathsf{B}}}; \bar{t})_\infty}{(\bar{t} q d \frac{x_{i+1,\mathsf{A}}}{x_{i,\mathsf{B}}}; \bar{t})_\infty \cdot (\bar{t} q d^{-1} \frac{x_{i-1,\mathsf{A}}}{x_{i,\mathsf{B}}}; \bar{t})_\infty} \cdot \prod_{k>0} \prod_{i \in [n]} \prod_{r=1}^N \bar{\psi}^-_i (\bar{t}^k q^{1/2} x_{i,r}),$$

*where $\bar{t} = t\gamma^{-1}$, $\tilde{\Omega} = \frac{1}{2\pi\sqrt{-1}} \cdot (a_{ij} \ln(\bar{t}))^{n-1}_{i,j=1}$, and*

$$\mathbf{x} = (x_1, \ldots, x_{n-1}) \quad \text{with} \quad x_i = \frac{1}{2\pi\sqrt{-1}} \ln \left( u_i^{-1} \bar{t}^{\delta_{pi}} \psi_{i,0} \prod_{r=1}^N \frac{x_{i-1,r} x_{i+1,r}}{x_{i,r}^2} \right).$$

*In the above products, we take all $x_{i,r}$'s to the left and all $\psi_{i,k}$'s to the right.*

This theorem follows directly from Propositions 3.12, 3.13 and Theorem 3.16 by applying the following technical lemma:

**Lemma 3.15** *(a) We have*

$$\varphi'(aa', bb') = \varphi'(a,b) \cdot \varphi'(a',b')$$

*for any elements* $a \in \ddot{U}_{q,d}^{'>}$, $a' \in \ddot{U}_{q,d}^{'\geq} \cap \ddot{U}_{q,d}^{'0}$, $b \in \ddot{U}_{q,d}^{'<}$, $b' \in \ddot{U}_{q,d}^{'\leq} \cap \ddot{U}_{q,d}^{'0}$.

*(b) For any* $k_i, k_i' \in \mathbb{N}$, $a_i, a_i' \in \mathbb{Z}$ *(with* $i \in [n]$*), and* $A, B, A', B' \in \mathbb{Z}$, *we have:*

$$\varphi'\left( \prod_{i \in [n]}^{r \leq k_i} \bar{\psi}_i^-(z_{i,r}) \prod_{i \in [n]} \psi_{i,0}^{a_i} \gamma^{A/2} q^{Bd_1}, \prod_{j \in [n]}^{r' \leq k_j'} \bar{\psi}_j^+(w_{j,r'}) \prod_{j \in [n]} \psi_{j,0}^{a_j'} \gamma^{A'/2} q^{B'd_1} \right) =$$

$$q^{-\frac{1}{2}A'B - \frac{1}{2}AB' + \sum_{i,j \in [n]} a_i a_j' a_{ij}} \cdot \prod_{i \in [n]}^{j \in [n]} \prod_{1 \leq r \leq k_i}^{1 \leq r' \leq k_j'} \frac{w_{j,r'} - q^{a_{ij}} d^{m_{ij}} z_{i,r}}{w_{j,r'} - q^{-a_{ij}} d^{m_{ij}} z_{i,r}} \, .$$

*(c) For any* $k \in \mathbb{N}$, $i_1, \dots, i_k \in I$, $d_1, \dots, d_k \in \mathbb{Z}$, $X \in \ddot{U}_{q,d}^{'>}$, *we have:*

$$\varphi'(X, f_{i_1, -d_1} \cdots f_{i_k, -d_k}) =$$

$$(q - q^{-1})^{-k} \int_{|z_1| \gg \cdots \gg |z_k|} \frac{\Psi(X)(z_1, \dots, z_k) z_1^{-d_1} \cdots z_k^{-d_k}}{\prod_{1 \leq a < b \leq k} \zeta_{i_a, i_b}(z_a/z_b)} \prod_{r=1}^k \frac{dz_r}{2\pi\sqrt{-1} z_r} \, , \qquad (3.101)$$

*where each variable* $z_a$ *is plugged into an argument of color* $i_a$ *of the function* $\Psi(X)$, *and the pairing vanishes if* $\deg(X) \neq (1_{i_1} + \cdots + 1_{i_k}, d_1 + \cdots + d_k) \in \mathbb{N}^{[n]} \times \mathbb{Z}$.

*(d) Using the notations of part (c), we have:*

$$\varphi'(X, f_{i_1}(z_1) \cdots f_{i_k}(z_k)) = (q - q^{-1})^{-k} \frac{\Psi(X)(z_1, \dots, z_k)}{\prod_{1 \leq a < b \leq k} \zeta_{i_a, i_b}(z_a/z_b)} \qquad (3.102)$$

*with the right-hand side of (3.102) expanded in the domain* $|z_1| \gg \cdots \gg |z_k|$.

*Proof* The equalities in parts (a, b) follow immediately from the formulas (3.7), (H1), and (P1)–(P6). Parts (c, d) are due to [20, Proposition 3.9].                    □

### 3.5.3  Transfer matrices and Bethe subalgebras

In this section, we recall a standard construction of Bethe subalgebras in any quasi-triangular Hopf algebra and discuss modifications needed for the *formal* versions.

Let $A$ be a quasitriangular Hopf algebra with a coproduct $\Delta$ and a universal $R$-matrix $R \in A \otimes A$ satisfying (3.9), and let $x \in A$ be a group-like element:

$$\Delta(x) = x \otimes x \, . \qquad (3.103)$$

For any finite dimensional $A$-representation $(V, \pi_V)$, we define the corresponding *L-operator* and the *transfer matrix* via:

$$L_V(x) := (\mathrm{Id} \otimes \pi_V)((1 \otimes x)R) \in A \otimes \mathrm{End}(V),$$
$$T_V(x) := (\mathrm{Id} \otimes \mathrm{Tr}_V)(L_V(x)) \in A.$$
(3.104)

The following is standard:

**Lemma 3.16** *(a)* $T_{V_1 \otimes V_2}(x) = T_{V_2}(x) \cdot T_{V_1}(x)$.

*(b)* $T_V(x) = T_{V_1}(x) + T_{V_2}(x)$ *for any short exact sequence* $0 \to V_1 \to V \to V_2 \to 0$.

*Proof* Part (a) follows from (3.103) and $(\mathrm{Id} \otimes \Delta)(R) = R^{13}R^{12}$ of (3.9):

$$T_{V_1 \otimes V_2}(x) = (\mathrm{Id} \otimes \mathrm{Tr}_{V_1 \otimes V_2})(1 \otimes \pi_{V_1 \otimes V_2})((1 \otimes x)R) =$$

$$(\mathrm{Id} \otimes \mathrm{Tr}_{V_1} \otimes \mathrm{Tr}_{V_2})(\mathrm{Id} \otimes \pi_{V_1} \otimes \pi_{V_2})\Big((1 \otimes x \otimes x)(R^{13}R^{12})\Big) = T_{V_2}(x) \cdot T_{V_1}(x).$$

Part (b) is obvious from the definition (3.104). □

**Lemma 3.17** *The transfer matrices* $\{T_V(x)\}_V$ *pairwise commute:*

$$T_{V_1}(x) \cdot T_{V_2}(x) = T_{V_2}(x) \cdot T_{V_1}(x).$$
(3.105)

*Proof* Define $R_{V_1, V_2} := (\pi_{V_1} \otimes \pi_{V_2})(R) \in \mathrm{End}(V_1) \otimes \mathrm{End}(V_2)$ and consider the permutation operator $P \colon V_1 \otimes V_2 \to V_2 \otimes V_1$ acting via $v_1 \otimes v_2 \mapsto v_2 \otimes v_1$. Then:

$$P R_{V_1, V_2}(\pi_{V_1} \otimes \pi_{V_2})(\Delta(y)) = (\pi_{V_2} \otimes \pi_{V_1})(\Delta(y)) P R_{V_1, V_2} \quad \forall y \in A,$$

due to $R\Delta(y) = \Delta^{\mathrm{op}}(y)R$ of (3.9). This implies $T_{V_1 \otimes V_2}(x) = T_{V_2 \otimes V_1}(x)$ and hence the claimed commutativity $[T_{V_1}(x), T_{V_2}(x)] = 0$, due to Lemma 3.16(a). □

As an immediate consequence of Lemmas 3.16–3.17, we get that the assignment $V \mapsto T_V(x)$ gives rise to a ring homomorphism from the Grothendieck group of the category of finite dimensional $A$-modules onto a commutative subalgebra $\mathcal{B}_x$ of $A$.

**Definition 3.4** These commutative subalgebras $\mathcal{B}_x$ of $A$ are called *Bethe subalgebras*.

In the present setup, we need a generalization of this construction to a quasitriangular topological Hopf algebra $A$ with $\Delta \colon A \to A \widehat{\otimes} A$, $R \in A \widehat{\otimes} A$, $x \in A^\wedge$ being in appropriate completions and a representation $V$ being not necessarily finite dimensional but for which the trace (3.104) is a well-defined element of the completion $A^\wedge$. Then, the above construction provides a commutative Bethe subalgebra $\mathcal{B}_x$ of $A^\wedge$.

*Remark 3.15* We note that such a treatment for quantum affinized algebras is usually accomplished by introducing an additional spectral parameter $z$. To this end, one replaces a representation $V$ with $V(z)$, the twist of $V$ by the algebra automorphism

$$\mathrm{Adj}(z^{-d_1}) \colon e_{i,r} \mapsto z^r e_{i,r}, \; f_{i,r} \mapsto z^r f_{i,r}, \; \psi_{i,r} \mapsto z^r \psi_{i,r}, \; \gamma^{1/2} \mapsto \gamma^{1/2}, \; q^{d_1} \mapsto q^{d_1}.$$

Thus, (3.104) will get replaced with $L_V(x; z)$ and $T_V(x; z)$ that are series in $z$. We refer the interested reader to [14, §3.1] for more details in the case of $A = U_q(\widehat{\mathfrak{g}})$, $x = q^{2\rho}$.

### 3.5.4 Bethe algebra incarnation of $\mathcal{A}(s_0, \ldots, s_{n-1})$

Combining (3.96) and (3.104) allows us to identify $X_{p,\underline{c}}^{u,t}$ with the transfer matrices:

$$X_{p,\underline{c}}^{u,t} = T_{\rho_{1,\underline{c}}^p}\left(u_1^{-\bar{\Lambda}_1} \cdots u_{n-1}^{-\bar{\Lambda}_{n-1}} t^{d_1}\right) \cdot \prod_{j=1}^{n-1} u_j^{\langle \bar{\Lambda}_j, \bar{\Lambda}_p \rangle}. \qquad (3.106)$$

Furthermore, the other two families of elements can be recovered as limits of (3.106):

- $X_{p,\underline{c}}^{u}$ are the $t^0$-coefficients of $X_{p,\underline{c}}^{u,t}$;

- $X_{p,\underline{c}}^{0}$ are the further limits of $X_{p,\underline{c}}^{u}$ as $u_1, \ldots, u_{n-1} \to 0$.

Combining (3.100) with Theorem 3.17 and the commutativity of the Bethe subalgebras $\mathcal{B}_x$, see Lemma 3.17, we obtain:

**Proposition 3.14** *(a) The elements* $\{\Gamma_{p;N}^{u,t}\}_{p \in [n]}^{N \geq 1}$ *pairwise commute.*

*(b) The elements* $\{\Gamma_{p;N}^{u}\}_{p \in [n]}^{N \geq 1}$ *pairwise commute.*

*(c) The elements* $\{\Gamma_{p;N}^{0}\}_{p \in [n]}^{N \geq 1}$ *pairwise commute.*

Due to Theorem 3.17(b), the elements $\Gamma_{p;N}^{u}$ have the same form as the generators of the subalgebra $\mathcal{A}(s_0, \ldots, s_{n-1})$ from Section 3.2.2 with the parameters

$$s_i = u_i \cdot q^{\bar{\Lambda}_{i+1} - 2\bar{\Lambda}_i + \bar{\Lambda}_{i-1}} \in \mathbb{C}^\times \cdot e^{\tilde{\rho}} \quad \forall\, i \in [n], \qquad (3.107)$$

where $u_0 := \frac{1}{u_1 \cdots u_{n-1}}$. As $\{e^h\}_{h \in \tilde{\rho}}$ commute with $\bigoplus_k S_{k\delta}$, we see that $\{s_i\}$ of (3.107) can be treated as formal parameters with $s_0 \cdots s_{n-1} = 1$ and being generic (3.37). This immediately implies Lemma 3.6 (thus completing our proof of Theorem 3.7):

**Corollary 3.4** *The elements* $\{F_k^\mu(\underline{s})\}_{k \geq 1}^{\mu \in \mathbb{C}}$ *of Lemma 3.2 pairwise commute.*

Let us recall the horizontal subalgebra $U_q^h(L\mathfrak{gl}_n)$ of $\ddot{U}_{q,d}'(\mathfrak{sl}_n)$ from Remark 3.5(a). The subspace $M(p)_n$ from Lemma 3.13 is clearly $U_q^h(L\mathfrak{gl}_n)$-invariant and is isomorphic to the $p$-th fundamental representation of the quantum loop algebra $U_q(L\mathfrak{gl}_n)$. Realizing $U_q(L\mathfrak{gl}_n)$ as the Drinfeld double (but now using the Drinfeld-Jimbo triangular decomposition!), we thus obtain the Bethe algebra realization [12] of $\mathcal{A}(\underline{s})$:

**Theorem 3.18** *The Bethe subalgebra* $\mathcal{B}_x$ *of* $U_q(L\mathfrak{gl}_n)$ *corresponding to the group-like element* $x = u_1^{-\bar{\Lambda}_1} \cdots u_{n-1}^{-\bar{\Lambda}_{n-1}}$ *and the category of finite dimensional* $U_q(L\mathfrak{gl}_n)$-*modules is naturally identified with* $\mathcal{A}(\{u_i \cdot q^{\bar{\Lambda}_{i+1} - 2\bar{\Lambda}_i + \bar{\Lambda}_{i-1}}\}_{i \in [n]})$ *from Definition 3.1.*

*Remark 3.16* (a) According to Theorem 3.8, the elements $\{\Gamma_{p;N}^{0}\}_{p \in [n]}^{N \geq 1}$ featuring in Theorem 3.17(a) generate the commutative subalgebra $\mathcal{A}^h \subset S$ from Definition 3.2.

(b) In contrast, the commutative subalgebras generated by $\{\Gamma_{p;N}^{u,t}\}_{p \in [n]}^{N \geq 1}$ can be rather viewed as one-parameter deformations of the algebras $\mathcal{A}(\underline{s})$.

*Remark 3.17* Applying the same analysis to $\ddot{U}'_{q_1,q_2,q_3}(\mathfrak{gl}_1)$-modules $\rho_{1,c}$ from Proposition 2.11, one obtains the transfer matrix realization of the elements $K_k$ from (2.37) and their one-parameter deformations. Explicitly, we have:

- The *top-to-top* functional $\phi^0_c \colon \ddot{U}'^{\le}_{q_1,q_2,q_3} \to \mathbb{C}$ given by $y \mapsto \langle v_0 | \rho_{1,c}(y) | v_0 \rangle$ corresponds to the pairing with

$$\Gamma^0_c = \sum_{N=0}^{\infty} c^{-N} q^{-N(N-1)} \cdot K_N(x_1, \ldots, x_N) \cdot q^{-d_1} \,.$$

- The *graded trace* functional $\phi^t_c \colon \ddot{U}'^{\le}_{q_1,q_2,q_3} \to \mathbb{C}[[t]]$ given by $y \mapsto \mathrm{Tr}_F(\rho_{1,c}(y)t^d)$ corresponds to the pairing with

$$\Gamma^t_c = \sum_{N=0}^{\infty} \frac{c^{-N} q^{-N(N-1)}}{(t; t)_\infty} \cdot \Gamma^t_N \cdot q^{-d_1} \,,$$

where

$$\Gamma^t_N = K_N(x_1, \ldots, x_N) \cdot \prod_{A,B=1}^{N} \frac{(t\frac{x_A}{x_B}; t)_\infty \cdot (tq^2 \frac{x_A}{x_B}; t)_\infty}{(tqd^{-1}\frac{x_A}{x_B}; t)_\infty \cdot (tqd\frac{x_A}{x_B}; t)_\infty} \cdot \prod_{k>0} \prod_{r=1}^{N} \tilde{\psi}^-(t^k q^{1/2} x_r) \,.$$

Similarly to (3.106), it can be realized as the transfer matrix $T_{\rho_{1,c}}(t^d)$ of (3.104).

### 3.5.5 Shuffle formulas for Miki's twist of Cartan generators

For any $r > 0$, let $\{\widetilde{h}_{i,-r}\}_{i \in [n]}$ be a basis of $\mathrm{span}_{\mathbb{C}}\{h_{j,-r}\}_{j \in [n]}$ dual to $\{h_{i,r}\}_{i \in [n]}$ with respect to the pairing $\underline{\varphi}'$ of Theorem 3.4: $\varphi'(\widetilde{h}_{i,-r}, h_{j,s}) = \delta_{ij}\delta_{rs}$ for any $i, j \in [n]$, $r, s > 0$. In particular, $\widetilde{h}_{0,-r}$ is a multiple of $h^\vee_{-r}$ from (3.25). The following result was proved in [23, Theorem 4.5] (while a shorter proof was recently provided in [25]):

$$X^0_{p,\underline{c}} = \exp\left(\sum_{r \ge 1} \frac{[r]_q}{r} (qu)^r \varpi(\widetilde{h}_{p,-r})\right) \cdot q^{\hat{\Lambda}_p} q^{-d_1} \tag{3.108}$$

with $u = \frac{(-1)^{(n-2)(n-3)/2}}{qd^{n/2}c_0 \cdots c_{n-1}}$ as in Theorem 3.15. Comparing this with the formula for $\Gamma^0_{p,\underline{c}} = \Psi^{\ge}(X^0_{p,\underline{c}})$ from Theorem 3.17(a), we finally obtain:

**Theorem 3.19** *The following equality holds:*

$$\exp\left(\sum_{r \ge 1} c^{-r} \cdot \frac{d^{-rn/2}[r]_q}{r} \cdot \Psi\left(\varpi(\widetilde{h}_{p,-r})\right)\right) = \\ \sum_{N \ge 0} c^{-N} \cdot (q^{-2} - 1)^{nN} (-q^n d^{-n/2})^{N^2} \cdot \Gamma^0_{p;N} \,. \tag{3.109}$$

# References

1. J. Beck: Braid group action and quantum affine algebras. Commun. Math. Phys. **165**, no. 3, 555–568 (1994)

2. V. Chari, N. Jing: Realization of level one representations of $U_q(\widehat{\mathfrak{g}})$ at a root of unity. Duke Math. J. **108**, no. 1, 183–197 (2001)

3. I. Damiani: From the Drinfeld realization to the Drinfeld-Jimbo presentation of affine quantum algebras: injectivity. Publ. Res. Inst. Math. Sci. **51**, no. 1, 131–171 (2015)

4. J. Ding, I. Frenkel: Isomorphism of two realizations of quantum affine algebra $U_q(\widehat{\mathfrak{gl}(n)})$. Commun. Math. Phys. **156**, no. 2, 277–300 (1993)

5. J. Ding, K. Iohara: Generalization of Drinfeld quantum affine algebras. Lett. Math. Phys. **41**, no. 2, 181–193 (1997)

6. V. Drinfeld: Quantum groups. Proceedings of the International Congress of Mathematicians **1**, no. 2 (Berkeley, Calif., 1986), 798–820, Amer. Math. Soc., Providence, RI (1987)

7. V. Drinfeld: A New realization of Yangians and quantized affine algebras. Sov. Math. Dokl. **36**, no. 2, 212–216 (1988)

8. B. Feigin, K. Hashizume, A. Hoshino, J. Shiraishi, S. Yanagida: A commutative algebra on degenerate $\mathbb{CP}^1$ and Macdonald polynomials. J. Math. Phys. **50**, no. 9, Paper No. 095215 (2009)

9. B. Feigin, M. Jimbo, T. Miwa, E. Mukhin: Representations of quantum toroidal $\mathfrak{gl}_n$. J. Algebra **380**, 78–108 (2013)

10. B. Feigin, M. Jimbo, T. Miwa, E. Mukhin: Quantum toroidal $\mathfrak{gl}_1$ and Bethe ansatz. J. Phys. A **48**, no. 24, Paper No. 244001 (2015)

11. B. Feigin, T. Kojima, J. Shiraishi, H. Watanabe: The integrals of motion for the deformed $W$-algebra $W_{q,t}(\widehat{\mathfrak{sl}}_N)$. Proceedings of Representation Theory 2006, Atami, Japan, 102–114 (2006)

12. B. Feigin, A. Tsymbaliuk: Bethe subalgebras of $U_q(\widehat{\mathfrak{gl}}_n)$ via shuffle algebras. Selecta Math. (N. S.) **22**, no. 2, 979–1011 (2016)

13. I. Frenkel, N. Jing: Vertex representations of quantum affine algebras. Proc. Nat. Acad. Sci. U.S.A. **85**, no. 24, 9373–9377 (1988)

14. E. Frenkel, N. Reshetikhin: The $q$-characters of representations of quantum affine algebras and deformations of $\mathcal{W}$-algebras. Recent developments in quantum affine algebras and related topics (Raleigh, NC), 163–205, Contemp. Math. **248**, Amer. Math. Soc., Providence, RI (1999)

15. M. Jimbo: A $q$-analogue of $U(\mathfrak{gl}(N+1))$, Hecke algebra and the Yang-Baxter equation. Lett. Math. Phys. **11**, no. 3, 247–252 (1986)

16. N. Jing: On Drinfeld realization of quantum affine algebras. In Monster and Lie Algebras, Columbus, OH **7**, 195–206 (1996)

17. C. Kassel, M. Rosso, V. Turaev: Quantum groups and knot invariants. Panoramas et Synthèses [Panoramas and Syntheses], 5. Société Mathématique de France, Paris (1997)

18. K. Miki: Toroidal braid group action and an automorphism of toroidal algebra $U_q(\mathfrak{sl}_{n+1,\mathrm{tor}})$ ($n \geq 2$). Lett. Math. Phys. **47**, no. 4, 365–378 (1999)

19. K. Miki: Representations of quantum toroidal algebra $U_q(\mathfrak{sl}_{n+1,\mathrm{tor}})$ ($n \geq 2$). J. Math. Phys. **41**, no. 10, 7079–7098 (2000)

20. A. Neguț: Quantum toroidal and shuffle algebras. Adv. Math. **372**, Paper No. 107288 (2020)

21. A. Neguț, A. Tsymbaliuk: Quantum loop groups and shuffle algebras via Lyndon words. Preprint, Ar$\chi$iv:2102.11269 (2021).

22. Y. Saito: Quantum toroidal algebras and their vertex representations. Publ. Res. Inst. Math. Sci. **34**, no. 2, 155–177 (1998)

23. A. Tsymbaliuk: Several realizations of Fock modules for toroidal $\ddot{U}_{q,d}(\mathfrak{sl}_n)$. Algebr. Represent. Theory **22**, no. 1, 177–209 (2019)

24. M. Varagnolo, E. Vasserot: Schur duality in the toroidal setting. Commun. Math. Phys. **182**, no. 2, 469–483 (1996)

25. J. Wen: Wreath Macdonald polynomials as eigenstates. Preprint, Ar$\chi$iv:1904.05015 (2019)